\documentclass[a4paper,12pt]{article}
\usepackage{amsmath,amsthm,verbatim,amssymb,amsfonts,amscd, graphicx}
\usepackage[all]{xy}
\usepackage{xcolor}
\usepackage{geometry}
\usepackage{hyperref}
\usepackage{cleveref}
\usepackage{enumitem}
\usepackage{tocloft}
\usepackage{indentfirst}
\usepackage{caption}
\usepackage{tikz}
\usetikzlibrary{arrows.meta,positioning,fit,backgrounds}
\definecolor{darkblue}{rgb}{0,0,0.3}
\definecolor{urlblue}{rgb}{0,0,0.7}
\hypersetup{
	colorlinks=true,
	citecolor=blue,
	linkcolor=darkblue,
	urlcolor=urlblue,
}

\newcommand{\RR}{\mathbb{R}}
\newcommand{\ZZ}{\mathbb{Z}}

\DeclareMathOperator{\Ree}{Re}

\DeclareMathOperator{\id}{id}

\let\div\undefined

\DeclareMathOperator{\div}{div}

\DeclareMathOperator{\Lip}{Lip}
\DeclareMathOperator{\loc}{loc}

\newcommand{\D}{\nabla}
\newcommand{\p}{\partial}
\renewcommand{\H}{\mathcal{H}}
\renewcommand{\L}{\mathcal{L}}
\renewcommand{\th}{\theta}
\newcommand{\metric}[2]{\langle#1\,,\,#2\rangle}

\renewcommand{\bar}{\overline}
\renewcommand{\tilde}{\widetilde}
\renewcommand{\epsilon}{\varepsilon}
\renewcommand{\leq}{\leqslant}
\renewcommand{\geq}{\geqslant}

\newcommand{\updatetag}[2]{}

\newtheorem{theorem}{Theorem}[section]
\newtheorem{lemma}[theorem]{Lemma}
\newtheorem{cor}[theorem]{Corollary}

\newtheorem{defn}[theorem]{Definition}
\numberwithin{equation}{section}
\theoremstyle{definition}

\theoremstyle{definition}
\newtheorem{remark}[theorem]{Remark}
\theoremstyle{definition}
\newtheorem{example}[theorem]{Example}
\theoremstyle{definition}
\newtheorem{claim}{Claim}

\newtheorem{claimheatmain}{Claim}
\newtheorem{claimheatreg}{Claim}

\renewcommand{\SS}{\mathbb{S}}

\newcommand{\bg}{\bar{g}}
\newcommand{\bnu}{\bar{\nu}}
\newcommand{\bw}{\bar{w}}
\newcommand{\cD}{\mathcal{D}}

\newcommand{\cJ}{\mathcal{J}}

\newcommand{\oT}{{\overline{T}}}
\newcommand{\oth}{\overline{\th}}

\newcommand{\uT}{{\underline{T}}}
\newcommand{\uth}{\underline{\th}}

\newcommand{\rM}{\mathring{M}}
\newcommand{\rth}{\mathring{\th}}

\newcommand{\sg}{\mathsf{g}}
\newcommand{\tsg}{\tilde{\mathsf{g}}}
\newcommand{\tg}{\tilde{g}}
\newcommand{\tw}{\tilde{w}}
\newcommand{\tZ}{\tilde{Z}}
\DeclareMathOperator{\BV}{BV}
\DeclareMathOperator{\curl}{curl}
\DeclareMathOperator{\Int}{Int}
\DeclareMathOperator{\graphh}{graph}

\newcommand{\Th}{\Theta}
\newcommand{\bmetric}[2]{{\big\langle{#1}\,,\,{#2}\big\rangle}}
\newcommand{\Bmetric}[2]{{\Big\langle{#1}\,,\,{#2}\Big\rangle}}
\let\P\undefined
\newcommand{\P}[1]{{P\big(#1\big)}}
\newcommand{\Ps}[1]{{P(#1)}}
\DeclareMathOperator{\Arg}{Arg}

\DeclareMathOperator{\Crit}{Crit}
\DeclareMathOperator{\dom}{dom}

\DeclareMathOperator{\esssup}{esssup}
\DeclareMathOperator{\osc}{osc}
\DeclareMathOperator{\reg}{reg}
\DeclareMathOperator{\ridge}{ridge}

\newcommand{\varth}{\vartheta}
\newcommand{\const}{\text{const}}
\let\Int\undefined
\DeclareMathOperator{\Int}{Int}
\def\Xint#1{\mathchoice
	{\XXint\displaystyle\textstyle{#1}}%
	{\XXint\textstyle\scriptstyle{#1}}%
	{\XXint\scriptstyle\scriptscriptstyle{#1}}%
	{\XXint\scriptscriptstyle\scriptscriptstyle{#1}}%
	\!\int}
\def\XXint#1#2#3{{\setbox0=\hbox{$#1{#2#3}{\int}$}
		\vcenter{\hbox{$#2#3$}}\kern-.5\wd0}}

\def\dashint{\Xint-}

\newcommand{\Phragmenlindelof}{{Phragm\'en-Lindel\"of }}
\newcommand{\termI}{\operatorname{I}}
\newcommand{\termII}{\operatorname{II}}

\newcommand{\arrowangle}[1]{\rotatebox[origin=c]{#1}{$\to$}}
\newcommand{\side}{\text{side}}


\begin{document}
	
\title{Infinity-harmonic functions and inverse mean curvature flow clusters}
\author{Kai Xu}
\date{}
\maketitle

\begin{abstract}
	An $\infty$-harmonic function is a viscosity solution of $\D^2 u(\D u,\D u)=0$, or equivalently, an absolute minimizer of $\|\D u\|_{L^\infty}$. We prove a variety of new structural and regularity results in two dimensions, including:
	\begin{itemize}[nosep]
		\item $\infty$-harmonic functions in domains of $\RR^2$ are $C^{1,1/3}$.
		\item Critical points are isolated, and at each critical point, the solution has a unique quasiradial blow-up.
		\item Entire solutions with polynomial growth have unique quasiradial blow-downs, and are determined by their Fourier modes at infinity.
	\end{itemize}
	These results are consequences of a new theory relating $\infty$-harmonic functions to inverse mean curvature flow (IMCF) clusters -- which are piecewise weak solutions of IMCF with common obstacle-type boundary conditions on the interfaces (a simple example is an embedded family of cuspidal curves evolving by inverse curvature). This connection arises as the $p\to\infty$ limit of the classical duality between $p$-harmonic and $q$-harmonic functions in $\RR^2$, where $\frac1p+\frac1q=1$.
\end{abstract}

\tableofcontents
	
\section{Introduction}\label{sec:intro}

% MSC 2020:
	% 35J60, 53E10; 53B65, 35J70, 35B08

% Figures: (42)
	% intro: (5)
		% fig_intro_imcf.eps
		% fig_intro_deg2.eps
		% fig_intro_quasiradial.eps
		% fig_intro_angle_param.eps
		% fig_intro_angle_splitting.eps
	% prelim: (2)
		% fig_prelim_models.eps
		% fig_prelim_solitons.eps
	% cluster: (4)
		% fig_cluster_excep_pt.eps
		% fig_cluster_mixed1.eps
		% fig_cluster_Yt.eps
		% fig_cluster_tZt.eps
	% ex: (7)
		% fig_entire_e1.eps
		% fig_ex_potential_square_2.eps
		% fig_ex_splitting.eps
		% fig_entire_31.eps
		% fig_entire_33.eps
		% fig_entire_32.eps
		% fig_ex_gluing.eps
	% approx: (6)
		% fig_approx_imcf.eps
		% fig_approx_mixed.eps
		% fig_approx_many_valleys.eps
		% fig_approx_entire32.eps
		% fig_approx_deg.eps
		% fig_approx_nondeg.eps
	% heat: (10)
		% fig_heat_noncontinuous.eps
		% fig_heat_chi_sigma.eps
		% fig_heat_no_model_ic.eps
		% fig_heat_1b_nondeg.eps
		% fig_heat_curve_comp.eps
		% fig_heat_1b_nondeg.eps [appeared twice]
		% fig_heat_1b_res.eps
		% fig_heat_2a_info.eps
		% fig_heat_2a_jump.eps
		% fig_heat_curve_comp_2.eps
	% ereg: (7)
		% fig_ereg_intro1.eps
		% fig_ereg_mixed.eps
		% fig_ereg_flatridge.eps
		% fig_ereg_Ek.eps
		% fig_ereg_fine_struc.eps
		% fig_ereg_fine_struc2.eps
		% fig_ereg_s1s4.eps
	% lvlset: (1)
		% fig_lvlset_impossible.eps
	% closed: (0)
	% crit: (0)
	% entire: (0)
	% kzz_es: (0)

% Unused references:
	% Bhattacharya_2005a, Bojarski-Iwaniec_1984, Crandall-Evans_2001, Crandall-Gunnarsson-Wang_2007, DeBlassie-Smits_2021, Drucker-Williams_2009, Crandall-Evans_2001, Hong-Zhao_2018, Juutinen-Lindqvist-Manfredi_2001, Steinerberger_2020

For $p\in(1,\infty)$, we say that $u$ is a $p$-harmonic function if it is a distributional solution of the equation
\begin{equation}\label{eq-intro:p-Laplacian}
	\Delta_pu:=\div\big(|\D u|^{p-2}\D u\big)=0,
\end{equation}
equivalently, if $u$ minimizes the $L^p$ Dirichlet energy $\|\D u\|_{L^p}$ among compactly supported perturbations. Let $\Omega$ be a simply-connected domain in $\RR^2$. Given a $p$-harmonic function $u$ in $\Omega$, one may construct a $q$-harmonic function $v$, where $\frac1p+\frac1q=1$, such that
\begin{equation}\label{eq-intro:pq_duality}
	\D v=|\D u|^{p-2}\D^\perp u,\qquad\text{or}\qquad\D u=-|\D v|^{q-2}\D^\perp v.
\end{equation}
In this paper, the symbol ``$\perp$'' means rotation of vectors by $90^\circ$ counterclockwise.

\vspace{3pt}

As $p\to\infty$, $p$-harmonic functions are known to converge to $\infty$-harmonic functions. The study of this class of objects was pioneered by Aronsson \cite{Aronsson_1967, Aronsson_1968} in the 1960s. It can be understood from the following perspectives. Expanding \eqref{eq-intro:p-Laplacian}, we find
\[\Delta_p u=|\D u|^{p-2}\Delta u+(p-2)|\D u|^{p-4}\D^2 u(\D u,\D u),\]
and the dominant term is $\D^2 u(\D u,\D u)$ as $p\to\infty$. The operator $\D^2 u(\D u,\D u)$ is called the \textit{$\infty$-Laplacian} and is denoted by $\Delta_\infty$. Next, formally taking $p\to\infty$ of the $L^p$ Dirichlet energy, we obtain the $L^\infty$ Dirichlet energy $\|\D u\|_{L^\infty}$. In summary, we find the following candidate definitions of $\infty$-harmonicity:
\begin{enumerate}[label={(\roman*)}, nosep]
	\item A viscosity solution of $\D^2 u(\D u,\D u)=0$.
	\item A minimizer of $\|\D u\|_{L^\infty}$ among compactly supported perturbations.
	\item The $C^0_{\loc}$ limit of a sequence of $p$-harmonic functions, as $p\to\infty$.
\end{enumerate}
The following \textit{comparison with cones} property, based on maximum principle, is also useful:
\begin{enumerate}[label={(\roman*)}, nosep]
	\setcounter{enumi}{3}
	\item A function $u$ so that, for any cone function $c(x)=a+b|x-x_0|$ with $u\geq c$ or $u\leq c$ on $\p(V\setminus\{x_0\})$ for a precompact $V$, it holds $u\geq c$ or $u\leq c$ in $V$.
\end{enumerate}
The four formulations are all equivalent, as was shown in a series of works by Bhattacharya-DiBenedetto-Manfredi \cite{Bhattacharya-DiBenedetto-Manfredi_1989}, Jensen \cite{Jensen_1993}, and Crandall-Evans-Gariepy \cite{Crandall-Evans-Gariepy_2001}. Initially, Aronsson \cite{Aronsson_1967, Aronsson_1968} worked with definitions (i)(ii) and only considered $C^2$ solutions. The use of viscosity solutions in (i) was first introduced in \cite{Bhattacharya-DiBenedetto-Manfredi_1989}. As we shall see, $\infty$-harmonic functions are typically non-$C^2$. Comparison with cones played a central role in the more recent developments \cite{Crandall_2004, Evans-Savin_2008, Savin_2005}. Later, Peres-Schramm-Sheffield-Wilson \cite{Peres-Schramm-Sheffield-Wilson_2009} studied the $\infty$-Laplacian on metric spaces and discovered the connection to the random tug-of-war game, and Daskalopoulos-Uhlenbeck \cite{Daskalopoulos-Uhlenbeck_2022, Daskalopoulos-Uhlenbeck_2024a, Daskalopoulos-Uhlenbeck_2024b} studied $\infty$-harmonic maps in relation with Thurston's asymmetric metric on Teichm\"uller spaces.

\vspace{3pt}

The regularity of $\infty$-harmonic functions remains a central question. In $\RR^2$, the best known result is $C^{1,\alpha}$ for some $\alpha>0$, by Evans-Savin \cite{Evans-Savin_2008}; see also Savin \cite{Savin_2005}, Zhang-Zhou \cite{Zhang-Zhou_2020}. Koch-Zhang-Zhou \cite{Koch-Zhang-Zhou_2019} also showed that $|\D u|\in W^{1,2}_{\loc}$. Much more is known for $\infty$-harmonic potentials in convex rings \cite{Lindgren-Lindqvist_2019, Lindgren-Lindqvist_2021, Peng-Zhang-Zhou_2025}. In $\RR^n$ ($n\geq3$), the best known regularity is everywhere differentiability, by Evans-Smart \cite{Evans-Smart_2011a, Evans-Smart_2011b}. Many standard techniques fail since $\Delta_\infty$ is degenerate elliptic. The fine structure of $\infty$-harmonic functions remains largely unknown; it is one of the main objectives of this paper.

In $\RR^2$, the solution
\begin{equation}\label{eq-intro:Aronsson_example}
	u_0=|x|^{4/3}-|y|^{4/3}
\end{equation}
was discovered by Aronsson \cite{Aronsson_1984} in 1984, and was the first known nonsmooth $\infty$-harmonic function. For each $d\in\ZZ_{\geq2}$, there exists a solution in $\RR^2$ of the form
\begin{equation}
	u=r^{\frac{d^2}{2d-1}}\varphi(\th),
\end{equation}
which we call a \textit{degree $d$ quasiradial solution}. This includes \eqref{eq-intro:Aronsson_example} as the $d=2$ instance; see Example \ref{ex-ex:quasiradial}. These examples are of $C^{1,1/3}$ regularity and no better.

Since Aronsson's work, it remained open whether one should expect the sharp $C^{1,1/3}$ regularity in general. We answer this question affirmatively.

\begin{theorem}\label{thm-intro:reg}
	Suppose $u$ is $\infty$-harmonic in a domain $\Omega\subset\RR^2$. Then $u\in C^{1,1/3}_{\loc}(\Omega)$.
\end{theorem}

The main component of the proof is the following $\epsilon$-regularity theorem: when $u$ is almost linear, we prove a sharp $C^{1,1/3}$ estimate of the type
\begin{equation}\label{eq-intro:ereg}
	\|u-x\|_{C^0}\leq\delta^2\quad\Rightarrow\quad\big\|\log|\D u|\big\|_{C^{0,2/3}}\leq C\delta^{4/3},\quad\Big\|\frac{\D u}{|\D u|}-\p_x\Big\|_{C^{0,1/3}}\leq C\delta^{2/3}.
\end{equation}
See Theorems \ref{thm-ereg:main_infhar} and \ref{thm-ereg:es_combined} for the precise statements. It is noteworthy that $|\D u|$ gains more regularity than the vector part of $\D u$.

The consequences of the main theory include a structural result near critical points:

\begin{theorem}\label{thm-intro:iso_crit}
	Suppose $u$ is nonconstant and $\infty$-harmonic in a domain $\Omega\subset\RR^2$. Then $\Crit(u)$ is discrete in $\Omega$. For each $z_0\in\Crit(u)$, there exists $d\in\ZZ_{\geq2}$ so that the family of functions
	\[\lambda^{-\frac{d^2}{2d-1}}\big[u(\lambda z+z_0)-u(z_0)\big]\]
	converges in $C^1_{\loc}(\RR^2)$ to a degree $d$ quasiradial solution, as $\lambda\to0$.
\end{theorem}

We also study entire $\infty$-harmonic functions with polynomial growth:

\begin{theorem}\label{thm-intro:entire}
	Suppose $u$ is a nonlinear entire $\infty$-harmonic function in $\RR^2$, with
	\[|u|\leq C\big(1+|z|^n\big)\qquad\text{in}\ \ \RR^2,\]
	for some $C,n>0$. Then there exists $d\in\ZZ_{\geq2}$ so that the family of functions
	\[\lambda^{-\frac{d^2}{2d-1}}u(\lambda z)\]
	converges in $C^1_{\loc}(\RR^2)$ to a degree $d$ quasiradial solution, as $\lambda\to\infty$.
\end{theorem}

Let us call $d$ the \textit{degree} of a polynomial growth solution. We also investigate low-degree solutions: in Theorem \ref{thm-entire:deg_2}, we show that $u_0=|x|^{4/3}-|y|^{4/3}$ is the unique degree 2 solution, up to rigid motion and scaling. The family of degree 3 solutions is constructed in Subsection \ref{subsec:entire3}. Further structural properties of the solutions in Theorem \ref{thm-intro:iso_crit}, \ref{thm-intro:entire} are discussed later in this introduction. The $d=1$ case was classified by Savin \cite{Savin_2005} by showing that entire solutions with linear growth must be linear. See also \cite{Dong-Peng-Zhang-Zhou_2024} for an integral version of linear Liouville theorem.

\vspace{3pt}

Let us turn to the $q\to1$ side of the duality \eqref{eq-intro:pq_duality}. It is known that $q$-harmonic functions converge to 1-harmonic functions, or functions of least gradient \cite{Gorny-Mazon}, once suitably normalized (for instance, one may multiply $v$ by a constant so that $\int|\D v|=1$). A renormalization is necessary since $\D v=|\D u|^{p-2}\D^\perp u$ does not converge in general. Taking $p\to\infty$ and $q\to1$ simultaneously, this yields a duality between $\infty$-harmonic and 1-harmonic functions. This duality was studied by Daskalopoulos-Uhlenbeck \cite{Daskalopoulos-Uhlenbeck_2022}, see also \cite{Backus_2024, Daskalopoulos-Uhlenbeck_2024a, Han_2023, Moser_2026}: roughly stated, the result is that the maximal gradient set
\begin{equation}\label{eq-intro:max_grad_set}
	S=\Big\{x\in\Omega:|\D u(x)|=\sup_\Omega|\D u|\Big\}
\end{equation}
of the $\infty$-harmonic function $u$ is a geodesic lamination, and the dual 1-harmonic function $v\in\BV$ has $Dv$ supported on $S$. From an analytic point of view, this duality only reveals information about $u$ on the maximal gradient set. In the study of regularity, a finer duality is perhaps needed.

\vspace{3pt}

A key idea of R. Moser \cite{Moser_2022}, built on the earlier work \cite{Moser_2007}, is to consider the transformation
\begin{equation}\label{eq-intro:moser_renormalization}
	w=(1-q)\log v,
\end{equation}
and the $q\to1$ limit of $w$. Using $\Delta_qv=0$, one may verify
\begin{equation}\label{eq-intro:p_IMCF_eq}
	\Delta_q w=|\D w|^q.
\end{equation}
Hence, the limit should formally solve the equation
\begin{equation}\label{eq-intro:imcf_eq}
	\div\Big(\frac{\D w}{|\D w|}\Big)=\Delta_1w=|\D w|.
\end{equation}

The equation \eqref{eq-intro:imcf_eq} is known as the level set form of the \textit{inverse mean curvature flow} (IMCF). If $w$ is a smooth solution of \eqref{eq-intro:imcf_eq}, then the level sets $\gamma_t=\{w=t\}$ evolve with normal speed equal to $1/(\text{curvature of $\gamma_t$})$ everywhere.

\begin{figure}[ht]
	\centering
	\captionsetup{width=0.85\textwidth}
	\includegraphics{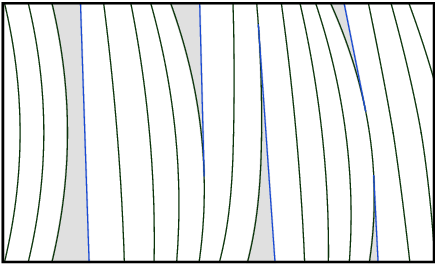}
	\caption{A right-evolving weak IMCF.}\label{fig-intro:imcf}
\end{figure}

The study of IMCF arose from its connections to scalar curvature via the monotonicity of Hawking mass. A milestone in the development was the introduction of the weak IMCF by Huisken-Ilmanen \cite{Huisken-Ilmanen_2001}. The weak IMCF is a variational form of \eqref{eq-intro:imcf_eq}. The function $w$ is no longer smooth. The evolving curves, now $\gamma_t=\p\{w<t\}$, may jump forward instantaneously. A jump corresponds to a set $\Int\big(\{w=t\}\big)$, see the grey regions in Figure \ref{fig-intro:imcf}. The portion of the curve after each jump, namely the set $\p\{w\leq t\}\setminus\p\{w<t\}$, is always a line segment (blue segments in Figure \ref{fig-intro:imcf}; see also Remark \ref{rmk-prelim:imcf_properties}\ref{item-prelim:imcf_jump_minsurf}). We refer to Subsection \ref{subsec:imcf} or \cite{Huisken-Ilmanen_2001, Xu_thesis} for more discussions.

The renormalization \eqref{eq-intro:moser_renormalization} was first introduced by R. Moser in \cite{Moser_2007}, where he showed the following: whenever a sequence of $q$-harmonic functions $v_q$ satisfies
\begin{equation}\label{eq-intro:moser_conv_to_imcf}
	\lim_{q\searrow1}(1-q)\log v_q=w\qquad\text{in}\ \ C^0_{\loc}(\Omega),
\end{equation}
then $w$ is a weak IMCF in $\Omega$. This major connection was further developed by subsequent works including \cite{Agostiniani-Borghini-Mazzieri_2026, AMMO_2025, Benatti-Fogagnolo-Mazzieri_2025, BMRSX_2025, Benatti-Pluda-Pozzetta_2024, Kotschwar-Ni_2009, Mari-Rigoli-Setti_2022}. As a model example, the function $v_q=|x|^{-\frac{n-q}{q-1}}$ is $q$-harmonic in $\RR^n\setminus\{0\}$, and its renormalized limit \eqref{eq-intro:moser_conv_to_imcf} in $C^0_{\loc}(\RR^n)$ is the spherically symmetric IMCF $w=(n-1)\log|x|$.

\vspace{3pt}

To summarize: given an $\infty$-harmonic function $u$, we have considered an approximating sequence of $p$-harmonic functions $u_p\to u$ as $p\to\infty$, and then the procedure
\begin{equation}\label{eq-intro:approx_naive}
	\begin{aligned}
		& \Delta_p u_p=0\quad\to\quad \D v_p=|\D u_p|^{p-2}\D^\perp u_p \\
		&\hspace{144pt} \quad\to\quad w_p=(1-q)\log v_p\quad\to\quad w=\lim_{q\searrow1}w_p.
	\end{aligned}
\end{equation}

In \cite{Moser_2022}, R. Moser applied \eqref{eq-intro:approx_naive} to a subclass of $\infty$-harmonic functions called oriented solutions. An oriented solution $u$ has no critical points, and $|\D u|$ is monotone in a uniform direction along all level sets of $u$. Applying \eqref{eq-intro:approx_naive} locally to an oriented solution, it is shown that the limit function is $w=-\log|\D u|$. Thus, \cite{Moser_2022} shows that:
\[\text{if $u$ is an oriented $\infty$-harmonic function, then $-\log|\D u|$ is a weak IMCF.}\]
The arguments of \cite{Moser_2022} inspired us in particular in the proof of Theorem \ref{thm-approx:global}. The reason that $w=-\log|\D u|$ for oriented solutions lies at the core of the proof; we refer to the detailed discussion in Section \ref{sec:approx}.

In fact, more is known from \cite{Moser_2022}: denoting $\nu=\pm\D^\perp u/|\D u|$ (where the sign depends on the orientation), it is shown that $w=-\log|\D u|$ satisfies
\begin{equation}\label{eq-intro:calibrated_imcf}
	|\nu|\leq1,\qquad \metric{\nu}{\D w}=|\D w|\ \ \text{a.e.},\qquad\div(\nu)=|\D w|\ \ \text{distributionally.}
\end{equation}
A function $w$ satisfying \eqref{eq-intro:calibrated_imcf} for some measurable vector field $\nu$ is called a \textit{calibrated IMCF}, and $\nu$ is called a \textit{calibration} of $w$. A calibrated IMCF is always a weak IMCF. The condition $|\nu|\leq1$ is trivial in the context of this paragraph, but is substantial in general.

Condition \eqref{eq-intro:calibrated_imcf} is another weak form of the IMCF equation \eqref{eq-intro:imcf_eq}: indeed, the calibration $\nu$ is a weak replacement of $\D w/|\D w|$ in \eqref{eq-intro:imcf_eq}. In calibrated IMCFs, $\nu$ is shown to be the outer unit normal of the level sets $\gamma_t=\p\{w<t\}$. Thus, each $\gamma_t$ is ``calibrated'' by $\nu$ similarly to the traditional sense. It is a general principle that calibrated submanifolds are area-minimizing. This explains many of the area-minimizing properties in IMCF: for example, in a calibrated IMCF, each sub-level set $\{w<t\}$ minimizes the perimeter-like energy
\[J_u(E)=|\p E|-\int_E|\D u|\]
among compact perturbations (this is in fact a defining property of weak IMCFs). See Subsection \ref{subsec:imcf} for more discussions. In this paper, we mostly work with calibrated IMCFs.

\vspace{3pt}

\begin{figure}[ht]
	\centering
	\includegraphics{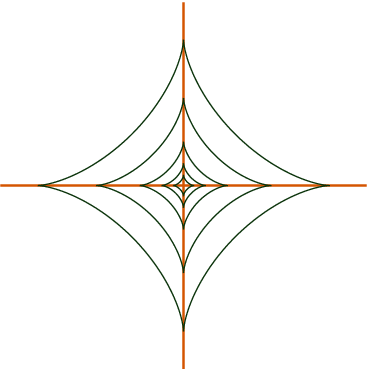}
	\caption{Level sets of $|\D u_0|$, where $u_0=|x|^{4/3}-|y|^{4/3}$.}\label{fig-intro:deg2}
\end{figure}

We turn our attention back to oriented solutions. The orientedness condition is non-generic. The solution $u_0=|x|^{4/3}-|y|^{4/3}$ is not oriented: it is oriented when restricted to each quadrant, but not across. Therefore, $-\log|\D u_0|$ is a weak IMCF (in fact, a smooth IMCF in this case) in each quadrant. It is natural to investigate the behavior on the $x$- and $y$-axes. Figure \ref{fig-intro:deg2} shows the shape of the level sets
\[\gamma_t=\big\{\!-\log|\D u_0|=t\big\}=\big\{|\D u_0|=e^{-t}\big\},\]
where we see cusps forming on the interfaces. More precisely: $\{\gamma_t\}$ is a self-shrinking solution of the IMCF consisting of 4 identical copies of hypocycloids.

As a useful note for future discussions: the edges of $\gamma_t$ are \textit{streamlines} of $u_0$ (meaning integral lines of the gradient flow) in addition to being level sets of $|\D u_0|$. Indeed, as $u_0$ is smooth in each quadrant, we may write
\[0=\Delta_\infty u_0=\D^2 u_0(\D u_0,\D u_0)=\frac12\bmetric{\D u_0}{\D|\D u_0|^2},\]
hence the level sets of $|\D u_0|$ are streamlines of $u_0$.

The structure of cusps and streamlines are general: for an $\infty$-harmonic function $u$, each set $\p\{|\D u|<e^{-t}\}$ is a union of concave curves with cusps at their endpoints, and all edges are streamlines of $u$ (Theorem \ref{thm-lvlset:global}).

\vspace{3pt}

The boundary condition on the interfaces can be understood in two natural ways: either as an IMCF of curves with cusps, or as a continuously glued piecewise IMCF where the level sets are tangent to the interfaces. The latter tangential condition is known as the \textit{outer obstacle} condition, initially introduced in \cite{Xu_2024_obstacle}. For calibrated IMCFs, i.e. solutions of \eqref{eq-intro:calibrated_imcf}, the outer obstacle condition for a region $D$ naturally states
\begin{equation}\label{eq-intro:bd_tangency}
	\nu=\nu_D\qquad\text{on}\ \ \p D
\end{equation}
(in a suitable sense), where $\nu_D$ denotes the outer unit normal of $D$. We refer the reader to \cite{Xu_2024_obstacle} for more discussions on the weak outer obstacle condition. This paper mainly involves calibrated solutions, and the main existence result in \cite{Xu_2024_obstacle} is not used.

\vspace{3pt}

The notion of \textit{simple IMCF cluster} refers to the structures of calibrated piecewise IMCF with outer obstacles. A simple IMCF cluster consists of a continuous function $w$ and vector field $\nu$, and a partition of $\Omega$ into $C^1$ regions $D_i$, so that each $w|_{D_i}$ is a weak IMCF with outer obstacle $\p D_i$ and is calibrated by $\pm\nu$ (Definition \ref{def-cluster:simple}). The sign of $\nu$ in the calibration corresponds to the orientations of $D_i$. The interfaces, i.e. the common boundaries of $D_i$, are called \textit{ridges}.

\begin{figure}[ht]
	\centering
	\includegraphics{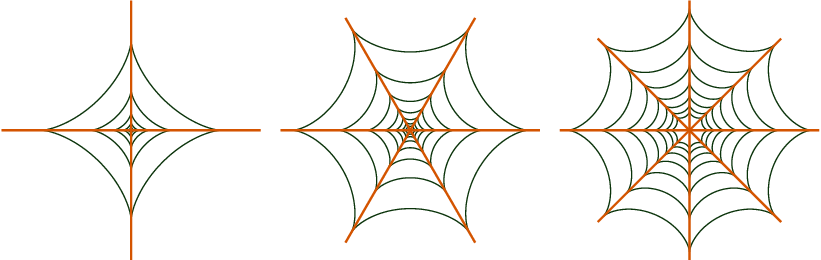}
	\caption{Quasiradial solutions with degree $2,3,4$.}\label{fig-intro:quasiradial}
\end{figure}

For $u_0=|x|^{4/3}-|y|^{4/3}$, the function $w=-\log|\D u_0|$ is a simple IMCF cluster in $\RR^2\setminus\{0\}$, with the four quadrants being the natural partition and the $x$, $y$-axes being the ridges. For a degree $d$ quasiradial solution, $w=-\log|\D u|$ is a simple IMCF cluster in $\RR^2\setminus\{0\}$ under the division into $2d$ sectors with angles $\pi/d$ (Figure \ref{fig-intro:quasiradial}).

\vspace{3pt}

It remains to be seen how the structure of ridges appears in the procedure \eqref{eq-intro:approx_naive}. In Moser's transformation \eqref{eq-intro:moser_renormalization}, the logarithm requires the positivity of $v$. However, as $v$ comes from integrating $|\D u|^{p-2}\D^\perp u$, it need not have a definite sign. The ridges in fact arise from the sign changes. We observe the general philosophy
\begin{equation}\label{eq-intro:limit_obstacle}
	\begin{aligned}
		& \text{The $q\to1$ limit of the Dirichlet condition for $q$-harmonic functions} \\
		&\hspace{120pt}\text{is the outer obstacle condition for the weak IMCF.}
	\end{aligned}
\end{equation}
This idea first appeared in \cite{Xu_2024_obstacle} and was formally addressed in \cite[Theorem 1.10]{BMRSX_2025}. It is shown in \cite{BMRSX_2025} that for smooth domains $E_0\Subset\Omega\Subset\RR^n$, the $q$-capacitary potentials $v$ with
\[\left\{\begin{aligned}
	& \Delta_qv=0\qquad\text{in}\ \ \Omega\setminus\bar{E_0} \\
	& v|_{\p E_0}=1 \\
	& v|_{\p\Omega}=0
\end{aligned}\right.\]
have $(1-q)\log v$ converging to a weak IMCF with outer obstacle $\p\Omega$, as $q\to1$.

In our present context, we have $q$-harmonic functions $v_p$ obtained from integrating $|\D u_p|^{p-2}\D^\perp u_p$. Then naturally, $|v_p|$ is a positive $q$-harmonic function in both $\{v_p>0\}$ and $\{v_p<0\}$, with common Dirichlet conditions on $\{v_p=0\}$ (which moves as $p$ changes). We consider the limiting procedures:
\begin{enumerate}[label={(\arabic*)}, nosep]
	\item The set $\{v_p=0\}$ subsequentially converges to a collection of $C^1$ curves $\{\gamma_j\}$, as $p\to\infty$ and $q\to1$.
	\item Take $\{D_i\}$ to be the connected components of $\Omega\setminus\cup\gamma_j$.
	\item Take a subsequential limit $w_i=\lim_{q\to1}(1-q)\log|v_p|$ in $C^0_{\loc}(D_i)$, for each $i$.
\end{enumerate}
Then, we show that $w_i$ is a weak IMCF in each $D_i$, and satisfies the outer obstacle condition, in accordance with \eqref{eq-intro:limit_obstacle}. This finally yields a simple IMCF cluster. The limit curves $\{\gamma_j\}$ become the ridges in the resulting cluster.

To summarize, we now consider the process
\begin{equation}\label{eq-intro:approx_main}
	\begin{aligned}
		\Delta_\infty u=0\quad\xleftarrow[\,p\to\infty]{C^1_{\loc}}\quad\Delta_p u_p=0\quad&\to\quad \D v_p=|\D u_p|^{p-2}\D^\perp u_p \\
		\quad&\to\quad w_p=(1-q)\log|v_p| \\
		\quad&\to\quad w=\lim_{q\searrow1}w_p\text{ is a simple IMCF cluster}.
	\end{aligned}
\end{equation}
We refer to Section \ref{sec:approx} for the full statement, proofs, and further discussion about this process. Notice that the convergence $u_p\to u$ holds in $C^1_{\loc}$: this follows from Theorem \ref{thm-equi:C1_equicont} and is new as well. The convergence of gradients is important to us.

\vspace{3pt}

The $C^{1,1/3}$ regularity in Theorem \ref{thm-intro:reg} relies on the optimal regularity of simple IMCF clusters, in particular, the regularity near ridges. A central tool therein is the use of \textit{support functions}. The support function of a $C^1$ curve $\gamma$ is
\begin{equation}\label{eq-intro:spt_function}
	h(\gamma(s))=\metric{\gamma(s)}{\nu(s)},
\end{equation}
where $\nu(s)=-\gamma'(s)^\perp$ is the orientation-compatible unit normal vector along $\gamma$.

For a smooth solution of IMCF by strictly convex curves $\gamma_t$, we may parametrize all curves by the angle of the tangent vectors. Then, the support function satisfies the evolution equation
\begin{equation}\label{eq-intro:heat}
	\frac{\p h}{\p t}=\frac{\p^2 h}{\p\th^2}+h.
\end{equation}
Indeed, $\p_th$ is the normal speed, and $\p_{\th\th}h+h$ is the inverse of the curvature (which is a crucial observation). This equation was observed by Chow-Tsai \cite{Chow-Tsai_1996} and Urbas \cite{Urbas_1999}.

We prove that \eqref{eq-intro:heat} holds in a suitable sense for weak IMCFs and for simple IMCF clusters as well. This is carried out in Section \ref{sec:heat}. Note that, on a union of concave curves with cusps at endpoints, there is a continuous angle parameter through all cusps.
\begin{figure}[ht]
	\centering
	\includegraphics{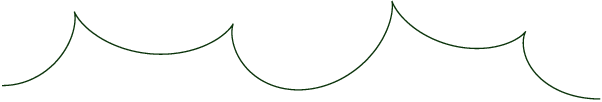}
	\begin{picture}(0,0)
		% Midpoints
		\put(-273.4,15){\arrowangle{43}}
		\put(-228,21.6){\arrowangle{175}}
		\put(-142.2,4.8){\arrowangle{16}}
		\put(-82,25.8){\arrowangle{167}}
		\put(-32,2.7){\arrowangle{-33}}
		
		% arrows at the four cusps
		\put(-257.2,28){\arrowangle{107}}
		\put(-190.5,27){\arrowangle{-114}}
		\put(-104.4,33){\arrowangle{104}}
		\put(-51.0,24.0){\arrowangle{-119}}
	\end{picture}
	\refstepcounter{figure}\label{fig-intro:angle_param}
\end{figure}
Under this parametrization, the edges alternate between being concave and convex. In terms of support functions, the inverse curvature $\p_{\th\th}h+h$ changes sign, and the zeros of $\p_{\th\th}h+h$ are the angles at the cusps. A main component of Theorem \ref{thm-heat:main} is that \eqref{eq-intro:heat} holds across the cusps.

For instance, in a degree $d$ quasiradial solution, the level sets $\gamma_t=\{|\D u|=e^{-t}\}$ have support functions given by
\begin{equation}\label{eq-intro:h_quasiradial}
	h(t,\th)=\exp\Big(\frac{1-2d}{(d-1)^2}t\Big)\sin\frac{-d\th}{d-1},\qquad\th\in\SS^1\big(2(d-1)\pi\big),
\end{equation}
where $\SS^1(l)$ denotes the circle with length $l$. Observe that $\p_{\th\th}h+h$ changes sign $2d$ times, corresponding to the $2d$ cusps in $\gamma_t$. See Example \ref{ex-ex:quasiradial} for details. It is notable that \eqref{eq-intro:h_quasiradial} is a parabolic scaling of the degree $d$ Fourier mode.

\vspace{3pt}

If $h$ solves \eqref{eq-intro:heat}, then the signed inverse curvature $\kappa^{-1}=\p_{\th\th}h+h$ also solves the same equation $\p_t\kappa^{-1}=\p_{\th\th}\kappa^{-1}+\kappa^{-1}$. At a point $(t_0,\th_0)$ where $\kappa^{-1}$ changes sign (corresponding to a cusp), one has the asymptotics
\begin{equation}\label{eq-intro:hopf}
	\kappa^{-1}\approx c(\th-\th_0)
\end{equation}
due to the nondegeneracy of roots. Then, switching to a length parameter $s$, we have
\begin{equation}\label{eq-intro:hopf2}
	\frac{d\th}{ds}=\kappa\approx\frac1{c(\th-\th_0)}\qquad\Leftrightarrow\qquad \frac{d\big[(\th-\th_0)^2\big]}{ds}\approx\frac2c.
\end{equation}
This shows that the tangential angle changes in a $C^{0,1/2}$ manner in length parameter, and hence, the edges of $\gamma_t$ have $C^{1,1/2}$ regularity up to the cusps. This regularity is sharp, since in quasiradial solutions, the edges of $\gamma_t$ are hypocycloids, hence are all $C^{1,1/2}$.

The following principle applies to the edges of $\gamma_t$, after a suitable formulation:
\begin{equation}\label{eq-intro:reg_joint_normal}
	\begin{aligned}
		&\text{If a disjoint family of curves has uniform $C^{1,\alpha}$ regularity,} \\
		&\hspace{64pt}\text{then the collection of all their normal vectors is $C^{0,\alpha/(1+\alpha)}$.}
	\end{aligned}
\end{equation}
See Lemma \ref{lemma-prelim:joint_cont_of_normals}, \ref{lemma-prelim:two_graphs}. The H\"older exponent $\frac{\alpha}{1+\alpha}$ is sharp as well. Finally, since
\[\frac{1/2}{1+1/2}=\frac13,\]
the $C^{1,1/2}$ regularity of $\gamma_t$ turns into the $C^{0,1/3}$ regularity of the joint normal vector field. As $\gamma_t$ are streamlines of $u$, this explains the source of $C^{1,1/3}$ regularity in Theorem \ref{thm-intro:reg}.

These ideas are used in Section \ref{sec:ereg} to prove the $\epsilon$-regularity statement \eqref{eq-intro:ereg}. We mention that, hidden in the technicality, the Sturmian principle is our main tool from parabolic theory. The Sturmian principle states that the number of roots of a solution of the heat equation is non-increasing in time; see for example \cite{Angenent_1988, Angenent-Fiedler_1988}. This is useful in showing a uniform version of \eqref{eq-intro:hopf}. See Section \ref{sec:ereg} for further discussion on the techniques.

\vspace{3pt}

Conversely, the equation \eqref{eq-intro:heat} can be used to construct simple IMCF clusters. Given an abstract function $h(t,\th)$ solving $\p_th=\p_{\th\th}h+h$, one may reconstruct a family of curves $\gamma_t$ (with cusps) that solve the IMCF, so that each $\gamma_t$ is parametrized by the angle $\th$ and has support function $h(t,\cdot)$. The exact formula is
\begin{equation}\label{eq-intro:recovery_formula}
	\gamma_t(\th)=h(t,\th)\nu_\th+\p_\th h(t,\th)\tau_\th,
\end{equation}
where in this paper we denote
\begin{equation}\label{eq-intro:nu_tau_th}
	\nu_\th=(\sin\th,-\cos\th)=e^{i(\th-\pi/2)},\qquad\tau_\th=(\cos\th,\sin\th)=e^{i\th}.
\end{equation}
Once $\{\gamma_t\}$ is shown to be an embedded family, this yields a simple IMCF cluster, with the ridges being the trajectories of cusps. From such a simple cluster, using Theorem \ref{thm-ex:reconstruction}, one may then construct an $\infty$-harmonic function $u$ so that
\begin{equation}\label{eq-intro:reconstruction}
	\D u(\gamma_t(\th))=e^{-t}\tau_\th.
\end{equation}
The formula \eqref{eq-intro:reconstruction} is chosen so that $\gamma_t=\{-\log|\D u|=t\}$ and $\gamma_t$ are streamlines of $u$ (since $\tau_\th$ is tangent to $\gamma_t$ at $\gamma_t(\th)$). This is consistent with the previous discussions: recall that for an $\infty$-harmonic function $u$, all edges of $\p\{|\D u|<e^{-t}\}$ are streamlines.

\vspace{3pt}

In Subsection \ref{subsec:entire3}, we construct $\infty$-harmonic functions that are asymptotic to degree 3 quasiradial solutions at infinity. To our knowledge, this is the first construction of non-quasiradial degree 3 solutions. The idea is to perturb the quasiradial support function \eqref{eq-intro:h_quasiradial} (with $d=3$) by adding lower order terms, namely, to consider
\begin{equation}\label{eq-intro:h_degree_3_pre}
	h(t,\th)=\sum_{j=0}^3\exp\Big(t-\frac{j^2}4t\Big)\Big(a_j\sin\frac{j\th}2+b_j\cos\frac{j\th}2\Big),\qquad\th\in\SS^1(4\pi).
\end{equation}
Then, we apply the above procedure \eqref{eq-intro:recovery_formula}\,--\,\eqref{eq-intro:reconstruction} to this $h$.

By an alternating length identity (see Subsection \ref{subsec:entire3} for details), the constant term must vanish for the construction to succeed. By translation symmetry, we may drop the $j=2$ term. Eventually, it suffices to consider
\begin{equation}\label{eq-intro:h_degree_3}
	h(t,\th)=e^{-5t/4}\Big(a\sin\frac32\th+b\cos\frac32\th\Big)+e^{3t/4}\Big(c\sin\frac\th2+d\cos\frac\th2\Big),\quad\th\in\SS^1(4\pi).
\end{equation}
Modulo rigid motions and scalings (thus reducing to the case $a=1, b=0$), we obtain a 2-dimensional family of solutions parametrized by $c,d$.

\begin{figure}[ht]
	\centering
	\includegraphics{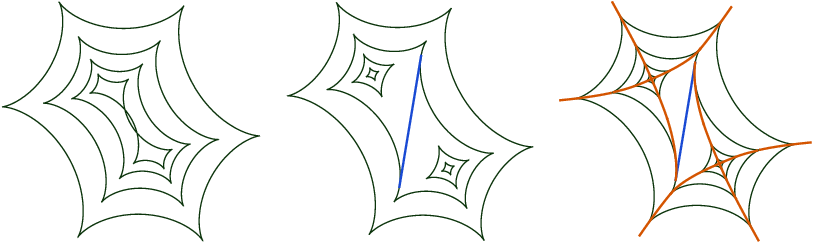}
	\caption{Splitting of polygon and mixed IMCF cluster.}\label{fig-intro:splitting}
\end{figure}

The major obstacle in the construction is to ensure the embeddedness of $\{\gamma_t\}$. In fact, the $\gamma_t$ arising from \eqref{eq-intro:recovery_formula} \eqref{eq-intro:h_degree_3} does not stay embedded for all $t\in\RR$ (Figure \ref{fig-intro:splitting}, left), except for the exact quasiradial case $c=d=0$. A new ingredient -- \textit{polygon splitting} -- is needed for the construction.

In reality, we first build $\gamma_t$ from \eqref{eq-intro:recovery_formula} \eqref{eq-intro:h_degree_3} for all sufficiently negative $t$. Denote by $Z_t$ the interior of $\gamma_t$ (it is combinatorially a 6-gon). We find a suitable time $T$ and use a suitable segment $\sigma$ to cut $Z_T$ into two 4-gons. Then, we evolve the two 4-gons individually by IMCF with cusps (Figure \ref{fig-intro:splitting}, middle). Finally, we construct the $\infty$-harmonic function $u$ by asking that all $\gamma_t$ are level sets of $|\D u|$ as well as streamlines of $u$, and additionally $-\log|\D u|\equiv T$ on $\sigma$.

We call the segment $\sigma$ a \textit{valley}. The solution locally forms two IMCFs evolving away from $\sigma$ (as opposed to moving toward a ridge). An IMCF cluster with both ridges and valleys is an instance of a \textit{mixed IMCF cluster} (Definition \ref{def-cluster:mixed}). The mixed cluster associated to our construction is shown in Figure \ref{fig-intro:splitting} (right); see also Figures \ref{fig-ex:entire_31}, \ref{fig-ex:entire_33}, \ref{fig-ex:entire_32} for other degree 3 solutions.

\begin{remark}
	In all IMCF cluster diagrams in this paper, orange curves are ridges and blue segments are valleys. Curves without colors are $\gamma_t$. All curves are streamlines of $u$, and all valleys and $\gamma_t$ are level sets of $|\D u|$ (ridges need not be level sets).
\end{remark}

The role of valleys is perhaps more apparent from an analytic point of view. Recalling $|\D u|=e^{-t}$ on $\gamma_t$, and comparing the direction of evolution, we find that $|\D u|$ becomes smaller when moving away from a valley. Thus, a valley is the maximal gradient set in its neighborhood. Also, it is mentioned below \eqref{eq-intro:max_grad_set} that the maximal gradient set of an $\infty$-harmonic function is a disjoint union of line segments. This explains that valleys must be line segments.

A mixed IMCF cluster also allows ``poles'' where the solution is quasiradial-like (e.g., there are two poles in Figure \ref{fig-intro:splitting}, right). The poles correspond to critical points in the resulting $\infty$-harmonic functions. See Sections \ref{sec:cluster} and \ref{sec:ex} for details.

\vspace{3pt}

In Section \ref{sec:closed}, we prove a converse of the above construction. Suppose $u$ is $\infty$-harmonic in a domain $\Omega$, and $F$ is a connected component of $\{|\D u|<e^{-T}\}$ with $F\Subset\Omega$. Then the curves $\gamma_t=\{|\D u|=e^{-t}\}\cap F$ follow the IMCF-and-splitting pattern for $t>T$:
\begin{enumerate}[label={(\arabic*)}, nosep]
	\item For all but finitely many $t$, the curves $\gamma_t$ evolve smoothly by IMCF with cusps.
	\item At each exceptional time, some connected components of $\gamma_t$ are split into polygonal curves with fewer edges. The cutting segments form the maximal gradient set in the interior of $\gamma_t$.
\end{enumerate}
Each isolated critical point is contained in a precompact connected component of $\big\{|\D u|<e^{-T}\big\}$, for some $T\gg1$. Thus, the above pattern applies. By the finiteness of exceptional times, there is no splitting beyond time $T'$ for some $T'>T$. Then, by investigating the $t\to\infty$ limit of the support function (which solves $\p_th=\p_{\th\th}h+h$ for $t\in(T',\infty)$), we eventually prove the asymptotics in Theorem \ref{thm-intro:iso_crit}.

\vspace{3pt}

In Section \ref{sec:crit}, we prove the discreteness of critical set in Theorem \ref{thm-intro:iso_crit}. The $\infty$-Laplacian does not have the unique continuation property: there are distinct (nonconstant) $\infty$-harmonic functions that coincide in an open subset \cite[p.21]{Lindqvist_note}. However, by Theorem \ref{thm-crit:main}, a nonconstant $\infty$-harmonic function cannot be constant in an open subset. This settles a question asked by Lindqvist \cite[p.22]{Lindqvist_note}.

\vspace{3pt}

In Section \ref{sec:entire}, we prove Theorem \ref{thm-intro:entire}. The key is to show that for an entire solution $u$ with polynomial growth, the curves $\gamma_t=\{|\D u|=e^{-t}\}$ follow the IMCF-and-splitting pattern for all $t\in\RR$, with finitely many splitting times in total. Thus, there is $\uT\in\RR$ so that $\gamma_t$ do not split, thus evolve as smooth IMCF with cusps, for $t\in(-\infty,\uT)$. Correspondingly, the support functions of $\gamma_t$ solve $\p_th=\p_{\th\th}h+h$ for $t\in(-\infty,\uT)$. An analytic lemma (Theorem \ref{thm-entire:ancient_sol}) is invoked to show that $h$ has finitely many Fourier frequencies. More specifically, it must take the form
\begin{equation}\label{eq-intro:entire_h_fourier}
	h(t,\th)=\sum_{j=1}^d\exp\Big(t-\frac{j^2t}{(d-1)^2}\Big)\Big(a_j\cos\frac{j\th}{d-1}+b_j\sin\frac{j\th}{d-1}\Big),
\end{equation}
for some $d\in\ZZ_{\geq2}$, where $\th\in\SS^1\big(2(d-1)\pi\big)$. The constant term and terms with frequency above $d$ must all vanish. The dominant term has $j=d$ and again corresponds to degree $d$ quasiradial solutions (see \eqref{eq-intro:h_quasiradial}). Then, Theorem \ref{thm-intro:entire} follows by taking the rescaled limit as $t\to-\infty$, with $d$ being the degree of $u$.

By Jensen's maximum principle, two $\infty$-harmonic functions in $\RR^2$ must coincide if they coincide in $\RR^2\setminus B(R)$ for some $R$. Hence, each choice of parameters $\{a_j,b_j\}$ in \eqref{eq-intro:entire_h_fourier} gives rise to a unique solution (if it exists). The moduli space of degree $d$ entire solutions is thus a subspace of $\RR^{2d-2}$ (where $a_j,b_j$ are the free parameters; by translation symmetry, the $j=d-1$ term can always be removed).

It is natural to ask whether any parameter yields a solution. The main difficulty is the determination of the splitting times and angles. The degree 3 case involves only one splitting with a unique combinatorial type (one 6-gon into two 4-gons), and existing necessary conditions already imply a unique choice (see Subsection \ref{subsec:entire3}). When $d\geq4$, the existence of combinatorially different splitting types poses a complication. How to determine the canonical splitting among the many different ones, given that the resulting $\infty$-harmonic function is unique?

\subsection{Structure of this paper; remarks and discussions}

The main content of each section has been outlined above and is also suggested by its title. Each section is organized around a main theorem, with the relevant main ideas and discussions included at the front. The proofs in different sections are mostly independent. The dependencies among sections are summarized below. In the diagram, ``[1.x]'' means Theorem 1.x.

\begin{figure}[ht]
	\centering
	\begin{tikzpicture}[
		sec/.style={font=\small},     % <-- no circles
		dep/.style={->},
		node distance=4mm and 4mm
		]
		
		% --- Nodes ---
		\node[sec] (prelim) {\ref{sec:prelim}};
		\node[sec, right=of prelim] (cluster) {\ref{sec:cluster}};
		\node[sec, below=of cluster] (equi) {\ref{sec:kzz_es},\,\ref{sec:equi}};
		\node[sec, right=of cluster] (approx) {\ref{sec:approx}};
		\node[sec, above=of approx] (ex) {\ref{sec:ex}};
		\node[sec, right=of approx] (heat) {\ref{sec:heat}};
		\node[sec, right=of heat] (lvlset) {\ref{sec:lvlset}};
		\node[sec, below=of lvlset] (ereg) {\ref{sec:ereg}};
		\node[sec, right=of lvlset] (closed) {\ref{sec:closed}};
		\node[sec, right=of closed] (crit) {\ref{sec:crit}};
		\node[sec, right=of crit] (entire) {\ref{sec:entire}};
		
		\node[sec, below=of crit] (th-reg) {[\ref{thm-intro:reg}]};
		\node[sec, above=of closed] (th-crit) {[\ref{thm-intro:iso_crit}]};
		\node[sec, above=of entire] (th-entire) {[\ref{thm-intro:entire}]};
		
		% --- Arrows ---
		\draw[dep] (prelim) -- (cluster);
		\draw[dep] (cluster) -- (approx);
		\draw[dep] (cluster) -- (ex);
		\draw[dep] (equi) -- (approx);
		\draw[dep] (approx) -- (heat);
		\draw[dep] (heat) -- (ereg);
		\draw[dep] (heat) -- (lvlset);
		\draw[dep] (lvlset) -- (closed);
		\draw[dep] (closed) -- (crit);
		\draw[dep] (crit) -- (entire);
		
		\draw[dep] (ereg) -- (th-reg);
		\draw[dep] (closed) -- (th-reg);
		\draw[dep] (crit) -- (th-reg);
		\draw[dep] (closed) -- (th-crit);
		\draw[dep] (crit) -- (th-crit);
		\draw[dep] (entire) -- (th-entire);
	\end{tikzpicture}
\end{figure}

A quick route towards understanding the main theorems might be to first read the introductions (while temporarily skip the proofs) of Sections \ref{sec:cluster}, \ref{sec:approx}, \ref{sec:ereg}, then to have a glance of the examples in Section \ref{sec:ex}, then to proceed with Sections \ref{sec:closed}, \ref{sec:crit}, \ref{sec:entire}. Sections \ref{sec:equi}, \ref{sec:heat}, \ref{sec:lvlset} may be treated as tools for the first read.

We end this introduction with a list of remarks of possible interest.

\begin{remark} {\ }
	
	\begin{enumerate}[label={(\roman*)}, topsep=1pt, itemsep=-0.5ex]
		\item It is natural to ask whether $\infty$-harmonic functions satisfy a $C^{1,1/3}$ estimate
		\begin{equation}
			\|u\|_{C^{1,1/3}(B(1))}\leq C\|u\|_{L^\infty(B(2))}.
		\end{equation}
		This seems to require new ingredients, for instance, a quantitative version of Theorem \ref{thm-crit:main} which relates the distribution of critical points and the growth of $u$.
		\item Several integral regularity conjectures were also posed for planar $\infty$-harmonic functions (for example, see \cite{Koch-Zhang-Zhou_2019}): one may ask if it always hold
		\begin{equation}\label{eq-intro:integral_reg_conj}
			u\in W^{2,3/2-\epsilon}_{\loc},\qquad|\D u|\in W^{1,3-\epsilon}_{\loc},\qquad\log|\D u|\in\operatorname{BMO}_{\loc}.
		\end{equation}
		These are sharply modelled on $u_0=|x|^{4/3}-|y|^{4/3}$ as well. We remark that the integral regularity for the vector part of simple IMCF clusters has not been established.
		\item According to Iwaniec-Manfredi \cite{Iwaniec-Manfredi_1989}, $p$-harmonic functions ($1<p<\infty$) in the plane are $C^{k,\alpha}$, where
		\[k+\alpha=\frac16\Big(7+\frac1{p-1}+\sqrt{1+\frac{14}{p-1}+\frac1{(p-1)^2}}\Big).\]
		Observe that
		\[\lim_{p\to\infty}(k+\alpha)=\frac43.\]
		Is it true that $p$-harmonic functions ($2\leq p\leq\infty$) are uniformly $C^{1,1/3}$? The proof in \cite{Iwaniec-Manfredi_1989} combines the $C^{k,\alpha}$ regularity at critical points and the smoothness away from critical points (due to strict ellipticity). The former may pass well to $p=\infty$, but the latter does not. The ridges in an $\infty$-harmonic function, where the optimality of $C^{1,1/3}$ also appears, are smoothed out in the approximating $p$-harmonic functions.
		\item The ridges are in natural correspondence with the ``attracting streamlines'' in the setting of Lindgren-Lindqvist \cite{Lindgren-Lindqvist_2021}. See precisely \cite[p.4]{Lindgren-Lindqvist_2021}.
		\item A very recent work of Brustad \cite{Brustad_2026} showed the following: if $u$ has no critical points and the map $z\mapsto\D u(z)$ has a Lipschitz inverse, then $u\in C^{1,1/3}(\Omega)$. We include a discussion in Remark \ref{rmk-approx:du_noninj}.
		\item Quasiradial solutions are the only known examples in $\RR^2$ with only one critical point, to the author's knowledge. Example \ref{ex-ex:entire_e1} is the only known solution in $\RR^2$ without critical points. Do there exist other examples?
	\end{enumerate}
\end{remark}

\textbf{Tool and computational resource disclosure.} The author thanks ChatGPT 5.5 Pro for suggesting improvements on the exposition and for finding mathematical and linguistic typos.

\vspace{3pt}

\textbf{Acknowledgements.} The author would like to thank Richard Bamler for many useful discussions, and Erik Lindgren for clarifications on the reference \cite{Lindgren-Lindqvist_2021}. The author also thanks Peking University and Yuguang Shi for their hospitality, as this work was partially completed during his visit.

%\newpage

\section{Preliminaries and notations}\label{sec:prelim}

The following symbols usually have a fixed meaning:

\begin{itemize}[topsep=0.5ex, itemsep=-0.6ex]
	\item $z,x,y$: a point in $\RR^2$ and its coordinates.
	\item $u$: an $\infty$-harmonic function.
	\item $\Omega$: a domain (i.e. connected open set) in $\RR^2$.
	\item $w$: an IMCF or IMCF cluster.
\end{itemize}

\noindent The following are our standard notations:

\begin{itemize}[topsep=0.5ex, itemsep=-0.6ex]
	\item $[z_1,z_2]$: the closed line segment connecting $z_1,z_2\in\RR^2$.
	\item $\SS^1(r)$: the circle with total length $r$; and $\SS^1=\SS^1(2\pi)$.
	\item $B(z,r)$, $Q(z,r)$: the open disk and square with center $z$ and radius $r$.
	\item $B(r)=B(0,r)$, and $Q(r)=Q(0,r)=(-r,r)\times(-r,r)$.
	\item $Q_\delta(r)=(-r,r)\times(-\delta r,\delta r)$.
	\item $N(A,r)$: the $r$-neighborhood of $A$, namely, $N(A,r)=\{z\in\RR^2:d(z,A)<r\}$.
	\item $\nu_D$: the outer unit normal of a region $D$, whenever $D$ has sufficient regularity.
	\item Streamline: a \textit{streamline} of a $C^1$ function $u$ is a gradient flow line: namely, it is a $C^1$ curve $\gamma$ satisfying $\gamma'(t)=\D u(\gamma(t))$.
	\item $\nu^\perp$ means $J\nu$, namely, the counterclockwise $90^\circ$ rotation of $\nu$.
	\item A degenerate root of a $C^1$ function $f(t,\th)$ is a point where $f(t,\th)=\p_\th f(t,\th)=0$.
\end{itemize}

\noindent We assume the following conventions:

\begin{itemize}[topsep=0.5ex, itemsep=-0.6ex]
	\item The curl of a vector field $X=(X_1,X_2)$ is defined as
	\[\curl(X)=\p_xX_2-\p_yX_1.\]
	\item A ``line segment'' means either a line segment or point or ray or line (this is particularly assumed in Section \ref{sec:heat}, \ref{sec:ereg}). A ``finite line segment'' means a line segment with finite length. A ``nontrivial line segment'' means a line segment that is not a point.
	\item A function $f$ is convex (or concave) if $f''\geq0$ (or $f''\leq0$) in the barrier sense.
	\item The closures and boundaries of all sets are relative to $\RR^2$. The relative boundary of a set $A$ in a domain $\Omega$ is written as $\p A\cap\Omega$.
	\item We say that $f_i\to f$ in $C^0_{\loc}(\Omega)$, if for all $K\Subset\Omega$ it holds
	\[\dom(f_i)\supset K\ \ \forall\,i\gg1,\qquad f_i\to f\ \ \text{in}\ \ C^0(K).\]
\end{itemize}

\subsection{\texorpdfstring{$\infty$}{∞}-Harmonic functions}

We define $\infty$-harmonic functions as viscosity solutions of $\D^2 u(\D u,\D u)=0$:

\begin{defn}\label{def-prelim:infty_har}
	Assume $\Omega$ is a domain in $\RR^n$, and $u\in C^0(\Omega)$, and $z\in\Omega$. We say that:
	\begin{itemize}[nosep]
		\item $\varphi\in C^2$ touches $u$ from above (resp. from below) at $z$, if $\varphi(z)=u(z)$ and $\varphi\geq u$ (resp. $\varphi\leq u$) in a small neighborhood of $z$;
		\item $u$ is $\infty$-subharmonic at $z$, if for all $\varphi\in C^2$ that touches $u$ from above at $z$, we have $\D^2\varphi(\D\varphi,\D\varphi)(z)\geq0$;
		\item $u$ is $\infty$-superharmonic at $z$, if for all $\varphi\in C^2$ that touches $u$ from below at $z$, we have $\D^2\varphi(\D\varphi,\D\varphi)(z)\leq0$;
		\item $u$ is $\infty$-harmonic at $z$, if it is both $\infty$-subharmonic and $\infty$-superharmonic at $z$;
		\item $u$ is $\infty$-harmonic $(\Delta_\infty u=0)$ in $\Omega$, if it is $\infty$-harmonic at every point in $\Omega$.
	\end{itemize}
\end{defn}

A main property is that it is the $C^0_{\loc}$ limit of $p$-harmonic functions as $p\to\infty$:

\begin{lemma}[{\cite{Bhattacharya-DiBenedetto-Manfredi_1989, Jensen_1993, Lindqvist-Manfredi_1995}}]\label{lemma-prelim:p-approx}
	Suppose $\Omega\subset\RR^n$ and $u\in C^0(\Omega)$. Then $u$ is $\infty$-harmonic in $\Omega$ if and only if $u$ is the $C^0_{\loc}(\Omega)$ limit of a sequence of $p_i$-harmonic functions, with $p_i\to\infty$.
\end{lemma}

Infinity-harmonicity is also equivalent to the \textit{AMLE} (absolutely minimizing Lipschitz extension) condition and the \textit{comparison with cones} condition. To state the former, denote
\[\Lip(u,K)=\sup_{x,y\in K}\frac{|u(x)-u(y)|}{|x-y|},\qquad F_\infty(u,K)=\esssup_K|\D u|.\]
The AMLE condition states that for all $V\Subset\Omega$ and all functions $v\in C^0(\bar V)$ with $u|_{\p V}=v|_{\p V}$, it holds $F_\infty(u,V)\leq F_\infty(v,V)$. This is (nontrivially) equivalent to $\Lip(u,\bar V)\leq\Lip(v,\bar V)$ for all $V\Subset\Omega$ and all $v\in C^0$ with $u|_{\p V}=v|_{\p V}$, \cite[Theorem 6.1]{Crandall_2004}. The AMLE property is not directly used in this paper. The comparison with cone property states that for all $V\Subset\Omega$ and cone functions of the form $c(x)=a+b|x-x_0|$, we have
\[u\leq c\text{ (resp. $u\geq c$) on $\p(V\setminus\{x_0\})$}\qquad\Rightarrow\qquad u\leq c\text{ (resp. $u\geq c$) in $V$.}\]
The equivalence of these properties and $\infty$-harmonicity were shown in \cite{Bhattacharya-DiBenedetto-Manfredi_1989, Crandall_2004, Crandall-Evans-Gariepy_2001, Jensen_1993}. For a detailed introduction of the fundamental theory, we refer to \cite{Araujo-Urbano_2023, Aronsson-Crandall-Juutinen_2001, Crandall_2004, Lindqvist_note, Wang}.

In addition to Lemma \ref{lemma-prelim:p-approx}, we will also need the following facts (most of which follow from comparison with cones). The gradient estimate is stated as follows:

\begin{lemma}[{\cite[Corollary 2.6]{Crandall-Evans-Gariepy_2001}}]\label{lemma-prelim:grad_est}
	If $\Delta_\infty u=0$ in $\Omega\subset\RR^n$, then $u\in\Lip_{\loc}(\Omega)$, and
	\[|\D u(x)|\leq\frac{2\sup_\Omega|u|}{d(x,\p\Omega)},\qquad\forall\,x\in\Omega.\]
\end{lemma}

In fact, it was shown by Evans-Smart \cite{Evans-Smart_2011a, Evans-Smart_2011b} that $u$ is everywhere differentiable (in all dimensions). The maximal gradient set of an $\infty$-harmonic function has a special structure: the following is a consequence of \cite[Theorem 2.2]{Savin-Wang-Yu_2008} or \cite[Proposition 6.2]{Crandall_2004}.

\begin{lemma}\label{lemma-prelim:interior_max_gradient}
	If $\Delta_\infty u=0$ and $|\D u|\leq L$ in $\Omega$, then the set $\{|\D u|=L\}$ is a disjoint union of nontrivial line segments which are streamlines of $u$ and have no endpoints in $\Omega$.
\end{lemma}

This lemma plays a crucial role in Section \ref{sec:closed}: these are the line segments that cut the concave polygons. The lemma implies that $\{x\in\Omega:|\D u(x)|=\sup_\Omega|\D u|\}$ cannot be precompact in $\Omega$ if nonempty. Hence, the maximum principle holds for $|\D u|$:

\begin{lemma}\label{lemma-prelim:max_principle_du}
	If $\Delta_\infty u=0$ in $\Omega\subset\RR^n$, and $K\Subset\Omega$, then $\sup_K|\D u|=\sup_{\p K}|\D u|$.
\end{lemma}

We recall Jensen's weak maximum principle:

\begin{lemma}[{\cite{Armstrong-Smart_2010, Jensen_1993}}]\label{lemma-prelim:jensen_max_principle}
	If $\Omega\Subset\RR^n$ and $u,v\in C^0(\bar\Omega)$ are respectively $\infty$-subharmonic and $\infty$-superharmonic in $\Omega$, then
	\[\max_{\bar\Omega}(u-v)=\max_{\p\Omega}(u-v).\]
\end{lemma}

In Section \ref{sec:ex}, we need a removable singularity theorem in Savin-Wang-Yu \cite{Savin-Wang-Yu_2008}:

\begin{lemma}\label{lemma-prelim:removable_sing}
	
	If $u\in C^1(\Omega)$ is $\infty$-harmonic in $\Omega\setminus\{z_0\}$, then $u$ is $\infty$-harmonic in $\Omega$.
\end{lemma}

In dimension 2, knowing the $C^{1,\alpha}$ regularity \cite{Evans-Savin_2008, Savin_2005} (in fact, knowing the $C^1$ equicontinuity is sufficient) is important to us throughout the paper:

\begin{lemma}[{\cite[Theorem 3]{Savin_2005}}]
	The family of $\infty$-harmonic functions in $B(2)\subset\RR^2$ that are bounded by 1 is $C^1$-equicontinuous in $B(1)$.
\end{lemma}

\subsection{Planar curves, concave polygons, support functions}\label{subsec:curves}

For a parametrized $C^1$ planar curve $\gamma$, we denote by $\tau=\gamma'$ its velocity vector. Whenever having fixed a unit speed parametrization of $\gamma$, it induces a choice of normal vector field $\nu=-\tau^\perp$. Conversely, whenever having chosen a normal vector field $\nu$, we arrange the length parametrization so that $\tau=\nu^\perp$.

\paragraph{Convex curves, angle parameter and support function.}\label{para-prelim:support_func} {\ }

Let $\gamma$ be an embedded open $C^1$ curve in $\RR^2$, and $\nu$ be a continuous unit normal vector on $\gamma$. We say that $\gamma$ is \textit{convex} (or \textit{strictly convex}) with respect to $\nu$, if for each $z\in\gamma$ there is $r>0$ and a rigid motion that sends $\gamma\cap B(z,r)$ to a convex (or strictly convex) graph and sends $\nu$ to $-\p_y$.

For an angle $\th\in\RR$ or $\th\in\SS^1(2k\pi)$, and $k\in\ZZ_+$, we fix the notation
\begin{equation}\label{eq-prelim:directional_vec}
	\nu_\th=(\sin\th,-\cos\th),\qquad\tau_\th=(\cos\th,\sin\th).
\end{equation}
Thus $\tau_\th=\nu_\th^\perp$.

Let $\gamma$ be an open convex $C^1$ curve with respect to a choice of normal vector $\nu$. Let $s$ be a compatible arclength parametrization (so that $\gamma'=\nu^\perp$). An \textit{angle parameter} of $\gamma$, usually denoted by $\varth$, is a continuous lift of $\nu^\perp$ to $\RR$. If $\gamma$ is the graph of a convex function $f$, then
\[\varth(x,f(x))=\arctan f'(x)\in(-\pi/2,\pi/2)\]
is a choice of angle parameter. If $\gamma$ is the graph of a concave function $f$ (it is convex with respect to the upward pointing normal and left-moving length parameter), then
\[\varth(x,f(x))=\pi+\arctan f'(x)\in(\pi/2,3\pi/2)\]
is a choice of angle parameter.

For the induced length parameter $s$, the function $s\mapsto\varth(\gamma(s))$ is non-decreasing. If $\varth(\gamma(s_1))=\varth(\gamma(s_2))$ for some $s_1<s_2$, then $\gamma|_{[s_1,s_2]}$ is a line segment.

Consider the following function on $\gamma$:
\[\bar h(\gamma(s))=\metric{\gamma(s)}{\nu(\gamma(s))}.\]
As each preimage $\varth^{-1}(\th)$ is a line segment orthogonal to $\nu$, the function $\bar h$ is constant on each $\varth^{-1}(\th)$ and hence descends to a function $h$ on the image $J=\varth(\gamma)$. We call $h$ the \textit{support function} of $\gamma$. Clearly,
\begin{equation}\label{eq-prelim:def_support_func}
	h(\th)=\metric{z}{\nu(z)}=\metric{z}{\nu_\th},\qquad\forall\,\th\in J,\ z\in\varth^{-1}(\th).
\end{equation}

When we translate $\gamma$ by a vector $a\p_x+b\p_y$, the angle parameter stays unchanged, while the support function becomes $h(\th)+a\sin\th-b\cos\th$. If we rotate $\gamma$ counterclockwise by an angle $\omega$ and denote the resulting curve by $\gamma_\omega$, then the angle parameter translates by $\omega$, and the support functions are related by
\[h_{\gamma_\omega}(\th)=h_\gamma(\th-\omega).\]

If $\gamma$ is smooth and strictly convex, then $\varth$ is a diffeomorphism onto its image, with $\frac{d\varth}{ds}=\kappa$, where $\kappa$ is the curvature of $\gamma$. Then we may also parametrize $\gamma$ by the angle. Differentiating \eqref{eq-prelim:def_support_func} with respect to $\th$, we have
\[\frac{dh}{d\th}=\Bmetric{\frac{d\gamma}{d\th}}{\nu_\th}+\metric{\gamma(\th)}{\tau_\th}=\metric{\gamma(\th)}{\tau_\th}.\]
Combining $h(\th)=\metric{\gamma(\th)}{\nu_\th}$, we have the recovery formula
\begin{equation}\label{eq-prelim:recovery_formula_smooth}
	\gamma(\th)=h\nu_\th+\frac{dh}{d\th}\tau_\th.
\end{equation}
We may further differentiate
\[\frac{d^2h}{d\th^2}=\Bmetric{\frac{d\gamma}{d\th}}{\tau_\th}-\metric{\gamma}{\nu_\th}=\kappa^{-1}-h(\th).\]
This gives the curvature formula
\begin{equation}\label{eq-prelim:curvature_formula}
	\frac1\kappa=\frac{d^2h}{d\th^2}+h.
\end{equation}

\paragraph{Geometric properties of support function.} {\ }

Let $\gamma$ be a convex $C^{1,1}$ curve with curvature bounded by $K$. Denote $J=\varth(\gamma)$. Then
\begin{equation}\label{eq-prelim:viscosity_convexity}
	\frac{d^2h}{d\th^2}+h\geq\frac1K\qquad\text{in the barrier sense in $\Int(\varth(\gamma))$.}
\end{equation}
To verify \eqref{eq-prelim:viscosity_convexity}, we may perform a rotation and assume that the reference point is $\th=0\in\Int(\varth(\gamma))$. Set $\sigma=\varth^{-1}(0)$. As $0\in\Int(J)$, notice that $\sigma$ has finite length and lies in the interior of $\gamma$. Since \eqref{eq-prelim:viscosity_convexity} is invariant under $h\mapsto h+a\sin\th-b\cos\th$, we may translate in $\RR^2$ and assume $\sigma=[-l,l]\times\{0\}$. Thus $\gamma$ is a convex graph in a neighborhood of $\sigma$. Denote $P=[-l,l]\times[1/K,\infty)$. By our convexity and curvature assumption, there is $\oth\ll1$ so that
\[\varth^{-1}\big((-\oth,\oth)\big)\subset\{y\geq0\}\setminus N(P,1/K).\]
Hence, the support function of $\gamma$ is no smaller than the support function of $\p N(P,1/K)$ in $(-\oth,\oth)$. But the latter function is calculated to be
\[\tilde h(\th)=K^{-1}+\sqrt{K^{-2}+l^2}\,\sin\Big(|\th|-\arctan(1/Kl)\Big),\]
hence
\[\tilde h+\frac{d^2\tilde h}{d\th^2}=\frac1K+2l\delta_0,\]
where $\delta_0$ is the Dirac delta centered at $\th=0$. This implies \eqref{eq-prelim:viscosity_convexity} as desired.

Suppose $\gamma$ is a $C^1$ convex curve, and $0\in\gamma$, $\varth(0)=0$, $0\in\Int(\varth(\gamma))$. Denote $\varth^{-1}(0)=[a,b]\times\{0\}$. Then in a small neighborhood of $\varth^{-1}(0)$, we may write $\gamma$ as the graph of a convex function $f:(a-r,b+r)\to\RR$ with $f|_{[a,b]}=0$ and $f|_{(a-r,a)\cup(b,b+r)}>0$. For each $x\in(b,b+r)$, the corresponding point $z=(x,f(x))$ satisfies
\begin{equation}\label{eq-prelim:h_convex_graph}
	\varth(z)=\arctan f'(x),\qquad\nu(z)=\frac{(f'(x),-1)}{\sqrt{1+f'(x)^2}},\qquad h(z)=\frac{xf'(x)-f(x)}{\sqrt{1+f'(x)^2}}.
\end{equation}
Thus
\begin{equation}\label{eq-prelim:h_convex_graph_2}
	h'_+(0)=\lim_{\th\searrow0}\frac{h(\th)}{\th}=\lim_{x\searrow b}\frac{h(z)}{\varth(z)}=b.
\end{equation}
Likewise, it can be shown that $h'_-(0)=a$. It follows that
\[\varth^{-1}(0)=[a,b]\times\{0\}=\big[h'_-(0),h'_+(0)\big]\times\{0\}.\]
Taking symmetry into account, we have in general
\begin{equation}\label{eq-prelim:recovery_formula_weak}
	\varth^{-1}(\th)=\Big\{h(\th)\nu_\th+c\tau_\th:c\in\big[h'_-(\th),h'_+(\th)\big]\Big\},\qquad\forall\,\th\in\Int(J).
\end{equation}
This is understood as the weak version of \eqref{eq-prelim:recovery_formula_smooth}. As a consequence, $h$ is differentiable at a point $\th\in\Int(J)$ if and only if $\varth^{-1}(\th)$ is a point.

Suppose $\gamma$ is a complete $C^1$ convex curve, and $\gamma$ is not a line, and $\gamma$ is a ray at both ends. Then $J=\varth(\gamma)=[\th_1,\th_2]$ for some $\th_1\leq\th_2$. If $\th_1=0$ and $\varth^{-1}(0)$ is contained in the $x$ axis, then arguing similarly as in \eqref{eq-prelim:h_convex_graph} \eqref{eq-prelim:h_convex_graph_2}, we have
\[\varth^{-1}(0)=\big(\!-\infty,h'_+(0)\big]\times\{0\}.\]
By symmetry, this implies the general recovery formula for the ray ends
\begin{equation}\label{eq-prelim:recovery_formula_endpt}
	\begin{aligned}
		\varth^{-1}(\th_1) &= \Big\{h(\th_1)\nu_{\th_1}+c\tau_{\th_1}:c\leq h'_+(\th_1)\Big\}, \\
		\varth^{-1}(\th_2) &= \Big\{h(\th_2)\nu_{\th_2}+c\tau_{\th_2}:c\geq h'_-(\th_2)\Big\}.
	\end{aligned}
\end{equation}

Conversely, given a smooth function $h$, we may construct its profile curve
\begin{equation}\label{eq-prelim:curve_recovering}
	\gamma(\th)=h(\th)\nu_\th+h'(\th)\tau_\th.
\end{equation}
The resulting curve $\gamma$ is parametrized by angle $\th$ and has support function $h$. We have
\[\begin{aligned}
	\frac{d\gamma}{d\th} &= \big(h'\nu_\th+h\tau_\th\big)+\big(h''\tau_\th-h'\nu_\th\big)
	= (h''+h)\tau_\th.
\end{aligned}\]
Thus $\gamma(\th)$ is a smoothly immersed curve as long as $h+h''\ne0$. It is convex with respect to $\nu_\th$ if $h+h''>0$, and is concave with respect to $\nu_\th$ if $h+h''<0$.

\paragraph{Concave polygons}\label{para-prelim:concave_polygon} {\ }

Let $E$ be a domain in an ambient domain $\Omega\subset\RR^2$. We say that $E$ is a \textit{concave polygon} in $\Omega$, if for each $z\in\p E\cap\Omega$, there is a rigid motion sending $z$ to $0$, after which we have:
\begin{itemize}[nosep]
	\item either
	\[E\cap Q(r)=\Big\{x\in(-r,r),\ y<g(x)\Big\}\]
	for a $C^1$ convex function $g:(-r,r)\to(-r/2,r/2)$ with $g(0)=g'(0)=0$;
	\item or
	\[E\cap Q(r)=\Big\{x\in(0,r),\ g_1(x)<y<g_2(x)\Big\}\]
	for a $C^1$ concave function $g_1:[0,r)\to(-r/2,0]$ and $C^1$ convex function $g_2:[0,r)\to[0,r/2)$ with $g_1(0)=g_2(0)=g'_1(0)=g'_2(0)=0$ and $g_1<g_2$ on $(0,r)$,
\end{itemize}
for some small $r>0$. Hence, $\p E$ is a union of concave curves (with respect to the outer unit normal), and $E$ is cusplike at the common endpoints of edges.

Let $\gamma$ be a connected component of $\p E$, with the counterclockwise orientation. Define the normal line map $\hat\nu:\gamma\to\SS^1(\pi)$ by $\hat\nu(z)=\text{(the unit normal vector of $E$ at $z$, mod $\pi$)}$. Looking at the local models, we find that this map is well-defined, continuous, and locally \textit{non-increasing} on $\gamma$ (due to the concavity). We are primarily concerned with two cases:
\begin{enumerate}[label={(\roman*)}, nosep]
	\item $E$ is simply-connected, $E\Subset\Omega$, and $\gamma=\p E$ contains $N$ edges where $N$ is even;
	\item $E$ is simply-connected, $E\not\Subset\Omega$, and $\gamma$ is a connected component of $\p E$.
\end{enumerate}
For case (i), $\gamma$ is a simple loop and $E$ lies inside $\gamma$. By Gauss-Bonnet, the total curvature of $\gamma$ is $-(N-2)\pi$, thus the map $\hat\nu$ has degree $-(N-2)$. We may lift $\hat\nu$ to a map
\[\nu:\gamma\to\SS^1.\]
Then $\nu$ is a continuous normal vector field on $\gamma$, and $\deg(\nu)=-\frac{N-2}2$. Denote $\tau=\nu^\perp$, which is a continuous tangential vector field on $\gamma$. In many cases, the edges of $\p E$ will be streamlines of a function $u$. In this case, we usually choose our lift so that $\tau=\D u/|\D u|$. The edges of $\gamma$ alternate between being concave and convex with respect to $\nu$, and $\nu=\pm\nu_E$ alternatively ($\nu_E$ always denotes the outer unit normal of $E$). Next, we may lift $\nu^\perp$ to a map
\[\varth:\gamma\to\SS^1\big((N-2)\pi\big).\]
Thus $\varth$ is a degree $-1$ non-increasing map. We call $\varth$ an \textit{angle parameter} of $\gamma$. For case (ii), we likewise define the lifted maps
\[\nu:\gamma\to\SS^1,\qquad\varth:\gamma\to\RR,\]
with $\varth$ being non-increasing as well. In either case, we denote $J=\varth(\gamma)$.

In both cases, each preimage $\varth^{-1}(\th)$ is a line segment orthogonal to $\nu_\th$. Hence, the map
\[\bar h(z)=\metric{z}{\nu(z)},\qquad z\in\gamma,\]
descends to a map
\[h(\th)=\metric{z}{\nu(z)}\big|_{z\in\varth^{-1}(\th)},\qquad\forall\,\th\in J,\qquad\text{where}\ \ J=\varth(\gamma),\]
which we call the \textit{support function} of $\gamma$. As $\nu(z)=\nu_\th$ for any $z\in\varth^{-1}(\th)$, we may also write
\[h(\th)=\metric{z}{\nu_\th}\big|_{z\in\varth^{-1}(\th)},\qquad\forall\,\th\in J.\]
It is easy to see that $h$ does not depend on the choice of $\varth$ once $\nu$ is fixed, and switching between the two choices of $\nu$ changes the sign of $h$. (Reminder: if $\p E$ has an odd number of edges, then the lift $\nu$ does not exist, hence $h$ is not defined.)

Let $\{\gamma_i\}$ be the edges of $\gamma$, and suppose that the interior of all $\gamma_i$ are $C^{1,1}$. We may arrange the indices so that $\gamma_i$ contacts $\gamma_{i-1}$ and $\gamma_{i+1}$. Then $\varth(\gamma_i)$ is an interval for each $i$. Denote it by $J_i$; then $\bar J=\bigcup\bar{J_i}$. Then note that
\[\pm\Big(\frac{d^2h}{d\th^2}+h\Big)>0\qquad\text{in}\ \ \Int(J_i),\]
due to \eqref{eq-prelim:viscosity_convexity}. Here, the sign changes alternatively in $i$. (Note: some $J_i$ may be a point.)

Suppose $\th\in\Int(J)$, and denote $\sigma=\varth^{-1}(\th)$. Then we have $\sigma\Subset\Omega$, by the monotonicity of $\varth$ and the connectedness of $\gamma$. With a rigid motion and a replacement $\varth\to\varth+2k\pi$, we may assume
\[\sigma=[x_1,x_2]\times\{0\}\quad(x_1\leq 0\leq x_2),\qquad\nu(0)=-\p_y,\qquad\varth(0)=0.\]
Then, the local shape of $E$ near $\sigma$ falls into one of the cases in Figure \ref{fig-prelim:models} below: shaded region represents $E$, and horizontal segment represents $\varth^{-1}(0)$. The second row is the ``orientation reversal'' of the first row. If $\sigma$ is a point, then cases 2,6 and cases 3,7 are of the same type, and cases 4,8 do not exist. In cases 1,5 we have $0\in\Int(J_i)$ for some $i$, and in cases 2,3,6,7 we have $0\in\p J_i\cap\p J_{i+1}$, and in cases 4,8 we have $J_i=\{0\}$.

\begin{figure}[ht]
	\centering
	\includegraphics{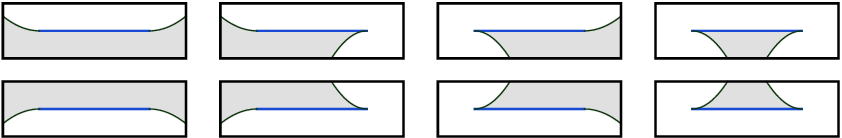}
	\caption{Local models near $\varth^{-1}(0)$.}\label{fig-prelim:models}
\end{figure}

Table \ref{table-prelim:models} collects the information of $h$ at $\th=0$. The items $\chi_\pm$ mean the sign of $h''+h$ for $\th<0$ resp. $\th>0$. The calculation of $h'_\pm(0)$ follows in the same way as \eqref{eq-prelim:h_convex_graph} \eqref{eq-prelim:h_convex_graph_2}.

\begin{table}[ht]
	\centering
	\begin{tabular}{|c|c|c|c|c|c|c|c|c|}
		\hline
		Model \# & 1 & 2 & 3 & 4 & 5 & 6 & 7 & 8 \\
		\hline
		$\chi_-$ & $+$ & $+$ & $-$ & $-$ & $-$ & $+$ & $-$ & $+$ \\
		\hline
		$\chi_+$ & $+$ & $-$ & $+$ & $-$ & $-$ & $-$ & $+$ & $+$ \\
		\hline
		$h'_-(0)$ & $x_1$ & $x_1$ & $x_1$ & $x_1$ & $x_2$ & $x_2$ & $x_2$ & $x_2$ \\
		\hline
		$h'_+(0)$ & $x_2$ & $x_2$ & $x_2$ & $x_2$ & $x_1$ & $x_1$ & $x_1$ & $x_1$ \\
		\hline
	\end{tabular}
	\caption{Information on the local models.}\label{table-prelim:models}
\end{table}

In particular, note that $h$ is left and right differentiable at $\th=0$, and
\[\sigma=\big[h'_-(0),h'_+(0)\big]\times\{0\},\]
where we denote $[a,b]=[b,a]$ if $a>b$. In view of symmetry, this shows that $h$ is everywhere left and right differentiable in $\Int(J)$, with
\begin{equation}\label{eq-prelim:recovery_formula_polygon}
	\varth^{-1}(\th)=\Big\{h(\th)\nu_\th+c\tau_\th:c\in\big[h'_-(\th),h'_+(\th)\big]\Big\},\qquad\th\in\Int(J).
\end{equation}
As a consequence, $h$ is differentiable at $\th\in\Int(J)$ if and only if $\varth^{-1}(\th)$ is a point, and in this case we have
\begin{equation}\label{eq-prelim:recovery_formula_point}
	\varth^{-1}(\th)=h(\th)\nu_\th+h'(\th)\tau_\th.
\end{equation}
Again, notice that \eqref{eq-prelim:recovery_formula_polygon} \eqref{eq-prelim:recovery_formula_point} are the polygonal versions of \eqref{eq-prelim:recovery_formula_smooth}.

\paragraph{Profile polygons.} {\ }

Suppose $h\in C^\infty(J)$ for $J=\SS^1(2k\pi)$ or $J\subset\RR$, so that $h''+h$ only has nondegenerate roots. Then
\begin{equation}\label{eq-prelim:profile_curve_2}
	\gamma(\th)=h(\th)\nu_\th+h'(\th)\tau_\th
\end{equation}
is an immersed concave polygonal curve (namely, it is a union of concave\,/\,convex curves with cusps at common endpoints), with angle parameter $\th$ and support function $h$. Each cusp corresponds to a root of $h''+h$. Indeed, suppose $\th=0$ is a nondegenerate root of $h''+h$. We may add a suitable multiple of $\sin\th$ and $\cos\th$ to $h$ (which is equivalent to a translation), and assume that $h(0)=h'(0)=0$. Hence, we obtain $h''(\th)=a\th+O(\th^2)$ near 0. By a scaling, we may assume $a=6$. Under these reductions, the Taylor expansion of $h$ becomes
\[h(\th)=\th^3+O(\th^4),\qquad h'(\th)=3\th^2+O(\th^3).\]
Then \eqref{eq-prelim:profile_curve_2} gives
\[\begin{aligned}
	\gamma(\th) &= \big(\th^3+O(\th^4)\big)(\sin\th,-\cos\th)+\big(3\th^2+O(\th^3)\big)(\cos\th,\sin\th) \\
	&= \big(3\th^2+O(\th^3),2\th^3+O(\th^4)\big).
\end{aligned}\]
Hence, the curve $\gamma$ is asymptotic to a canonical cubic cusp $y^2=4x^3/27$ near the origin.

Suppose $\gamma$ is a loop, namely, $J=\SS^1(2k\pi)$ for some $k\in\ZZ_{\geq1}$. Also, suppose that $\gamma$ has an even number of cusps. Let $\{\th_i\}$ be the roots of $h''+h$, and $J_i=[\th_i,\th_{i+1}]$, and $\gamma_i=\gamma(J_i)$. Let us assume that $h''+h>0$ in $J_0$, thus $(-1)^i(h''+h)>0$ in each $J_i$. Then we have the alternating length formula
\begin{equation}\label{eq-prelim:alternating_length}
	\int_{\SS^1(2k\pi)}h(\th)\,d\th=\sum_i(-1)^i|\gamma_i|.
\end{equation}
Indeed, if $s$ is a length parametrization, then $\big|\frac{d\th}{ds}\big|=|\kappa|=\frac1{|h''+h|}$ in each $\Int(\gamma_i)$, then
\[\begin{aligned}
	\sum(-1)^i|\gamma_i| &= \sum\int_{\gamma_i}(-1)^i\,ds=\sum\int_{\th_i}^{\th_{i+1}}(-1)^i|h''+h|\,d\th=\int_{\SS^1(2k\pi)}(h''+h)\,d\th \\
	&= \int_{\SS^1(2k\pi)}h\,d\th.
\end{aligned}\]

\subsection{Inverse mean curvature flow in \texorpdfstring{$\RR^2$}{R²}}\label{subsec:imcf}

\paragraph{Smooth IMCF in $\RR^2$ (or inverse curve shortening flow).} {\ }

A smooth family of strictly convex curves $\{\gamma_t\subset\Omega\}$ is a solution of the \textit{inverse mean curvature flow} (IMCF) if the normal speed is everywhere equal to $1/\kappa$. A function $w\in C^\infty(\Omega)$ with nonvanishing gradient is called a \textit{smooth IMCF} if it satisfies $\div\big(\frac{\D w}{|\D w|}\big)=|\D w|$. The latter is known as the level set form of IMCF. The transformation between these two formulations is given by $\gamma_t=\{w=t\}$: indeed, for a smooth function $w$ with nonvanishing gradient, the quantity $\div\big(\frac{\D w}{|\D w|}\big)$ is equal to the mean curvature of level sets, while $|\D w|$ is equal to $1/\text{(normal speed of movement)}$.

Suppose that $\{\gamma_t\}$ is smoothly parametrized by angles. Let $h(t,\th)$ be the (collection of) support functions of $\gamma_t$. Then recall from \eqref{eq-prelim:recovery_formula_smooth} that
\[\gamma_t(\th)=h(t,\th)\nu_\th+\p_\th h(t,\th)\tau_\th.\]
Calculating its speed:
\[\frac{\p\gamma_t}{\p t}=\frac{\p h}{\p t}(t,\th)\nu_\th+\frac{\p^2h}{\p t\,\p \th}(t,\th)\tau_\th.\]
The normal velocity is thus $\p_t h$. As the curvature of $\gamma_t$ is $1/(\p_{\th\th}h+h)$, the IMCF equation is therefore
\begin{equation}
	\frac{\p h}{\p t}=\frac{\p^2h}{\p\th^2}+h.
\end{equation}
This calculation can also be found in \cite{Chow-Tsai_1996, Urbas_1999}.

Conversely, suppose $h(t,\th)$ is a smooth function in a domain $R\subset\RR^2$, with $\p_th=\p_{\th\th}h+h>0$. Consider the map
\begin{equation}\label{eq-prelim:aux1}
	F(t,\th)=h(t,\th)\nu_\th+\p_\th h(t,\th)\tau_\th.
\end{equation}
We may compute
\begin{equation}\label{eq-prelim:aux2}
	\p_t F=(\p_th)\nu_\th+(\p_{t\th}h)\tau_\th,\qquad
	\p_\th F=(\p_{\th\th}h+h)\tau_\th,
\end{equation}
thus $\det(DF)=(\p_th)(\p_{\th\th}h+h)>0$. Therefore, $F$ is a local diffeomorphism, and the family of (possibly non-embedded) curves $\gamma_t=F\big(\{t\}\times(\th_1,\th_2)\big)$ form a smooth solution of IMCF. Whenever $F$ is known to be injective, it is then a diffeomorphism onto its image, and thus $\{\gamma_t\}$ form an embedded smooth IMCF in $F(R)$.

A particular case is $R=(t_1,t_2)\times(\th_1,\th_2)$ with $\th_2-\th_1<\pi$. Indeed, after a rotation, we may assume $0<\th_1<\th_2<\pi$, then $DF$ has $(2,2)$-entry $(\p_{\th\th}h+h)\sin\th>0$ and determinant $(\p_th)(\p_{\th\th}h+h)>0$. Then the Gale-Nikaid\^o theorem \cite[Theorem 4]{Gale-Nikaido_1965} implies that $F$ is injective, hence a diffeomorphism from $R$ to $F(R)$.

\paragraph{The weak IMCF.} {\ }

We call $w\in\Lip_{\loc}(\Omega)$ a \textit{weak IMCF} in $\Omega$, if every sub-level set $\{w<t\}$ is a local minimizer of the energy
\begin{equation}\label{eq-prelim:weak_imcf_energy}
	J_w^K(E):=\Ps{E;K}-\int_{E\cap K}|\D w|,
\end{equation}
meaning that $J_w^K(\{w<t\})\leq J_w^K(E)$ for all $t\in\RR$, all sets $E$ with locally finite perimeter, and all $K\Subset\Omega$, so that $E\Delta\{w<t\}\Subset K$. Here, $\Ps{E;K}$ denotes the perimeter of a set $E$ in $K$, see \cite{Maggi}. In most places, the objects have sufficient regularity, thus $\Ps{E;K}$ is the area of $\p E\cap K$.

Given a weak IMCF $w$, the sets $\p\{w<t\}$ are viewed as the evolving hypersurfaces. A hypersurface could instantly move forward at a time $t$ if $\{w=t\}$ has nonempty interior. A connected component of $\Int\big(\{w=t\}\big)$ is often called a jump. The weak IMCF can be (heuristically) understood as a combination of the $1/H$-evolution and jumps.

We refer the reader to \cite{Huisken-Ilmanen_2001, Lee, Xu_thesis} for general introductions of the weak IMCF and how the energy \eqref{eq-prelim:weak_imcf_energy} was motivated. The introduction in this subsection is specialized to providing preliminary materials for this paper. For the same reason, we specialize the ambient space to domains in $\RR^2$. As a note to experts, the perspective of strictly outward minimizing hulls (see \cite[Property 1.4]{Huisken-Ilmanen_2001} or \cite[Subsection 1.2]{Xu_thesis}) is less important for us, since the weak IMCF encountered in this paper are mostly non-proper. Nevertheless, the following geometric observation is important: the sets $\p\{w\leq t\}\setminus\p\{w<t\}$, namely the hypersurface reached after each jump, are minimal surfaces (i.e. line segments).

The perspective of calibration is more central. We call $w\in\Lip_{\loc}(\Omega)$ a \textit{calibrated weak IMCF}, if there is a measurable vector field $\nu$ (called a \textit{calibration} of $w$) that satisfies
\begin{equation}\label{eq-prelim:def_calibrated_imcf}
	|\nu|\leq1\quad\text{and}\quad\nu\cdot\D w=|\D w|\ \ \text{a.e.},\qquad\div(\nu)=|\D w|\ \ \text{distributionally.}
\end{equation}
Note that this implies
\begin{equation}\label{eq-prelim:div_identity}
	\div(e^{-w}\nu)=0\qquad\text{distributionally.}
\end{equation}
A smooth IMCF is always a calibrated weak IMCF: indeed, if $w\in C^\infty$ satisfies $\div\big(\frac{\D w}{|\D w|}\big)=|\D w|$, then $\nu=\D w/|\D w|$ serves as a calibration. A calibrated weak IMCF is always a weak IMCF; see the end of this subsection for the proof. Modulo analytic subtlety, this can be understood as follows: denoting $E_t=\{w<t\}$, we have $\nu=\nu_{E_t}$ on $\p E_t$, so for any competitor set $E$ the divergence theorem gives
\[\begin{aligned}
	J_w^K(E)-J_w^K(E_t) &= \Ps{E;K}-\Ps{E_t;K}-\int(\chi_E-\chi_{E_t})|\D w| \nonumber\\
	&\geq \int_{\p E\cap K}\metric{\nu_E}{\nu}-\int_{\p E_t\cap K}\metric{\nu_{E_t}}{\nu}-\int(\chi_E-\chi_{E_t})|\D w| \\
	&= 0.
\end{aligned}\]

The notion of calibrated weak IMCF first appeared in \cite[p.\,391]{Huisken-Ilmanen_2001}, and also played a central role in various $L^1$ variational problems, see for instance \cite{Andreu-Ballester-Caselles-Mazon_2000, Gorny-Mazon, Mazon-Leon_2015} and the references therein. We further say that $w$ is
\begin{itemize}[nosep]
	\item \textit{continuously calibrated}, if there is a continuous vector field $\nu$ satisfying \eqref{eq-prelim:def_calibrated_imcf};
	\item \textit{unit-calibrated}, if there is a continuous $\nu$ satisfying \eqref{eq-prelim:def_calibrated_imcf} with $|\nu|=1$ everywhere.
\end{itemize}
We include a list of useful properties of the weak IMCF. For the reader's convenience, we include the nontrivial proofs in the end of this subsection.
\begin{remark}\label{rmk-prelim:imcf_properties}
	Suppose $w$ is a weak IMCF in a domain $\Omega$.
	\begin{enumerate}[label={(\roman*)}, nosep]
		\item The regions $\{w<t\}$, $\{w\leq t\}$ are locally outward perimeter-minimizing \cite[Fact 1.2.9, 1.2.11]{Xu_thesis}, hence have mean convex boundaries.
		\item\label{item-prelim:imcf_gamma_t} The boundaries $\gamma_t=\p\{w<t\}\cap\Omega$, $\gamma_t^+=\p\{w\leq t\}\cap\Omega$ are $C^{1,1}$ with curvature almost everywhere $|\D w|$ for almost every $t$. See \cite[(1.12)]{Huisken-Ilmanen_2001}. As $w\in\Lip_{\loc}(\Omega)$, the curves $\gamma_t,\gamma_t^+$ are locally uniformly $C^{1,1}$. Observe that
		\[\bigcup_{s<t}\{w<s\}=\bigcup_{s<t}\{w\leq s\}=\{w<t\},\qquad\bigcap_{s>t}\{w<s\}=\bigcap_{s>t}\{w\leq s\}=\{w\leq t\}.\]
		Hence
		\begin{equation}\label{eq-prelim:convergence_of_lvlset}
			\lim_{s\nearrow t}\gamma_s=\lim_{s\nearrow t}\gamma_s^+=\gamma_t,\qquad\lim_{s\searrow t}\gamma_s=\lim_{s\searrow t}\gamma_s^+=\gamma_t^+,
		\end{equation}
		all in $C^1_{\loc}$ graphical sense.
		\item\label{item-prelim:imcf_orthogonality} If $w$ is continuously calibrated, then $\nu$ is the outer unit normal of $\gamma_t,\gamma_t^+$, for all $t$.
		\item\label{item-prelim:imcf_jump_minsurf} $\gamma_t^+\setminus\gamma_t$ is a union of line segments, for all $t$.
		\item\label{item-prelim:imcf_line_foliation} The regions $\Int\big(\{w=t\}\big)$ are called jumps. If $w$ is unit-calibrated, then each jump is foliated by line segments which are orthogonal to $\nu$.
		\item\label{item-prelim:imcf_extension} The following fact is useful in Section \ref{sec:heat}, \ref{sec:ereg}: suppose $f:[-r,r]\to(-l,l)$ is a $C^1$ convex function, and $w$ is a weak IMCF in the region
		\[D=\Big\{|x|<r,\ f(x)<y<l\Big\},\]
		so that $w\in\Lip(\bar D)$ and $w\leq T$ and $w|_{\graphh(f)}\equiv T$. Then the new function
		\[w'=\left\{\begin{aligned}
			& w\qquad\text{in}\ \ D \\
			& T\qquad\text{in}\ \ \big((-r,r)\times(-l,l)\big)\setminus D
		\end{aligned}\right.\]
		is a weak IMCF in $(-r,r)\times(-l,l)$.
	\end{enumerate}
\end{remark}

The following lemma is useful. Its proof is postponed to the end of this subsection.

\begin{lemma}\label{lemma-prelim:joint_cont_of_normals}
	Suppose $\alpha\in(0,1]$, and $\{E_i\}_{i\in I}$ is a family of open subsets of $\Omega$, so that:
	\begin{enumerate}[label={(\roman*)}, nosep]
		\item the $C^{1,\alpha}$ norms of $\p E_i$ are locally uniformly bounded in $\Omega$;
		\item for any $i,j\in I$ we have either $E_i\subset E_j$ or $E_j\subset E_i$.
	\end{enumerate}
	Then the outer unit normals of $E_i$ glue together to give a well-defined $C^{0,\alpha/(1+\alpha)}_{\loc}$ map
	\begin{equation}\label{eq-prelim:joint_unit_normal}
		\nu:\bigcup_i\p E_i\to\SS^1.
	\end{equation}
\end{lemma}
The H\"older exponent $\alpha/(1+\alpha)$ is sharp: for example, taking $\alpha=1$ and
\[\Omega=\RR^2,\qquad E_1=\{y<0\},\qquad E_2=\{y<x^2\},\]
the upward pointing unit normal of $\p E_1\cup\p E_2$ is exactly $C^{0,1/2}$ in $\RR^2$.

For a weak IMCF $w$, we consider
\begin{equation}\label{eq-prelim:imcf_Gamma0}
	\Gamma_0=\Big(\bigcup_{t\in\RR}\gamma_t\Big)\cup\Big(\bigcup_{t\in\RR}\gamma_t^+\Big).
\end{equation}
Then Lemma \ref{lemma-prelim:joint_cont_of_normals} yields a well-defined $C^{0,1/2}_{\loc}$ normal vector field
\begin{equation}\label{eq-prelim:imcf_nu}
	\nu_0:\Gamma_0\to\SS^1.
\end{equation}
When $w$ is continuously calibrated, by Remark \ref{rmk-prelim:imcf_properties}\ref{item-prelim:imcf_orthogonality} we have $\nu_0=\nu|_{\Gamma_0}$, and hence, the calibration $\nu$ is always $C^{0,1/2}$ on $\Gamma_0$. It is generally impossible to obtain regularity of $\nu$ outside $\Gamma_0$. If $w$ is known to be unit-calibrated, then Remark \ref{rmk-prelim:imcf_properties}\ref{item-prelim:imcf_line_foliation} and a finer use of Lemma \ref{lemma-prelim:joint_cont_of_normals} is possible to yield $\nu\in C^{0,1/2}_{\loc}(\Omega)$. This observation is not used anywhere in this paper, but it implies the following: if $u$ is an oriented $\infty$-harmonic function in $\Omega\subset\RR^2$ (in the sense of Moser \cite{Moser_2022}), then $u\in C^{1,1/2}_{\loc}(\Omega)$.

\paragraph{IMCF with outer obstacle.} {\ }

Let $D$ be a $C^1$ open set in an ambient domain $\Omega$. We call $w\in C^0(\bar D\cap\Omega)\cap\Lip_{\loc}(D)$ a \textit{continuously calibrated weak IMCF} in $D$ with \textit{outer obstacle} $\p D\cap\Omega$, if there is a continuous vector field $\nu$ on $\bar{D}\cap\Omega$, such that
\begin{equation}\label{eq-prelim:def_imcfoo_1}
	|\nu|\leq1,\qquad\metric{\nu}{\D w}=|\D w|,\qquad\div(\nu)=|\D w|\qquad\text{in}\ \ D,
\end{equation}
and
\begin{equation}\label{eq-prelim:def_imcfoo_2}
	\nu\cdot\nu_D\equiv1\qquad\text{on}\ \ \p D\cap\Omega.
\end{equation}
We often call \eqref{eq-prelim:def_imcfoo_2} the \textit{boundary tangency} or \textit{outer obstacle} condition. In the smooth case, this asks each hypersurface to stay tangent to $\p D$.

A solution $w$ of \eqref{eq-prelim:def_imcfoo_1} \eqref{eq-prelim:def_imcfoo_2} is a weak IMCF with outer obstacle $\p D\cap\Omega$ in the sense of \cite[Definition 3.1]{Xu_2024_obstacle}. This is due to \cite[Lemma 3.19]{Xu_2024_obstacle}. However, it is not known whether all solutions in \cite[Definition 3.1]{Xu_2024_obstacle} can be calibrated. As another difference, we assume $w$ here to be continuous up to $\p D$, while this is not assumed in \cite{Xu_2024_obstacle}. We choose to narrow the definition since all solutions in this paper are continuously calibrated and continuous up to $\p D$. We refer interested readers to \cite{Xu_thesis, Xu_2024_obstacle} for the general theory.

The solitons of IMCF in $\RR^2$ turn out to satisfy the boundary tangency condition in suitable domains. This is in particular the case for the shrinking soliton in (i) below which models quasiradial $\infty$-harmonic functions.

\begin{example}[Solitons]\label{ex-prelim:soliton} {\ }
	
	The translating, shrinking and expanding solitons of IMCF include cycloids, hypocycloids and epicycloids respectively, and other unbounded spiraling curves, see \cite{Urbas_1999} and \cite{Castro-Lerma_2017, Drugan-Lee-Wheeler_2016}. We are primarily interested in the classical curves.
	
	\begin{figure}[ht]
		\centering
		\includegraphics{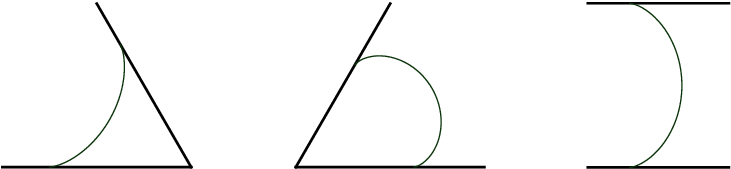}
		\begin{picture}(0,0)
			% Cycloid:
			\put(-55,-8){\arrowangle{0}}
			\put(-41,-11){$\th=0$}
			\put(-57,83){\arrowangle{180}}
			\put(-86,86){$\th=\pi$}
			% Hypocycloid:
			\put(-334,-8.5){\arrowangle{0}}
			\put(-320,-11.5){$\th=0$}
			\put(-297,59){\arrowangle{120}}
			\put(-294,68){$\th=\pi-\alpha$}
			% Epicycloid:
			\put(-160,-8.5){\arrowangle{0}}
			\put(-146,-11.5){$\th=0$}
			\put(-197,49){\arrowangle{240}}
			\put(-246,40){$\th=\pi+\alpha$}
		\end{picture}
		\vspace{9pt}
		\caption{Solitons of IMCF and their angle parameters.}\label{fig-prelim:solitons}
	\end{figure}
	
	\begin{enumerate}[label={(\roman*)}, topsep=1pt, itemsep=-0.5ex]
		\item For each $\alpha\in(0,\pi)$, consider the sector $D=\{\pi-\alpha<\Arg(z)<\pi\}$. There is a hypocycloid $\gamma_0$ contained in this sector and contacts the two boundary rays tangentially: $\gamma_0$ is the trajectory of the point $(-1,0)$ as we roll the disk $B(-1+\alpha/2\pi,\alpha/2\pi)$ clockwise against the unit circle. With an explicit calculation (see the end of this subsection), $\gamma_0$ has angle parameter ranging in $(0,\pi-\alpha)$ and has support function
		\begin{equation}\label{eq-prelim:single_hypocycloid}
			h(\th)=(1-\alpha/\pi)\sin\Big(\frac{-\th}{1-\alpha/\pi}\Big),\qquad\th\in(0,\pi-\alpha).
		\end{equation}
		Then, (after a scaling that removes the $(1-\alpha/\pi)$ factor,) consider
		\begin{equation}\label{eq-prelim:hypocycloid_h}
			h(t,\th)=\exp\Big(t-\frac{t}{(1-\alpha/\pi)^2}\Big)\sin\Big(\frac{-\th}{1-\alpha/\pi}\Big),\qquad t\in\RR,\ \th\in(0,\pi-\alpha).
		\end{equation}
		It is easy to see that $\p_th=\p_{\th\th}h+h>0$ in $\RR\times(0,\pi-\alpha)$. Denoting
		\[c=\frac{1}{(1-\alpha/\pi)^2}-1>0,\]
		the curves $\gamma_t=e^{-ct}\gamma_0$ solve the IMCF and have support functions $h(t,\cdot)$, up to a scaling. Another calculation (which leads to the same result) can be found in \cite[Section 3.1]{Xu_thesis}. This soliton (more precisely, the function $w$ so that $\{w=t\}=\gamma_t$) is then a smooth IMCF in $D$ that satisfies the boundary tangency condition on $\p D\setminus\{0\}$. In our language introduced above, setting the ambient domain as $\Omega=\RR^2\setminus\{0\}$, then $w$ is a continuously calibrated IMCF in $D$ with outer obstacle $\p D\cap\Omega$. Note: the origin needs to be removed since $w\to+\infty$ as $z\to0$.
		
		\item For each $\alpha\in(0,2\pi)$, consider the sector $\{0<\Arg(z)<\alpha\}$ and the epicycloid $\gamma_0$ supported in this sector. Now $\gamma_0$ is obtained as the trajectory of $(1,0)$ as we roll the exterior disk $B(1+\alpha/2\pi,\alpha/2\pi)$ counterclockwise against the unit disk. Through a similar calculation, the family of support functions
		\begin{equation}\label{eq-prelim:epicycloid_h}
			h(t,\th)=\exp\Big(t-\frac{t}{(1+\alpha/\pi)^2}\Big)\sin\Big(\frac{\th}{1+\alpha/\pi}\Big),\qquad\,t\in\RR,\ \th\in(0,\pi+\alpha)
		\end{equation}
		generate a self-expanding solution $\gamma_t=e^{ct}\gamma_0$, where $c=1-\frac{1}{(\alpha/\pi+1)^2}$.
		
		\item For each $\lambda>0$, consider the cycloid $\gamma_0$ obtained as the trajectory of $0$ by rolling $B(\lambda/2\pi,\lambda/2\pi)$ upward against the $y$-axis. Then $\gamma_0$ is calculated to have support function
		\begin{equation}\label{eq-prelim:single_cycloid}
			h(\th)=\frac\lambda\pi(\sin\th-\th\cos\th),\qquad\th\in(0,\pi).
		\end{equation}
		This induces a family of support functions
		\begin{equation}\label{eq-prelim:cycloid_h}
			h(t,\th)=\frac\lambda\pi(1+2t)\sin\th-\frac\lambda\pi\th\cos\th,\qquad t\in\RR,\ \th\in(0,\pi).
		\end{equation}
		One may verify that $\p_th=\p_{\th\th}h+h>0$. Then the self-translating solution
		\[\gamma_t=\gamma_0+\Big(\frac{2\lambda}\pi t,0\Big)\]
		has support function given by $h(t,\cdot)$.
	\end{enumerate}
\end{example}

We will need the following properties of solutions of \eqref{eq-prelim:def_imcfoo_1} \eqref{eq-prelim:def_imcfoo_2}:

\begin{remark}
	Let $w,\nu$ be as in \eqref{eq-prelim:def_imcfoo_1} \eqref{eq-prelim:def_imcfoo_2}, with the notations introduced there.
	\begin{enumerate}[label={(\roman*)}, nosep]
		\item Each sub-level set $\{w<t\}$, $\{w\leq t\}$ has $C^1$ boundary in $\Omega$. Their relative boundaries in $\Omega$ are orthogonal to $\nu$ everywhere, and contact $\p D\cap\Omega$ tangentially. Note: these sub-level sets are only mean convex in the interior $D$.
		\item The following is useful in Section \ref{sec:crit}. Denote $E_t=\{w<t\}$, then for any $C^1$ domain $F\Subset\Omega$ with $\H^{n-1}(\p E_t\cap\p F)=0$, we have
		\begin{equation}\label{eq-prelim:excess_ineq}
			|\p E_t\cap F|\leq e^{t-\inf_D(w)}|\p F\cap E_t|.
		\end{equation}
		Here $|\cdot|$ denotes the area. Indeed, using $\div(e^{-w}\nu)=0$ in $E_t\cap F$ we have
		\[\begin{aligned}
			0 &= \int_{\p E_t\cap F}\metric{e^{-w}\nu}{\nu_{E_t}}+\int_{\p F\cap E_t}\metric{e^{-w}\nu}{\nu_F},
		\end{aligned}\]
		Using $\nu=\nu_{E_t}$, $w\leq t$ on $\p E_t\cap F$, and $|\metric{\nu}{\nu_F}|\leq1$, $w\geq\inf_D(w)$ on $\p F\cap E_t$:
		\[e^{-t}|\p E_t\cap F|\leq\int_{\p E_t\cap F}\metric{e^{-w}\nu}{\nu_{E_t}}=-\int_{\p F\cap E_t}\metric{e^{-w}\nu}{\nu_F}\leq e^{-\inf_D(w)}|\p F\cap E_t|,\]
		as desired. In fact, \eqref{eq-prelim:excess_ineq} holds without calibration, see \cite[Lemma 4.2]{Xu_2024_obstacle}.
	\end{enumerate}
\end{remark}

\paragraph{Maximum principles.} {\ }

The following results are useful in Section \ref{sec:heat}. A function $w\in C^\infty(\Omega)$ with nonvanishing gradient is called a smooth \textit{subsolution} (resp. \textit{supersolution}) of IMCF if
\[\div\Big(\frac{\D w}{|\D w|}\Big)\geq|\D w|,\qquad\text{resp.}\qquad\div\Big(\frac{\D w}{|\D w|}\Big)\leq|\D w|.\]
This is equivalent to its level sets $\gamma_t=\{w=t\}$ evolving with normal speed at least $1/\kappa$ (resp. at most $1/\kappa$) everywhere.

A function $w\in\Lip_{\loc}(\Omega)$ is called a weak \textit{subsolution} (resp. \textit{supersolution}) of IMCF in a domain $\Omega$, if the energy \eqref{eq-prelim:weak_imcf_energy} is locally outward minimized (resp. inward minimized), namely, if $J_w^K(\{w<t\})\leq J_w^K(F)$ for all competitors $F\supset\{w<t\}$ (resp. $F\subset\{w<t\}$) with $\{w<t\}\Delta F\Subset K\Subset\Omega$.

\begin{remark}\label{rmk-prelim:max_principles} {\ }
	
	\begin{enumerate}[label={(\roman*)}, nosep]
		\item\label{item-prelim:max_prin_smooth} A smooth sub\,/\,supersolution of IMCF is also a weak sub\,/\,supersolution.
		\item\label{item-prelim:max_from_h_to_sol} Suppose $h(t,\th)$ is smooth in a box $R=(t_1,t_2)\times(\th_1,\th_2)$, with $\p_t h\geq\p_{\th\th}h+h>0$ (resp. $\p_{\th\th}h+h\geq\p_t h>0$) and $\th_2-\th_1<\pi$. Set $F(t,\th)=h(t,\th)\nu_\th+\p_\th h(t,\th)\tau_\th$. Arguing as in \eqref{eq-prelim:aux1} \eqref{eq-prelim:aux2}, the curves $\gamma_t=F\big(\{t\}\times(\th_1,\th_2)\big)$ form a smooth subsolution (resp. supersolution) of IMCF. Setting $w=t\circ F^{-1}$, namely, setting $w=t$ on each $\gamma_t$, then $w$ is a smooth subsolution (or supersolution) of IMCF.
		\item\label{item-prelim:max_prin_truncate} If $w$ is a weak subsolution (resp. supersolution) of IMCF in a domain, then so is $\min\{w,T\}$ for any $T\in\RR$. See \cite[Remark 2.1.4(iv)]{Xu_2024_obstacle}.
		\item\label{item-prelim:max_principle} If $w_1,w_2$ are a supersolution and subsolution of IMCF in $\Omega$, respectively, and $\{w_1<w_2\}\Subset\Omega$, then $w_1\geq w_2$ in $\Omega$. See \cite[Theorem 2.2]{Huisken-Ilmanen_2001}
		\item\label{item-prelim:max_principle_imcfoo} If $w_1$ is a continuously calibrated IMCF in $D$ with outer obstacle $\p D\cap\Omega$, and $w_2$ is a weak subsolution of IMCF in $D$ (with no boundary condition assumed), and
		\[\{w_1<w_2\}\Subset\Omega,\]
		then $w_1\geq w_2$ in $D$. Notice: we a priori allow $w_1<w_2$ near $\p D\cap\Omega$. This follows by combining Theorem 3.8 and Theorem 3.21 in \cite{Xu_2024_obstacle}.
	\end{enumerate}
\end{remark}

\paragraph{$q$-Harmonic functions as $q\to1$.} {\ }

Recall that the $q$-Laplacian is
\[\Delta_q v=\div\big(|\D v|^{q-2}\D v\big),\]
and that $v$ is $q$-harmonic if $\Delta_q v=0$ distributionally. Suppose $q\leq2$, and $v_q$ is a positive $q$-harmonic function in a domain $\Omega\subset\RR^n$. Take $w_q=(1-q)\log v_q$. It may be verified that $w_q$ satisfies
\begin{equation}\label{eq-prelim:p-imcf}
	\Delta_q w_q=|\D w_q|^q,
\end{equation}
see \cite{Moser_2007}. We have the interior gradient estimate \cite{Kotschwar-Ni_2009}, see also \cite{Moser_2007, Mari-Rigoli-Setti_2022}
\begin{equation}\label{eq-prelim:kotschwar_ni}
	\sup_{\Omega'}|\D w_q|\leq\frac{C(n)}{d(\Omega',\p\Omega)},\qquad\forall\,\Omega'\Subset\Omega.
\end{equation}
The following general observation was due to Moser \cite{Moser_2007, Moser_2008}: if $v_q$ is a sequence of positive $q$-harmonic functions, and $w_q:=(1-q)\log v_q$ converges in $C^0_{\loc}(\Omega)$ to a function $w\in\Lip_{\loc}(\Omega)$ as $q\to1$, then $w$ is a weak IMCF in $\Omega$.

Along with the convergence $w_q\to w$, it is also possible to show that the vector fields $\nu_q=|\D w_q|^{q-2}\D w_q$ converge weakly in $L^1_{\loc}(\Omega)$ to a calibration of $w$. This requires an adaptation of the arguments in \cite{Andreu-Ballester-Caselles-Mazon_2000, Mazon-Leon_2015}.

In the recent work \cite{BMRSX_2025} the following is shown. Suppose $E_0\Subset\Omega\Subset\RR^n$ are smooth domains. Consider the $q$-capacitary potentials $v_q$ satisfying
\[\left\{\begin{aligned}
	& \Delta_q v_q=0\qquad\text{in}\ \ \Omega\setminus\bar{E_0}, \\
	& v_q|_{\p E_0}=1,\qquad v_q|_{\p\Omega}=0.
\end{aligned}\right.\]
Then $w_q:=(1-q)\log v_q$ converges in $C^0_{\loc}(\Omega\setminus E_0)$ to the weak IMCF in $\Omega$ with initial value $E_0$ and outer obstacle $\p\Omega$ (see \cite[Definition 3.6]{Xu_2024_obstacle} for its meaning).

These observations partially inspire the main theorem of Section \ref{sec:approx} and its proof. In this paper, we have that $v_q$ are the conjugate of a $p$-harmonic function $u_p$, and $u_p$ that approximates a given $\infty$-harmonic function. The $C^1$ equicontinuity theorem in Section \ref{sec:equi} will imply that $|\D v_q|^{q-2}\D v_q$ converges in $C^0_{\loc}$ to a continuous calibration of $w$. This significantly simplifies the analysis when $q\to1$. Eventually, we do not need adaptations of the convergence arguments in \cite{Andreu-Ballester-Caselles-Mazon_2000, BMRSX_2025, Moser_2008, Mazon-Leon_2015}.

\paragraph{Proofs of facts and lemmas.} {\ }

We end this subsection with the proof of the facts mentioned above. In the proofs we employ the notations $E_t=\{w<t\}$, $E_t^+=\{w\leq t\}$ and $\gamma_t=\p E_t\cap\Omega$, $\gamma_t^+=\p E_t^+\cap\Omega$. It is well-known that each $E_t^+$ also locally minimizes the energy \eqref{eq-prelim:weak_imcf_energy} \cite[Fact 1.2.10]{Xu_thesis}.

\begin{proof}[A calibrated IMCF is a weak IMCF] {\ }
	
	Let $w\in\Lip_{\loc}(\Omega)$ be a calibrated weak IMCF. Then for any $v\in\Lip_{\loc}(\Omega)$ and set $K$ so that $\{v\ne w\}\Subset K\Subset\Omega$, we test $\div(\nu)=|\D w|$ with the function $v-w$, and obtain
	\[\int(v-w)|\D w|=-\int\metric{\nu}{\D v-\D w}\geq\int\big(|\D w|-|\D v|\big).\]
	Hence
	\[\int_K|\D w|+w|\D w|\leq\int_K|\D v|+v|\D w|.\]
	The well-known equivalence \cite[Lemma 1.1]{Huisken-Ilmanen_2001} then implies that $w$ is a weak IMCF.
\end{proof}

\begin{proof}[Proof of Remark \ref{rmk-prelim:imcf_properties}\ref{item-prelim:imcf_orthogonality}] {\ }
	
	For any $K\Subset\Omega$, the coarea formula and $\metric{\nu}{\D w}=|\D w|$ gives
	\[\int_{-\infty}^\infty\Ps{E_t;K}\,dt=\int_K|\D w|=\int_K\metric{\nu}{\D w}=\int_{-\infty}^\infty\int_{\p E_t\cap K}\Bmetric{\nu}{\frac{\D w}{|\D w|}}\,dt.\]
	Hence, $\nu=\frac{\D w}{|\D w|}$ almost everywhere on $\p E_t$ for a.e. $t$. On the other hand, since $w\in\Lip_{\loc}$, it is well-known that $\nu_{E_t}=\frac{\D w}{|\D w|}$ almost everywhere on $\p E_t$ for a.e. $t$. Hence
	\[\nu=\nu_{E_t}\text{ a.e. on $\p E_t$, for a.e. $t$.}\]
	Using the continuity of $\nu$ and regularity of $\p E_t$ and \eqref{eq-prelim:convergence_of_lvlset}, the full result follows.
\end{proof}

\begin{proof}[Proof of Remark \ref{rmk-prelim:imcf_properties}\ref{item-prelim:imcf_jump_minsurf}] {\ }
	
	Assume $0\in\gamma_t^+\setminus\gamma_t$. Then there is $r\ll1$ so that $w\geq t$ in $B(r)$ and $w\equiv t$ in $E_t^+\cap B(r)$. For any $E$ with $E\Delta E_t^+\Subset K\Subset B(r)$, we have
	\[\Ps{E;K}=J_w^K(E)+\int_{E\cap K}|\D w|\geq J_w^K(E)\geq J_w^K(E_t^+)=\P{E_t^+;K}.\]
	Hence, $E_t^+$ is perimeter-minimizing in $B(r)$.
\end{proof}

\begin{proof}[Proof of Remark \ref{rmk-prelim:imcf_properties}\ref{item-prelim:imcf_line_foliation}] {\ }
	
	 In each open disk $B\subset\Int\big(\{w=t\}\big)$, we have $0=\div(\nu)=\curl(\nu^\perp)$ distributionally, hence there is a function $f\in C^1(B)$ with $\D f=\nu^\perp$. Since $|\D f|\equiv1$, the streamlines of $f$ (which are integral curves of $\nu^\perp$) must all be line segments. The global statement follows from a continuation argument.
\end{proof}

\begin{proof}[Proof of Remark \ref{rmk-prelim:imcf_properties}\ref{item-prelim:imcf_extension}] {\ }
	
	Denote $R=(-r,r)\times(-l,l)$ and $E_t=\{w<t\}\subset D\subset R$. Notice that $\{w'<t\}=E_t$ when $t\leq T$, and $\{w'<t\}=R$ when $t>T$. Thus, it suffices to show that
	\[J_{w'}^K(E_t)\leq J_{w'}^K(F),\qquad\forall\,t\leq T,\ \forall\,F,\text{ so that }E_t\Delta F\Subset K\Subset R.\]
	We first assume $t<T$. Denote $D_\epsilon=\{|x|<r,\,f(x)+\epsilon<y<l\}$. By continuity, we have $E_t\subset D_\epsilon$ for some $\epsilon>0$. Since the latter set is outward perimeter-minimizing in $R$, we obtain
	\[\begin{aligned}
		J_{w'}^K(F) &= \Ps{F;K}-\int_{F\cap K}|\D w'|\geq\P{F\cap D_\epsilon;K}-\int_{F\cap K}|\D w'| \\
		&= J_w^{K\cap D_{\epsilon/2}}(F\cap D_\epsilon)-\int_{(F\setminus D_\epsilon)\cap K}|\D w'|.
	\end{aligned}\]
	Using the energy minimizing of $E_t$ in $D$, we have
	\[J_{w'}^K(F)\geq J_w^{K\cap D_{\epsilon/2}}(E_t)-\int_{(F\setminus D_\epsilon)\cap K}|\D w'|.\]
	Taking $\epsilon\to0$, it follows that
	\[J_{w'}^K(F)\geq J_w^{K\cap D}(E_t)=J_{w'}^K(E_t),\]
	as desired. Then, taking $t\nearrow T$ shows that $E_T$ minimizes $J_{w'}$ as well.
\end{proof}

\begin{proof}[Derivation of \eqref{eq-prelim:single_hypocycloid}] {\ }
	
	Set $r=\alpha/2\pi$. The given hypocycloid has the well-known parametrization
	\[\gamma(s)=-(1-r)e^{-is}-re^{\frac{1-r}ris},\]
	with parameter $s=\pi-\Arg(z)$ where $z$ is the point where rolling circle with the unit circle contact. The velocity is then computed as
	\[\gamma'(s)=(1-r)ie^{-is}\big[1-e^{is/r}\big].\]
	Therefore, the angle parameter is given by
	\[\varth(s)=\Arg\gamma'(s)=\frac\pi2-s-\frac{\pi-s/r}2=\frac{1-2r}{2r}s.\]
	Finally, switching to the angle parametrization, we obtain (notice $\nu_\th=e^{i(\th-\pi/2)}$)
	\[\begin{aligned}
		h(\th) &= \metric{\gamma(s(\th))}{e^{i(\th-\pi/2)}}
		= \Ree\Big[\Big(\!-(1-r)e^{-\frac{2r}{1-2r}i\th}-re^{\frac{1-r}{1-2r}2i\th}\Big)\cdot e^{i(\pi/2-\th)}\Big] \\
		&= \Ree\Big[\!-(1-r)ie^{-\frac{i\th}{1-2r}}-ire^{\frac{i\th}{1-2r}}\Big]
		= (1-2r)\sin\Big(\frac{-\th}{1-2r}\Big).
	\end{aligned}\]
	The result follows by inserting $r=\alpha/2\pi$. Note that
	\[s\in(0,\alpha)\qquad\Rightarrow\qquad\th\in\Big(0,\frac{1-2r}{2r}\alpha\Big)=(0,\pi-\alpha). \qedhere\]
\end{proof}

\begin{proof}[Derivation of \eqref{eq-prelim:single_cycloid}] {\ }
	
	Set $r=\lambda/2\pi$. Then the cycloid can be parametrized as
	\[\gamma(s)=r+is-re^{is/r}.\]
	Its velocity and angle parameter is
	\[\gamma'(s)=i(1-e^{is/r}),\qquad\varth(s)=\Arg\gamma'(s)=s/2r.\]
	Hence, the support function is
	\[h(\th)=\Ree\Big[\Big(r+2ri\th-re^{2i\th}\Big)\cdot e^{i(\pi/2-\th)}\Big]=2r(\sin\th-\th\cos\th).\qedhere\]
\end{proof}

\begin{proof}[Proof of Remark \ref{rmk-prelim:max_principles}\ref{item-prelim:max_prin_smooth}] {\ }
	
	Suppose $w$ is a smooth subsolution, namely, $\div\big(\frac{\D w}{|\D w|}\big)\geq|\D w|$. Denote $E_t=\{w<t\}$. For any competitor $F\supset E_t$ and $K$ as stated, the divergence theorem gives
	\[\begin{aligned}
		J_w^K(E_t)
		&= \Ps{E_t;K}-\int_{E_t\cap K}|\D w|
		= \int_{\p E_t\cap K}\Bmetric{\frac{\D w}{|\D w|}}{\nu_{E_t}}-\int_{E_t\cap K}|\D w| \\
		&= \int_{\p^*F\cap K}\Bmetric{\frac{\D w}{|\D w|}}{\nu_F}-\int_{F\setminus E_t}\div\Big(\frac{\D w}{|\D w|}\Big)-\int_{E_t\cap K}|\D w| \\
		&\leq \Ps{F;K}-\int_{F\cap K}|\D w|
		= J_w^K(F),
	\end{aligned}\]
	as desired. The case of supersolution is similar.
\end{proof}

\begin{proof}[Proof of Lemma \ref{lemma-prelim:joint_cont_of_normals}] {\ }

	Clearly, the contact of any $\p E_i$ with $\p E_j$ is tangential. Hence, the outer unit normals of $\p E_i$ agree on overlaps, hence the map \eqref{eq-prelim:joint_unit_normal} is defined. We first state a local version of the regularity, which will also be useful later.
	
	\begin{lemma}\label{lemma-prelim:two_graphs}
		Suppose $f,g$ are functions on $(-1,1)$, with $f'=g'$ whenever $f=g$, and
		\[\sup_{(-1,1)}|f'|\leq1,\qquad\sup_{(-1,1)}|g'|\leq1,\]
		and
		\[\|f'\|_{C^{0,\alpha}(-1,1)}\leq\mu,\qquad\|g'\|_{C^{0,\alpha}(-1,1)}\leq\mu,\]
		for some $\alpha\in(0,1]$ and $\mu>0$. Let $\nu$ be the (well-defined) upper normal vector field of $\graphh(f)\cup\graphh(g)$. For $s\leq1$, denote $S=\big(\graphh(f)\cup\graphh(g)\big)\cap\{|x|\leq 1-s\}$. Then
		\[\|\nu\|_{C^{0,\alpha/(1+\alpha)}(S)}\leq 16\mu^{\frac1{1+\alpha}}s^{-\alpha}.\]
	\end{lemma}
	
	Back to the main lemma. It suffices to show that for any $z\in\Omega$, there is $r>0$ so that $\nu\in C^{\alpha/(1+\alpha)}\big(B(z,r)\cap(\cup\,\p E_i)\big)$. We may assume $z=0$. Choosing $r$ sufficiently small and magnifying the picture by $r^{-1}$, we may assume that the $C^{1,\alpha}$ norm of all $\p E_i$ are bounded by $e^{-100}$ in $B(100)$, and aim ourselves at showing $\nu\in C^{0,\alpha/(1+\alpha)}\big(B(1)\cap(\cup\,\p E_i)\big)$.
	
	Fix an $i$ with $\p E_i\cap B(1)\ne\emptyset$. After a rotation, we may assume that $\p E_i\cap B(100)$ is a $e^{-90}$-Lipschitz subgraph. Hence
	\[\{y<-2\}\cap B(100)\subset E_i\cap B(100)\subset\{y<2\}\cap B(100).\]
	For any other $j$ with $\p E_j\cap B(1)\ne\emptyset$, there is $\th_j\in\SS^1$ so that
	\[\big\{\metric{z}{\nu_{\th_j}}<-2\big\}\cap B(100)\subset E_j\cap B(100)\subset\big\{\metric{z}{\nu_{\th_j}}<2\big\}\cap B(100).\]
	Since $E_i\subset E_j$ or $E_j\subset E_i$, these imply $|\th_j-\pi|<1/10$. So $\p E_j\cap B(100)$ is a $1/5$-Lipschitz graph (without rotation). Now for any $j,k$ with $\p E_j\cap B(1)\ne\emptyset$ and $\p E_k\cap B(1)\ne\emptyset$, we apply Lemma \ref{lemma-prelim:two_graphs} to obtain
	\[\|\nu\|_{C^{0,\alpha/(1+\alpha)}}\leq e^{10}\qquad\text{on}\ \ (\p E_j\cup\p E_k)\cap B(1).\]
	This implies the lemma since $j,k$ are arbitrary.
\end{proof}

\begin{proof}[Proof of Lemma \ref{lemma-prelim:two_graphs}]
	Note that
	\begin{equation}\label{eq-prelim:two_graphs_nu}
		\nu=\frac{(-f',1)}{\sqrt{1+(f')^2}}\text{\ \ \ on\ \ }\graphh(f),\qquad
		\nu=\frac{(-g',1)}{\sqrt{1+(g')^2}}\text{\ \ \ on\ \ }\graphh(g),
	\end{equation}
	and the map $a\to(-a,1)/(1+a^2)^{1/2}$ is 1-Lipschitz. It suffices to show that
	\begin{equation}\label{eq-prelim:two_graphs_aux1}
		|f'(x)-g'(x)|\leq 4\mu^{\frac1{1+\alpha}}s^{-\alpha}|f(x)-g(x)|^{\frac{\alpha}{1+\alpha}},\qquad x\in[-1+s,1-s].
	\end{equation}
	Indeed, once \eqref{eq-prelim:two_graphs_aux1} is shown, for all $x,y\in[-1+s,1-s]$ we have as desired
	\[\begin{aligned}
		\big|\nu(x,f(x))-\nu(y,g(y))\big| &\leq |f'(x)-g'(y)|
		\leq |f'(x)-g'(x)|+|g'(x)-g'(y)| \\
		&\leq 4\mu^{\frac1{1+\alpha}}s^{-\alpha}|f(x)-g(x)|^{\frac{\alpha}{1+\alpha}}+2\big(\mu|x-y|^\alpha\big)^{\frac1{1+\alpha}} \\
		&\leq 4\mu^{\frac1{1+\alpha}}s^{-\alpha}\Big[|f(x)-g(y)|^{\frac\alpha{1+\alpha}}+|g(y)-g(x)|^{\frac\alpha{1+\alpha}}\Big] \\
		&\qquad +2\mu^{\frac1{1+\alpha}}s^{-\alpha}|x-y|^{\frac\alpha{1+\alpha}} \\
		&\leq 8\mu^{\frac1{1+\alpha}}s^{-\alpha}\Big[|f(x)-g(y)|^{\frac\alpha{1+\alpha}}+|x-y|^{\frac\alpha{1+\alpha}}\Big] \\
		&\leq 16\mu^{\frac1{1+\alpha}}s^{-\alpha}\big|(x,f(x))-(y,g(y))\big|^{\frac{\alpha}{1+\alpha}}.
	\end{aligned}\]
	In the second line, we used $|g'(x)-g'(y)|\leq\min\big\{\mu|x-y|^\alpha,2\big\}$.
	
	As \eqref{eq-prelim:two_graphs_aux1} is linear, we may replace $f\to f-g$ and $g\to0$ (thus $f'=0$ whenever $f=0$). It suffices to show
	\[\|f'\|_{C^{0,\alpha}(-1,1)}\leq2\mu\qquad\Rightarrow\qquad|f'|\leq 4\mu^{\frac1{1+\alpha}}|f|^{\frac{\alpha}{1+\alpha}}s^{-\alpha}\ \ \text{in}\ \ [-1+s,1-s].\]
	Fix $x\in[-1+s,1-s]$. By changing $x\to-x$ and $f\mapsto-f$, we may assume $f(x)>0$ and $f'(x)>0$. If $f(x)\geq\mu s^{1+\alpha}$, then using $|f'|\leq2\mu$ we have
	\[f'(x)\leq2\mu\leq2\mu\big[f(x)\mu^{-1}s^{-(1+\alpha)}\big]^{\frac{\alpha}{1+\alpha}}=2\mu^{\frac1{1+\alpha}}f(x)^{\frac{\alpha}{1+\alpha}}s^{-\alpha}.\]
	Now assume $f(x)<\mu s^{1+\alpha}$. Set $\Delta x=[f(x)/\mu]^{1/(1+\alpha)}$, thus $x-\Delta x\in(-1,1)$. If $f(x-\Delta x)<0$, then there is $y\in(x-\Delta x,x)$ with $f'(y)=f(y)=0$, hence
	\[f'(x)\leq2\mu(x-y)^\alpha\leq2\mu(\Delta x)^\alpha=2\mu^{\frac1{1+\alpha}}f(x)^{\frac{\alpha}{1+\alpha}}.\]
	If $f(x-\Delta x)\geq0$, then
	\[\begin{aligned}
		0 &\leq f(x-\Delta x) = f(x)-\int_{x-\Delta x}^x f'(t)\,dt
		\leq f(x)-f'(x)\Delta x+\int_{x-\Delta x}^x2\mu(x-t)^\alpha\,dt \\
		&= f(x)-f'(x)f(x)^{\frac1{1+\alpha}}\mu^{-\frac1{1+\alpha}}+\frac2{1+\alpha}f(x),
	\end{aligned}\]
	hence
	\[f'(x)\leq 4\mu^{\frac1{1+\alpha}}f(x)^{\frac{\alpha}{1+\alpha}}.\]
	In either case, the result follows.
\end{proof}

\subsection{Sturmian theory for 1D heat equation}\label{subsec:sturmian}

Let $f(t,\th)$ be a function that solves the heat equation $\p_tf=\p_{\th\th}f$. Suppose $(t_0,\th_0)$ is a degenerate root of $f$, meaning that $f(t_0,\th_0)=\p_\th f(t_0,\th_0)=0$. Let $k$ be the order of this root, meaning that $f(t_0,\th_0)=\p_\th f(t_0,\th_0)=\cdots\p_\th^{k-1}f(t_0,\th_0)=0$ but $\p_\th^kf(t_0,\th_0)\ne0$. The following is shown by Angenent \cite{Angenent_1988} and also Angenent-Fiedler \cite{Angenent-Fiedler_1988}: for any $\epsilon>0$, there are $\mu,\delta<\epsilon$ so that
\begin{enumerate}[label={(\roman*)}, nosep]
	\item $f$ has no roots on the sides $[t_0-\mu,t_0+\mu]\times\{\th_0\pm\delta\}$, and $\th_0$ is the only root of $f$ in $\{t_0\}\times[\th_0-\delta,\th_0+\delta]$;
	\item $f(t_0+\tau,\cdot)$ has $n\in\{0,1\}$ nondegenerate roots in $(\th_0-\delta,\th_0+\delta)$ for all $0<\tau\leq\mu$, where $n=0$ if $k$ is even and $n=1$ if $k$ is odd;
	\item $f(t_0-\tau,\cdot)$ has $k$ nondegenerate roots in $(\th_0-\delta,\th_0+\delta)$ for all $0<\tau\leq\mu$.
\end{enumerate}
In particular, (ii) (iii) implies that the number of roots of $f(t,\cdot)$ locally drops by at least 2 whenever we encounter a degenerate root.

Item (i) follows from \cite[Theorem A, B]{Angenent_1988}, and from the same source it is known that $f(t_0+\tau,\cdot)$ has at most one root and $f(t_0-\tau,\cdot)$ has at least two roots. Also, from this one knows that the number of roots must locally drop at a degenerate root. The full item (iii) follows from the asymptotic analysis in \cite[p.566]{Angenent-Fiedler_1988}. The full item (ii) follows from the maximum principle: if $k$ is even, then $f$ has the same sign in the set
\[\p_0P:=\Big(\{t_0\}\times[\th_0-\delta,\th_0+\delta]\Big)\cup\Big([t_0,t_0+\mu]\times\{\th_0\pm\delta\}\Big),\]
which implies $f>0$ in $(t_0,t_0+\mu]\times(\th_0-\delta,\th_0+\delta)$. If $k$ is odd, then $f$ does not change sign in the following regions respectively: \[\begin{aligned}
	& \p_-P:=\{t_0\}\times[\th_0-\delta,\th_0)\cup[t_0,t_0+\mu]\times\{\th_0-\delta\}, \\
	& \p_+P:=\{t_0\}\times(\th_0,\th_0+\delta]\cup[t_0,t_0+\mu]\times\{\th_0+\delta\}.
\end{aligned}\]
If $f(t+\tau,\cdot)$ has more than one root in $(\th_0-\delta,\th_0+\delta)$ for some $\tau\in(0,\mu]$, then either $\{f>0\}$ or $\{f<0\}$ would have a connected component not intersecting $\p_0P$, contradicting the strong maximum principle.

Globally, suppose $f$ solves the heat equation on $(t_1,t_2)\times\SS^1$. Then the number of roots of $f(t,\cdot)$ is nonincreasing in $t$, and decreases by at least 2 whenever a degenerate root appears at time $t$. If the number of roots stays constant in $(t_3,t_4)\subset(t_1,t_2)$, then all the roots in $(t_3,t_4)\times\SS^1$ are nondegenerate.

If a function $f$ solves $\p_t f=\p_{\th\th}f+f$, then $e^{-t}f$ solves the heat equation. Therefore, the above theory holds for solutions of $\p_t f=\p_{\th\th}f+f$ as well.

For a function $0\not\equiv f\in C^0(\SS^1)$, recall its Fourier coefficients
\[a_0=\frac1{2\pi}\int_{\SS^1}f,\quad b_0=0,\quad a_k=\frac1\pi\int_{\SS^1}f\cos(k\th),\quad b_k=\frac1\pi\int_{\SS^1}f\sin(k\th).\]
The Sturm-Hurwitz theorem states that if $f$ satisfies $a_k=b_k=0$ for all $k\leq n$, then $f$ has at least $2n+2$ roots in $\SS^1$. See also Section \ref{sec:entire} for relevant discussions. In particular, if $f$ is smooth with $\int_{\SS^1}f=0$, then integrating by parts gives
\[\int_{\SS^1}(\p_{\th\th}f+f)=\int_{\SS^1}(\p_{\th\th}f+f)\cos\th=\int_{\SS^1}(\p_{\th\th}f+f)\sin\th=0,\]
hence $\p_{\th\th}f+f$ has at least 4 roots in $\SS^1$.

%\newpage

\section{IMCF clusters}\label{sec:cluster}

We start with defining simple IMCF clusters:

\begin{defn}\label{def-cluster:simple}
	We call $\big(w,\nu,\{D_i\},\{\chi_i\}\big)$ a simple IMCF cluster in a domain $\Omega\subset\RR^2$, where $w$ is a function and $\nu$ is a vector field in $\Omega$, if the following hold:
	\begin{enumerate}[label={(\roman*)}, nosep]
		\item $\{D_i\}$ are disjoint $C^1$ domains in $\Omega$, and $\Omega\setminus\bigcup D_i$ is the disjoint union of countably many $C^1$ curves $\gamma_j$ without endpoints in $\Omega$ (we call $\gamma_j$ ridges).
		\item $\chi_i\in\{\pm1\}$ for each $i$ (it is understood that $\chi_i$ is a choice of orientation);
		\item $w,\nu$ are continuous with $|\nu|\leq1$ in $\Omega$;
		\item $w|_{D_i}\in\Lip_{\loc}(D_i)$, and
		\begin{equation}\label{eq-cluster:def_calibration}
			\chi_i\nu\cdot\D w=|\D w|\quad\text{a.e.,}\qquad \div(\chi_i\nu)=|\D w|\quad\text{weakly\ \ in $D_i$,}
		\end{equation}
		and
		\begin{equation}\label{eq-cluster:def_ridge}
			\chi_i\nu=\nu_{D_i}\qquad\text{on}\ \ \p D_i.
		\end{equation}
	\end{enumerate}
	We say that the simple IMCF cluster is unit-calibrated, if $|\nu|=1$ everywhere.
\end{defn}

\begin{remark}\label{rmk-cluster:simple}
	The following facts are trivial but may help clarify the definition:
	\begin{enumerate}[label={(\roman*)}, nosep]
		\item Each $\p D_i$ is a union of ridges.
		\item A continuously calibrated weak IMCF is trivially a simple cluster, by setting $D_1=\Omega$, $\chi_1=1$ and $\{\gamma_j\}=\emptyset$. On the other hand, we remind that not all weak IMCFs are known to be continuously calibrated.
		\item $w|_{D_i}$ is a continuously calibrated weak IMCF with outer obstacle $\p D_i\cap\Omega$. The calibration is given by $\chi_i\nu$.
		\item Adjacent regions have the opposite orientations, due to \eqref{eq-cluster:def_ridge}.
		\item If we change $\nu\to-\nu$ and switch all orientations $\chi_i$ to $-\chi_i$, then the resulting data remains a simple IMCF cluster.
	\end{enumerate}
\end{remark}

\begin{remark}\label{rmk-cluster:simple_div_free}
	Since $e^{-w}\nu$ is divergence-free in each $D_i$ and is continuous across all $\gamma_j$, it follows that
	\begin{equation}\label{eq-cluster:div_free}
		\div(e^{-w}\nu)=0\qquad\text{distributionally in }\Omega.
	\end{equation}
\end{remark}

Then, we define mixed clusters which allow valleys in addition to ridges:

\begin{defn}\label{def-cluster:mixed}
	We say that $(w,\nu)$ is a mixed IMCF cluster in a domain $\Omega$, if there are points $\{z_i\}$ and disjoint relatively closed line segments $\sigma_j$ in $\Omega$ so that:
	\begin{enumerate}[label={(\roman*)}, nosep]
		\item $w:\Omega\to\RR\cup\{+\infty\}$ is continuous, and $e^{-w}\nu$ is a continuous vector field in $\Omega$;
		\item $\{w=+\infty\}\subset\{z_i\}$, and $w,\nu$ are constant with $\nu\equiv\pm(\sigma'_j)^\perp$ on each $\sigma_j$;
		\item Each compact subset of $\Omega$ meets finitely many $z_i$ and $\sigma_j$;
		\item Denote $A=\{z_i\}\cup(\cup\sigma_j)$. Then $(w,\nu)$ is a simple IMCF cluster in $\Omega\setminus A$.
	\end{enumerate}
	We say that $w$ is unit-calibrated if $|\nu|=1$ everywhere in $\{w<\infty\}$.
\end{defn}

In (iv), we mean that there exists a division of $\Omega\setminus A$ into regions $D_i$, and a choice of orientations $\chi_i$, so that $\big(w,\nu,\{D_i\},\{\chi_i\})$ is a simple IMCF cluster in $\Omega\setminus A$. Such a division need not be unique.

\begin{figure}[ht]
	\centering
	\includegraphics{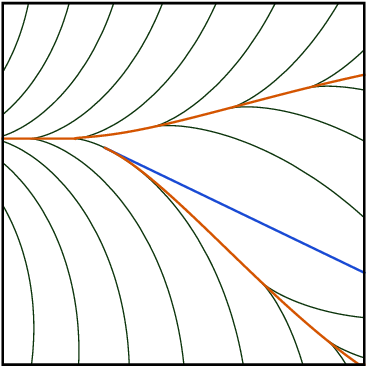}
	\caption{A local picture when a valley and a ridge contact.}\label{fig-cluster:mixed1}
\end{figure}

See Section \ref{sec:ex} for examples of mixed IMCF clusters and their associated $\infty$-harmonic functions. An exceptional segment $\sigma$ is called a \textit{valley} if $w$ forms two IMCFs evolving away from $\sigma$ in its neighborhood. More precisely, assuming that $\sigma$ lies in the $x$-axis and $\nu|_\sigma=-\p_y$ after a rigid motion, we say that $\sigma$ is a valley if for all points $(x,0)\in\Int(\sigma)$, there is $r>0$ so that $(x-r,x+r)\times(-r,0)$ and $(x-r,x+r)\times(0,r)$ lie in two regions $D_1,D_2$ in the simple cluster in (iv), so that
\[\chi_1=1,\qquad\chi_2=-1.\]
In the examples in Section \ref{sec:ex}, all exceptional segments are valleys. In general, an exceptional segment could be a valley or other trivial objects (e.g. a straight portion in the level set of a weak IMCF that is intentionally removed). A valley and a ridge (in $\Omega\setminus A$) may contact: Figure \ref{fig-cluster:mixed1} shows a possible local model.

\begin{figure}[ht]
	\centering
	\includegraphics{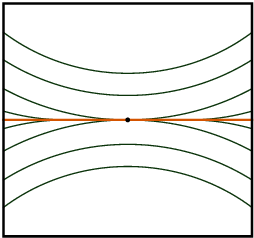}
	\caption{Exceptional point in mixed cluster.}\label{fig-cluster:excep_pt}
\end{figure}

A nontrivial exceptional point $z\in\{z_i\}$ contains two cases: near $z$, the solution $w$ either is asymptotic to the quasiradial solution in Example \ref{ex-ex:quasiradial}, or consists of two IMCFs colliding together as in Figure \ref{fig-cluster:excep_pt}. The second type of exceptional point may also be viewed as a degenerate valley. Trivial cases (e.g. a point removed from a simple cluster) may occur as well. For simplicity, we choose not to formally classify the types of exceptional segments and points in Definition \ref{def-cluster:mixed}.

We prove several auxiliary lemmas on simple IMCF clusters. In the lemmas, we denote
\[Q_\delta(r)=(-r,r)\times(-\delta r,\delta r)\]
for $\delta\leq1$. Note that $Q_\delta(r)=Q(r)$ if $\delta=1$.

\begin{lemma}\label{lemma-cluster:unique_ridge_in_square}
	Let $\big(w,\nu,\{D_i\},\{\chi_i\}\big)$ be a simple cluster in $\Omega\supset Q_\delta(r)$ with $\delta\leq1$, so that
	\begin{equation}\label{eq-cluster:almost_vertical}
		|\nu+\p_y|<e^{-10}\delta\qquad\text{in}\ \ Q_\delta(r).
	\end{equation}
	If $0\in\cup\gamma_j$, then $(\cup\gamma_j)\cap Q_\delta(r)$ is a connected $e^{-9}\delta$-Lipschitz graph.
\end{lemma}
\begin{proof}
	\eqref{eq-cluster:def_ridge} \eqref{eq-cluster:almost_vertical} imply that $(\cup\gamma_j)\cap Q_\delta(r)$ is a disjoint union of $e^{-9}\delta$-Lipschitz graphs whose boundaries are on $\p Q_\delta(r)$, and one of these graphs passes through 0. If the lemma is false, then there is a connected component of $Q_\delta(r)\setminus\cup\gamma_j$, called $D$, with the following property: $\p D\cap Q_\delta(r)$ contains two $C^1$ curves $\gamma_1,\gamma_2$ so that $\nu_D$ points almost upwards on $\gamma_1$ and points almost downwards on $\gamma_2$. If $D$ is positively oriented, then this violates \eqref{eq-cluster:def_ridge} \eqref{eq-cluster:almost_vertical} on $\gamma_1$. If $D$ is negatively oriented, then this violates \eqref{eq-cluster:def_ridge} \eqref{eq-cluster:almost_vertical} on $\gamma_2$.
\end{proof}

\begin{cor}\label{cor-cluster:loc_finiteness}
	In a simple cluster, each compact set intersects finitely many ridges.
\end{cor}
\begin{proof}
	This follows from Lemma \ref{lemma-cluster:unique_ridge_in_square}, the continuity of $\nu$, and a covering argument.
\end{proof}

Throughout this paper, we reserve the following notations:
\[Y_t=\{w>t\},\qquad\tZ_t=\{w\geq t\},\qquad Z_t=\Int(\tZ_t).\]
If $w$ is a weak IMCF, then it is known that $\tZ_t=\bar{Z_t}$, and $Y_t,Z_t$ are mean concave.

\begin{figure}[ht]
	\centering
	\includegraphics{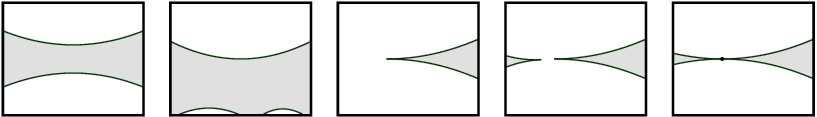}
	\captionsetup{width=0.85\textwidth}
	\caption{Possible local shapes of $Y_t$ or $Z_t$. The fifth picture prevents $Y_t,Z_t$ from being a concave polygon.}\label{fig-cluster:Yt}
\end{figure}

\begin{figure}[ht]
	\centering
	\includegraphics{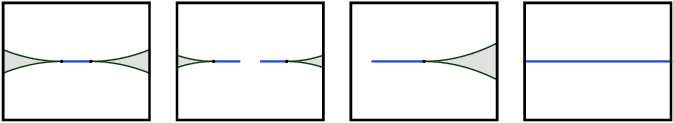}
	\caption{Possible local shapes $\tZ_t$. Blue segments represent $\tZ_t\setminus\bar{Z_t}$.}\label{fig-cluster:tZt}
\end{figure}

The following lemma concerns the shape of these sets in a simple cluster. Roughly saying, each $Y_t,Z_t$ is a ``concave polygon'' except that some vertices could stick together (Figure \ref{fig-cluster:Yt}(v)), and each $\tZ_t$ is a ``closed concave polygon'' but possibly with line segments attached to the vertices (Figure \ref{fig-cluster:tZt}). In any case where $-r<a=b<r$ in \ref{item-cluster:Yt_Zt_model} could be ruled out, it follows that $Y_t$ or $Z_t$ is a disjoint union of concave polygons. One notable outcome is that $\p Y_t$ and $\p Z_t$ must have a cusp when they travel across regions.

Using the strategies developed in Section \ref{sec:heat}, it can be shown that Figure \ref{fig-cluster:tZt}(i)(ii)(iii) cannot happen. More precisely, for any segment $\sigma\subset\p\tZ_t$ it can be shown that $\sigma\not\Subset\Omega$. As it is not used in this paper, we omit the formal statement.

\begin{lemma}\label{lemma-cluster:local_model}
	Suppose $\big(w,\nu,\{D_i\},\{\chi_i\}\big)$ is a simple IMCF cluster in a domain $\Omega$, with $\{\gamma_j\}$ the associated collection of ridges. Then the following hold for each $t\in\RR$.
	\begin{enumerate}[label={(\roman*)}, nosep]
		\item\label{item-cluster:model} Suppose $0\in\cup\gamma_j$ and $\Omega\supset Q_\delta(r)$ for a $\delta\leq1$, and
		\begin{equation}\label{eq-cluster:nu_almost_vertical}
			|\nu+\p_y|<e^{-10}\delta\qquad\text{in}\ \ Q_\delta(r).
		\end{equation}
		Then $(\cup\gamma_j)\cap Q_\delta(r)$ is the graph of a $e^{-9}\delta$-Lipschitz function $f$ on $(-r,r)$ with $f(0)=0$. In $Q_\delta(r)$ we have:
		\begin{enumerate}[label={(i\alph*)}, nosep]
			\item \label{item-cluster:Yt_Zt_model} Let $F$ denote either $Y_t$ or $Z_t$. Then there is a set $I$ of the form
			\[I=(-r,r)\qquad\text{or}\qquad I=(-r,a)\cup(b,r)\quad(\text{where }-r\leq a\leq b\leq r),\]
			and functions $g_1,g_2:\bar I\cap(-r,r)\to[-\delta r,\delta r]$, with $g_1<g_2$ and $g_1\leq f\leq g_2$ in $I$, so that
			\begin{equation}\label{eq-cluster:Yt_Zt_model}
				F\cap Q_\delta(r)=\Big\{x\in I,\ g_1(x)<y<g_2(x)\Big\}.
			\end{equation}
			If $F=Y_t$, then $g_1<f<g_2$ in $I$. More properties of $g_1,g_2$ are found in \ref{item-cluster:g1_g2}.
			\item\label{item-cluster:tZt_model} There is a set
			\[I=(-r,r)\qquad\text{or}\qquad I=(-r,a]\cup[b,r)\quad(\text{where} -r\leq a<b\leq r),\]
			and two functions $g_1,g_2:I\to[-\delta r,\delta r]$, with $g_1\leq f\leq g_2$, so that
			\begin{equation}\label{eq-cluster:tilde_Zt_model}
				\tZ_t\cap Q_\delta(r)=\Big\{x\in I,\ g_1(x)\leq y\leq g_2(x)\Big\}.
			\end{equation}
			\item\label{item-cluster:g1_g2} The functions $g_1,g_2$ in \ref{item-cluster:Yt_Zt_model} \ref{item-cluster:tZt_model} satisfy the following common properties:
			\begin{itemize}[nosep]
				\item\label{item-cluster:iiia} If $I=(-r,r)$, then $g_1$ is concave $($resp. $g_2$ is convex$)$ in each connected component of $\{g_1>-\delta r\}$ $($resp. $\{g_2<\delta r\})$;
				\item\label{item-cluster:iiib} If $I\ne(-r,r)$, then $g_1$ is concave and $g_2$ is convex in each connected component of $\bar I\cap(-r,r)$;
				\item\label{item-cluster:iiic} $g_1=g_2=f$ on $\p I\cap(-r,r)$, and $g'_1=g'_2=f'$ whenever $g_1=g_2$;
				\item\label{item-cluster:iiid} $\graphh(g_1)$ and $\graphh(g_2)$ have $\nu$ as unit normal vector.
			\end{itemize}
		\end{enumerate}
		\item\label{item-cluster:tZt_setminus_Zt} $\tZ_t\setminus\bar Z_t$ is contained in $\cup\gamma_j$ and is a disjoint union of nontrivial line segments.
		\item\label{item-cluster:no_interior_segment} Let $\sigma$ be a connected component of $\tZ_t\setminus\bar Z_t$. If $\sigma\Subset\Omega$, then $\bar\sigma\cap\p Z_t\ne\emptyset$.
		\item\label{item-cluster:monotonicity_on_ridge} For any $z\in\gamma_j$, either
		\begin{itemize}[nosep]
			\item $w$ is monotone along $\gamma_j$ in a small neighborhood of $z$, or
			\item $z$ is a strict local minimum of $w|_{\gamma_j}$ $($i.e. $w>w(z)$ in $U\cap\gamma_j\setminus\{z\}$ for some neighborhood $U\ni z$$)$.
		\end{itemize}
		The second type of points are discrete in $\Omega$.
	\end{enumerate}
\end{lemma}

\begin{proof}
	Fix $j$ so that $0\in\gamma_j$. By Lemma \ref{lemma-cluster:unique_ridge_in_square}, the set $(\cup\gamma_j)\cap Q_\delta(r)$ is connected and hence equals $\gamma_j\cap Q_\delta(r)$. We may write $\gamma_j\cap Q_\delta(r)$ as the graph of a $e^{-9}\delta$-Lipschitz function $f:(-r,r)\to(-e^{-9}\delta r,e^{-9}\delta r)$ with $f(0)=0$. Denote
	\[D_1=Q_\delta(r)\cap\{y<f(x)\},\qquad D_2=Q_\delta(r)\cap\{y>f(x)\}.\]
	By \eqref{eq-cluster:nu_almost_vertical} and that $\gamma_j$ is a ridge, we have $\chi_1=-1$ and $\chi_2=1$, and the cluster consists of an upward evolving IMCF in $D_1$ and a downward evolving IMCF in $D_2$. Thus
	\begin{equation}\label{eq-cluster:w_monotonicity}
		w(x,y)\text{ is }\left\{\begin{aligned}
			& \text{nondecreasing in $y$ for $y\leq f(x)$,} \\
			& \text{nonincreasing in $y$ for $y\geq f(x)$,}
		\end{aligned}\right.\qquad\forall\,x\in(-r,r),
	\end{equation}
	and that $y=f(x)$ attains a maximum of $y\mapsto w(x,y)$ in $(-\delta r,\delta r)$, for all $x\in(-r,r)$.
	
	\vspace{3pt}
	
	Recall $Y_t=\{w>t\}$. By \eqref{eq-cluster:w_monotonicity}, there is a function $\bg_1:(-r,r)\to[-\delta r,\delta r)$ with $\bg_1\leq f$, so that
	\[Y_t\cap D_1=\Big\{\!-r<x<r,\ \bg_1(x)<y<f(x)\Big\}.\]
	Since $w$ is a weak IMCF with outer obstacle in $D_1$ and is calibrated by $-\nu$, it follows that $\bg_1$ is concave in each connected component of $\{-\delta r<\bg_1<f\}$, and $\graphh(\bg_1)\cap Q_\delta(r)$ have $\nu$ as unit normal vector. For brevity, throughout this proof, we say that a function $g$ is $\nu$-orthogonal if $\nu$ is the unit normal vector of $\graphh(g)$. The $\nu$-orthogonality and \eqref{eq-cluster:nu_almost_vertical} imply that $\bg_1$ is $e^{-9}\delta$-Lipschitz. Similarly, there is a function $\bg_2:(-r,r)\to(-\delta r,\delta r]$ with $\bg_2\geq f$, so that
	\[Y_t\cap D_2=\Big\{\!-r<x<r,\ f(x)<y<\bg_2(x)\Big\},\]
	and $\bg_2$ has similar convexity properties as $\bg_1$. Furthermore, notice that $\bg_1(x)<f(x)<\bg_2(x)$ whenever $\bg_1(x)<\bg_2(x)$: this follows from \eqref{eq-cluster:w_monotonicity} and $Y_t=\{w>t\}$. And obviously, we have $\bg_1(x)=f(x)=\bg_2(x)$ whenever $\bg_1(x)=\bg_2(x)$.
	
	Set $I=\big\{x:\bg_1(x)<\bg_2(x)\big\}\subset(-r,r)$. Then
	\[Y_t\cap Q_\delta(r)=\Big\{x\in I,\ \bg_1(x)<y<\bg_2(x)\Big\}.\]
	We claim that either $I=(-r,r)$, or $I=(-r,a)\cup(b,r)$ for some $-r\leq a\leq b\leq r$. If this is not true, then some connected component $(c,d)$ of $I$ would satisfy $-r<c<d<r$. Recall that we have obtained:
	\begin{enumerate}[label={(\arabic*)}, nosep]
		\item $\bg_1(c)=f(c)=\bg_2(c)$, $\bg_1(d)=f(d)=\bg_2(d)$;
		\item $\bg_1<f<\bg_2$ in $(c,d)$;
		\item $\bg_1,\bg_2$ are $e^{-9}\delta$-Lipschitz and $f:(-r,r)\to(-e^{-9}\delta r,e^{-9}\delta r)$;
		\item $\bg_1$ (resp. $\bg_2$) is concave (resp. convex) in each connected component of $\{-\delta r<\bg_1<f\}$ (resp. $\{f<\bg_2<\delta r\}$).
	\end{enumerate}
	Combining (1)(3) we know that $\bg_1>-\delta r$ and $\bg_2<\delta r$ in $[c,d]$, then (4) implies that $\bg_1$ is concave and $\bg_2$ is convex in $(c,d)$ hence in $[c,d]$ as well. Combining this and (1)(2), we obtain a contradiction with the following elementary lemma, as desired.
	
	\begin{lemma}\label{lemma-cluster:elementary}
		There does not exist a concave function $f_1$ and convex function $f_2$ on an interval $J$, so that $f_1<f_2$ in $\Int(J)$ and $f_1=f_2$ on $\p J$.
	\end{lemma}
	
	Set $g_1=\bg_1|_{\bar I\cap(-r,r)}$ and $g_2=\bg_2|_{\bar I\cap(-r,r)}$. If $I=(-r,r)$, then note that $g_1<f<g_2$ in $(-r,r)$. The convexity property \ref{item-cluster:g1_g2} follows from (4). If $I\ne(-r,r)$, then clearly $g_1=g_2=f$ on $\p I\cap(-r,r)$. Since $\bg_1,\bg_2$ are $e^{-9}\delta$-Lipschitz, they never reach $\pm\delta r$. The convexity follows from (4) as well. Finally, the fact $g'_1=g'_2=f'$ whenever $g_1=g_2$ follows from the $\nu$-orthogonality of their graphs. This proves items \ref{item-cluster:Yt_Zt_model}\ref{item-cluster:g1_g2} for $F=Y_t$.
	
	\vspace{3pt}
	
	Then we analyze $\tZ_t$ using $\tZ_t=\bigcap_{s<t}Y_s$. If the ridge $\gamma_j\cap Q_\delta(r)$ is contained in $\tZ_t$, then for all $s<t$ there are functions $g_{1,s}:(-r,r)\to[-\delta r,\delta r)$ and $g_{2,s}:(-r,r)\to(-\delta r,\delta r]$ with $g_{1,s}<f<g_{2,s}$, so that
	\[Y_s\cap Q_\delta(r)=\Big\{x\in(-r,r),\ g_{1,s}(x)<y<g_{2,s}(x)\Big\}.\]
	Moreover, $g_{1,s}$ is concave and $\nu$-orthogonal in each connected component of $\{g_{1,s}>-\delta r\}$, and similarly for $g_{2,s}$. Taking an ascending\,/\,descending limit
	\[g_1=\lim_{s\nearrow t}g_{1,s},\qquad g_2=\lim_{s\nearrow t}g_{2,s},\]
	we find
	\[\tZ_t\cap Q_\delta(r)=\Big\{x\in(-r,r),\ g_1(x)\leq y\leq g_2(x)\Big\}.\]
	The convexity and $\nu$-orthogonality of $g_{1,s},g_{2,s}$ pass to $g_1,g_2$, thus proving \ref{item-cluster:tZt_model} \ref{item-cluster:g1_g2} for this case (the fact that $g'_1=f'=g'_2$ whenever $g_1=g_2$ follows from $\nu$-orthogonality). Note that $g_{1,s}<f<g_{2,s}$ passes to $g_1\leq f\leq g_2$, and equality may hold somewhere.
	
	Assume $\gamma_j\cap Q_\delta(r)\not\subset\tZ_t$. Pick $x_0\in(-r,r)$ so that $w(x_0,f(x_0))<t$. Then for all $s<t$ close enough to $t$, there are intervals $I_s=(-r,a_s)\cup(b_s,r)$ with $-r\leq a_s<x_0<b_s\leq r$, and functions $g_{1,s}:\bar{I_s}\cap(-r,r)\to[-\delta r,\delta r]$ and $g_{2,s}:\bar{I_s}\cap(-r,r)\to[-\delta r,\delta r]$ as in items \ref{item-cluster:Yt_Zt_model} \ref{item-cluster:g1_g2}, so that
	\[Y_s\cap Q_\delta(r)=\Big\{x\in I_s,\ g_{1,s}(x)<y<g_{2,s}(x)\Big\}.\]
	Recall that $g_{1,s},g_{2,s}$ are concave\,/\,convex in each interval in their domains.
	
	Notice the strict monotonicity (which comes from \eqref{eq-cluster:w_monotonicity} and the continuity of $w$):
	\begin{itemize}[nosep]
		\item $a_s<a_{s'}$ whenever $t>s>s'$ and $a_{s'}>-r$,
		\item $g_{1,s}(x)>g_{1,s'}(x)$ whenever $t>s>s'$ and $x\in(-r,a_s]$.
	\end{itemize}
	Similar facts hold for $b_s$, $g_{2,s}$ as well. So the following monotone limits are defined:
	\[a=\lim_{s\nearrow t}a_s,\qquad b=\lim_{s\nearrow t}b_s,\]
	and
	\[g_1=\lim_{s\nearrow t}g_{1,s},\qquad g_2=\lim_{s\nearrow t}g_{2,s}\qquad\text{in}\ \ C^1_{\loc}((-r,a])\ \ \text{and}\ \ C^1_{\loc}([b,r)).\]
	The $C^1$ convergence is due to $\nu$-orthogonality. It follows that
	\[\tZ_t\cap Q_\delta(r)=\bigcap_{s<t}\big(Y_s\cap Q_\delta(r)\big)=\Big\{x\in I,\ g_1(x)\leq y\leq g_2(x)\Big\},\]
	where $I:=(-r,a]\cup[b,r)$. Note that $a<x_0<b$ since $(x_0,f(x_0))\notin\tZ_t$. The convexity and $\nu$-orthogonality of $g_1,g_2$ follows from taking the limit. This proves \ref{item-cluster:tZt_model} \ref{item-cluster:g1_g2} for $\tZ_t$.
	
	\vspace{3pt}
	
	Next, consider $Z_t=\Int(\tZ_t)$. By \ref{item-cluster:tZt_model}, there is a set $\tilde I=(-r,r)$ or $\tilde I=(-r,\tilde a]\cup[\tilde b,r)$, and functions $\tg_1,\tg_2:\tilde I\to[-\delta r,\delta r]$ with properties stated in \ref{item-cluster:tZt_model} \ref{item-cluster:g1_g2}, so that
	\begin{equation}\label{eq-cluster:tZt_local}
		\tZ_t\cap Q_\delta(r)=\Big\{x\in\tilde I,\ \tg_1(x)\leq y\leq\tg_2(x)\Big\}.
	\end{equation}
	Letting
	\[I=\Big\{x\in\tilde I: \tg_1(x)<\tg_2(x)\Big\},\]
	it follows that
	\begin{equation}\label{eq-cluster:Zt_local}
		Z_t\cap Q_\delta(r)=\Big\{x\in I,\ \tg_1(x)<y<\tg_2(x)\Big\}.
	\end{equation}
	
	If $I=(-r,r)$, then setting $g_i=\tg_i$ already proves \ref{item-cluster:Yt_Zt_model}\ref{item-cluster:g1_g2} with $F=Z_t$. Now suppose $I\ne(-r,r)$. This implies $\tg_1=\tg_2$ hence $\tg_1=\tg_2=f$ somewhere. Thus $\tg_1,\tg_2$ do not reach $\pm\delta r$ due to $e^{-9}\delta$-Lipschitzness, and hence $\tg_1$ is concave and $\tg_2$ is convex in each connected component of $\tilde I$. The following fact is elementary:
	\begin{lemma}\label{lemma-cluster:two_convex_func}
		Suppose $f_1$ is concave and $f_2$ is convex in an interval, with $f_1\leq f_2$. Then $\{f_1=f_2\}$ is a subinterval, on which $f_1,f_2$ are both linear.
	\end{lemma}
	This implies that
	\begin{itemize}[nosep]
		\item if $\tilde I=(-r,r)$, then $\tilde I\setminus I$ is a subinterval;
		\item if $\tilde I\ne(-r,r)$, then $I$ has the form $I=(-r,a)\cup(b,r)$ with $a\leq\tilde a$ and $b\geq\tilde b$.
	\end{itemize}
	In either case, $I$ has the form as in \ref{item-cluster:Yt_Zt_model}. Then set $g_1=\tg_1|_{\bar I\cap(-r,r)}$ and $g_2=\tg_2|_{\bar I\cap(-r,r)}$. The convexity, $\nu$-orthogonality and endpoint tangency for $\tg_1,\tg_2$ pass to their restrictions. This shows \ref{item-cluster:Yt_Zt_model} \ref{item-cluster:g1_g2} for $F=Z_t$.
	
	\vspace{3pt}
	
	Furthermore, from \eqref{eq-cluster:tZt_local} \eqref{eq-cluster:Zt_local} and Lemma \ref{lemma-cluster:two_convex_func} we find
	\[(\tilde Z_t\setminus\bar{Z_t})\cap Q_\delta(r)=\Big\{x\in\tilde I\setminus\bar I,\ y=\tg_1(x)=\tg_2(x)=f(x)\Big\}.\]
	Hence, $(\tilde Z_t\setminus\bar{Z_t})\cap Q_\delta(r)$ is either empty or is a union of at most two nontrivial (open or half-open) line segments contained in $\graphh(f)$. This also implies that each segment is orthogonal to $\nu$. Since $\nu$ is continuous, this proves item \ref{item-cluster:tZt_setminus_Zt} by a covering argument.
	
	\vspace{3pt}
		
	Next, we prove \ref{item-cluster:no_interior_segment}. Suppose there exists a $\sigma$ as stated, so that $\bar\sigma\cap\p Z_t=\emptyset$. This implies $\bar\sigma\cap\bar{Z_t}=\emptyset$. By \ref{item-cluster:tZt_setminus_Zt} we know that $\sigma$ is nontrivial, and $\sigma\subset\gamma_j$ for some ridge $\gamma_j$. We may assume $\bar\sigma=[-l,l]\times\{0\}$ for some $l>0$, and $\nu=-\p_y$ on $\sigma$. Choose a small $r>0$ so that
	\[|\nu+\p_y|<e^{-10}\qquad\text{in}\ \ P:=(-l-r,l+r)\times(-r,r).\]
	So there is a $e^{-9}$-Lipschitz function $f:(-l-r,l+r)\to(-e^{-9}r,e^{-9}r)$ with $\gamma_j\cap P=\graphh(f)$. Recall that $\bar\sigma\cap\bar{Z_t}=\emptyset$, and from the proof of \ref{item-cluster:tZt_setminus_Zt} above, that $\tZ_t\setminus\bar{Z_t}$ is a discrete union of line segments in $\cup\gamma_j$ (i.e. any compact set intersects finitely many segments). Hence, for some sufficiently small $r$, we have $\tZ_t\cap\p P=\emptyset$. Then it follows that, for some $s<t$ close enough to $t$, it holds $Y_s\cap P\Subset P$.
	
	Arguing similarly as in \ref{item-cluster:Yt_Zt_model}, there are functions $\bg_1:(-l-r,l+r)\to(-r,r)$ and $\bg_2:(-l-r,l+r)$, with $\bg_1\leq f\leq\bg_2$, so that
	\[Y_s\cap P=\Big\{x\in(-l-r,l+r),\ \bg_1(x)<y<\bg_2(x)\Big\}.\]
	Note that $Y_s\cap P\ne\emptyset$ since $\sigma\subset\tZ_t\subset Y_s$. Then, a connected component of $\{\bg_1<\bg_2\}$ would take the form $(c,d)$ with $-l-r<c<d<l+r$, contradicting Lemma \ref{lemma-cluster:elementary}.
	
	\vspace{3pt}
	
	Finally, we prove \ref{item-cluster:monotonicity_on_ridge}. After a rotation, we may assume $z=0$ and $\exists\,r\ll1$ so that $|\nu+\p_y|<e^{-10}$ in $Q(r)$. Let $f:(-r,r)\to(-e^{-9}r,e^{-9}r)$ be such that $\gamma_j\cap Q(r)=\graphh(f)$. Call $w'(x)=w(x,f(x))$. By Item \ref{item-cluster:Yt_Zt_model}, for all $z\in\gamma_j\cap Q(r)$ and all squares $Q(z,s)\subset Q(r)$, we have $Y_t\cap Q(z,s)\not\Subset Q(z,s)$ if the former is nonempty. Hence
	\[\max_I(w')=\max_{\p I}(w'),\qquad\forall\text{ interval }I\subset(-r,r),\]
	or
	\begin{equation}\label{eq-cluster:aux1}
		w'(x_2)\leq\max\big\{w'(x_1),w'(x_3)\big\},\qquad\forall\,-r<x_1<x_2<x_3<r.
	\end{equation}
	It is then elementary to verify that: there can be at most one strict local minimum of $w'$ in $(-r,r)$, and $w$ is monotone in a neighborhood of each non strict local minimum.
\end{proof}

%\newpage

\section{From IMCF clusters to \texorpdfstring{$\infty$}{∞}-harmonic functions}\label{sec:ex}

Mixed IMCF clusters can be used to produce $\infty$-harmonic functions:

\begin{theorem}\label{thm-ex:reconstruction}
	Suppose $\Omega$ is simply-connected, and $(w,\nu)$ is a unit-calibrated mixed IMCF cluster in $\Omega$. Then there is an $\infty$-harmonic function $u\in C^1(\Omega)$ so that $\D u=e^{-w}\nu^\perp$.
\end{theorem}

We start with collecting several known examples and describing their IMCF cluster counterparts. Then in Subsection \ref{subsec:entire3}, we construct all degree 3 entire solutions in $\RR^2$ (meaning that they are asymptotic to degree 3 quasiradial solutions at $\infty$). In all figures, orange curves are ridges and blue segments are valleys.

\begin{example}[elementary examples]\label{ex-ex:classical} {\ }
	
	\begin{enumerate}[label={(\arabic*)}, nosep]
		\item The linear function $u(z)=x$ is $\infty$-harmonic in $\RR^2$. It corresponds to a constant weak IMCF with constant calibration.
		\item The cone function $u(z)=|z|$ is $\infty$-harmonic in $\RR^2\setminus\{0\}$. It corresponds to the constant weak IMCF with calibration $\nu=-z^\perp/|z|$.
		\item The angle function $u(re^{i\th})=\th$ is $\infty$-harmonic in $\RR^2$ with any ray removed (or, it is an $\infty$-harmonic map from $\RR^2\setminus\{0\}$ to $\SS^1$). It corresponds to the smooth IMCF $w=-\log|\D u|=\log|z|$ whose level sets are $\gamma_t=\{|z|=e^t\}$.
	\end{enumerate}
\end{example}

\begin{example}[quasiradial solutions]\label{ex-ex:quasiradial} {\ }
	
	For each $k=\frac{d^2}{2d-1}$, $d\in\ZZ_{\geq2}$, there is a (unique up to scaling and rotation) $\infty$-harmonic function of the form
	\begin{equation}\label{eq-ex:quasiradial}
		u(z)=r^k\phi(\th),
	\end{equation}
	where $\phi\in C^1(\SS^1)$ satisfies $\phi(\th+\pi/d)=-\phi(\th)$. We call it a \textit{degree $d$ quasiradial solution.} It was constructed by Aronsson in \cite{Aronsson_1984} (see also \cite{Aronsson_1986}) by converting the $\infty$-Laplacian into polar coordinates, and solving the corresponding second-order ODE \cite[(6)]{Aronsson_1986}
	\begin{equation}\label{eq-ex:quasiradial_ode}
		(\phi')^2\phi''+(2k-1)k\phi(\phi')^2+(k-1)k^3\phi^3=0.
	\end{equation}
	The generic solutions (see \cite[Section 5]{Aronsson_1984}) are given by the implicit representation formula
	\begin{equation}\label{eq-ex:quasiradial_formula}
		\phi=\frac Ck\Big(1-\frac1k\cos^2\omega\Big)^{\frac{k-1}2}\cos\omega,\qquad\th=\th_0+\int_0^\omega\frac{\sin^2x}{k-\cos^2x}\,dx,
	\end{equation}
	and $\phi(\th)$ is smooth in $(\th_0,\th_0+\pi/d)$ and $C^1$ in $[\th_0,\th_0+\pi/d]$, with
	\[\phi(\th_0+\pi/d)=-\phi(\th_0),\qquad\phi'(\th_0)=\phi'(\th_0+\pi/d)=0.\]
	Here, $C>0$, $\th_0\in\SS^1$ are generic parameters. Then, the solution is extended periodically to $\SS^1$. For the case $d=2$, this coincides (up to scaling) with the well-known explicit solution
	\begin{equation}\label{eq-ex:x4/3-y4/3}
		u=|x|^{4/3}-|y|^{4/3}.
	\end{equation}
	
	Let us consider the IMCF counterparts (see also Figure \ref{fig-intro:quasiradial}). Let $\{\bar\gamma_t\}$ be the shrinking hypocycloid supported in the sector $D_0=\big\{\pi-\pi/d<\operatorname{Arg}(z)<\pi\big\}$. See Example \ref{ex-prelim:soliton}(i), where it is set $\alpha=\pi/d$. It has support function
	\begin{equation}\label{eq-ex:quasiradial_building_piece}
		h_{\bar\gamma_t}(\th)=\exp\Big(\frac{1-2d}{(d-1)^2}t\Big)\sin\frac{-d\th}{d-1},\qquad t\in\RR,\ \th\in\Big(0,\frac{d-1}{d}\pi\Big).
	\end{equation}
	Then, consider $2d$ identical copies with each rotated by angle $\{-k\pi/d\}_{0\leq k<2d}$. Set
	\[\gamma_t=\bigcup_{0\leq k<2d}e^{-k\pi i/d}\bar\gamma_t,\qquad D_k=e^{-k\pi i/d}D_0,\]
	and set the function
	\[w|_{\gamma_t}=t,\qquad w(0)=+\infty,\]
	and set the orientation and vector
	\[\chi_k=(-1)^k,\qquad \nu=\chi_k\,\frac{\D w}{|\D w|}\ \ \text{in}\ \ D_k,\qquad \nu=(-1)^k\frac{z^\perp}{|z|}\ \ \text{in}\ \ -e^{-k\pi i/d}\,\RR_+,\]
	for each $k\in\{0,1,\cdots,2d-1\}$. Then, these data form a mixed IMCF cluster in $\RR^2$, with $z=0$ being an exceptional point, and $\{D_k\}$ be the natural partition of $\RR^2\setminus\{0\}$.
	
	Let $u$ be the resulting $\infty$-harmonic function by applying Theorem \ref{thm-ex:reconstruction}. Note that $\D u=e^{-t}\nu^\perp$ on $\gamma_t$. The homogeneity of $\{\gamma_t\}$ implies
	\[\D u(\lambda z)=\lambda^{\frac{(d-1)^2}{2d-1}}\D u(z),\qquad\forall\,\lambda\in(0,\infty),\]
	hence $u$ is $\frac{(d-1)^2}{2d-1}+1=\frac{d^2}{2d-1}$-homogeneous, hence a degree $d$ quasiradial solution.
	
	Denote
	\begin{equation}\label{eq-ex:quasiradial_joint_h}
		h(t,\th)=\exp\Big(\frac{1-2d}{(d-1)^2}t\Big)\sin\frac{-d\th}{d-1},\qquad\th\in\SS^1\big(2(d-1)\pi\big).
	\end{equation}
	Let us show that the support function of $\gamma_t$ is $h(t,\cdot)$. By Example \ref{ex-prelim:soliton}(i), this is true for the first copy of $\bar\gamma_t$, corresponding to $\th\in(0,\frac{d-1}d\pi)$. In the $k$-th copy of $\gamma_t$, considering the rotation and change of sign of the normal vector, the support function is given by
	\[(-1)^kh\Big(t,\th+\frac{k\pi}d\Big),\qquad\th\in\Big(\!-\frac{k\pi}d,\frac{d-1}d\pi-\frac{k\pi}d\Big).\]
	Changing variable $\th'=\th+k\pi$, this is the same as
	\[(-1)^kh\Big(t,\th'-\frac{d-1}dk\pi\Big),\qquad\th'\in\Big(\frac{d-1}dk\pi,\frac{d-1}d(k+1)\pi\Big).\]
	Noting that
	\[h\Big(t,\th'-\frac{d-1}dk\pi\Big)=(-1)^kh(t,\th'),\]
	this exactly fits as the $k$-th portion of \eqref{eq-ex:quasiradial_joint_h}, as desired.
	
	Thus, by rotating and scaling, for each $(a,b)\in\RR^2\setminus\{0\}$, there is a quasiradial solution $u$ so that $Y_t=\{|\D u|<e^{-t}\}$ has support function
	\[h(t,\th)=\exp\Big(\frac{1-2d}{(d-1)^2}t\Big)\Big[a\sin\frac{d\th}{d-1}+b\cos\frac{d\th}{d-1}\Big],\]
	for $t\in\RR$ and $\th\in\SS^1\big(2(d-1)\pi\big)$. By the recovery formula \eqref{eq-prelim:recovery_formula_polygon} and the fact that $\p Y_t$ are streamlines of $u$, we have the general gradient formula
	\[\D u\Big(h(t,\th)\nu_\th+\p_\th h(t,\th)\tau_\th\Big)=e^{-t}\tau_\th.\]
	
	Note: for a non-integer $d>1$, the formula \eqref{eq-ex:quasiradial_formula} still gives a $\frac{d^2}{2d-1}$-homogeneous $\infty$-harmonic function in an open sector with angle $\pi/d$, but this solution cannot be extended periodically to $\RR^2$. For $d\in\ZZ_{<0}$, one may also construct a $\frac{d^2}{2d-1}$-homogeneous $\infty$-harmonic function in $\RR^2\setminus\{0\}$, with $z=0$ being a nonremovable singularity. The corresponding IMCF cluster is made of $2|d|$ copies of the expanding soliton supported in a sector with angle $\pi/|d|$. We refer to \cite{Bhattacharya_2004, Bhattacharya_2005b, Savin-Wang-Yu_2008} for general studies on isolated singularities.
\end{example}

\begin{example}\label{ex-ex:entire_e1}
	There is an entire solution in $\RR^2$ of the form
	\begin{equation}\label{eq-ex:transl_sol}
		u=e^{-x}\phi(y),
	\end{equation}
	with $\phi(y)$ a $C^1$ periodic function that solves
	\[(\phi')^2\phi''+2\phi(\phi')^2+\phi^3=0.\]
	This solution corresponds to a simple IMCF cluster obtained by tessellating the translating cycloids. See Figure \ref{fig-ex:entire_e1}. Indeed, up to a scaling, the level sets $\gamma_t$ satisfy $\gamma_t=\gamma_s+(t-s,0)$. Hence
	\[\D u(x+\lambda,y)=e^{-\lambda}\D u(x,y),\]
	hence $u$ is of the form \eqref{eq-ex:transl_sol}.
	\begin{figure}[ht]
		\centering
		\includegraphics{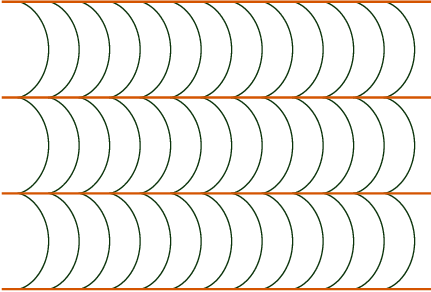}
		\caption{Translating soliton and solution $u=e^{-x}\phi(y)$.}\label{fig-ex:entire_e1}
	\end{figure}
\end{example}

\begin{example}[$\infty$-harmonic potential in the square]\label{ex-ex:potential_square} {\ }
	
	Let $u$ be the unique continuous function in $\bar{Q(1)}$, so that $u$ is $\infty$-harmonic in $Q(1)\setminus\{0\}$, and $u(0)=1$, $u|_{\p Q(1)}=0$. See \cite{Brustad_2022, Brustad-Lindgren-Lindqvist_2026, Lindgren-Lindqvist_2019, Lindgren-Lindqvist_2021} for related works. It is known that (see \cite{Lindgren-Lindqvist_2021}) $u$ is 1-Lipschitz in $Q(1)$, and is linear with slope 1 on the $x$- and $y$-axes. Thus, $u$ can be constructed from odd reflection by solving the following boundary value problem in $Q=(0,2)\times(0,2)$:
	\begin{equation}\label{eq-ex:square_new_eq}
		\left\{\begin{aligned}
			& \Delta_\infty u=0\qquad\text{in}\ \ Q, \\
			& \begin{aligned}
				& u(s,0)=1-s,\qquad u(s,2)=s-1, \\
				& u(0,s)=1-s,\qquad u(2,s)=s-1,
			\end{aligned}\qquad\forall\,s\in[0,2].
		\end{aligned}\right.
	\end{equation}
	To solve this boundary value problem, we view $\p Q$ as a degenerate concave 4-gon in the sense that adjacent edges do not meet tangentially. Consider the initial support function
	\[h_0(\th)=\left\{\begin{aligned}
		& 0\qquad(0\leq\th\leq\pi/2), \\
		& \!-2\cos\th\qquad(\pi/2\leq\th\leq\pi), \\
		& 2\sin\th-2\cos\th\qquad(\pi\leq\th\leq3\pi/2), \\
		& 2\sin\th\qquad(3\pi/2\leq\th\leq2\pi).
	\end{aligned}\right.\]
	Then let $h(t,\th)$ solve $\p_t h=\p_{\th\th}h+h$ on $[0,\infty)\times\SS^1$, with initial value $h_0$. Using similar argument as in Subsection \ref{subsec:entire3} below, it may be shown that the profile curves of $h$ generate a mixed IMCF cluster in $Q$, so that $(1,1)$ is an exceptional point and $\{y=x\}$, $\{y=2-x\}$ are ridges. See Figure \ref{fig-ex:potential_square}: the $x$- and $y$-axes are valleys, and $\{x=\pm y\}$ are ridges.
	\begin{figure}[ht]
		\centering
		\includegraphics{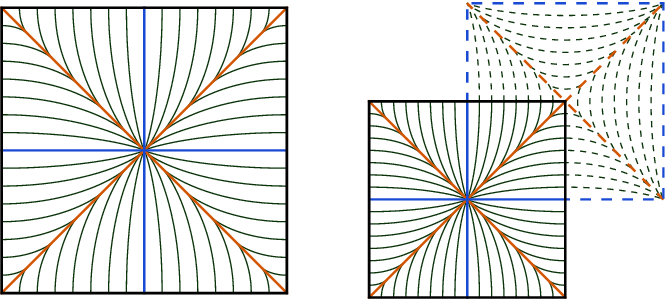}
		\caption{$\infty$-harmonic potential in $Q(1)$, and its odd reflection.}\label{fig-ex:potential_square}
	\end{figure}
\end{example}

We now prove Theorem \ref{thm-ex:reconstruction}. The following is a useful criterion for showing $\infty$-harmonicity. Recall that a streamline is an ascending gradient flow line.

\begin{lemma}\label{lemma-ex:linear_criterion}
	Suppose $u\in C^1(\Omega)$, and $z\in\Omega$, and $\D u(z)\ne0$. Suppose $\gamma:(-\epsilon,\epsilon)\to\Omega$ is a streamline of $u$ parametrized by arclength, with $\gamma(0)=z$.
	\begin{enumerate}[label={(\roman*)}, nosep]
		\item If $t\mapsto u(\gamma(t))$ is convex for $t\in[0,\epsilon)$, then $u$ is $\infty$-subharmonic at $z$.
		\item If $t\mapsto u(\gamma(t))$ is concave for $t\in(-\epsilon,0]$, then $u$ is $\infty$-superharmonic at $z$.
	\end{enumerate}
\end{lemma}
\begin{proof}
	We may assume $z=0$, $u(0)=0$ and $\D u(0)=\p_x$. Let $\varphi\in C^2$ be a test function that touches $u$ at $0$; hence $\varphi(0)=0$ and $\D\varphi(0)=\p_x$. We may Taylor expand $\varphi$:
	\[\varphi(x,y)=x+ax^2+bx y+cy^2+o\big(x^2+y^2\big).\]
	Note that $\D^2\varphi(\D\varphi,\D\varphi)(0)=2a$. We need to show that if $\varphi$ touches $u$ from above as in (i), then $a\geq0$, and if $\varphi$ touches $u$ from below as in (ii), then $a\leq0$.
	
	Let $\gamma$ be the given streamline. Then our conditions imply $u(\gamma(t))\geq t$ for $t\geq0$ in (i), or $u(\gamma(t))\leq t$ for $t\leq0$ in (ii). Since $\gamma$ is a $C^1$ curve, there is a continuous function $\th:(-\epsilon,\epsilon)\to\RR$, with $\th(0)=0$, so that $\gamma'(t)=e^{i\th(t)}$. Then
	\[\gamma(t)=\int_0^t\big(\!\cos\th(s),\sin\th(s)\big)\,ds=\Big(\int_0^t\cos\th(s)\,ds,\,o(t)\Big).\]
	The Taylor expansion of $\varphi$ gives
	\[\begin{aligned}
		\varphi(\gamma(t)) &= \int_0^t\cos\th(s)\,ds+at^2+o(t^2),\quad\text{which is}\quad
		\left\{\begin{aligned}
			&\leq t+at^2+o(t^2)\quad(t\geq0), \\
			&\geq t+at^2+o(t^2)\quad(t\leq0).
		\end{aligned}\right.
	\end{aligned}\]
	If $\varphi$ touches from above, then $t+at^2+o(t^2)\geq\varphi(\gamma(t))\geq u(\gamma(t))\geq t$ for all $t\geq0$, hence $a\geq0$. If $\varphi$ touches from below, then $t+at^2+o(t^2)\leq\varphi(\gamma(t))\leq u(\gamma(t))\leq t$ for all $t\leq0$, hence $a\leq0$, as desired.
\end{proof}

\begin{proof}[Proof of Theorem \ref{thm-ex:reconstruction}] {\ }
	
	Let $\{z_i\}$, $\{\sigma_j\}$ be as in Definition \ref{def-cluster:mixed}. Let $\sigma'_j$ be the interior of the segment $\sigma_j$. Denote
	\[A=\{z_i\}\cup(\cup\sigma_j),\qquad \Omega'=\Omega\setminus A,\qquad\Omega''=\Omega'\cup(\cup\sigma'_j).\]
	Note that $\Omega\setminus\Omega''$ is a discrete set, by Definition \ref{def-cluster:mixed}(iii). Recall that $(w,\nu)$ fits into a simple cluster in $\Omega'$. Let $\{D_i\}$ be any compatible choice of partition of $\Omega'$, with $\{\gamma_j\subset\Omega'\}$ be the corresponding set of ridges and $\{\chi_i\}$ the set of orientations.
	
	By Remark \ref{rmk-cluster:simple_div_free}, we have $\curl(e^{-w}\nu^\perp)=\div(e^{-w}\nu)=0$ distributionally in $\Omega'$. Then since $e^{-w}\nu^\perp$ is continuous near each $\sigma_j$, we have $\curl(e^{-w}\nu^\perp)=0$ distributionally in $\Omega''$. Finally, since $\Omega\setminus\Omega''$ is discrete, we have $\curl(e^{-w}\nu^\perp)=0$ distributionally in the entire $\Omega$. By the simply-connectedness of $\Omega$, there is $u\in C^1(\Omega)$ so that $\D u=e^{-w}\nu^\perp$.
	
	It remains to verify that $u$ is $\infty$-harmonic at every point $z\in\Omega$.
	
	\vspace{3pt}
	\textbf{Case 1:} $z\in D_i$ for some $i$. Recall that $w$ is a weak IMCF in $D_i$ and is calibrated by the unit vector field $\chi_i\nu$. If $z\in\p\{w<t\}$ or $z\in\p\{w\leq t\}$ for some $t$, then $\p\{w<t\}$ or $\p\{w\leq t\}$ is a streamline of $u$ on which $u$ is linear (under the length parametrization). Then we may apply Lemma \ref{lemma-ex:linear_criterion} as desired. If $z\in\Int\big(\{w=t\}\big)$ for some $t$, then $\nu^\perp$ gives a line foliation near $z$ (see Remark \ref{rmk-prelim:imcf_properties}\ref{item-prelim:imcf_line_foliation}), and each line is a streamline of $u$ on which $u$ is linear. The result follows as well.
	
	\vspace{3pt}
	\textbf{Case 2:} $z\in\sigma'_j$ for an exceptional curve $\sigma_j$. Then $\sigma'_j$ is a streamline of $u$ on which $u$ is linear, by Definition \ref{def-cluster:mixed}(ii), and the result follows easily from Lemma \ref{lemma-ex:linear_criterion}.
	
	\vspace{3pt}
	\textbf{Case 3:} $z$ lies on a ridge $\gamma$, and $w|_\gamma$ is monotone in a neighborhood of $z$. Then after a rigid motion and scaling and orientation reversing, we may assume
	\begin{itemize}[nosep]
		\item $\Omega'\supset Q(1)$, $z=0$, $u(0)=0$, $w(0)=0$, $\nu(0)=-\p_y$ (so $\D u(0)=\p_x$),
		\item $|\nu+\p_y|<e^{-10}$ and $|w|<e^{-10}$ in $Q(1)$,
		\item $\gamma\cap Q(1)=\graphh(f)$ for a $e^{-9}$-Lipschitz function $f:(-1,1)\to(-e^{-9},e^{-9})$,
		\item $w(x,f(x))$ is nondecreasing in $x$ for $x\in(-1,1)$,
		\item $Q(1)=D_1\cup\gamma\cup D_2$, where $D_1$ lies below $\gamma$, and $\chi_1=-1$, $\chi_2=1$.
	\end{itemize}
	Parametrize $\gamma$ by arclength, so that $\gamma(0)=0$ and $\gamma'(t)=\nu^\perp(\gamma(t))$. Notice that $\frac d{dt}u(\gamma(t))=\metric{\D u(\gamma(t))}{\gamma'(t)}=e^{-w(\gamma(t))}$ is non-increasing in $t$. Hence $u(\gamma(t))$ is concave in $t$ for $-\frac12\leq t\leq\frac12$, and by Lemma \ref{lemma-ex:linear_criterion}(ii), it follows that $u$ is $\infty$-superharmonic at $0$.
	
	Next, we show that $u$ is $\infty$-subharmonic at $0$. Suppose there exists $\varphi\in C^2(Q(1))$ that touches $u$ from above at $0$. Then $\varphi(0)=0$, $\D\varphi(0)=\p_x$. Since $u\leq\varphi\in C^2$, there is a constant $C>1$ so that
	\begin{equation}\label{eq-ex:reconst_aux1}
		u(0,y)\leq Cy^2\qquad\forall\,y\in(-1,1).
	\end{equation}
	We claim that $w(x,0)=0$ for all $x\in[-\frac1{8C},0]$. Once this is shown, then $|\p_x u(x,0)|\leq e^{-w(x,0)}=1$, hence $\varphi(x,0)\geq u(x,0)\geq x$ for all $x\in[-\frac1{8C},0]$, hence $\p_x^2\varphi(0)\geq0$, as desired.
	
	It remains to prove our claim. For each $x\in[-1,1]$, note that $w(x,y)$ is nondecreasing in $y$ when $y\leq f(x)$, and is nonincreasing in $y$ when $y\geq f(x)$. In addition, we have assumed that $w|_\gamma$ is nondecreasing in $x$. Hence, $w(x,0)\leq0$ for all $x\leq0$. If our claim is not true, then there exists $t<0$ and $x_0\in(-1/8C,0)$ so that $[x_0,0]\times\{0\}\subset\bar{Y_t}$ while $(x_0,0)\in\p Y_t$. By Lemma \ref{lemma-cluster:local_model}\ref{item-cluster:Yt_Zt_model} and the monotonicity of $w|_{\gamma\cap Q(1)}$, there is a subset $I=(b,1)$ with $b\leq x_0$ and functions $g_1,g_2:[b,1)\to\RR$, so that
	\[Y_t\cap Q(1)=\Big\{x\in I,\ g_1(x)<y<g_2(x)\Big\}.\]
	So either $g_1(x_0)=0$ or $g_2(x_0)=0$. We may assume $g_1(x_0)=0$ (same argument if $g_2(x_0)=0$). So $g_1$ is a concave $e^{-9}$-Lipschitz function on $[b,1)$. Denote $\sigma=\graphh(g_1|_{[x_0,0]})$ and $y_0=g_1(0)$. Notice that $|y_0|\leq|x_0|$. Recall $\div(e^{-w}\nu)=0$, and observe that
	\[\metric{\p_x}{e^{-w}\nu}=\metric{\p_y}{e^{-w}\nu^\perp}=\p_y u.\]
	The divergence theorem and \eqref{eq-ex:reconst_aux1} and $w\geq t$ on $[x_0,0]\times\{0\}$ gives
	\[\begin{aligned}
		0 &= \int_\sigma\metric{-(\sigma')^\perp}{e^{-w}\nu}+\int_{[x_0,0]\times\{0\}}\metric{\p_y}{e^{-w}\nu}+\int_{\{0\}\times[0,y_0]}\metric{\p_x}{e^{-w}\nu} \\
		&= e^{-t}|\sigma|+\int_{[x_0,0]\times\{0\}}\metric{\p_y}{e^{-w}\nu}+\big[u(0,0)-u(0,y_0)\big] \\
		&\geq e^{-t}\sqrt{x_0^2+y_0^2}-e^{-t}|x_0|-Cy_0^2
		\geq \frac{e^{-t}y_0^2}{2\sqrt{x_0^2+y_0^2}}-Cy_0^2
		\geq \Big(\frac1{4|x_0|}-C\Big)y_0^2.
	\end{aligned}\]
	Hence $y_0=0$, and then $0\in\p Y_t$, then $t=w(0)=0$, contradicting our hypothesis $t<0$.
	
	\vspace{3pt}
	
	\textbf{Case 4:} $z$ lies on a ridge $\gamma$, and $w|_\gamma$ is not monotone near $z$. By Lemma \ref{lemma-cluster:local_model}\ref{item-cluster:monotonicity_on_ridge}, these points are discrete in $\Omega'$. Then Lemma \ref{lemma-prelim:removable_sing} implies that $u$ is $\infty$-harmonic at $z$ as well.
	
	\vspace{3pt}
	
	In summary, we have shown that $u$ is $\infty$-harmonic at all points in $\Omega'$ and in all $\sigma'_j$. The remaining set of points is discrete in $\Omega$, hence Lemma \ref{lemma-prelim:removable_sing} applies to show that $u$ is $\infty$-harmonic in $\Omega$.
\end{proof}

\subsection{Degree 3 entire solutions}\label{subsec:entire3}

For each $a,b,c,d\in\RR$ with $(a,b)\ne(0,0)$ and $(c,d)\ne(0,0)$, we construct a mixed IMCF cluster in $\RR^2$ whose level set $\gamma_t$ for all $t\ll-1$ has support function
\begin{equation}\label{eq-ex:entire_3_asymp}
	h(t,\th)=e^{-5t/4}\Big(a\sin\frac32\th+b\cos\frac32\th\Big)+e^{3t/4}\Big(c\sin\frac\th2+d\cos\frac\th2\Big),\quad\th\in\SS^1(4\pi).
\end{equation}
The $e^{-5t/4}$ terms dominate as $t\to-\infty$, thus the resulting $\infty$-harmonic function is asymptotic to a degree 3 quasiradial solution at $\infty$. It also has two critical points at which it is asymptotic to degree 2 quasiradial solutions (see the proof of Theorem \ref{thm-closed:main}\ref{item-closed:main_asymp} and \ref{thm-entire:poly_growth}).

As was mentioned in the introduction, the support function \eqref{eq-ex:entire_3_asymp} is intended as a perturbation of the quasiradial model by adding lower order terms. A generic perturbation should take the form
\[h(t,\th)=\sum_{j=0}^3\exp\Big(t-\frac{j^2}4t\Big)\Big(a_j\sin\frac{j\th}2+b_j\cos\frac{j\th}2\Big),\qquad\th\in\SS^1(4\pi).\]
To see how this is reduced to \eqref{eq-ex:entire_3_asymp}: first, noting that $h\mapsto h+a\sin\th-b\cos\th$ corresponds to translations, we may assume that the $j=2$ term vanishes. We claim that the constant term must also vanish. By \eqref{eq-prelim:alternating_length}, the quantity $\int_{\SS^1(4\pi)}h(t,\th)\,d\th$ is the alternating length of $\gamma_t$: namely, denoting by $\alpha_1,\cdots,\alpha_6$ the edges of $\gamma_t$, we have
\[\int_{\SS^1(4\pi)}h(t,\th)\,d\th=\sum_i(-1)^i|\alpha_i|.\]
However, each $\alpha_i$ is a streamline of the resulting $\infty$-harmonic function $u$, with $|\D u|=e^{-t}$ there. Noting that $\int_{\gamma_t}du=0$ since $du$ is exact, and that $du$ has alternating orientations on $\alpha_i$, we obtain $\sum(-1)^i|\alpha_i|=0$, hence $\int_{\SS^1(4\pi)}h=0$, hence $a_0=b_0=0$.

\begin{figure}[ht]
	\centering
	\includegraphics{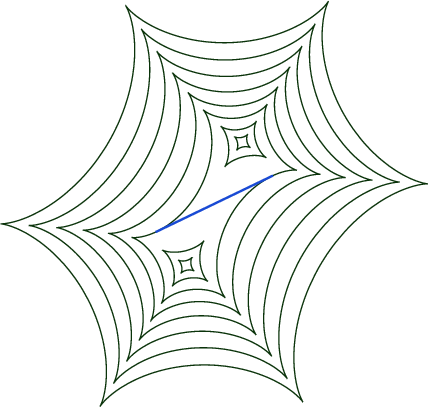}
	\begin{picture}(0,0)
		\put(-136,88){\arrowangle{25.7}}
		\put(-145,94){$\th_0$}
		\put(-88,102){\arrowangle{25.7}}
		\put(-79,95){$\th_0+2\pi$}
	\end{picture}
	\caption{Splitting of 6-gon into two 4-gons.}\label{fig-ex:splitting}
\end{figure}

Now, let $h$ be as in \eqref{eq-ex:entire_3_asymp}. Then $(\p_{\th\th}h+h)(t,\cdot)$ has 6 roots for $t\ll-1$, and $\gamma_t$ has 6 vertices. The solution cannot be constructed by running \eqref{eq-ex:entire_3_asymp} for all times, since otherwise $\gamma_t$ would stop being embedded at a certain $t$. To produce an IMCF cluster, we split the 6-gon into two 4-gons at an appropriate time by cutting along a line segment (see Figure \ref{fig-ex:splitting}).
Then we evolve each resulting 4-gon by IMCF with cusps. Let $T$ be the splitting time, and $\th_0$ be the angle of direction of the dividing segment. This segment clearly must be tangent to $\gamma_T$ at both endpoints. Translating to the language of $h$, this means 
\[h(T,\th_0)=h(T,\th_0+2\pi).\]
Also, each edge of $\gamma_t$ is a streamline of $u$ on which $|\D u|=e^{-t}$. By the same alternating length reasoning as above, for each $t>T$ we must have
\[\int_{\SS^1}h(t,\th)\,d\th=0.\]
Letting $t\searrow T$, we obtain the second constraint
\begin{equation}\label{eq-ex:constraint2}
	\int_{\th_0}^{\th_0+2\pi}h(T,\th)\,d\th=0.
\end{equation}
In the proof, we will see that there exists a unique pair $(T,\th_0)$ with these properties. It is possibly surprising that the alternating length argument only serves to inspire \eqref{eq-ex:entire_3_asymp} and \eqref{eq-ex:constraint2} in our initial discussion. It does not appear in the actual construction below.

\paragraph{Formal construction.} {\ }
	
Let $h$ be as in \eqref{eq-ex:entire_3_asymp}, with $(a,b)\ne(0,0)$ and $(c,d)\ne(0,0)$. Observe $\p_t h=\p_{\th\th}h+h$. It is helpful to see that
\begin{equation}\label{eq-ex:odd_symmetry}
	h(t,\th+2\pi)=-h(t,\th).
\end{equation}
Let $T\in\RR$ and $\th_0\in\SS^1(4\pi)$ be determined by
\begin{equation}\label{eq-ex:split_time}
	h(T,\th_0)=h(T,\th_0+2\pi)=0,\qquad\int_{\th_0}^{\th_0+2\pi}h(T,y)\,dy=0.
\end{equation}
Such $(T,\th_0)$ exists uniquely up to $\th_0\leftrightarrow\th_0+2\pi$: indeed, consider
\[H(t,\th)=\int_\th^{\th+2\pi}h(t,y)\,dy.\]
So \eqref{eq-ex:split_time} is equivalent to $H(T,\th_0)=\p_\th H(T,\th_0)=0$, namely, $(T,\th_0)$ is a degenerate root of $H$. Note that $H(t,\cdot)$ has 2 roots for $t\gg1$ and has 6 roots for $t\ll-1$, due to the specific form of $h$ given by \eqref{eq-ex:entire_3_asymp}. By the Sturmian theory, the number of roots decrease by at least 2 whenever a degenerate root appears. By \eqref{eq-ex:odd_symmetry}, degenerate roots always come in pair $(t,\th)\leftrightarrow(t,\th+2\pi)$. So $H$ has a unique degenerate root (modulo $\th_0\leftrightarrow\th_0+2\pi$) in $\RR\times\SS^1(4\pi)$, as desired.

Let $\SS^1_1$ and $\SS^1_2$ denote two identical copies of $\SS^1$. Identify $\SS^1_1$ as $[\th_0,\th_0+2\pi]$ with the endpoints glued, then identify $\SS^1_2$ as $[\th_0+2\pi,\th_0+4\pi]$ with endpoints glued. Denote
\[M_0=(-\infty,T]\times\SS^1(4\pi),\qquad M_1=[T,\infty)\times\SS^1_1,\qquad M_2=[T,\infty)\times\SS^1_2.\]
Let $M$ be the topological manifold obtained by gluing $M_0$ with $M_1$, $M_2$ in the natural way. The inclusion $M_0\to M$ is not injective since $(T,\th_0)$ and $(T,\th_0+2\pi)$ are identified. However, $h(T,\th_0)=h(T,\th_0+2\pi)=0$ ensures that $h$ is a continuous function on $M\cap\{t\leq T\}$. Then, we evolve $h$ by $\p_t h=\p_{\th\th}h+h$ individually in $M_1,M_2$ from the initial data $h|_{t=T}$. For notational clarity, we denote
\[h_0=h|_{M_0},\qquad h_1=h|_{M_1},\qquad h_2=h|_{M_2}.\]

The second condition of \eqref{eq-ex:split_time} is preserved in $\SS^1_1$ and $\SS^1_2$ respectively, hence
\begin{equation}\label{eq-ex:integral_zero}
	\int_{\SS^1_1}h_1(t,\cdot)=\int_{\SS^1_2}h_2(t,\cdot)=0,\qquad\forall\,t\geq T.
\end{equation}
As $t\to\infty$, the function $h$ has Fourier expansions in $M_1,M_2$ of the form
\begin{equation}\label{eq-ex:asymp_h}
	\begin{aligned}
		& h_1=x_1\sin\th-y_1\cos\th+O\big(e^{-3t}\big)\qquad\text{in}\ \ M_1, \\
		& h_2=x_2\sin\th-y_2\cos\th+O\big(e^{-3t}\big)\qquad\text{in}\ \ M_2,
	\end{aligned}
\end{equation}
for some $x_1,x_2,y_1,y_2\in\RR$. Denote $z_1=(x_1,y_1)$, $z_2=(x_2,y_2)$. Due to the odd symmetry of $h$, we have $z_1=-z_2$. It is shown later that $z_1,z_2\ne0$. Then denote
\[\{z_3,z_4\}=\big\{\!\pm\p_\th h_0(T,\th_0)\tau_{\th_0}\big\},\]
and let $\sigma$ be the closed line segment connecting $z_3,z_4$. We set $z_3=z_4=0$ and $\sigma=\{0\}$ if $\p_\th h_0(T,\th_0)=0$. Finally, define the function
\begin{equation}\label{eq-ex:profile_curves_1}
	F(t,\th)=h(t,\th)\nu_\th+\p_\th h(t,\th)\tau_\th
\end{equation}
on $\rM:=M\setminus\{(T,\th_0)\}$, where we recall
\begin{equation}\label{eq-ex:nu_th}
	\nu_\th=(\sin\th,-\cos\th),\qquad\tau_\th=(\cos\th,\sin\th).
\end{equation}
Note: $h,\nu_\th,\tau_\th$ are continuous in $M$, while $\p_\th h$ is only continuous in $\rM$.

\begin{lemma}\label{lemma-ex:produce_cluster}
	Let $M,\rM,h,F$ be defined as above. Then $z_1,z_2\ne0$ and $z_1,z_2\notin\sigma$, and $F$ is a homeomorphism from $\rM$ to $\RR^2\setminus\big(\{z_1,z_2\}\cup\sigma\big)$. Define the function and vector
	\[w(F(t,\th))=t,\qquad\nu(F(t,\th))=\nu_\th,\]
	and define $w(z_1)=w(z_2)=+\infty$, $w|_\sigma=T$, $\nu|_\sigma=\nu_{\th_0}$. Then $(w,\nu)$ is a unit-calibrated mixed IMCF cluster in $\RR^2$, with exceptional points and segments given by $\{z_1,z_2,\sigma\}$. Consequently, the function $u$ given by $\D u=e^{-w}\nu^\perp$ is $\infty$-harmonic in $\RR^2$.
\end{lemma}

\begin{figure}[ht!]
	\centering
	\includegraphics{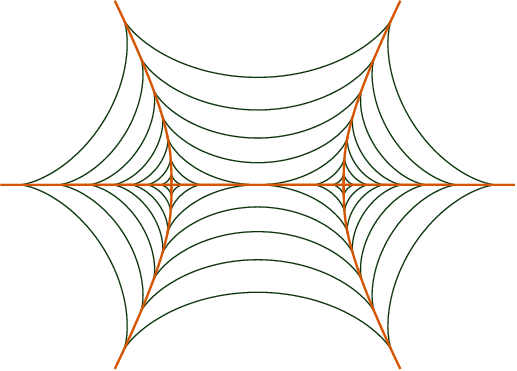}
	\caption{The solution with data $(a,b,c,d)=(1,0,1,0)$.}\label{fig-ex:entire_31}
\end{figure}

\begin{figure}[ht!]
	\centering
	\includegraphics{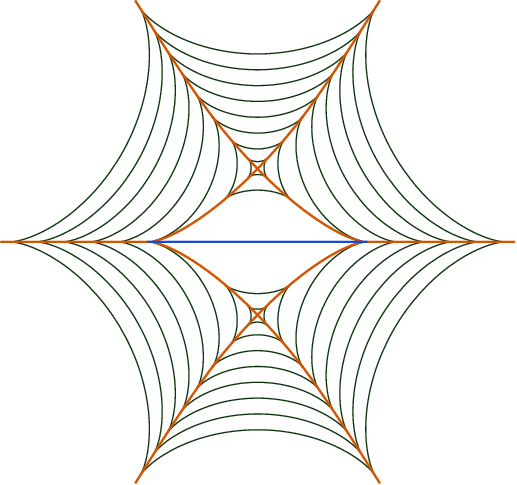}
	\caption{The solution with data $(a,b,c,d)=(1,0,-1,0)$.}\label{fig-ex:entire_33}
\end{figure}

\begin{figure}[ht!]
	\centering
	\includegraphics{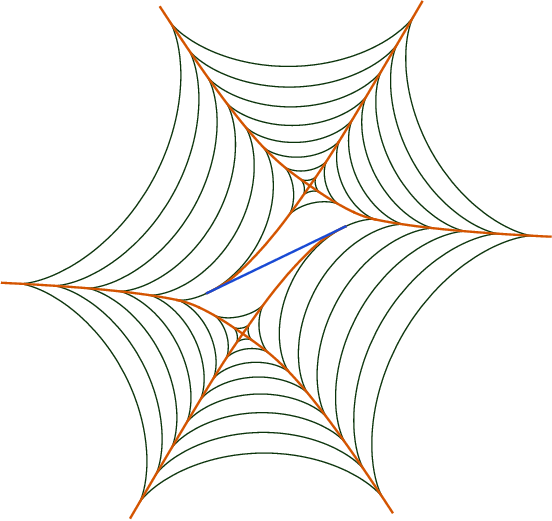}
	\caption{The solution so that $T=0$ and $\th_0=\pi/7$.}\label{fig-ex:entire_32}
\end{figure}

Figure \ref{fig-ex:entire_31}, \ref{fig-ex:entire_33} and \ref{fig-ex:entire_32} respectively describe the solutions with data
\[(a,b,c,d)=(1,0,1,0),\qquad(a,b,c,d)=(1,0,-1,0)\]
and
\[a=1,\quad b=0,\quad c=\frac13\Big(\!\sin\frac{3\pi}{14}-2\cos\frac\pi7\Big),\quad d=-\frac43\cos^2\frac\pi{14}\sin\frac\pi7.\]
In Figure \ref{fig-ex:entire_31}, the segment $\sigma$ degenerates to a point. The parameters for Figure \ref{fig-ex:entire_32} are chosen so that $T=0$ and $\th_0=\pi/7$. Blue segments in the figures are valleys, while red curves are ridges.

\begin{proof}[Proof of Lemma \ref{lemma-ex:produce_cluster}] {\ }
	
	For a point $z=(x,y)\in\RR^2$, we know from \eqref{eq-ex:profile_curves_1} \eqref{eq-ex:nu_th} that $F(t,\th)=z$ is equivalent to
	\[h(t,\th)=x\sin\th-y\cos\th,\qquad \p_\th h(t,\th)=x\cos\th+y\sin\th,\]
	hence equivalent to
	\[(t,\th)\text{\ \ is a degenerate root of\ \ }g_z:=h-x\sin\th+y\cos\th.\]
	
	\vspace{3pt}
	
	\textbf{Step 1.} Injectivity and image of $F$. Let $z=(x,y)\in\RR^2$ be given. We first collect some information about the function $g_z=h-x\sin\th+y\cos\th$. Observe that $g_z(t,\cdot)$ has 6 nondegenerate roots for all $t\ll-1$, and for large positive time we have:
	\begin{itemize}[nosep]
		\item $g_z(t,\cdot)$ has 2 roots in $\SS^1_1$ for all $t\gg1$ if $z\ne z_1$, by the asymptotics \eqref{eq-ex:asymp_h};
		\item $g_{z_1}(t,\cdot)$ has at least 4 roots in $\SS^1_1$ for all $t>T$, by the Sturm-Hurwitz theorem and the fact $\int_{\SS^1}g_{z_1}(t,\th)\,d\th=\int_{\SS^1}g_{z_1}(t,\th)\sin\th\,d\th=\int_{\SS^1}g_{z_1}(t,\th)\cos\th\,d\th=0$.
	\end{itemize}
	Similar fact holds for $\SS^1_2$. Hence, $g_z$ has 4 roots for all $t\gg1$ if $z\notin\{z_1,z_2\}$, and has $\geq6$ roots for all $t\gg1$ if $z=z_1\ne z_2$ or $z=z_2\ne z_1$, and has $\geq8$ roots if $z=z_1=z_2$.
	
	The following lemma keeps track of the change of roots at the saddle point:
	
	\begin{lemma}\label{lemma-ex:sturmian_at_junction}
		Suppose $(T,\th_0)$ is a root of $g_z$. For $\mu,\delta>0$, denote
		\[\begin{aligned}
			P_0 &= [T-\mu,T]\times\Big([\th_0-\delta,\th_0+\delta]\cup[\th_0+2\pi-\delta,\th_0+2\pi+\delta]\Big)\subset M_0, \\
			P_1 &= [T,T+\mu]\times[\th_0-\delta,\th_0+\delta]\subset M_1, \\
			P_2 &= [T,T+\mu]\times[\th_0+2\pi-\delta,\th_0+2\pi+\delta]\subset M_2.
		\end{aligned}\]
		Call $P=P_0\cup P_1\cup P_2$. Let $P_{0,t}$ be the $t$-slice of $P_0$, similarly for $P_{1,t}$, $P_{2,t}$. Denote
		\begin{equation}\label{eq-ex:p_side}
			\begin{aligned}
				\p_\side P &= \Big([T-\mu,T]\times\{\th_0\pm\delta,\th_0+2\pi\pm\delta\}\Big) \\
				&\hspace{60pt} \cup\Big([T,T+\mu]\times\{\th_0\pm\delta\}\Big)\cup\Big([T,T+\mu]\times\{\th_0+2\pi\pm\delta\}\Big)
			\end{aligned}
		\end{equation}
		For any $\epsilon>0$, there are $\mu,\delta<\epsilon$ so that $g_z$ has no root in $\p_\side P$, and the following hold:
		\begin{enumerate}[label={(\roman*)}, nosep]
			\item If $|z|<|\p_\th h(T,\th_0)|$, then $g_z$ has two nondegenerate roots in $P_{0,t}$ for each $t\in[T-\mu,T]$, and $g_z$ has no root in $P_{1,t},P_{2,t}$ for all $t\in(T,T+\mu]$.
			\item If $|z|>|\p_\th h(T,\th_0)|$, then $g_z$ has two nondegenerate roots in $P_{0,t}$ for all $t\in[T-\mu,T]$, and has one nondegenerate root in $P_{1,t},P_{2,t}$ respectively for each $t\in(T,T+\mu]$.
			\item If $|z|=|\p_\th h(T,\th_0)|$, then $g_z$ has $n\geq3$ nondegenerate roots in $P_{0,t}$ for all $t\in[T-\mu,T)$, and has $n'\leq n-2$ nondegenerate roots in $P_{1,t}\cup P_{2,t}$ for all $t\in(T,T+\mu]$.
		\end{enumerate}
	\end{lemma}
	
	The proof is postponed later. Combined with the standard Sturmian theory, it follows that the number of roots of $g_z(t,\cdot)$ is non-increasing in $t$, and whenever a degenerate root in $\rM$ or type (i)(iii) root at $(T,\th_0)$ appears, the number of roots drop by at least 2.
	
	Using this, it quickly follows that $z_1\ne z_2$, hence $z_1,z_2\ne0$ (since $z_2=-z_1$). If $z\in\{z_1,z_2\}$, then $g_z(t,\cdot)$ has 6 roots for all $t\ll-1$ and in total $\geq6$ roots for all $t\gg1$. Therefore, $g_z$ has 6 roots for all $t\in\RR$, and hence no degenerate roots can appear in $\rM$. This shows $z_1,z_2\notin F(\rM)$. If $z\in\sigma$, then $z=c\p_\th h(T,\th_0)\tau_{\th_0}$ for some $c\in[-1,1]$, hence
	\[g_z(T,\th_0)=h(T,\th_0)-\metric{z}{\nu_{\th_0}}=0,\]
	hence $(T,\th_0)$ is a root of $g_z$. Note that $|z|\leq|\p_\th h_0(T,\th_0)|$, thus $(T,\th_0)$ is a type (i)(iii) root as in Lemma \ref{lemma-ex:sturmian_at_junction}, and the number of roots drop by at least 2 at $(T,\th_0)$. This forces $g_z$ to have 4 roots for $t\gg1$, and the saddle point $(T,\th_0)$ occupies the unique decrease of number of roots. Hence $g_z$ has no degenerate root in $\rM$. We conclude that $\sigma\cap F(\rM)=\emptyset$ and $z_1,z_2\notin\sigma$. Denoting
	\[\Omega=\RR^2\setminus\big(\{z_1,z_2\}\cup\sigma\big),\]
	we have shown that $F(\rM)\subset\Omega$.
	
	If $z\notin\{z_1,z_2\}\cup\sigma$, then $g_z$ has 6 roots for all $t\ll-1$ and 4 roots for all $t\gg1$, and $(T,\th_0)$ either is not a root or is a type (ii) root of $g_z$. Therefore, there must be exactly one degenerate root of $g_z$ in $\rM$. This shows that $F$ maps $\rM$ bijectively to $\Omega$.
	
	\vspace{3pt}
	
	\textbf{Step 2.} Continuity of $F^{-1},w,\nu$. By the invariance of domain, note that $F:\rM\to\Omega$ is a homeomorphism. Then,  $(w,\nu)$ is continuous in $\Omega$ since $(w,\nu)=(\id,\th\mapsto\nu_\th)\circ F^{-1}$.
	
	Assume $z\in\sigma$. Then $(T,\th_0)$ is a type (i)(iii) root of $g_z$ as in Lemma \ref{lemma-ex:sturmian_at_junction}, by Step 1 above. For any given $\epsilon$, we use Lemma \ref{lemma-ex:sturmian_at_junction} to obtain a ``box'' $P=P_0\cup P_1\cup P_2$, so that $g_z$ has no roots in $\p_\side P$, and has $n\geq2$ nondegenerate roots in $\p P\cap\{t=T-\mu\}$, and has $n'\leq n-2$ nondegenerate roots in $\p P\cap\{t=T+\mu\}$. By continuity, $g_{z'}$ have the same pattern of roots in $\p P$ for all $|z'-z|<\rho$, for some $\rho=\rho(\epsilon)$. So the number of roots of $g_{z'}$ goes through a decrease in $P$, which shows that there is either a degenerate root in $P\cap\rM$ or a type (i)(iii) root at $(T,\th_0)$. In the former case, if the degenerate root is $(t',\th')$, then $w(z')=t'$ and $\nu(z')=\nu_{\th'}$ and hence
	\[|w(z')-w(z)|=|t'-T|<\epsilon,\qquad |\nu(z')-\nu(z)|\leq|\th'-\th_0|<\epsilon.\]
	In the latter case, we have $z'\in\sigma$ hence $w(z')=w(z)$, $\nu(z')=\nu(z)$ by construction. In either case, this shows the continuity of $w,\nu$ at $z$.
	
	Assume $z=z_1$ (same argument for $z=z_2$). Let $\epsilon>0$ be given. Recall that $g_{z_1}(1/\epsilon,\cdot)$ has 4 nondegenerate roots in $\SS^1_1$ (Sturm-Hurwitz implies that it has $\geq4$, but it can only have 4 roots due to the global behavior of roots). By continuity, there is $\rho>0$ so that $g_{z'}(1/\epsilon,\cdot)$ also has 4 nondegenerate roots in $\SS^1_1$, for all $z'\in B(z_1,\rho)\setminus\{z_1\}$. It follows that the unique degenerate root of $g_{z'}$ must occur in $\SS^1_1\cap\{t>1/\epsilon\}$ for all $z'\in B(z_1,\rho)\setminus\{z_1\}$. This shows $w(z')>1/\epsilon$ for all $z'\in B(z,\rho)$. Hence $w$ is continuous at $z_1$.
	
	To summarize, we have shown that $w$ is a continuous function from $\RR^2$ to $\RR\cup\{+\infty\}$, and $\nu$ is continuous in $\{w<\infty\}$.
	
	\vspace{3pt}
	\textbf{Step 3.} Finally, we show that $(w,\nu)$ produces a simple IMCF cluster in $\Omega$. Once this is done, Definition \ref{def-cluster:mixed} implies that $(w,\nu)$ is a mixed cluster in $\RR^2$ (unit calibration is obvious by our construction). Denote
	\[\kappa^{-1}=\p_{\th\th}h+h=\p_th\qquad\text{in}\ \ \rM.\]
	Here we view $\kappa^{-1}$ only as a notation. We first show that all roots of $\kappa^{-1}$ in $\rM$ are simple. Note that $\kappa^{-1}$ has 6 roots for $t\ll-1$, and by the Sturm-Hurwitz theorem and the fact
	\[\int_{\SS^1}(\p_{\th\th}f+f)\sin\th=\int_{\SS^1}(\p_{\th\th}f+f)\cos\th=0,\qquad\forall\,f\in C^\infty(\SS^1),\]
	$\kappa^{-1}$ has $\geq8$ roots for $t\gg1$. The change of roots at the saddle point is as follows:
	
	\begin{lemma}\label{lemma-ex:curvature_at_junction}
		Adopt the notations in Lemma \ref{lemma-ex:sturmian_at_junction}. Then there exists $\mu,\delta>0$ so that $\kappa^{-1}$ has no roots in $\p_\side P$, and has $n\geq0$ roots in $P_{0,t}$ for each $t\in[T-\mu,T)$, and has $n'\leq n+2$ roots in $P_{1,t}\cup P_{2,t}$ for each $t\in(T,T+\mu]$.
	\end{lemma}
	
	We postpone the proof later. Since $\p_t\kappa^{-1}=\p_{\th\th}\kappa^{-1}+\kappa^{-1}$, the Sturmian principle and this lemma imply that $\kappa^{-1}$ cannot have degenerate roots in $\rM$, as desired.
	
	Let $\Gamma=\big\{(t,\th)\in\rM:\kappa^{-1}(t,\th)=0\big\}$. By the nondegeneracy of roots, $\Gamma$ is a disjoint union of smooth curves $\{\Gamma_j\}$ in $\rM$, with each $\Gamma_j$ of the form $\{\th=g(t)\}$. The implicit function theorem gives
	\[g'(t)=-\frac{\p_t\kappa^{-1}(t,g(t))}{\p_\th\kappa^{-1}(t,g(t))}\qquad\text{on}\ \ \Gamma_j.\]
	The image of $\Gamma_j$ under $F$ is given by
	\[\gamma_j(t)=h(t,g(t))\nu_{g(t)}+\p_\th h(t,g(t))\tau_{g(t)}.\]
	Its velocity vector is
	\[\begin{aligned}
		\gamma'_j(t) &= \p_t h(t,g(t))\nu_{g(t)}+\p_\th h(t,g(t))g'(t)\nu_{g(t)}+h(t,g(t))g'(t)\tau_{g(t)} \\
		&\qquad +\p_{t\th}h(t,g(t))\tau_{g(t)}+\p_{\th\th}h(t,g(t))g'(t)\tau_{g(t)}-\p_\th h(t,g(t))g'(t)\nu_{g(t)} \\
		&= \p_\th\kappa^{-1}(t,g(t))\cdot\tau_{g(t)}\ne0.
	\end{aligned}\]
	Hence, each $\gamma_j$ is an embedded $C^1$ curve in $\Omega$, with $\gamma'_j\perp\nu$ (recall $\nu(F(t,\th))=\nu_\th\perp\tau_\th$).
	
	Let $\{V_i\}$ be the connected components of $\rM\setminus\bigcup\Gamma_j$, and let $D_i=F(V_i)$. Set $\chi_i=1$ if $\kappa^{-1}>0$ in $V_i$, and $\chi_i=-1$ otherwise. Definition \ref{def-cluster:simple}(i)(ii)(iii) and \eqref{eq-cluster:def_calibration} are readily verified by our construction, and it remains to verify the ridge condition \eqref{eq-cluster:def_ridge}, namely $\chi_i\nu=\nu_{D_i}$ on each $\p D_i$.
	
	Consider any point $z\in\p D_i\cap\p D_j$. Let $(t_1,\th_1)\in\rM$ be the preimage of $z$. By symmetry, it suffices to consider the case $t_1=\th_1=0$. So $(0,0)\in\Gamma_k$ for some $k$. Let $\epsilon$ be small enough so that $\{0\}\times(-\epsilon,0)\subset V_i$ and $\{0\}\times(0,\epsilon)\subset V_j$. Consider the Taylor expansion of $h$ at $\th=0$, in which we used $(\p_{\th\th}h+h)(0,0)=0$:
	\[h(0,\th)=a+b\th-\frac12a\th^2+c\th^3+O(\th^4),\qquad \p_\th h(0,\th)=b-a\th+3c\th^2+O(\th^3).\]
	Then
	\begin{equation}\label{eq-ex:taylor}
		\begin{aligned}
			F(0,\th) &= h(0,\th)(\sin\th,-\cos\th)+\p_\th h(0,\th)(\cos\th,\sin\th) \\
			&= \Big(b+\frac12\big(b+6c\big)\th^2+O(\th^3),\ -a+\frac13\big(b+6c\big)\th^3+O(\th^4)\Big).
		\end{aligned}
	\end{equation}
	Observe that
	\begin{equation}\label{eq-ex:sign}
		b+6c=\p_\th\kappa^{-1}(0,0)\left\{\begin{aligned}
			& >0\qquad\text{if $\chi_i=-1$, $\chi_j=1$,} \\
			& <0\qquad\text{if $\chi_i=1$, $\chi_j=-1$.}
		\end{aligned}\right.
	\end{equation}
	Denote $\gamma_\pm=F\big(\{0\}\times\pm(0,\epsilon)\big)$. If $\chi_i=-1$ and $\chi_j=1$, then \eqref{eq-ex:taylor} \eqref{eq-ex:sign} imply that $\gamma_-$ is the graph of $y=-a-C(x-b)^{3/2}+o(|x-b|^{3/2})$ for $0\leq x-b\ll1$ and $\gamma_+$ is the graph of $y=-a+C(x-b)^{3/2}+o(|x-b|^{3/2})$ for $0\leq x-b\ll1$, where
	\[C=\frac{2\sqrt 2}{3\sqrt{b+6c}}>0.\]
	As $\gamma_-\in D_i$ and $\gamma_+\in D_j$, this shows $\nu_{D_i}(z)=\p_y$ and $\nu_{D_j}(z)=-\p_y$. Then recalling $\nu(z)=\nu(F(0,0))=\nu_0=-\p_y$, it follows that $\nu_{D_i}=\chi_i\nu$ and $\nu_{D_j}=\chi_j\nu$ at $z$. For $\chi_i=1$ and $\chi_j=-1$ one may argue in the same manner. This completes the construction.
\end{proof}

\begin{proof}[Proof of Lemma \ref{lemma-ex:sturmian_at_junction}] {\ }
	
	By a rotation in $\RR^2$, we may assume $\th_0=0$. Recall
	\[z=(x,y),\qquad g_z=h-x\sin\th+y\cos\th.\]
	Note $h(T,0)=h(T,2\pi)=0$ by construction. As $(T,0)$ is a root of $g_z$, we have $(x\sin\th-y\cos\th)|_{\th=0}=0$, thus $y=0$, $|z|=|x|$. Notice that $g_z(T,2\pi)=0$ as well. (Note: $h$ is odd and $x\sin\th-y\cos\th$ is even under the symmetry $\th\to\th+2\pi$.) Let us summarize
	\begin{equation}\label{eq-ex:conditions}
		h(T,0)=h(T,2\pi)=g_z(T,0)=g_z(T,2\pi)=0,\qquad y=0,\qquad |z|=|x|.
	\end{equation}
	
	We further denote
	\[P_0^-=[T-\mu,T]\times[-\delta,\delta],\qquad P_0^+=[T-\mu,T]\times[2\pi-\delta,2\pi+\delta].\]
	So $P=P_0^-\cup P_0^+\cup P_1\cup P_2$. Then denote $g_{z,0}=g_z|_{P_0}$, $g_{z,1}=g_z|_{P_1}$, $g_{z,2}=g_z|_{P_2}$. The local gluing pattern for $P$ is shown in Figure \ref{fig-ex:gluing} (edges with the same color are glued).
	
	\begin{figure}[ht]
		\centering
		\includegraphics{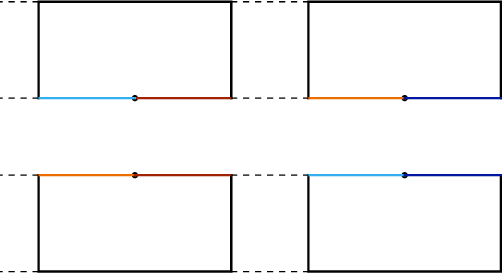}
		\begin{picture}(0,0)
			\put(-187,20){$P_0^-$}
			\put(-57,20){$P_0^+$}
			\put(-187,103){$P_1$}
			\put(-57,103){$P_2$}
			\put(-279,-2){$T-\mu$}
			\put(-257,43){$T$}
			\put(-257,81){$T$}
			\put(-279,127){$T+\mu$}
			\put(-20,51){$2\pi+\delta$}
			\put(-57,51){$2\pi$}
			\put(-113,51){$2\pi-\delta$}
			\put(-138,51){$\delta$}
			\put(-183.5,51){$0$}
			\put(-234,51){$-\delta$}
			\put(-8,72){$\delta$}
			\put(-54,72){$0$}
			\put(-104,72){$-\delta$}
			\put(-138,72){$\delta$}
			\put(-183.5,72){$0$}
			\put(-234,72){$-\delta$}
		\end{picture}
		\caption{Coordinates in $P$.}\label{fig-ex:gluing}
	\end{figure}
	
	Since solutions of heat equation have isolated roots, we may choose $\delta<\epsilon$ so that
	\begin{equation}\label{eq-ex:no_sign_change}
		g_{z,0}\ne0\text{\ \ everywhere in\ \ }\Big\{t=T,\,\th\in\big([-\delta,\delta]\cup[2\pi-\delta,2\pi+\delta]\big)\setminus\{0,2\pi\}\Big\}.
	\end{equation}
	By continuity, there is $\mu<\epsilon$ so that $g_z$ has no roots in $\p_\side P$. Here, $\epsilon$ is the small constant given in Lemma \ref{lemma-ex:sturmian_at_junction}. The standard maximum principle (see Subsection \ref{subsec:sturmian}) implies that $g_{z,1}$ has $n_1\in\{0,1\}$ nondegenerate roots in $P_{1,t}$ for all $t\in(T,T+\mu]$, with $n_1=1$ if $g_{z,1}(T,\pm\delta)$ has the opposite sign, and $n_1=0$ otherwise. Similar fact holds for $g_{z,2}$.
	
	Expand $g_{z,0}$ as Taylor polynomial near $\th=0,2\pi$ using \eqref{eq-ex:conditions} and the symmetry of $h$:
	\begin{align}
		g_{z,0}(T,\th) &= (a-x)\th+b\th^2+O(\th^3), \label{eq-ex:gz_taylor1} \\
		g_{z,0}(T,\th+2\pi) &= (-a-x)\th-b\th^2+O(\th^3), \label{eq-ex:gz_taylor2}
	\end{align}
	where $a=\p_\th h(T,0)=-\p_\th h(T,2\pi)$ and $b=\frac12\p_{\th\th} h(T,0)=-\frac12\p_{\th\th}h(T,2\pi)$. Then $g_{z,1}$, $g_{z,2}$ have the following expansions at $(T,0)$ respectively:
	\begin{equation}\label{eq-ex:taylor_Pt1}
		g_{z,1}(T,\th)=\left\{\begin{aligned}
			& (a-x)\th+b\th^2+O(\th^3)\qquad(\th\geq0), \\
			& (-a-x)\th-b\th^2+O(\th^3)\qquad(\th\leq0).
		\end{aligned}\right.,
	\end{equation}
	and
	\begin{equation}\label{eq-ex:taylor_Pt2}
		g_{z,2}(T,\th)=\left\{\begin{aligned}
			& (-a-x)\th-b\th^2+O(\th^3)\qquad(\th\geq0), \\
			& (a-x)\th+b\th^2+O(\th^3)\qquad(\th\leq0).
		\end{aligned}\right..
	\end{equation}
	The no-sign-change condition \eqref{eq-ex:no_sign_change} implies that the signs of $g_{z,1}(T,\pm\delta)$ and $g_{z,2}(T,\pm\delta)$ coincide with the signs of the leading terms in \eqref{eq-ex:taylor_Pt1} \eqref{eq-ex:taylor_Pt2}.
	
	If $|a|>|x|$, then \eqref{eq-ex:taylor_Pt1} \eqref{eq-ex:taylor_Pt2} implies that $g_{z,1}(T,\pm\delta)$ have the same sign, and $g_{z,2}(T,\pm\delta)$ have the same sign. So $g_z$ has no root in $P_1,P_2$ when $t>T$. If $|a|<|x|$, then \eqref{eq-ex:taylor_Pt1} \eqref{eq-ex:taylor_Pt2} implies that $g_{z,1}(T,\pm\delta)$ have the opposite signs, and $g_{z,2}(T,\pm\delta)$ also have opposite signs, so $g_{z,1},g_{z,2}$ have respectively one root in $P_1,P_2$ when $t>T$. Also, $\th=0,2\pi$ are nondegenerate roots of $g_{z,0}$ if $|a|\ne|x|$. These imply Lemma \ref{lemma-ex:sturmian_at_junction}(i)(ii).
	
	If $|a|=|x|\ne0$, then $g_{z,0}$ has a degenerate root at $\th=0$ (if $a=x$) or at $\th=2\pi$ (if $a=-x$), due to \eqref{eq-ex:gz_taylor1} \eqref{eq-ex:gz_taylor2}. Therefore, $g_z$ has at least 3 roots in $P_{0,t}$ for $t\in[T-\mu,T)$. If $b>0$, then:
	\begin{itemize}[nosep]
		\item If $a=x>0$, then \eqref{eq-ex:taylor_Pt1} gives (for $0<\th\ll1$)
		\[g_{z,1}(T,\th)=b\th^2+O(\th^3)>0,\qquad g_{z,1}(T,-\th)=2a\th+O(\th^2)>0.\]
		\item If $a=-x<0$, then \eqref{eq-ex:taylor_Pt1} gives
		\[g_{z,1}(T,\th)=2a\th+O(\th^2)<0,\qquad g_{z,1}(T,-\th)=-b\th^2+O(\th^3)<0.\]
		\item If $a=x<0$, then \eqref{eq-ex:taylor_Pt2} gives
		\[g_{z,2}(T,\th)=-2a\th+O(\th^2)>0,\qquad g_{z,2}(T,-\th)=b\th^2+O(\th^3)>0.\]
		\item If $a=-x>0$, then \eqref{eq-ex:taylor_Pt2} gives
		\[g_{z,2}(T,\th)=-b\th^2+O(\th^3)<0,\qquad g_{z,2}(T,-\th)=-2a\th+O(\th^2)<0.\]
	\end{itemize}
	In either case, $g_{z,1}(T,\pm\delta)$ or $g_{z,2}(T,\pm\delta)$ have the same sign. Hence $g_z$ has at most one root in $P_{1,t}\cup P_{2,t}$ for all $t\in(T,T+\mu]$. The case $b<0$ is analogous. If $b=0$, then one of the roots $\th=0,2\pi$ has degree $\geq3$, so $g_{z,0}$ has at least 4 roots in $P_0$ for $t<T$. Finally, if $a=x=0$, then both roots of $g_{z,0}$ are degenerate, so $g_z$ has at least 4 roots in $P_0$ when $t<T$. So in either case, the number of roots for $t>T$ decreases by at least 2 compared to $t<T$. This proves Lemma \ref{lemma-ex:sturmian_at_junction}(iii).
\end{proof}

\begin{proof}[Proof of Lemma \ref{lemma-ex:curvature_at_junction}] {\ }
	
	We may assume $\th_0=0$. Similar to the proof of Lemma \ref{lemma-ex:sturmian_at_junction}, we may choose $\mu,\delta$ so that $\p_{\side}P$ contains no roots of $\kappa^{-1}$, and
	\[\kappa^{-1}\ne0\text{\ \ everywhere in\ \ }\Big\{t=T,\,\th\in\big([-\delta,\delta]\cup[2\pi-\delta,2\pi+\delta]\big)\setminus\{0,2\pi\}\Big\}.\]
	Denote $\kappa^{-1}_0=\kappa^{-1}|_{P_0}$, $\kappa^{-1}_1=\kappa^{-1}|_{P_1}$, $\kappa^{-1}_2=\kappa^{-1}|_{P_2}$. Notice that $\kappa^{-1}_0(t,\th+2\pi)=-\kappa^{-1}_0(t,\th)$. Recall $h(T,0)=0$. Denote $a=\p_\th h(T,0)$ and $b=\kappa_0^{-1}(T,0)$. Then the evolution of $\kappa^{-1}_1$, $\kappa^{-1}_2$ in $P_1,P_2$ have the initial data (see Figure \ref{fig-ex:gluing} for the coordinate change in gluing)
	\begin{equation}\label{eq-ex:iv_kappa_1}
		\kappa^{-1}_1(T,\th)=2a\delta_0+\left\{\begin{aligned}
			b+O(\th)\qquad\th>0, \\
			-b+O(\th)\qquad\th<0,
		\end{aligned}\right.
	\end{equation}
	and
	\begin{equation}\label{eq-ex:iv_kappa_2}
		\kappa^{-1}_2(T,\th)=-2a\delta_0+\left\{\begin{aligned}
			-b+O(\th)\qquad\th>0, \\
			b+O(\th)\qquad\th<0,
		\end{aligned}\right.
	\end{equation}
	where $\delta_0$ is the Dirac delta at $\th=0$.
	
	If $b\ne0$, then $\kappa^{-1}_1$, $\kappa^{-1}_2$ has 1 root in $P_{1,t}$, $P_{2,t}$ respectively for all $t\in(T,T+\mu]$. If $b=0$, then $\th=0,2\pi$ are roots of $\kappa_0^{-1}(T,\cdot)$, thus $\kappa^{-1}_0$ has at least 2 roots in $P_{0,t}$ for all $t\in[T-\mu,T)$; moreover, by \eqref{eq-ex:iv_kappa_1} \eqref{eq-ex:iv_kappa_2} we know that $\kappa_1^{-1},\kappa_2^{-1}$ has at most 2 roots in $P_{1,t},P_{2,t}$ respectively for all $t\in(T,T+\mu]$. In either case, the number of roots at most increases by 2, as desired.
\end{proof}

\section{\texorpdfstring{$C^1$}{C1} equicontinuity of \texorpdfstring{$p$}{p}-harmonic functions}\label{sec:equi}

In this section we prove the following:

\begin{theorem}\label{thm-equi:C1_equicont}
	The family
	\[\Big\{v\in\Lip_{\loc}(B(4)): \Delta_p v=0,\ p\in[3,\infty],\ \|v\|_{L^\infty}\leq1\Big\}\]
	is $C^1$ equibounded and equicontinuous in $B(1)$.
\end{theorem}

A major consequence is that the $p$-harmonic approximation, initially set up in $C^0_{\loc}$ in \cite{Bhattacharya-DiBenedetto-Manfredi_1989, Jensen_1993, Lindqvist-Manfredi_1995}, can automatically be upgraded to $C^1_{\loc}$:

\begin{theorem}\label{thm-equi:C1_approx}
	Suppose $u$ is $\infty$-harmonic in a domain $\Omega\subset\RR^2$, and $u_i$ is a sequence of $p_i$-harmonic functions with $p_i\to\infty$ and $u_i\to u$ in $C^0_{\loc}(\Omega)$. Then in fact
	\[u_i\to u\qquad\text{in}\ \ C^1_{\loc}(\Omega).\]
\end{theorem}

The following consequence is also useful:

\begin{cor}\label{cor-equi:min_principle}
	Suppose $u$ is $\infty$-harmonic in $\Omega\subset\RR^2$, and $K\Subset\Omega$, and $u$ has no critical points in $K$. Then $\inf_K|\D u|=\inf_{\p K}|\D u|$.
\end{cor}

The convergence $|\D u_i|\to|\D u|$ was previously shown by Lindgren-Lindqvist \cite[Section 3]{Lindgren-Lindqvist_2021} in the context of $\infty$-harmonic potentials in planar convex rings. Theorem \ref{thm-equi:C1_equicont} answers a question raised in \cite{Lindgren-Lindqvist_2021}. More recently, Peng-Zhang-Zhou \cite{Peng-Zhang-Zhou_2025} proved the full convergence $\D u_p\to\D u$ for $\infty$-harmonic potentials in convex rings in all dimensions.

The proof of Theorem \ref{thm-equi:C1_equicont} combines the following ingredients:
\begin{itemize}[nosep]
	\item integral estimates by Koch-Zhang-Zhou \cite{Koch-Zhang-Zhou_2019};
	\item pointwise flatness estimates by Evans-Smart \cite{Evans-Smart_2011b} (with a technical improvement);
	\item the convergence arguments in \cite[Theorem 1.4]{Koch-Zhang-Zhou_2019} and \cite[Theorem 7]{Lindgren-Lindqvist_2021};
	\item a new observation on the presence of critical points (see Case 2 in the proof below).
\end{itemize}
We need to adapt the estimates in \cite{Evans-Smart_2011b, Koch-Zhang-Zhou_2019} to the context of $p$-harmonic functions (since they were originally carried out for exponentially harmonic functions). We provide details in Appendix \ref{sec:kzz_es} for the reader's convenience; see also the appendix of \cite{Lindgren-Lindqvist_2021}.

\begin{proof}[Proof of Theorem \ref{thm-equi:C1_equicont}] {\ }
	
	Let us prove that the family of solutions of the equation
	\begin{equation}\label{eq-equi:reg_eq}
		\div\Big(\big(\epsilon^2+|\D v|^2\big)^{p/2-1}\D v\Big)=0,\qquad\text{where}\ p\in[3,\infty),\ \epsilon\in(0,1],
	\end{equation}
	with $|v|\leq1$ in $B(4)$, is $C^1$-equicontinuous in $B(1)$. Notice that the equation \eqref{eq-equi:reg_eq} is invariant under $v(z)\mapsto \lambda v(\lambda^{-1}z)$, but not invariant under $v(z)\mapsto \lambda v(z)$ if $\epsilon\ne0$.
	
	The following preliminary estimates are useful; see Appendix \ref{sec:kzz_es} for the proofs.
	
	\begin{lemma}\label{lemma-equi:estimates}
		There is a universal constant $C_1$ so that the following holds.
		
		\begin{enumerate}[label={(\roman*)}, nosep]
			\item If $v$ either solves \eqref{eq-equi:reg_eq} or is $p$-harmonic $(3\leq p\leq\infty)$ in a ball $B(z,4r)$, then
			\begin{equation}\label{eq-equi:grad_est}
				\sup_{B(z,3.5r)}|\D v|\leq C_1r^{-1}\sup_{B(z,4r)}|v|.
			\end{equation}
			\item If $v$ solves \eqref{eq-equi:reg_eq} in $B(z,4r)$, then
			\begin{equation}\label{eq-equi:kzz_sobolev}
				\int_{B(z,3r)}\big|\D|\D v|^2\big|^2\leq C_1r^{-2}\int_{B(z,3.5r)}\big(\epsilon^2+|\D v|^2\big)^2.
			\end{equation}
			\item If $v$ solves \eqref{eq-equi:reg_eq} in $B(z,4r)$, and $P$ is an affine function, then
			\begin{align}
				\dashint_{B(z,r)}\big(|\D v|^2-\metric{\D P}{\D v}\big)^2 &\leq C_1\Big(\dashint_{B(z,2r)}\big(\epsilon^2+|\D v|^2\big)^2\Big)^{\frac12} \label{eq-equi:kzz_flatness}\\
				&\quad\times \Big(\dashint_{B(z,2r)}\frac{(v-P)^4}{r^4}+\frac{(v-P)^2}{r^2}\big(|\D P|+|\D v|\big)^2\Big)^{\frac12}. \nonumber
			\end{align}
			\item If $0<\lambda\leq a$, and $v$ solves \eqref{eq-equi:reg_eq} and satisfies $|v-ay|\leq\lambda r$ in $B(4r)$, then
			\begin{equation}\label{eq-equi:es_flatness}
				|\D v|^2\leq a\p_yv+C_1\big(\epsilon p^{-1/2}+a\big)^{3/2}\lambda^{1/2}\qquad\text{in}\ \ B(r).
			\end{equation}
		\end{enumerate}
	\end{lemma}
	Item (i) implies that the $C^0$-bounded family of solutions of \eqref{eq-equi:reg_eq}, along with the family of $p$-harmonic functions ($3\leq p\leq\infty$), is equicontinuous and $C^1$-equibounded in $B(1)$.
	
	Let us show that the $C^1$ equicontinuity for \eqref{eq-equi:reg_eq} implies the main theorem: if Theorem \ref{thm-equi:C1_equicont} is false, then there exists a $\delta>0$, a sequence of $p_i$-harmonic functions $v_i$, and points $z_i,z'_i\in B(1)$, so that $|z_i-z'_i|\to0$ but $|\D v_i(z_i)-\D v_i(z'_i)|\geq\delta$. Since $\infty$-harmonic functions are uniformly $C^1$ \cite[Theorem 3]{Savin_2005}, we may remove finitely many terms and assume $p_i<\infty$ for all $i$. It is well-known that $p$-harmonic functions can be $C^{1,\alpha(p)}$ approximated by solutions of \eqref{eq-equi:reg_eq} (see e.g. \cite{Lewis_1983}). Hence, for each $i$ there exists an $\epsilon_i\in(0,1]$ and a solution $\tilde v_i$ of \eqref{eq-equi:reg_eq}, so that $\|\tilde v_i-v_i\|_{C^1(B(1))}\leq\delta/4$. Then we have $|\D\tilde v_i(z_i)-\D\tilde v_i(z'_i)|\geq\delta/2$, contradicting the hypothesized $C^1$ equicontinuity for \eqref{eq-equi:reg_eq}.
	
	The rest of the proof is thus devoted to the $C^1$-equicontinuity for solutions of \eqref{eq-equi:reg_eq}. Since solutions of \eqref{eq-equi:reg_eq} are smooth, all the computations will be in the classical sense.
	
	We prove by contradiction: suppose that there is a sequence of solutions $v_i$ of \eqref{eq-equi:reg_eq}, with $p=p_i\in[3,\infty)$ and $\epsilon=\epsilon_i\in(0,1]$, with $\|v_i\|_{L^\infty(B(4))}\leq1$, and there are points $z_i,z'_i\in B(1)$, with $|z_i-z'_i|\to0$, but $|\D v_i(z_i)-\D v_i(z'_i)|\geq\delta$ for some $\delta>0$. By the interior $C^{1,\alpha(p)}$ estimate for \eqref{eq-equi:reg_eq}, see for example \cite{Lewis_1983, Lieberman_1988}, we must have $p_i\to\infty$.
	
	Passing to a subsequence, we may assume $z_i,z'_i\to z_\infty\in\bar{B(1)}$, and $v_i\to u$ in $C^0_{\loc}(B(4))$ for some $u$. We replace each $v_i$ by the new function $z\mapsto 2v_i(\frac12(z+z_\infty))-2u(z_\infty)$. Then, to summarize, we have a sequence $v_i$ of solutions to \eqref{eq-equi:reg_eq} in $B(4)$, with
	\begin{equation}\label{eq-equi:hypo1}
		\|v_i\|_{L^\infty(B(4))}\leq4,\qquad z_i,z'_i\to0,\qquad |\D v_i(z_i)-\D v_i(z'_i)|\geq\delta,
	\end{equation}
	and
	\begin{equation}\label{eq-equi:hypo2}
		v_i\xrightarrow{C^0_{\loc}(B(4))}u,\qquad \|u\|_{L^\infty(B(4))}\leq4,\qquad u(0)=0.
	\end{equation}
	Following the argument in \cite[Theorem 1.22]{Jensen_1993}, we find that $u$ is $\infty$-harmonic in $B(4)$ (this is true even if $\epsilon_i$ does not converge to 0). Recall that planar $\infty$-harmonic functions are $C^1$: so there is a function $\mu:[0,1/2]\to\RR_+$, with $\lim_{r\to0^+}\mu(r)=0$, such that
	\begin{equation}\label{eq-equi:u_is_C1}
		\sup_{z\in B(r)}\big|u(z)-\metric{z}{\D u(0)}\big|\leq r\mu(r),\qquad\sup_{z\in B(r)}\big|\D u(z)-\D u(0)\big|\leq\mu(r),
	\end{equation}
	for all $r\leq1/2$.
	
	Using \eqref{eq-equi:grad_est} \eqref{eq-equi:kzz_sobolev} and passing to a subsequence, we may assume that $|\D v_i|^2$ converges in $L^2(B(3))$ to some function $f$. It is shown in \cite{Koch-Zhang-Zhou_2019} that $f=|\D u|^2$ a.e. in $B(3)$. We recall the proof for the reader's convenience. Let $z_0\in B(3)$ be a Lebesgue point of $f$. Since $u\in C^1$, for all $\lambda<1/100$ there is a small $r\in(0,1/8)$, such that
	\[\sup_{z\in B(z_0,4r)}\big|u(z)-u(z_0)-\metric{\D u(z_0)}{z-z_0}\big|\leq\lambda r.\]
	Hence, for all sufficiently large $i$ (depending on $\lambda,r$), it holds
	\[\sup_{z\in B(z_0,4r)}\big|v_i(z)-v_i(z_0)-\metric{\D u(z_0)}{z-z_0}\big|\leq 2\lambda r.\]
	Using the integral flatness estimate \eqref{eq-equi:kzz_flatness}, for sufficiently large $i$ it follows that
	\[\dashint_{B(z_0,r)}\big(|\D v_i|^2-\metric{\D u(z_0)}{\D v_i}\big)^2\leq C\lambda,\]
	where $C$ is a universal constant. Take $i\to\infty$ of this inequality. Since $|\D v_i|^2\to f$ in $L^2$ and $\D v_i\rightharpoonup\D u$ weakly in $L^2$, we obtain by lower semi-continuity
	\[\dashint_{B(z_0,r)}\big(f-\metric{\D u(z_0)}{\D u}\big)^2\leq C\lambda.\]
	Finally, taking $\lambda\to0$ and $r\to0$, we obtain $f(z_0)=|\D u(z_0)|^2$ as desired. To summarize, we have obtained
	\begin{equation}\label{eq-equi:conv_of_du2}
		|\D v_i|^2\xrightarrow{L^2(B(3))}|\D u|^2.
	\end{equation}
	
	\vspace{3pt}
	
	We are ready to derive a contradiction.
	
	\vspace{3pt}
	\textbf{Case 1:} $\D u(0)=0$. Choose $r\ll1$ so that $C_1\mu(r)<\delta/40$. Notice that $|u|\leq r\mu(r)$ in $B(r)$, hence $|v_i|\leq 2r\mu(r)$ in $B(r)$ for all large enough $i$. Then the gradient estimate \eqref{eq-equi:grad_est} implies
	\[|\D v_i|\leq 8C_1\mu(r)<\delta/5\qquad\text{in}\ \ B(r/4),\]
	for all large $i$, which contradicts our hypothesis \eqref{eq-equi:hypo1}.
	
	\vspace{3pt}
	\textbf{Case 2:} $\D u(0)\ne0$. With a rotation in $\RR^2$, we may assume $\D u(0)=a\p_y$, $a>0$. Choose a small $r$ so that
	\begin{equation}\label{eq-equi:cond_for_r_1}
		\left\{\begin{aligned}
			& 8\mu(4r)\leq\frac{a}{64C_1^2}, \\
			& a^2-\big(a-2\mu(r)\big)^2+4C_1a^{3/2}\big(8\mu(4r)\big)^{1/2}<(\delta/5)^2.
		\end{aligned}\right.\quad 
	\end{equation}
	Since $u\in C^1$, we may further decrease $r$ so that
	\begin{equation}\label{eq-equi:cond_for_r_2}
		|\D u|^2\geq\frac78a^2\qquad\text{in}\ \ B(4r).
	\end{equation}
	We fix our choice of $r$. By \eqref{eq-equi:u_is_C1} and $v_i\to u$ in $C^0_{\loc}(B(4))$, for all large enough $i$ we have
	\begin{equation}\label{eq-equi:flatness_vi}
		\operatorname{sup}_{B(4r)}|v_i-ay|\leq8\mu(4r)\cdot r.
	\end{equation}
	By \eqref{eq-equi:conv_of_du2}, there is a sequence $\eta_i\searrow0$ so that $\int_{B(3)}\big||\D v_i|^2-|\D u|^2\big|^2\leq\eta_i^2$ for all $i$. Denote
	\[F_i=\big|\D|\D v_i|^2\big|^2+\eta_i^{-1}\big||\D v_i|^2-|\D u|^2\big|^2.\]
	Disintegrating \eqref{eq-equi:kzz_sobolev} \eqref{eq-equi:conv_of_du2} on circles and using Fatou's lemma, we have
	\[\begin{aligned}
		\int_{r/3}^{2r/3}\Big(\liminf_{i\to\infty}\int_{\p B(s)}F_i\Big)\,ds
		\leq \liminf_{i\to\infty}\Big(\int_{r/3}^{2r/3}\int_{\p B(s)}F_i\Big)
		\leq \liminf_{i\to\infty}\int_{B(2)}F_i
		< \infty.
	\end{aligned}\]
	Hence, there exists $s\in[r/3,2r/3]$, a constant $C$, and a subsequence, such that
	\[\int_{\p B(s)}\big|\D|\D v_i|^2\big|^2\leq C\ \ \forall\,i,\qquad \int_{\p B(s)}\big||\D v_i|^2-|\D u|^2\big|^2\to0\ \ \text{as}\ i\to\infty.\]
	This implies
	\[|\D v_i|^2\xrightarrow{C^0(\p B(s))}|\D u|^2.\]
	Then by \eqref{eq-equi:cond_for_r_2}, it holds $|\D v_i|^2\geq\frac34a^2$ on $\p B(s)$ for all large $i$. Using the flatness estimate \eqref{eq-equi:es_flatness} in $B(r)$, together with the conditions \eqref{eq-equi:flatness_vi} \eqref{eq-equi:cond_for_r_1}, we obtain for all large $i$
	\[\p_yv_i\geq a^{-1}\Big(\frac34a^2-C_1\cdot2a^{3/2}\cdot\big(8\mu(4r)\big)^{1/2}\Big)\geq\frac a2\qquad\text{on}\ \ \p B(s).\]
	It is well-known that the map $z\mapsto\big(\p_x v_i(z),-\p_y v_i(z)\big)$ is quasiregular \cite[Proof of Theorem 1]{Manfredi_1988}. The maximum principle for quasiregular maps then implies $\p_yv_i\geq a/2$ in $B(s)$, so $v_i$ have no critical points in $B(s)$ (hence in $B(r/3)$ as well), for all large $i$.
	
	By the maximum principle for quasiregular maps, this implies that
	\[\max_L|\D v_i|=\max_{\p L}|\D v_i|,\qquad\min_L|\D v_i|=\min_{\p L}|\D v_i|,\qquad\forall\,\text{compact }L\subset B(r/3).\]
	By Lebesgue's lemma, see \cite[p.10]{Lindgren-Lindqvist_2021}, for all balls $B(z,\rho)\subset B(z,r/6)\subset B(r/3)$ we have
	\[\Big(\!\osc_{B(z,\rho)}|\D v_i|^2\Big)^2\cdot\log\frac {r/6}{\rho}\leq\pi\int_{B(z,r/6)}\big|\D|\D v_i|^2\big|^2\leq C,\]
	where $C$ is a universal constant. Hence, the functions $|\D v_i|^2$ are equicontinuous in $B(r/6)$. Together with \eqref{eq-equi:conv_of_du2}, we have $|\D v_i|^2\to|\D u|^2$ in $C^0(B(r/12))$.
	
	Recall by our condition \eqref{eq-equi:u_is_C1} that $|\D u-a\p_y|\leq\mu(r)$ in $B(r/12)$. Hence, for all large enough $i$ it holds
	\[\big(a-2\mu(r)\big)^2\leq|\D v_i|^2\leq\big(a+2\mu(r)\big)^2\qquad\text{in}\ \ B(r/12).\]
	Moreover, the flatness estimate \eqref{eq-equi:es_flatness} and conditions \eqref{eq-equi:flatness_vi} \eqref{eq-equi:cond_for_r_1} imply
	\[a\p_yv_i\geq |\D v_i|^2-2C_1a^{3/2}\big(8\mu(4r)\big)^{1/2}\qquad\text{in}\ \ B(r/12),\]
	for all large $i$. Hence, we obtain
	\[\begin{aligned}
		|\D v_i-a\p_y|^2 &= |\D v_i|^2+a^2-2a\p_yv_i \\
		&\leq a^2-\big(a-2\mu(r)\big)^2+4C_1a^{3/2}\big(8\mu(4r)\big)^{1/2}
	\end{aligned}\]
	in $B(r/12)$, for all large $i$. By our assumption \eqref{eq-equi:cond_for_r_1}, this implies
	\[|\D v_i-a\p_y|\leq\delta/5\qquad\text{in}\ \ B(r/12),\]
	for all large $i$, which contradicts \eqref{eq-equi:hypo1}.
\end{proof}

\begin{proof}[Proof of Theorem \ref{thm-equi:C1_approx}] {\ }
	
	It is generally known that if $u_i\to u$ in $C^0_{\loc}(\Omega)$ and $u_i$ is locally $C^1$-equibounded and equicontinuous in $\Omega$, then $u_i\to u$ in $C^1_{\loc}(\Omega)$.
\end{proof}

\begin{proof}[Proof of Corollary \ref{cor-equi:min_principle}] {\ }
	
	We may assume that $u$ has no critical points in $\bar K$. Hence there is $U\Supset K$ and $c>0$ so that
	\[|\D u|\geq c\qquad\text{in}\ \ \bar U.\]
	By Lemma \ref{lemma-prelim:p-approx} and Theorem \ref{thm-equi:C1_approx}, we may consider a sequence of $p$-harmonic functions $u_p$ that $C^1$-converges to $u$ (as $p\to\infty$) in a neighborhood of $\bar U$. Hence $|\D u_p|\geq c/2$ in $\bar U$ for all large enough $p$. Recall that $z\mapsto\frac{\p u_p}{\p z}$ is a quasiregular map. Hence, the map
	\[z\mapsto\frac1{\p u_p/\p z}\]
	is also quasiregular in $U$. Then the maximum principle implies
	\[\min_{\bar K}|\D u_p|=\min_{\p K}|\D u_p|.\]
	Taking $p\to\infty$, the result follows.
\end{proof}

%\newpage

\section{Limit of \texorpdfstring{$p$-$q$}{p-q} duality and simple clusters}\label{sec:approx}

Let $u$ be an $\infty$-harmonic function in a simply-connected domain $\Omega\subset\RR^2$. Denote
\[\Omega_0=\Omega\setminus\Crit(u).\]
Let $z_0\in\Omega$ be a basepoint, and $u_i$ be any sequence of $p_i$-harmonic functions, so that $p_i\to\infty$ and $u_i\to u$ in $C^0_{\loc}(\Omega)$. For the ease of notation, let us suppress the index $i$ and write ``$u_p$ are $p$-harmonic functions so that $p\to\infty$ and $u_p\to u$''. By Theorem \ref{thm-equi:C1_approx}, we in fact have $u_p\to u$ in $C^1_{\loc}(\Omega)$. Let $v_p$ be the conjugate $q$-harmonic function (where $q=\frac{p}{p-1}$), defined by
\begin{equation}\label{eq-approx:def_vp}
	\D v_p=|\D u_p|^{p-2}\D^\perp u_p,\qquad v_p(z_0)=0.
\end{equation}
It is easy to verify that $|\D u_p|^{p-2}\D^\perp u_p$ is indeed a gradient, and $v_p$ is indeed $q$-harmonic (see also \cite{Aronsson-Lindqvist_1988}). Following the idea of Moser \cite{Moser_2022}, we wish to take the renormalization
\[w_p=(1-q)\log v_p=-\frac{\log v_p}{p-1},\]
and take a subsequential limit as $p\to\infty$. The limit is a weak IMCF whenever it is defined, as was seen in \cite{Moser_2007, Moser_2008} (see also \eqref{eq-prelim:p-imcf}\,--\,\eqref{eq-prelim:kotschwar_ni}).

A major ``issue'' is that $v_p$ need not have a definite sign in $\Omega$ (in fact, $v_p$ must change sign since we have intentionally set $v_p(z_0)=0$). The actual limiting process goes roughly as follows. We first take a $C^1_{\loc}$ subsequential limit of $\{v_p=0\}$. These curves correspond to $\{\gamma_i\}$ in the main theorem below. The curves divide $\Omega$ into regions. In each region $D$, we individually take the renormalized limit
\begin{equation}\label{eq-approx:def_wp_aux}
	w_p=-\frac{\log|v_p|}{p-1},\qquad w=\lim_{p\to\infty}w_p\ \ \text{in}\ \ C^0_{\loc}(D).
\end{equation}
Finally, the limit solutions glue continuously across the interfaces, thus giving a piecewise IMCF. It will be shown that the interfaces are ridges. Namely, we obtain a simple IMCF cluster (Definition \ref{def-cluster:simple}). Note that the limit \eqref{eq-approx:def_wp_aux} is taken away from $\p D$: the Dirichlet condition $v_p=0$ $\leftrightarrow$ $w_p=+\infty$ on $\{v_p=0\}$ does not pass to the $p\to\infty$ limit. The function $w$ is bounded (in fact, is equal to $-\log|\D u|$) on all ridges $\gamma_i$.

Another layer of complexity arises from the presence of critical points: the set $\{v_p=0\}$ may not behave nicely in $\Crit(u)$. We will instead take a subsequential limit of $\{v_p=0\}\cap\Omega_0$ in $C^1_{\loc}(\Omega_0)$, and use (the closure of) the union of limit curves to divide $\Omega$ into regions. There are two types of regions: the regular ones, in which \eqref{eq-approx:def_wp_aux} works well, and the degenerate ones, in which $w_p\to+\infty$. The degenerate regions must all lie inside $\Crit(u)$, and the simple cluster will be objects in the regular regions.

Let us call $w$ the resulting simple cluster. The natural question is the relation between $w$ and the original $\infty$-harmonic function $u$. The key ingredient to connecting them is the following observation (inspired by similar estimates in \cite{Moser_2022}, see Lemma \ref{lemma-approx:vp_bounds}). Suppose $\sigma$ is a $C^1$ path joining $z_0$ and some $z\in\Omega$. Then
\[v_p(z)=v_p(z)-v_p(z_0)=\int_\sigma\metric{\D v_p}{\sigma'}=\int_\sigma|\D u_p|^{p-2}\metric{\D^\perp u_p}{\sigma'}.\]
If $\metric{\D^\perp u}{\sigma'}$ is uniformly positive or negative on $\sigma$, then the $1/(p-1)$-th power of this integral is roughly the $L^p$ norm of $|\D u_p|$ on $\sigma$, hence converges to the $L^\infty$ norm as $p\to\infty$. So we have
\[\lim_{p\to\infty}|v_p(z)|^{\frac1{p-1}}=\max_\sigma|\D u|,\]
which in view of \eqref{eq-approx:def_wp_aux} gives
\[w(z)=-\log\max_\sigma|\D u|.\]
Hence $w(z)=-\log|\D u|(z)$ whenever there exists a $\sigma$ for which $|\D u|(z)=\max_\sigma|\D u|$ and $\metric{\D^\perp u}{\sigma'}$ has a definite sign. This condition is highly restrictive (so $w\ne-\log|\D u|$ in general). However, see items \ref{item-approx:w_leq_-logdu}, \ref{item-approx:w_on_z0}\,$\sim$\,\ref{item-approx:w_on_path} and \ref{item-approx:lvlset} below, there are several important cases where $w=-\log|\D u|$.

The resulting cluster $w$ is generally nonunique. It depends on the choice of basepoint $z_0$ and also on the selected subsequence. By Theorem \ref{thm-approx:global}\ref{item-approx:w_on_nonconst}, the set $\{w\ne-\log|\D u|\}$ is made of jumps in the weak IMCFs. Therefore, $w$ may be viewed as a reduction of $u$ by showing partial information and erasing the remaining and converting into jumps. The full data on $u$ may be recovered by considering the approximations at all basepoints $z_0$. To prove the integral regularity conjectures in \eqref{eq-intro:integral_reg_conj}, one perhaps needs to investigate how the different IMCF clusters are patched together.

\begin{theorem}\label{thm-approx:global}
	Let $u$ be $\infty$-harmonic in a simply-connected domain $\Omega\subset\RR^2$. Let $z_0\in\Omega$, and $u_p$ be a sequence of $p$-harmonic functions, with $p\to\infty$ and $u_p\to u$ in $C^0_{\loc}(\Omega)$. Define $v_p$ by
	\begin{equation}
		\D v_p=|\D u_p|^{p-2}\D^\perp u_p,\qquad v_p(z_0)=0.
	\end{equation}
	Set $\Omega_0=\Omega\setminus\Crit(u)$. Then there is a continuous function $w:\Omega\to\RR\cup\{+\infty\}$ and a countable disjoint collection of embedded $C^1$ curves $\{\gamma_i\subset\Omega_0\}$, so that the following hold.
	\begin{enumerate}[label={(\roman*)}, nosep]
		\item\label{item-approx:gamma_i} Each $\gamma_i$ is a streamline of $u$, and has no endpoint in $\Omega_0$. Each compact set in $\Omega_0$ intersects finitely many $\gamma_i$.
		\item\label{item-approx:Dreg_Dcrit} Let $\cD$ be the set of connected components of $\Omega\setminus\bar{\cup\gamma_i}$. Denote
		\[\cD_{\Crit}=\Big\{D\in\cD: D\subset\Crit(u)\Big\},\qquad\cD_{\reg}=\cD\setminus\cD_{\Crit}.\]
		Then by passing to a subsequence of $p$, it holds:
		\begin{enumerate}[label={(\alph*)}, nosep]
			\item If $D\in\cD_{\Crit}$, then $w\equiv+\infty$ in $D$ and $|v_p|^{\frac1{p-1}}\to0$ in $C^0_{\loc}(D)$.
			\item If $D\in\cD_{\reg}$, then $w<+\infty$ everywhere in $D$, and
			\begin{equation}\label{eq-approx:main_convergence}
				-\frac{\log|v_p|}{p-1}\xrightarrow{p\to\infty} w\qquad\text{in}\ \ C^0_{\loc}(D).
			\end{equation}
		\end{enumerate}
		\item\label{item-approx:w_leq_-logdu} $w\leq-\log|\D u|$ everywhere in $\Omega$;
		\item\label{item-approx:Omega_reg} $\{w<+\infty\}$ coincides with the set
		\begin{equation}\label{eq-approx:main_omega_reg}
			\Omega_{\reg}:=\Big(\bigcup_{D\in\cD_{\reg}}D\Big)\cup\big(\!\cup\gamma_i\big).
		\end{equation}
		We have $\Omega_{\reg}\supset\Omega_0$, i.e. $\Omega\setminus\Omega_{\reg}\subset\Crit(u)$. Each $D\in\cD_{\reg}$ is a $C^1$ domain in $\Omega_{\reg}$.
		\item\label{item-approx:simple_cluster} Set $\nu=-e^w\D^\perp u$. For each $D\in\cD_{\reg}$, set $\chi_D$ to be the (well-defined) sign of $v_p(z_D)$ for all large enough $p$, for any fixed $z_D\in D$. Then the data
		\[\big(w,\nu,\cD_{\reg},\{\chi_D\}_{D\in\cD_{\reg}}\big)\]
		is a simple IMCF cluster in $\Omega_{\reg}$.
	\end{enumerate}
	Furthermore, $w$ has the following properties:
	\begin{enumerate}[label={(\roman*)}, nosep]
		\setcounter{enumi}{5}
		\item\label{item-approx:w_on_z0} $w=-\log|\D u|$ at the basepoint $z_0$;
		\item\label{item-approx:w_on_gamma} $w=-\log|\D u|$ on all $\gamma_i$;
		\item\label{item-approx:w_on_nonconst} $w=-\log|\D u|$ on the set
		\[\{w\ne\const\}:=\Big\{z\in\Omega:\text{ $w$ is non-constant in all neighborhoods of $z$}\Big\}.\]
		Consequently, $w=-\log|\D u|$ on $\p\{w<t\}\cap\Omega$ and $\p\{w\leq t\}\cap\Omega$, $\forall\,t\in\RR$.
		\item\label{item-approx:w_on_path} Suppose $z_1\in\cup\gamma_i$ or $z_1=z_0$, and $\sigma:[0,l]\to\Omega$ is a $C^1$ curve joining $z_1$ and another point $z_2\in\Omega$. Then
		\begin{equation}\label{eq-approx:main_w_on_path_1}
			w(z_2)\geq-\log\max_\sigma|\D u|.
		\end{equation}
		Denote
		\[M_\sigma=\Big\{t\in[0,l]:|\D u|(\sigma(t))=\max_\sigma|\D u|\Big\}.\]
		If there is $\chi\in\{\pm1\}$ so that $\chi\metric{\D^\perp u}{\sigma'}>0$ on $M_\sigma$, then $z_2$ lies in some $D\in\cD_{\reg}$, and
		\begin{equation}\label{eq-approx:main_w_on_path_2}
			w(z_2)=-\log\max_\sigma|\D u|,\qquad\chi=\chi_D.
		\end{equation}
		\item\label{item-approx:z0_on_gamma} If $z_0\notin\Crit(u)$, then $z_0\in\cup\gamma_i$.
		\item\label{item-approx:z0_in_crit} If $z_0\in\Int(\Crit(u))$, then $z_0\in D$ for some $D\in\cD_{\Crit}$.
	\end{enumerate}
	Regarding the sublevel sets of $|\D u|$, we have:
	\begin{enumerate}[label={(\roman*)}, nosep]
		\setcounter{enumi}{11}
		\item\label{item-approx:lvlset} Denote
		\begin{equation}\label{eq-approx:lvlset_notations}
			\begin{aligned}
				& Y_t=\big\{|\D u|<e^{-t}\big\},\qquad\tZ_t=\{|\D u|\leq e^{-t}\},\qquad Z_t=\Int(\tZ_t),\\
				&\hspace{24pt} Y'_t=\{w>t\},\qquad\tZ'_t=\{w\geq t\},\qquad Z'_t=\Int(\tZ'_t).
			\end{aligned}
		\end{equation}
		Let $F$ and $F'$ be the connected component of $Y_t$ and $Y'_t$ (or $Z_t$ and $Z'_t$, or $\tZ_t$ and $\tZ'_t$) that contains $z_0$, respectively. Then $F=F'$.
	\end{enumerate}
\end{theorem}

\begin{remark}\label{rmk-approx:misc}
	The following facts are obvious but may help clarify the statements.
	\begin{enumerate}[label={(\roman*)}, nosep]
		\item The basepoint $z_0$ may lie inside or outside $\Crit(u)$.
		\item When $\Crit(u)=\emptyset$, we have $\Omega=\Omega_{\reg}$ hence obtain a simple IMCF cluster in $\Omega$. In this case, the proof of Theorem \ref{thm-approx:global} can be significantly shortened.
		\item A regular region may contain critical points of $u$, and if this occurs, then the critical points must lie in $\Int\big(\{w=t\}\big)$ for some $t\in\RR$. This follows from item \ref{item-approx:w_on_nonconst}.
		\item Regarding \ref{item-approx:lvlset}: since $w\leq-\log|\D u|$, we clearly have $Y_t\supset Y'_t$, $\tZ_t\supset\tZ'_t$ and $Z_t\supset Z'_t$. Note that $F,F'$ might be empty.
		\item Eventually, Theorem \ref{thm-crit:main} implies that $\Int(\Crit(u))=\emptyset$, hence $\cD_{\Crit}=\emptyset$. However, Theorem \ref{thm-approx:global} is a step towards proving Theorem \ref{thm-crit:main}, and the proof seems not substantially simplifiable even assuming a priori $\Int(\Crit(u))=\emptyset$.
	\end{enumerate}
\end{remark}

It is useful to investigate Theorem \ref{thm-approx:global} in the following ``$\epsilon$-regularity'' case which appears in Section \ref{sec:ereg}. Denote $Q_\delta(r)=(-r,r)\times(-\delta r,\delta r)$. Suppose $\delta\leq1$, and $u$ is $\infty$-harmonic in $Q_\delta(4)$, in which
\begin{equation}\label{eq-approx:aux16}
	\Big|\frac{\D u}{|\D u|}-\p_x\Big|<e^{-10}\delta.
\end{equation}
In particular, $\Crit(u)=\emptyset$. Consider $z_0\in Q_\delta(3)$. Notice that
\[\Big|\frac{\D v_p}{|\D v_p|}-\p_y\Big|<e^{-10}\delta\qquad\text{in}\ \ Q_\delta(4-\epsilon),\qquad\forall\,p\geq p(\epsilon).\]
Hence $\{v_p=0\}\cap Q_\delta(4-\epsilon)$ consists of one almost horizontal curve passing through $z_0$, for all large $p$. As $\{\gamma_i\}$ arises from taking subsequential limits of $\{v_p=0\}\setminus\Crit(u)$, in the present scenario there is only one almost horizontal ridge, which we name by $\gamma$ (recall that ridges are streamlines of $u$, hence almost horizontal, by \eqref{eq-approx:aux16}). Let $f:(-4,4)\to(-4\delta,4\delta)$ be a $e^{-9}\delta$-Lipschitz function so that $\gamma=\graphh(f)$. So $\gamma$ divides $Q_\delta(4)$ into two regions
\[D_1=\{y<f(x)\},\qquad D_2=\{y>f(x)\}.\]
As $\D v_p=|\D u_p|^{p-2}\D^\perp u_p$ points upward, the choice of sign in Theorem \ref{thm-approx:global}\ref{item-approx:simple_cluster} gives $\chi_1=-1$, $\chi_2=1$, and $\nu=-e^w\D^\perp u$ points downward. Hence, the resulting simple cluster in \ref{item-approx:simple_cluster} contains an upward moving IMCF in $D_1$ and a downward moving IMCF in $D_2$. 

Item \ref{item-approx:w_on_path} gives a precise formula for $w$ as follows. For each $z=(x,y)\in Q_\delta(4)$, consider the vertical line segment $\sigma$ connecting $z$ to $(x,f(x))\in\gamma$. The sign condition in \ref{item-approx:w_on_path} clearly holds, since $\sigma$ is vertical and $\D^\perp u$ is almost vertical. Thus, we obtain the formula
\begin{equation}\label{eq-approx:ep_reg_w}
	w(x,y)=-\log\max\Big\{|\D u|(x,y'): y'\in[y,f(x)]\Big\}.
\end{equation}

\begin{example}\label{ex-approx:moser}
	Let $\delta=1$, and suppose $u$ is $\infty$-harmonic in $Q(4)$ so that
	\begin{equation}\label{eq-approx:moser_assump_modified}
		|\D u-\p_x|<e^{-20},\qquad |\D u|(x,y)\text{ is nondecreasing in $y$ for each $x$.}
	\end{equation}
	Take the basepoint $z_0=(0,-2)$. Then we have $\gamma\cap Q(1)=\emptyset$ and $Q(1)\subset D_2$, and hence $w$ is a weak IMCF in $Q(1)$. Combining \eqref{eq-approx:ep_reg_w} \eqref{eq-approx:moser_assump_modified}, we obtain
	\[w=-\log|\D u|,\quad|\nu|=1\qquad\text{in}\ \ Q(1).\]
	In summary, $-\log|\D u|$ is a weak IMCF in $Q(1)$ calibrated by $\nu=-\frac{\D^\perp u}{|\D u|}$, assuming \eqref{eq-approx:moser_assump_modified}.
	
	The same conclusion holds under a slightly different assumption
	\begin{equation}\label{eq-approx:moser_assump}
		|\D u|\text{ is nondecreasing in all upward integral lines of $\D^\perp u$}.
	\end{equation}
	Indeed, it suffices to apply \eqref{eq-approx:main_w_on_path_2} on upward integral lines of $\D^\perp u$ connecting $\gamma$ to each $z\in Q(1)$. The condition \eqref{eq-approx:moser_assump} is the local version of R. Moser's orientedness assumption in \cite{Moser_2022}, and our conclusion recovers \cite[Proposition 5]{Moser_2022}.
\end{example}

\begin{remark}\label{rmk-intro:why_imcf}
	We note a formal calculation: suppose $u\in C^\infty$ is $\infty$-harmonic in a domain $\Omega$, so that $|\D u|,|\D|\D u||$ are both nonvanishing. Then there is $\chi\in\{\pm1\}$ with
	\[\frac{\D|\D u|}{|\D|\D u||}=\chi\frac{\D^\perp u}{|\D u|}.\]
	Hence
	\[\begin{aligned}
			\div\Big(\frac{-\D\log|\D u|}{|\D\log|\D u||}\Big) &= -\div\Big(\frac{\D|\D u|}{|\D|\D u||}\Big)=-\chi\div\Big(\frac{\D^\perp u}{|\D u|}\Big) \\
			&= -\chi\frac{\div(\D^\perp u)}{|\D u|}+\chi\frac{\metric{\D|\D u|}{\D^\perp u}}{|\D u|^2}
			= \frac{|\D|\D u||}{|\D u|}=|\D\log|\D u||.
		\end{aligned}\]
	Hence $-\log|\D u|$ is a smooth IMCF. This calculation also appeared in Drucker-Williams \cite[Theorem 2.1]{Drucker-Williams_2009} and Koch-Zhang-Zhou \cite[Proposition 4.2]{Koch-Zhang-Zhou_2019} in slightly different forms.
\end{remark}

Suppose we have a unit-calibrated mixed IMCF cluster $(w_0,\nu_0)$ in $Q(4)$, so that the calibration $\nu_0$ satisfies $|\nu_0+\p_y|<e^{-20}$. Consider first using Theorem \ref{thm-ex:reconstruction} to obtain an $\infty$-harmonic function $U$, then using Theorem \ref{thm-approx:global} to obtain a simple IMCF cluster $w$. In general, we do not have $w=w_0$ even if $w_0$ is simple. The change can be described precisely by \eqref{eq-approx:ep_reg_w}: recall that $\D u=e^{-w_0}\nu_0^\perp$, and hence $|\D u|=e^{-w_0}$, then \eqref{eq-approx:ep_reg_w} implies
\begin{equation}\label{eq-approx:ep_reg_w_2}
	w(x,y)=\min\Big\{w_0(x,y'):y'\in[y,f(x)]\Big\}.
\end{equation}
We include several schematic examples.

\begin{example}\label{ex-approx:imcf}
	See Figure \ref{fig-approx:imcf} (left): consider a unit-calibrated weak IMCF in $Q(4)$ with a jump at $t=0$ represented by the grey region. The line foliation in $\{w_0=0\}$ are displayed as the blue segments. Then the above process (with basepoint $z_0$) produces the simple cluster in Figure \ref{fig-approx:imcf} (right): the ridge $\gamma$ is the integral curve of $\nu^\perp$ that passes through $z_0$ (in this picture it is the unique allowable curve, but in general $\gamma$ is not unique), and the function $w$ satisfies $w|_{D_1}\equiv0$ and $w|_{D_2}\equiv w_0$.
	\begin{figure}[ht!]
		\centering
		\includegraphics{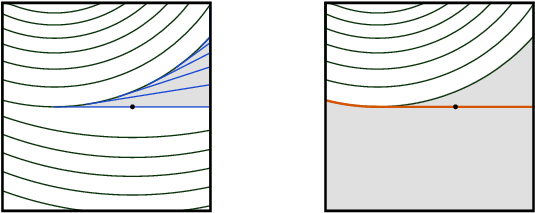}
		\begin{picture}(0,0)
			\put(-203,41){$z_0$}
			\put(-48,41){$z_0$}
			\put(-150,48){$\xrightarrow{\hspace{24pt}}$}
		\end{picture}
		\caption{A weak IMCF and its reduction to a simple cluster.}\label{fig-approx:imcf}
	\end{figure}
\end{example}

\begin{example}
	See Figure \ref{fig-approx:mixed}: consider a mixed IMCF cluster in $Q(4)$. The blue segment represents a valley. The above process with basepoint $z_0$ produces many possible simple clusters including the three examples shown on the right of the figure. In each possible reduction, a part of the original solution is erased and becomes a jump region.
	\begin{figure}[ht!]
		\centering
		\includegraphics{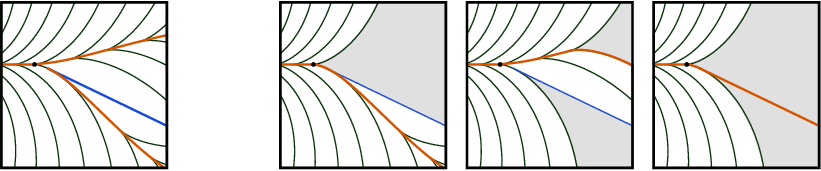}
		\begin{picture}(0,0)
			\put(-386,56){$z_0$}
			\put(-305,38){$\xrightarrow{\hspace{18pt}}$}
		\end{picture}
		\caption{A nonsimple IMCF cluster and its reductions to simple clusters.}\label{fig-approx:mixed}
	\end{figure}
\end{example}

\begin{example}
	Figure \ref{fig-approx:many_valleys} shows a mixed cluster. The valleys divide the cluster into pieces of simple clusters. The result of $p$-harmonic approximation in this picture is to select the ``simple component'' of the cluster that contains the basepoint (the components not containing the basepoint become jumps, shown as the grey regions).
	\begin{figure}[ht!]
		\centering
		\includegraphics{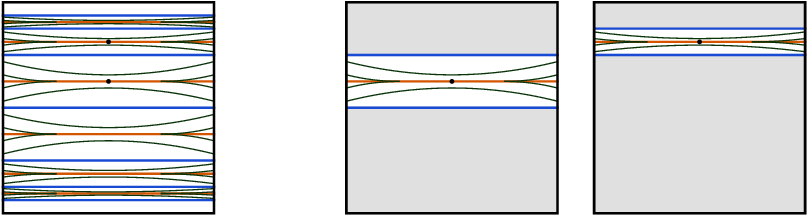}
		\begin{picture}(0,0)
			\put(-350,57){$z_1$}
			\put(-350,75){$z_2$}
			\put(-272,48.5){$\xrightarrow{\hspace{18pt}}$}
		\end{picture}
		\caption{A mixed cluster with many edges, and its reductions at different basepoints.}\label{fig-approx:many_valleys}
	\end{figure}
\end{example}

\begin{remark}\label{rmk-approx:du_noninj}
	If the mixed cluster contains valleys or jumps in the weak IMCF, Then $\D u$ (where $u$ is the resulting $\infty$-harmonic function) is non-injective. Hence, the map $z\mapsto\D u(z)$ is generically non-injective. See also the recent work \cite{Brustad_2026}, as well as Section \ref{sec:ereg} for discussion on jumps.
\end{remark}

\begin{example}
	Figure \ref{fig-approx:entire32} shows the mixed cluster in Figure \ref{fig-ex:entire_32} and possible reductions to simple clusters. Light grey regions represent jumps in the weak IMCF with outer obstacle. In the last cluster, notice that the region $\Omega_{\reg}$ does not contain the degree 2 vertex.
	\begin{figure}[ht!]
		\centering
		\includegraphics{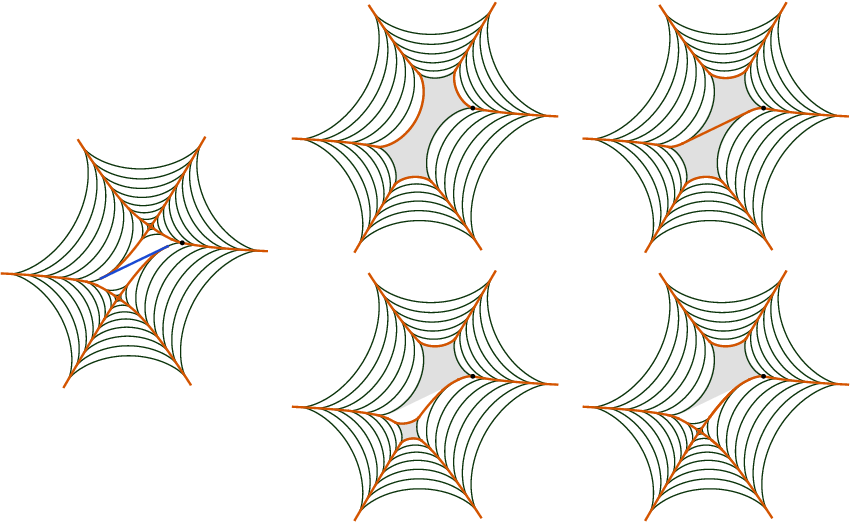}
		\begin{picture}(0,0)
			\put(-312,99){$z_0$}
			\put(-329,115){\rotatebox[origin=c]{120}{$\xrightarrow{\hspace{24pt}}$}}
		\end{picture}
		\caption{Possible clusters arising from Figure \ref{fig-ex:entire_32}.}\label{fig-approx:entire32}
	\end{figure}
\end{example}

Finally, we notice a possibly interesting phenomenon. Let us call a curve $\sigma$ transversal if $\metric{\D^\perp u}{\sigma'}$ has a definite sign on $\sigma$. By item \ref{item-approx:w_on_path}, a transversal curve cannot leave a region $D\in\cD_{\reg}$ once it enters $D$ through some curve $\gamma_i$. In particular, a regular level set that enters some $D\in\cD_{\reg}$ remains in $D$ for all time.

As a note to readers, we stated the strongest version of Theorem \ref{thm-approx:global}\ref{item-approx:w_on_path} though it is never used in full generality in this paper. For all places where it is used, $\metric{\D^\perp u}{\sigma'}$ in fact has a sign on the entire $\sigma$ (not only on the maximal gradient set).

\vspace{3pt}

The proof of Theorem \ref{thm-approx:global} is outlined as follows: we first construct the curves $\gamma_i$, then prove some useful preliminary bounds for $v_p$, then study the sets of regions $\cD_{\Crit}$ and $\cD_{\reg}$, then perform the limiting argument $-\log|v_p|/(p-1)\to w$ in each regular region, then prove item \ref{item-approx:w_on_path}, then extend $w$ continuously to all $\Omega$, and verify the remaining statements.

\begin{proof}[Proof of Theorem \ref{thm-approx:global}] {\ }
	
	We may assume $p\geq3$ in the sequence. By Theorem \ref{thm-equi:C1_equicont}, \ref{thm-equi:C1_approx}, we have
	\begin{equation}\label{eq-approx:equicont}
		\{u_p\},\{\D u_p\}\text{ are equicontinuous in $L$, $\forall\,L\Subset\Omega$},
	\end{equation}
	and
	\begin{equation}\label{eq-approx:C1}
		u_p\xrightarrow{C^1_{\loc}(\Omega)}u.
	\end{equation}
	Recall $\Omega_0=\Omega\setminus\Crit(u)$. Hence, for each $K\Subset\Omega_0$ there is a constant $c(K)>0$, so that
	\begin{equation}\label{eq-approx:unif_non_crit}
		|\D u_p|\geq c(K)\qquad\text{in $K$, $\forall\,p\in[p(K),\infty]$.}
	\end{equation}
	Here and for the remaining proof, the phrase ``$\,\forall\,p\geq p(\cdot)$ it holds $\cdots$'' means that ``there exists a constant $p_0$ depending on $(\cdot)$, such that for all $p\geq p_0$ it holds $\cdots$''.
	
	It is helpful to keep in mind that
	\begin{equation}\label{eq-approx:lvlset_streamline}
		\text{Regular level sets of $v_p$ are streamlines of $u_p$.}
	\end{equation}
	In the proof, all the closures and interiors are taken in $\RR^2$.
	
	\vspace{3pt}
	
	\noindent\textbf{Selection of $\gamma_i$.}
	
	Fix a sequence of smooth domains $K_n$ with $K_1\Subset K_2\Subset\cdots\Subset\Omega_0$ and $\bigcup K_n=\Omega_0$. If $z_0\notin\Crit(u)$, then we may arrange so that $z_0\in K_1$. For all $n\in\ZZ_{\geq2}$ and $p\geq p(n)$, let $\Gamma_{p,n}=\{\gamma_{i,p,n}\}_{i\geq1}$ be the set of curves in $\{v_p=0\}\cap K_n$ that intersect $\bar{K_{n-1}}$; thus tiny fragments near $\p K_n$ are removed. If $z_0\notin\Crit(u)$, then we may relabel the curves so that $z_0\in\gamma_{1,p,n}$ (recall that we set $v_p(z_0)=0$). By \eqref{eq-approx:unif_non_crit}, $\Gamma_{p,n}$ is a disjoint collection of streamlines of $u_p$, and the endpoints of each $\gamma_{i,p,n}$ lie on $\p K_n$. By \eqref{eq-approx:equicont} \eqref{eq-approx:unif_non_crit}, it follows that $\bigcup_{p\geq p(n)}\Gamma_{p,n}$ is a $C^1$-equicontinuous family of curves.
	
	By \eqref{eq-approx:equicont} \eqref{eq-approx:unif_non_crit} again, the following quantitative separation property holds:
	
	\setcounter{claim}{0}
	
	\begin{lemma}\label{lemma-approx:quant_sep}
		There is a radius $\delta_n>0$ so that the following holds: for all $z\in\bar{K_n}$ and $p\geq p(n)$, at most one connected component of $B(z,5\delta_n)\cap\{v_p=0\}$ can intersect $B(z,\delta_n)$, and any such component is the graph of a $1/100$-Lipschitz function after a rotation.
	\end{lemma}
	\begin{proof}
		Take $\delta_n<\frac1{40}d(K_n,\p K_{n+1})$ so that $|\D u_p(z_1)-\D u_p(z_2)|<e^{-10}c(K_{n+1})$ whenever $p\geq3$, $z_1,z_2\in K_{n+1}$, and $|z_1-z_2|<20\delta_n$. Here $c(K_{n+1})$ is the constant in \eqref{eq-approx:unif_non_crit}. Given any $z\in\bar{K_n}$, for all $p$ large enough so that $|\D u-\D u_p|<e^{-10}c(K_{n+1})$ in $K_{n+1}$, we have
		\[\begin{aligned}
			|\D u_p-\D u(z)| &\leq |\D u_p-\D u_p(z)|+|\D u_p(z)-\D u(z)| \\
			&< 2\cdot e^{-10}c(K_{n+1}) \\
			&< e^{-9}|\D u(z)|\qquad\text{in}\ \ B(z,10\delta_n).
		\end{aligned}\]
		This implies
		\[\Big|\frac{\D u_p}{|\D u_p|}-\frac{\D u(z)}{|\D u(z)|}\Big|<e^{-8}\qquad\Rightarrow\qquad\Big|\frac{\D v_p}{|\D v_p|}-\frac{\D^\perp u(z)}{|\D u(z)|}\Big|<e^{-8} \qquad\text{in}\ \ B(z,10\delta_n).\]
		This forbids two connected components of $\{v_p=0\}\cap B(z,5\delta_n)$ intersecting $B(z,\delta_n)$, and shows that any such connected component is almost parallel to $\frac{\D u(z)}{|\D u(z)|}$.
	\end{proof}
	
	As a consequence, for all $n\geq2$ and $p\geq p(n)$ we have:
	\begin{itemize}[nosep]
		\item There are at most $C(n)$ curves in $\Gamma_{p,n}$, each curve has length at most $C(n)$;
		\item These curves are at least $C(n)^{-1}$-separated in $K_{n-1}$.
	\end{itemize}
	By a general compactness theorem \cite[Theorem 6.3.11]{Simon}, we may select a diagonal subsequence $p_n$ so that $\Gamma_{p_n,n}$ converge as currents in $\Omega_0$. Denote by $\Gamma$ the limit. By the equicontinuity and separation properties above, $\Gamma$ is a locally finite union of multiplicity one $C^1$ curves in $\Omega_0$. Moreover, we have $\p\Gamma\subset\p\Omega_0$, and the convergence $\Gamma_{p_n,n}\to\Gamma$ is $C^1_{\loc}$-graphical with multiplicity one in $\Omega_0$.
	
	Let $\{\gamma_i\}$ be the connected components of $\Gamma$. Thus, each $\gamma_i$ is a streamline of $u$ (since $\Gamma_{p_n,n}$ are made of streamlines of $u_p$), and $\gamma_i$ is an open curve without endpoint in $\Omega_0$. Also, we have $z_0\in\cup\gamma_i$ if $z_0\notin\Crit(u)$, since $z_0\in\{v_p=0\}$ for all $p$.
	
	For our ease of notation, we suppress the index $n$ and write $p_n=p$, and keep in mind that we have chosen a specific subsequence of $p$. The exhaustion $\{K_n\}$ is only used in constructing $\{\gamma_i\}$. The following useful facts follow directly from the construction:
	
	\begin{lemma}\label{lemma-approx:approx_curves}
		$\Omega_0\setminus\cup\gamma_i$ is open. For each $z\in\gamma_i$ where $\gamma_i$ is a curve in $\Gamma$, there are curves $\gamma_p\subset\{v_p=0\}$ that converge to $\gamma_i$ in the $C^1$ graphical sense in a neighborhood of $z$.
	\end{lemma}
	
	\begin{lemma}\label{lemma-approx:no_sign_change_naive}
		If $K\Subset\Omega_0\setminus\cup\gamma_i$, then $\{v_p=0\}\cap K=\emptyset$ for all $p\geq p(K)$.
	\end{lemma}
	
	\vspace{3pt}
	
	\noindent\textbf{Effective bounds on $v_p$.}
	
	The following lemma is of central importance.
	
	\begin{lemma}\label{lemma-approx:vp_bounds}
		Suppose $\sigma:[0,l]\to\Omega$ is a $C^1$ curve with $z_1=\sigma(0)$ and $z_2=\sigma(l)$. Then
		\begin{equation}\label{eq-approx:bound1}
			\limsup_{p\to\infty}\Big|v_p(z_2)-v_p(z_1)\Big|^{\frac1{p-1}}\leq\max_\sigma|\D u|.
		\end{equation}
		Denote
		\[M_\sigma=\Big\{t\in[0,l]: |\D u|(\sigma(t))=\max_\sigma|\D u|\Big\}.\]
		If $\metric{\D^\perp u}{\sigma'}$ has a definite sign $\chi\in\{\pm1\}$ on $M_\sigma$, then
		\begin{equation}\label{eq-approx:bound2}
			\lim_{p\to\infty}\Big[\chi\cdot\big(v_p(z_2)-v_p(z_1)\big)\Big]^{\frac1{p-1}}=\max_\sigma|\D u|.
		\end{equation}
	\end{lemma}
	\begin{proof}
		Note that
		\begin{equation}\label{eq-approx:aux14}
			v_p(z_2)-v_p(z_1)=\int_0^l\metric{\D v_p}{\sigma'}\,dt=\int_0^l|\D u_p|^{p-2}\metric{\D^\perp u_p}{\sigma'}\,dt,
		\end{equation}
		and by the $C^1$ convergence we have
		\[\D u_p\to\D u\qquad\text{in $C^0(\sigma)$.}\]
		
		Denote $A=\max_\sigma|\D u|$. For \eqref{eq-approx:bound1}, we use $|\D u_p|\leq|\D u|+o_p(1)\leq A+o_p(1)$ and \eqref{eq-approx:aux14} to obtain
		\[\begin{aligned}
			\big|v_p(z_2)-v_p(z_1)\big| &\leq\big(A+o_p(1)\big)^{p-1}\cdot|\sigma|,
		\end{aligned}\]
		and the result follows since $|\sigma|^{1/(p-1)}=1+o_p(1)$. Here, $o_p(1)$ denotes certain quantities that converge to 0 uniformly on $\sigma$ as $p\to\infty$.
		
		Then we derive \eqref{eq-approx:bound2}. We may assume $A>0$, otherwise $\metric{\D^\perp u}{\sigma'}$ would be identically zero. Our sign condition implies $\chi\metric{\D^\perp u}{\sigma'}\geq b(1+|\D u|)$ on $M_\sigma$ for some constant $b>0$. For each $\delta>0$, denote
		\[M_\sigma^\delta:=\big(\text{the $\delta$-neighborhood of $M_\sigma$ in $[0,l]$}\big).\]
		By continuity, we may find $\delta_1>0$ so that $\chi\metric{\D^\perp u}{\sigma'}\geq\frac12b(1+|\D u|)$ in $M_\sigma^{\delta_1}$. Next, notice that
		\begin{equation}\label{eq-approx:aux1}
			\sup\Big\{|\D u|(\sigma(t)): t\in[0,l]\setminus M_\sigma^{\delta_1/2}\Big\}<(1-2\mu)A
		\end{equation}
		for some $\mu>0$. Suppose that we are given any $\epsilon\in(0,\mu)$. Choose $\delta_2\in(0,\delta_1/2)$ so that
		\[|\D u|\geq(1-\epsilon)A\qquad\text{in}\ \ M_\sigma^{\delta_2}.\]
		Notice that $M_\sigma^{\delta_2}$ has measure at least $\delta_2$. We now estimate
		\[\begin{aligned}
			\chi\big(v_p(z_2)-v_p(z_1)\big) &= \int_0^l|\D u_p|^{p-2}\cdot\chi\metric{\D^\perp u_p}{\sigma'}\,dt \\
			&= \int_{M_\sigma^{\delta_2}}|\D u_p|^{p-2}\cdot\chi\metric{\D^\perp u_p}{\sigma'}+\int_{M_\sigma^{\delta_1}\setminus M_\sigma^{\delta_2}}|\D u_p|^{p-2}\cdot\chi\metric{\D^\perp u_p}{\sigma'} \\
			&\qquad +\int_{[0,l]\setminus M_\sigma^{\delta_1}}|\D u_p|^{p-2}\cdot\chi\metric{\D^\perp u_p}{\sigma'} \\
			&= \operatorname{I}+\operatorname{II}+\operatorname{III}.
		\end{aligned}\]
		To bound I, we use $\chi\metric{\D^\perp u_p}{\sigma'}\geq\frac12b|\D u_p|$ in $M_\sigma^{\delta_1}$ (hence in $M_\sigma^{\delta_2}$) for all large $p$:
		\[\begin{aligned}
			\operatorname{I} \geq \int_{M_\sigma^{\delta_2}}|\D u_p|^{p-2}\cdot\frac12b|\D u_p|
			&\geq \frac12b|M_\sigma^{\delta_2}|\cdot\Big(\inf_{M_\sigma^{\delta_2}}|\D u_p|\Big)^{p-1} \\
			&\geq \frac{b\delta_2}{2}\Big((1-\epsilon)\big(1-o_p(1)\big)A\Big)^{p-1},\qquad\forall\,p\gg1,
		\end{aligned}\]
		where the $1-o_p(1)$ term originates from the error between $|\D u|$ and $|\D u_p|$. Next, we simply relax $\operatorname{II}\geq0$ for all large $p$. To estimate $\operatorname{III}$, we use \eqref{eq-approx:aux1} to obtain
		\[\begin{aligned}
			\operatorname{III} \geq -\Big(\max_{[0,l]\setminus M_\sigma^{\delta_1}}|\D u_p|\Big)^{p-1}\cdot|\sigma|
			\geq -\Big((1-2\mu)\big(1+o_p(1)\big)A\Big)^{p-1}\cdot|\sigma|.
		\end{aligned}\]
		Hence, in total we have
		\[\begin{aligned}
			&\quad \liminf_{p\to\infty}\Big[\chi\big(v_p(z_2)-v_p(z_1)\big)\Big]^{\frac1{p-1}}
			= \liminf_{p\to\infty}\Big(\!\operatorname{I}+\operatorname{II}+\operatorname{III}\!\Big)^{\frac1{p-1}} \\
			&\geq \liminf_{p\to\infty}\Big[\frac{b\delta_2}{2}\Big((1-\epsilon)\big(1-o_p(1)\big)A\Big)^{p-1}-|\sigma|\Big((1-2\mu)\big(1+o_p(1)\big)A\Big)^{p-1}\Big]^{\frac1{p-1}} \\
			&= (1-\epsilon)A,
		\end{aligned}\]
		since $1-\epsilon>1-2\mu$. By the arbitrariness of $\epsilon$, it follows that
		\[\liminf_{p\to\infty}\Big[\chi\big(v_p(z_2)-v_p(z_1)\big)\Big]^{\frac1{p-1}}\geq A.\]
		Combining this with \eqref{eq-approx:bound1}, we obtain \eqref{eq-approx:bound2} as desired.
	\end{proof}
	
	\vspace{3pt}
	\noindent\textbf{An upper bound of $w_p$.}
	
	Denote
	\begin{equation}\label{eq-approx:def_wp}
		w_p=-\frac{\log|v_p|}{p-1}
	\end{equation}
	whenever $v_p\ne0$. So $e^{-w_p}=|v_p|^{1/(p-1)}$. In this step we show that
	
	\begin{lemma}\label{lemma-approx:w_upper_bound_weak}
		For all $z\in\Omega_0\setminus\cup\gamma_i$ we have
		\begin{equation}\label{eq-approx:w_upper_bound_weak1}
			\limsup_{p\to\infty}w_p(z)\leq-\log|\D u|(z),
		\end{equation}
		or equivalently,
		\begin{equation}\label{eq-approx:w_upper_bound_weak2}
			\liminf_{p\to\infty}|v_p(z)|^{\frac1{p-1}}\geq|\D u|(z).
		\end{equation}
		In fact, this bound holds uniformly: for all $K\Subset\Omega_0\setminus\cup\gamma_i$, we have
		\begin{equation}\label{eq-approx:w_upper_bound_weak_uniform}
			\liminf_{p\to\infty}\Big[\inf_K\Big(|v_p|^{\frac1{p-1}}|\D u|^{-1}\Big)\Big]\geq1.
		\end{equation}
	\end{lemma}
	\begin{proof}
		For a fixed $z\in\Omega_0\setminus\cup\gamma_i$, consider the curve
		\[\sigma(t)=z+t\D^\perp u(z)\ \ (-\delta\leq t\leq\delta),\qquad\sigma_+=\sigma|_{[0,\delta]},\qquad\sigma_-=\sigma|_{[-\delta,0]}.\]
		Choose a small $\delta$ so that $\sigma\Subset\Omega_0\setminus\cup\gamma_i$ and $\metric{\D^\perp u}{\sigma'}>0$ on $\sigma$. So $v_p|_\sigma$ has a definite sign for all large $p$ (Lemma \ref{lemma-approx:no_sign_change_naive}). Take a subsequence that realizes the liminf in \eqref{eq-approx:w_upper_bound_weak2}. If the sign of $v_p|_\sigma$ is positive for infinitely many $p$ in this subsequence, then we may pass to a further subsequence $\{u_{p_k}\}$ and use \eqref{eq-approx:bound2} with $\chi=1$ to estimate
		\[\liminf_{p\to\infty}\Big(v_p(z)^{\frac1{p-1}}\Big)
			\geq \lim_{k\to\infty}\Big(v_{p_k}(z)-\overbrace{v_{p_k}(\sigma(-\delta))}^{\geq0}\Big)^{\frac1{p-1}}
			= \max_{\sigma_-}|\D u|.\]
		If this sign is negative for infinitely many $p$ in the subsequence, then we may pass to an analogous subsequence and use \eqref{eq-approx:bound2} with $\chi=1$ to estimate
		\[\liminf_{p\to\infty}\Big(\big(\!-v_p(z)\big)^{\frac1{p-1}}\Big)
			\geq \liminf_{k\to\infty}\Big(\overbrace{v_{p_k}(\sigma(\delta))}^{\leq0}-v_{p_k}(z)\Big)^{\frac1{p-1}}
			= \max_{\sigma_+}|\D u|.\]
		In combination, we obtain
		\[\liminf_{p\to\infty}|v_p(z)|^{\frac1{p-1}}\geq\min\Big\{\max_{\sigma_-}|\D u|,\max_{\sigma_+}|\D u|\Big\}.\]
		Taking $\delta\to0$, the inequalities \eqref{eq-approx:w_upper_bound_weak1} \eqref{eq-approx:w_upper_bound_weak2} follow.
		
		Finally, we prove \eqref{eq-approx:w_upper_bound_weak_uniform}. Note that $\{v_p=0\}\cap N(K,\epsilon)=\emptyset$ for some $\epsilon=\epsilon(K)>0$ and for all $p\geq p(K)$, by Lemma \ref{lemma-approx:no_sign_change_naive}. Then Kotschwar-Ni's interior gradient estimate \cite[Theorem 1.1]{Kotschwar-Ni_2009} implies
		\[\sup_K|\D w_p|\leq C(K)\qquad\forall\,K\Subset\Omega_0\setminus\cup\gamma_i,\ \forall\,p\geq p(K).\]
		This shows that the lower bound \eqref{eq-approx:w_upper_bound_weak2} is uniform in all compact sets in $\Omega_0\setminus\cup\gamma_i$.
	\end{proof}
	
	\vspace{3pt}
	
	\noindent\textbf{Critical and regular regions.}
	
	Recall that we have chosen a countable collection of streamlines $\{\gamma_i\subset\Omega_0\}$. Define $\cD,\cD_{\Crit},\cD_{\reg}$ as in item \ref{item-approx:Dreg_Dcrit} of the main theorem:
	\begin{equation}\label{eq-approx:def_D}
		\cD=\Big\{\text{connected components of $\Omega\setminus\bar{\cup\gamma_i}$}\Big\},
	\end{equation}
	and
	\begin{equation}\label{eq-approx:def_Dreg_Dcrit}
		\cD_{\Crit}=\Big\{D\in\cD: D\subset\Crit(u)\Big\},\qquad \cD_{\reg}=\cD\setminus\cD_{\Crit}.
	\end{equation}
	It is useful to keep in mind that $(\bar{\cup\gamma_i}\setminus\gamma_i)\cap\Omega\subset\p\Crit(u)$. Finally, set as in \ref{item-approx:Omega_reg}
	\begin{equation}\label{eq-approx:def_omega_reg_2}
		\Omega_{\reg}=\Big(\bigcup_{D\in\cD_{\reg}}D\Big)\cup\big(\!\cup\!\gamma_i\big).
	\end{equation}
	Note that $\Omega_{\reg}$ is open (since any $z\in\cup\gamma_i$ is adjacent to two regions in $\cD_{\reg}$), and $\Omega_{\reg}\supset\Omega_0$, which implies $\Omega\setminus\Omega_{\reg}\subset\Crit(u)$. In this step, we prove that:
	
	\begin{lemma}\label{lemma-approx:nondeg} {\ }
		
		\begin{enumerate}[label={(\roman*)}, nosep]
			\item For any $L\Subset\Omega\setminus\Omega_{\reg}$ and any $\epsilon>0$, there exists a $\delta>0$ such that
			\begin{equation}\label{eq-approx:unif_blow_up}
				|v_p(z)|^{\frac1{p-1}}\leq\epsilon\qquad\forall\,z\in N(L,\delta),\ p\geq p(\epsilon,L).
			\end{equation}
			\item If $D\in\cD_{\reg}$, then for each $L\Subset D$, it holds $\{v_p=0\}\cap L=\emptyset$ for all $p\geq p(L)$.
		\end{enumerate}
	\end{lemma}
	
	Item (ii) is a consequence of Lemma \ref{lemma-approx:no_sign_change_naive} if $L\Subset D\cap\Omega_0$, thus we need to extend this property when $D$ contains critical points. To prove (ii), we need the following Harnack inequality (applied with $q=\frac{p}{p-1}$ $\Leftrightarrow$ $q-1=\frac1{p-1}$), which follows from Kotschwar-Ni's gradient estimate \cite[Theorem 1.1]{Kotschwar-Ni_2009} or R. Moser's Harnack inequality \cite[Section 5]{Moser_2008}:
	
	\begin{theorem}\label{thm-approx:harnack}
		Let $K\Subset L\subset\RR^2$ be domains (i.e. connected open sets). Then for all $q\in(1,3/2]$ and positive $q$-harmonic function $v$ in $L$, it holds
		\[\Big(\frac{\sup_K(v)}{\inf_K(v)}\Big)^{q-1}\leq C(K,L).\]
	\end{theorem}
	
	\begin{proof}[Proof of Lemma \ref{lemma-approx:nondeg}] {\ }
		
		(i) Let $L,\epsilon$ be as stated. Choose a small $\rho=\rho(L)$ such that any point in $L$ can be connected to $z_0$ by a curve in $N(\Omega,-2\rho)\cap B(1/(2\rho))$ with length $\leq C(L)$. Here we denote
		\[N(\Omega,-\rho)=\big\{z\in\Omega: B(z,\rho)\Subset\Omega\big\}.\]
		Denote
		\[L'=\big(N(\Omega,-\rho)\cap B(1/\rho)\big)\setminus\Omega_{\reg}.\]
		We may further decrease $\rho$ and assume $L\Subset L'$. Since $L'\Subset\Crit(u)$, there is $\delta<\rho/4$ so that
		\begin{equation}\label{eq-approx:aux2}
			\sup_{N(L',4\delta)}|\D u_p|\leq\epsilon/2\qquad\forall\,p\geq p(\epsilon,\rho,L').
		\end{equation}
		Let $z\in N(L,\delta)$ be given. We find a point $z_1\in L$ with $|z-z_1|<\delta$, and then find a $C^1$ curve
		\[\sigma:[0,l]\to N(\Omega,-2\rho)\cap B(1/(2\rho)),\qquad|\sigma|\leq C(L),\]
		that joins $z_1$ to $z_0$, as was guaranteed above. Let $t$ be the largest time so that $\sigma([0,t])\subset\bar{N(L',\delta)}$. If $t=l$, then by \eqref{eq-approx:aux2} and direct estimates, we have
		\[\begin{aligned}
			|v_p(z)| &= |v_p(z)-v_p(z_1)+v_p(z_1)-v_p(z_0)| \\
			&\leq \int_{[z,z_1]+\sigma}|\D v_p|\leq\big(\delta+|\sigma|\big)\cdot(\epsilon/2)^{p-1}\leq\epsilon^{p-1},\qquad\forall p\geq p(\epsilon,\rho,L',\delta,L).
		\end{aligned}\]
		Suppose $t<l$. So $\sigma(t)\notin L'$, and the construction implies $\sigma(t)\in\Omega_{\reg}$. Then we observe that $\sigma([0,t])\cap(\bar{\cup\gamma_i})\ne\emptyset$: otherwise, $\sigma([0,t])$ would be contained in some $D\in\cD_{\reg}$, contradicting $\sigma(0)=z_1\in L$ and $L\cap\Omega_{\reg}=\emptyset$. Recall that $\cup\gamma_i$ is obtained by taking the limit of $\{v_p=0\}$ (see also Lemma \ref{lemma-approx:approx_curves}), so for all $p\geq p(\delta,\sigma)$ we have
		\[S\cap\{v_p=0\}\ne\emptyset,\quad\text{where}\quad S:=N\big(\sigma([0,t]),\delta\big).\]
		Take any point $z_2$ in this set. Then by \eqref{eq-approx:aux2} and $S\Subset N(L',4\delta)$ we have
		\[|v_p(z)|=|v_p(z)-v_p(z_2)|\leq\big(2\delta+|\sigma|\big)\cdot\sup_S|\D v_p|<\epsilon^{p-1},\qquad\forall\,p\geq p(\epsilon,\rho,L',\delta,L).\]
		In either case, the lemma follows.
		
		\begin{figure}[ht]
			\centering
			\includegraphics{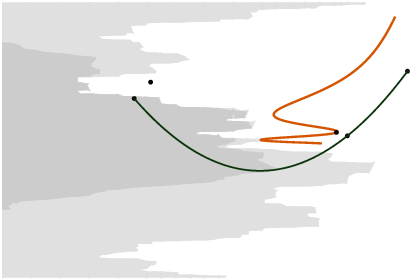}
			\begin{picture}(0,0)
				\put(-198,121){$L'$}
				\put(-198,100){$L$}
				\put(-127,90){$z$}
				\put(-146,78){$z_1$}
				\put(-45,77){$z_2$}
				\put(-7,103){$z_0$}
				\put(-115,52){$\sigma$}
				\put(-33,62){$\sigma(t)$}
			\end{picture}
			\caption{Objects in the proof of Lemma \ref{lemma-approx:nondeg}(i).}\label{fig-approx:deg}
		\end{figure}
		
		\vspace{3pt}
		
		(ii) Fix $D\in\cD_{\reg}$. By a covering argument, it suffices to show that for any $z\in D$, there exists a radius $\rho=\rho(z)>0$ so that
		\[\{v_p=0\}\cap B(z,\rho)=\emptyset\qquad\forall\,p\geq p(z).\]
		Let $U$ be the set of points in $D$ that has this property. Clearly, $U$ is open. Recall Lemma \ref{lemma-approx:no_sign_change_naive}: for any $K\Subset D\setminus\Crit(u)$, we have
		\begin{equation}\label{eq-approx:aux3}
			\{v_p=0\}\cap K=\emptyset\qquad\forall\,p\geq p(K).
		\end{equation}
		This implies $D\setminus\Crit(u)\subset U$, or $D\setminus U\subset\Crit(u)$.
		
		Suppose $U\ne D$. Then there is a smooth path $\sigma:[0,l]\to D$ and radius $r>0$, so that
		\[N(\sigma,4r)\subset D,\quad B(\sigma(0),2r)\subset D\setminus\Crit(u),\quad N(\sigma,2r)\subset U,\quad \bar{B(\sigma(l),2r)}\not\subset U.\]
		Denote $z_1=\sigma(0)$ and $z_2=\sigma(l)$. See Figure \ref{fig-approx:nondeg} for a depiction of these objects. Hence
		\[B(z_1,2r)\subset D\setminus\Crit(u),\qquad \p B(z_2,2r)\setminus U\ne\emptyset.\]
		
		\begin{figure}[ht]
			\centering
			\includegraphics{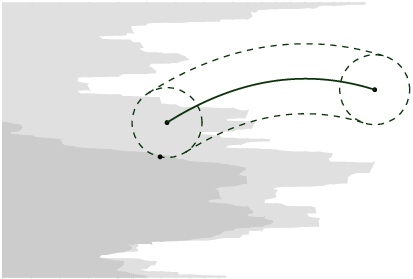}
			\begin{picture}(0,0)
					\put(-198,119){$\Crit(u)$}
					\put(-198,50){$D\setminus U$}
					\put(-19,93){$z_1$}
					\put(-135,77){$z_2$}
					\put(-131,49){$z_3$}
					\put(-77,98){$\sigma$}
					\put(-45,117){$N(\sigma,2r)$}
				\end{picture}
			\caption{Objects in the proof of Lemma \ref{lemma-approx:nondeg}(ii).}\label{fig-approx:nondeg}
		\end{figure}
		
		By Lemma \ref{lemma-approx:w_upper_bound_weak} and $B(z_1,r)\Subset D\setminus\Crit(u)\subset\Omega_0\setminus\cup\gamma_i$, there is a constant $c'>0$ so that
		\begin{equation}
			\inf_{B(z_1,r)}|v_p|^{\frac1{p-1}}\geq c',\qquad\forall\,p\geq p(z_1,r).
		\end{equation}
		Then, by our assumption $N(\sigma,2r)\subset U$, the function $v_p$ has a definite sign in $N(\sigma,3r/2)$ for all $p\geq p(\sigma,r)$. Applying Theorem \ref{thm-approx:harnack} with $\text{``$K$''}=N(\sigma,r)$ and $\text{``$L$''}=N(\sigma,3r/2)$, we have
		\begin{equation}\label{eq-approx:vp_inner_bd}
			\inf_{B(z_2,r)}|v_p|^{\frac1{p-1}}\geq c,\qquad\forall\,p\geq p(\sigma,r),
		\end{equation}
		for some other constant $c>0$.
		
		Next, fix $z_3\in\p B(z_2,2r)\setminus U$, so $z_3\in\Crit(u)$. According to the definition of $U$, we may find a subsequence of $p$ and points $\zeta_p\to z_3$ so that $v_p(\zeta_p)=0$.

		For each $s\in[r,2r)$, we denote $a_{s,p}=\min_{\bar{B(z_2,s)}}|v_p|^{1/(p-1)}$. We claim that
		\begin{equation}\label{eq-approx:liminf_asp}
			\lim_{s\to2r}\Big(\limsup_{p\to\infty}a_{s,p}\Big)=0,
		\end{equation}
		where the limsup is taken along our selected subsequence. To show this, consider $z_s=z_2+\frac s{2r}(z_3-z_2)$. Note that $z_s\in\p B(z_2,s)$ and $|z_s-z_3|=2r-s$. We estimate
		\[a_{s,p}^{p-1}\leq|v_p(z_s)|\leq|v_p(z_s)-v_p(z_3)|+|v_p(z_3)-v_p(\zeta_p)|.\]
		Then for all $p$ large enough so that $|\zeta_p-z_3|\leq2r-s$, we obtain
		\[a_{s,p}\leq\Big(2(2r-s)\cdot\sup_{B(z_3,2r-s)}|\D v_p|\Big)^{\frac1{p-1}}=\big(2(2r-s)\big)^{\frac1{p-1}}\cdot\sup_{B(z_3,2r-s)}|\D u_p|.\]
		Then \eqref{eq-approx:liminf_asp} follows by $z_3\in\Crit(u)$ and the uniform convergence $\D u_p\to\D u$.
		
		For each $p$ in our selected subsequence, consider the $q$-harmonic function
		\[\phi(\zeta)=\frac{a_{s,p}^{p-1}s^{p-2}-c^{p-1}r^{p-2}}{s^{p-2}-r^{p-2}}+\frac{c^{p-1}-a_{s,p}^{p-1}}{s^{p-2}-r^{p-2}}\Big(\frac{rs}{|\zeta-z_2|}\Big)^{p-2},\]
		where $c$ is the constant in \eqref{eq-approx:vp_inner_bd}. It can be verified that
		\[\phi|_{\p B(z_2,r)}=c^{p-1},\qquad\phi|_{\p B(z_2,s)}=a_{s,p}^{p-1}.\]
		Recall that $v_p$ is also $q$-harmonic and has no sign change in $B(z_2,s)$ for all $p\geq p(z_2,s)$. The maximum principle and \eqref{eq-approx:vp_inner_bd} then show that either $v_p\geq\phi$ or $v_p\leq-\phi$ in $B(z_2,s)\setminus B(z_2,r)$, for all $p\geq p(z_2,\sigma,s)$. For each $s,p$ we choose an extremal point $\zeta_{s,p}\in\p B(z_2,s)$ so that $a_{s,p}=|v_p(\zeta_{s,p})|^{1/(p-1)}$. The boundary gradient comparison then implies
		\begin{equation}\label{eq-approx:bd_grad_comp}
			|\D v_p(\zeta_{s,p})|\geq|\D\phi(\zeta_{s,p})|=(p-2)\frac{c^{p-1}-a_{s,p}^{p-1}}{s^{p-2}-r^{p-2}}\frac{(rs)^{p-2}}{s^{p-1}},\qquad\forall\,p\geq p(z_2,s).
		\end{equation}
		Taking a sequence $s_n\nearrow2r$ and passing to a further diagonal subsequence $p_n$, we have $\zeta_{s_n,p_n}\to z_\infty$ for some $z_\infty\in\p B(z_2,2r)$, and the limit in \eqref{eq-approx:liminf_asp} is attained (namely, $a_{s_n,p_n}\to0$), and \eqref{eq-approx:bd_grad_comp} holds along this subsequence. (Note: $z_\infty$ need not be $z_3$.) Hence
		\[|\D u(z_\infty)|=\lim_{n\to\infty}|\D u_{p_n}(\zeta_{s_n,p_n})|=\lim_{n\to\infty}\Big(|\D v_{p_n}(\zeta_{s_n,p_n})|^{\frac1{p_n-1}}\Big)\geq\frac c2.\]
		Thus $z_\infty\notin\Crit(u)$. On the other hand, by Lemma \ref{lemma-approx:w_upper_bound_weak}, there exists $\rho>0$ so that
		\[|v_p|^{\frac1{p-1}}\geq\frac12|\D u(z_\infty)|\quad\text{in}\ \ B(z_\infty,\rho),\ \ \forall\,p\geq p(z_\infty).\]
		This contradicts
		\[0=\lim_{n\to\infty}a_{s_n,p_n}=\lim_{n\to\infty}|v_{p_n}(\zeta_{s_n,p_n})|^{\frac1{p_n-1}}.\]
		The lemma follows as desired.
	\end{proof}
	
	\vspace{3pt}
	
	\noindent\textbf{Convergence of $v_p$ in regular regions.}
	
	Fix a basepoint $z_D\in D\setminus\Crit(u)$ for each $D\in\cD_{\reg}$. Passing to a diagonal subsequence of $p$, we may assume that $v_p(z_D)$ have the same sign for all sufficiently large $p$, for each $D$. Set $\chi_D=1$ if this sign is positive, and $\chi_D=-1$ if it is negative. In Lemma \ref{lemma-approx:nondeg}(ii), we have shown that for each connected $L\Subset D\in\cD_{\reg}$, the function $v_p|_L$ has a definite sign for all large $p$. Clearly, this sign agrees with $\chi_D$.
	
	Each $\p D\cap\Omega$ is the union of a subset of $\cup\gamma_i$ and a subset of $\p\Crit(u)$, since $\bar{\cup\gamma_i}\setminus\cup\gamma_i\subset\p\Crit(u)$. Note that if $D_1,D_2$ are adjacent (meaning that their boundaries contain the same portion of a curve $\gamma_i$), then $\chi_{D_1}=-\chi_{D_2}$. This is because the convergence $\Gamma_{p,n}\to\{\gamma_i\}$ is one-fold, and $v_p$ changes sign when crossing any curve in $\Gamma_{p,n}$.
	
	Each $D\in\cD_{\reg}$ is a $C^1$ domain in $\Omega_{\reg}$. Indeed, any $z\in\p D\cap\Omega_{\reg}$ must lie on some $\gamma_i$. Then, for a small $r>0$, the set $B(z,r)\setminus\cup\gamma_i$ is the union of two $C^1$ sub-graphs, and $D$ contains exactly one of them due to sign distinction.
	
	Fix $D\in\cD_{\reg}$, and consider any connected $P_1,P_2$ such that $z_D\in P_1\Subset P_2\Subset D$. For all large enough $p\geq p(P_2)$, we may set
	\[w_p=(1-q)\log|v_p|=-\frac{\log|v_p|}{p-1}=-\frac{\log(\chi_Dv_p)}{p-1}\qquad\text{in}\ \ P_2.\]
	By Kotschwar-Ni's interior gradient estimate \cite[Theorem 1.1]{Kotschwar-Ni_2009}, we have
	\begin{equation}\label{eq-approx:grad_est}
		|\D w_p|\leq\frac{C}{d(P_1,\p P_2)}\qquad\text{in}\ \ \bar{P_1}.
	\end{equation}
	To apply Arzel\`a-Ascoli, we need a pointwise bound for $w_p$. The upper bound follows from \eqref{eq-approx:w_upper_bound_weak1}:
	\begin{equation}\label{eq-approx:wp_upper_bound}
		\limsup_{p\to\infty}w_p(z_D)\leq-\log|\D u|(z_D)<\infty.
	\end{equation}
	
	To derive a lower bound, choose a $C^1$ path $\sigma\subset\Omega$ that joins $z_0$ and $z_D$. By \eqref{eq-approx:bound1} and $v_p(z_0)=0$, we have
	\begin{equation}\label{eq-approx:wp_lower_bound}
		\liminf_{p\to\infty}w_p(z_D)=-\limsup_{p\to\infty}\log\Big(|v_p(z_D)|^{\frac1{p-1}}\Big)\geq-\log\max_\sigma|\D u|>-\infty.
	\end{equation}
	
	Combining \eqref{eq-approx:grad_est} \eqref{eq-approx:wp_upper_bound} \eqref{eq-approx:wp_lower_bound} and the Arzel\`a-Ascoli theorem, and passing to a subsequence, we have $w_p\to w$ in $C^0(P_1)$ for some $w$. Letting $P_1,P_2\nearrow D$, and letting $D$ range in $\cD_{\reg}$, then finally passing to a diagonal subsequence, it follows that we can find functions $w_D\in C^0_{\loc}(D)$ for each $D\in\cD_{\reg}$, so that
	\[w_p\to w_D\qquad\text{in}\ \ C^0_{\loc}(D).\]
	Define the function $w$ on $\bigcup_{D\in\cD_{\reg}}D$ by setting $w|_D=w_D$. We have
	\begin{equation}\label{eq-approx:w_upper_bound_partial}
		w\leq-\log|\D u|\qquad\text{in}\ \ \bigcup_{D\in\cD_{\reg}}D.
	\end{equation}
	Indeed, if $z\in\Omega_0$ then this follows from \eqref{eq-approx:w_upper_bound_weak1}, and if $z\in\Crit(u)$ then this is obvious. For the remaining proof, we work with the chosen subsequence of $p$.
	
	\vspace{3pt}
	
	\noindent\textbf{Almost proof of item \ref{item-approx:w_on_path}.}
	
	The goal of this step is to prove the following two lemmas:
	
	\begin{lemma}\label{lemma-approx:weaker_item_9}
		Suppose $z_1\in\cup\gamma_i$ or $z_1=z_0$, and $\sigma:[0,l]\to\Omega$ is a $C^1$ curve joining $z_1$ and a point $z_2\in D\in\cD_{\reg}$. Denote
		\begin{equation}\label{eq-approx:def_A_Msigma}
			A=\max_\sigma|\D u|,\qquad M_\sigma=\Big\{t\in[0,l]: |\D u|(\sigma(t))=A\Big\}.
		\end{equation}
		Then $w(z_2)\geq-\log A$. If in addition, $\metric{\D^\perp u}{\sigma'}$ has a definite sign $\chi\in\{\pm1\}$ on $M_\sigma$, then this sign agrees with $\chi_D$, and we have $w(z_2)=-\log A$.
	\end{lemma}
	
	\begin{lemma}\label{lemma-approx:weaker_item_9_2}
		Suppose $z_1\in\cup\gamma_i$ or $z_1=z_0$, and $\sigma:[0,l]\to\Omega$ is a $C^1$ curve joining $z_1$ and a point $z_2\in\Omega\setminus\bigcup_{D\in\cD_{\reg}}D$. Define $A,M_\sigma$ as in \eqref{eq-approx:def_A_Msigma}. Then $\metric{\D^\perp u}{\sigma'}$ cannot have a definite sign on $M_\sigma$.
	\end{lemma}
	
	These almost imply item \ref{item-approx:w_on_path} of the main theorem; the only (minor) difference is that $w$ has not been defined outside $\cD_{\reg}$ yet.
	
	\begin{proof}[Proof of Lemma \ref{lemma-approx:weaker_item_9}] {\ }
		
		In the case $z_1=z_0$, Lemma \ref{lemma-approx:vp_bounds} and $v_p(z_0)=0$ together implies
		\[\left\{\begin{aligned}
			& \limsup_{p\to\infty}|v_p(z_2)|^{\frac1{p-1}}\leq A, \\
			& \lim_{p\to\infty}\big[\chi v_p(z_2)\big]^{\frac1{p-1}}=A\qquad\text{if the sign condition holds.}
		\end{aligned}\right.\]
		Also, recall that
		\begin{equation}\label{eq-approx:aux15}
			e^{-w}=\lim_{p\to\infty}|v_p|^{\frac1{p-1}}=\lim_{p\to\infty}\big(\chi_D v_p\big)^{\frac1{p-1}}.
		\end{equation}
		The result follows directly. We may now assume $z_1\in\gamma_i$ for some fixed $i$. Recall that $\gamma_i\subset\Omega_0$, so we always have $A>0$.
		
		We first prove $w(z_2)\geq-\log A$. Recall that $\gamma_i$ is a limit of subsets of $\{v_p=0\}$ (see Lemma \ref{lemma-approx:approx_curves}). For any $r>0$ and $p\geq p(z_1,r)$, we may find a point $z_{r,p}\in B(z_1,r)\cap\{v_p=0\}$. From this we may estimate
		\[|v_p(z_2)|\leq|v_p(z_2)-v_p(z_1)|+|v_p(z_1)-v_p(z_{r,p})|\leq\big(|\sigma|+r\big)\cdot\sup_{\sigma\cup B(z_1,r)}|\D u_p|^{p-1}.\]
		Taking $\frac1{p-1}$-th power of both sides and sending $p\to\infty$, we obtain
		\[e^{-w(z_2)}\leq\max\Big\{A,\sup_{B(z_1,r)}|\D u|\Big\}.\]
		Finally, taking $r\to0$, this implies $e^{-w(z_2)}\leq A$ as desired.
		
		For the rest of the proof, we assume that $\metric{\D^\perp u}{\sigma'}$ has a sign $\chi\in\{\pm1\}$ in $M_\sigma$, and we are aimed at showing $w(z_2)=-\log A$. We divide the proof into two cases.
		
		\vspace{3pt}
		
		\textbf{Case 1:} Suppose $0\notin M_\sigma$. Hence $|\D u|(z_1)<A$. Recall that there are points $\zeta_p\in\{v_p=0\}$ with $\zeta_p\to z_1$. And by the uniform convergence $\D u_p\to\D u$, there exists a $\delta>0$ so that
		\begin{equation}\label{eq-approx:aux4}
			\sup_{p>1/\delta}\sup_{B(z_1,\delta)}|\D u_p|\leq(1-\delta)A.
		\end{equation}
		Then we calculate
		\begin{equation}\label{eq-approx:aux5}
			\begin{aligned}
				e^{-w(z_2)}
				&= \lim_{p\to\infty}\big[\chi_Dv_p(z_2)\big]^{\frac1{p-1}} \\
				&= \lim_{p\to\infty}\Big[\chi_D\cdot\big(v_p(z_2)-v_p(z_1)+v_p(z_1)-v_p(\zeta_p)\big)\Big]^{\frac1{p-1}},
			\end{aligned}
		\end{equation}
		and by \eqref{eq-approx:bound2} we have
		\begin{equation}\label{eq-approx:aux6}
			\lim_{p\to\infty}\Big[\chi\cdot\big(v_p(z_2)-v_p(z_1)\big)\Big]^{\frac1{p-1}}=A,
		\end{equation}
		and by direct estimate using \eqref{eq-approx:aux4} we have
		\begin{equation}\label{eq-approx:aux7}
			\big|v_p(z_1)-v_p(\zeta_p)\big|^{\frac1{p-1}}\leq\Big(\delta\sup_{B(z_1,\delta)}|\D v_p|\Big)^{\frac1{p-1}}\leq\delta^{\frac1{p-1}}(1-\delta)A,\qquad\forall\,p\geq p(z_1,\delta).
		\end{equation}
		Combining \eqref{eq-approx:aux5}\,$\sim$\,\eqref{eq-approx:aux7}, we obtain $\chi=\chi_D$ and $e^{-w(z_2)}=A$, as desired.
		
		\vspace{3pt}
		
		\textbf{Case 2:} Suppose $0\in M_\sigma$. So $|\D u|(z_1)=A$, and $\sigma$ intersects $\gamma_i$ transversally at $z_1$ (since $\gamma'_i$ is parallel to $\D u$). Assume $\chi\metric{\D^\perp u}{\sigma'}\geq b>0$ on $M_\sigma$ for some constant $b>0$. Extend $\sigma$ to a $C^1$ curve on $[-\epsilon,l]$. We may choose a small enough $\epsilon$ and arrange the extension so that
		\begin{equation}\label{eq-approx:aux_sign}
			\chi\metric{\D^\perp u}{\sigma'}(\sigma(t))\geq b/2\qquad\forall\,t\in[-\epsilon,\epsilon].
		\end{equation}
		Recall by Lemma \ref{lemma-approx:approx_curves} that $\gamma_i$ is the $C^1_{\loc}$ graphical limit of curves $\gamma_p\subset\{v_p=0\}$ near $z_1$. So for large enough $p$, the intersection $\gamma_p\cap\sigma$ remains transversal near $z_1$. Let $\zeta_p=\sigma(t_p)$ be this unique intersection point. Clearly $t_p\to0$.
		
		For a constant $\delta$ with $|\delta|$ small, we need to bound $v_p(z_2)-v_p(\sigma(\delta))$. For this purpose, we consider the curve $\sigma_\delta=\sigma|_{[\delta,l]}$, and the associated objects
		\[A_\delta=\max_{\sigma_\delta}|\D u|,\qquad M_{\sigma,\delta}=\Big\{t\in[\delta,l]: |\D u|(\sigma(t))=A_\delta\Big\}.\]
		The following lemma is useful:
		\begin{lemma}\label{lemma-approx:stability}
			We have $\lim_{\delta\to0}A_\delta=A$. Furthermore, $\metric{\D^\perp u}{\sigma'}$ also has a definite sign on $M_{\sigma,\delta}$ as long as $|\delta|$ is small enough.
		\end{lemma}
		\begin{proof}
			The first statement is elementary. If the second statement fails, then there are $\delta_i\to0$ and $t_i\in M_{\sigma,\delta_i}$, so that $\chi\metric{\D^\perp u}{\sigma'}(t_i)\leq0$. Passing to a subsequence, we may assume $t_i\to t\in[0,l]$. Then $|\D u|(\sigma(t))=\lim_{i\to\infty}|\D u|(\sigma(t_i))=\lim_{i\to\infty}A_{\delta_i}=A$, hence $t\in M_\sigma$, hence $\chi\metric{\D^\perp u}{\sigma'}(t)>0$, contradiction.
		\end{proof}
		
		For each $0<\delta<\epsilon$ so that Lemma \ref{lemma-approx:stability} holds for $\pm\delta$, we use Lemma \ref{lemma-approx:vp_bounds} to obtain
		\begin{equation}\label{eq-approx:aux8}
			\begin{aligned}
				\lim_{p\to\infty}\Big[\chi\big(v_p(z_2)-v_p(\sigma(\delta))\big)\Big]^{\frac1{p-1}} &= A_\delta, \\
				\lim_{p\to\infty}\Big[\chi\big(v_p(z_2)-v_p(\sigma(-\delta))\big)\Big]^{\frac1{p-1}} &= A_{-\delta}.
			\end{aligned}
		\end{equation}
		Combining $v_p(\sigma(t_p))=0$ and \eqref{eq-approx:aux_sign}, we have
		\[\chi v_p(\sigma(\delta))>0,\qquad \chi v_p(\sigma(-\delta))<0,\qquad\forall\,p\geq p(\sigma,\delta).\]
		Then \eqref{eq-approx:aux8} implies that
		\[A_\delta\leq\liminf_{p\to\infty}\big[\chi v_p(z_2)\big]^{\frac1{p-1}}\leq\limsup_{p\to\infty}\big[\chi v_p(z_2)\big]^{\frac1{p-1}}\leq A_{-\delta}.\]
		Taking $\delta\to0$ of this inequality, it follows that
		\[\lim_{p\to\infty}\big[\chi v_p(z_2)\big]^{\frac1{p-1}}=A,\]
		which implies $e^{-w(z_2)}=A$ and $\chi=\chi_D$ (where $D\in\cD_{\reg}$ is the region containing $z_2$).
	\end{proof}
	
	\begin{proof}[Proof of Lemma \ref{lemma-approx:weaker_item_9_2}] {\ }
		
		Suppose otherwise that $\metric{\D^\perp u}{\sigma'}$ has a sign $\chi\in\{\pm1\}$ on $M_\sigma$.
		
		If $0\in M_\sigma$, then we must have $z_1\in\gamma_i$ for some $i$ (since in this case, $z_1=z_0$ implies $z_0\notin\Crit(u)$, hence $z_0\in\cup\gamma_i$), and $\sigma$ intersects $\gamma_i$ transversally (as $\gamma'_i$ is parallel to $\D u$). By Lemma \ref{lemma-approx:approx_curves} and the same extension argument in the proof of Lemma \ref{lemma-approx:weaker_item_9} above, there are $t_p\to0$ so that $v_p(\sigma(t_p))=0$. Integrating from $\sigma(t_p)$ to $\sigma(\delta)$ gives
		\begin{equation}\label{eq-approx:aux12}
			\chi v_p(\sigma(\delta))>0\qquad\forall\,\delta<\delta(\sigma),\,p>p(\sigma,\delta).
		\end{equation}
		If $0\notin M_\sigma$, then $|\D u(z_1)|<A$, so by the equicontinuity of $\D u_p$, there is $\epsilon>0$ so that
		\[\sup_{p>1/\epsilon}\sup_{B(z_1,\epsilon)}|\D u_p|\leq(1-\epsilon)A.\]
		Then, note that $B(z_1,\epsilon)\cap\{v_p=0\}\ne\emptyset$ for all $p\geq p(z_1,\epsilon)$, since $z_1\in\{z_0\}\cup(\cup\gamma_i)$. Hence
		\begin{equation}\label{eq-approx:aux13}
			|v_p(\sigma(\delta))|\leq2\epsilon(1-\epsilon)^{p-1}A^{p-1},\qquad\forall\,\delta<\epsilon,\,p>p(z_1,\epsilon).
		\end{equation}
		Combining both cases, we find that there is always an $\epsilon>0$ so that
		\begin{equation}\label{eq-approx:aux9}
			\chi v_p(\sigma(\delta))\geq-(1-\epsilon)^{p-1}A^{p-1},\qquad\forall\,\delta<\delta(\sigma,\epsilon),\,p>p(\sigma,\delta,\epsilon).
		\end{equation}
		
		Consider the other end of $\sigma$. Set
		\[T=\min\Big\{t\in[0,l]:\sigma|_{[t,l]}\subset\Omega\,\big\backslash\!\bigcup_{D\in\cD_{\reg}}D\Big\}.\]
		Then $\sigma([T,l])\subset(\cup\gamma_i)\cup\Crit(u)$, and hence $\metric{\D^\perp u}{\sigma'}=0$ on $(T,l]$, and hence $M_\sigma\subset[0,T]$. Also, $\sigma(T)$ lies on the boundary of $\bigcup_{D\in\cD_{\reg}}D$, hence lies on $\bar{\cup\gamma_i}$. Arguing analogously as above, we conclude that
		\begin{equation}\label{eq-approx:aux10}
			\chi v_p(\sigma(T-\delta))\leq(1-\epsilon)^{p-1}A^{p-1},\qquad\forall\,\delta<\delta(\sigma,\epsilon),\,p>p(\sigma,\delta,\epsilon),
		\end{equation}
		after possibly decreasing $\epsilon$. Denote
		\[A_\delta=\max_{\sigma([\delta,T-\delta])}|\D u|,\qquad M_{\sigma,\delta}=\Big\{t\in[\delta,T-\delta]:|\D u|(\sigma(t))=A_\delta\Big\}.\]
		An analogue of Lemma \ref{lemma-approx:stability} holds in this context as well: for all small enough $\delta$, we have $A_\delta>(1-\epsilon)A$, and $\metric{\D^\perp u}{\sigma'}$ has a sign in the set where $|\D u|$ attains maximum on $[\delta,T-\delta]$. Then, Lemma \ref{lemma-approx:vp_bounds} implies
		\begin{equation}\label{eq-approx:aux11}
			\lim_{p\to\infty}\Big[\chi\cdot\big(v_p(\sigma(T-\delta))-v_p(\sigma(\delta))\big)\Big]^{\frac1{p-1}}=A_\delta.
		\end{equation}
		Combining \eqref{eq-approx:aux9} \eqref{eq-approx:aux10} \eqref{eq-approx:aux11} and $A_\delta>(1-\epsilon)A$, we obtain
		\[0=\lim_{p\to\infty}\Big[\chi v_p(\sigma(\delta))+ \chi\Big(v_p(\sigma(T-\delta))-v_p(\sigma(\delta))\Big)-\chi v_p(\sigma(T-\delta))\Big]^{\frac1{p-1}}=A_\delta,\]
		which is a contradiction.
	\end{proof}
	
	\vspace{3pt}
	
	\noindent\textbf{Full extension of $w$.}
	
	We first extend $w$ continuously to $\cup\gamma_i$. Fix $z\in\gamma_i$. Up to a translation and rotation, we may assume $z=0$ and $\D u(0)=\p_x$. Let $\rho$ be small enough so that $Q(2\rho)\Subset\Omega_0$, and
	\[|\D u-\p_x|<e^{-10},\quad|\D^\perp u-\p_y|<e^{-10}\qquad\text{in}\ \ Q(2\rho),\]
	and that $\gamma_i\cap Q(\rho)$ is the graph $y=f(x)$ of a function $f:(-\rho,\rho)\to(-\rho/10,\rho/10)$.
	
	For each $r<\rho$ and $z_1=(x_1,y_1)\in Q(r)\setminus\gamma_i$, we apply Lemma \ref{lemma-approx:weaker_item_9} to the vertical segment $\sigma=\big\{(x_1,y): y\in[f(x_1),y_1]\}$ oriented so that $\sigma(0)=(x_1,f(x_1))$. Note that $\metric{\D^\perp u}{\sigma'}$ has a definite sign everywhere on $\sigma$ (which is positive iff $y_1>f(x_1)$). Lemma \ref{lemma-approx:weaker_item_9} gives
	\[w(z_1)=-\log\max_\sigma|\D u|.\]
	Hence
	\[-\log\sup_{Q(r)}|\D u|\leq\inf_{Q(r)\setminus\gamma_i}(w)\leq\sup_{Q(r)\setminus\gamma_i}(w)\leq-\log\inf_{Q(r)}|\D u|.\]
	It follows that $w$ is continuous at $z=0$ if we set $w=-\log|\D u|$ on $\cup\gamma_i$. This argument applies to all $z\in\cup\gamma_i$, hence $w$ can be extended continuously to the set
	\[\Big(\bigcup_{D\in\cD_{\reg}}D\Big)\cup\big(\cup\gamma_i\big)\]
	by setting $w=-\log|\D u|$ on $\cup\gamma_i$.
	This set is exactly $\Omega_{\reg}$ as appeared in item \ref{item-approx:Omega_reg} and \eqref{eq-approx:def_omega_reg_2}. Finally, we define $w\equiv+\infty$ on $\Omega\setminus\Omega_{\reg}$. By the $p\to\infty$ limit of Lemma \ref{lemma-approx:nondeg}(i), it follows that $w$ is a continuous function from $\Omega$ to $\RR\cup\{+\infty\}$. The upper bound \eqref{eq-approx:w_upper_bound_partial} clearly extends and yields
	\begin{equation}\label{eq-approx:w_upper_bound}
		w\leq-\log|\D u|\qquad\text{on}\ \ \Omega.
	\end{equation}
	
	\vspace{3pt}
	\noindent\textbf{Production of simple cluster.}
	
	Set $\nu=-e^w\D^\perp u$, which is a continuous vector field on $\Omega_{\reg}=\{w<+\infty\}$. By \eqref{eq-approx:w_upper_bound} it follows that $|\nu|\leq1$. Recall that in each $D\in\cD_{\reg}$ the function $w$ arises as the limit
	
	\[w=\lim_{p\to\infty}w_p=\lim_{p\to\infty}\Big(\!-\frac{\log(\chi_Dv_p)}{p-1}\Big)\qquad\text{in}\ \ C^0_{\loc}(D).\]
	
	Recall from \cite{Moser_2007} that $w_p$ is a distributional solution of
	\begin{equation}\label{eq-approx:p-imcf}
		\div\big(|\D w_p|^{q-2}\D w_p\big)=|\D w_p|^q,\qquad q=\frac{p}{p-1}.
	\end{equation}
	Next, we may calculate
	\begin{equation}
		|\D w_p|^{q-2}\D w_p=-(p-1)^{1-q}|v_p|^{1-q}|\D v_p|^{q-2}\chi_D\D v_p=-\chi_D(p-1)^{-\frac1{p-1}}e^{w_p}\D^\perp u_p,
	\end{equation}
	which taking $p\to\infty$ yields
	\begin{equation}\label{eq-approx:convergence_of_vector}
		|\D w_p|^{q-2}\D w_p\to-\chi_D e^w\D^\perp u=\chi_D\nu\qquad\text{in}\ \ C^0_{\loc}(D).
	\end{equation}
	In particular,
	\[|\D w_p|^{q-1}\to|\nu|\qquad\text{in}\ \ C^0_{\loc}(D).\]
	Testing \eqref{eq-approx:p-imcf} with the function $\varphi(w-w_p)$, where $\varphi\in\Lip_{\loc}(D)$ is nonnegative, we have
	\[\begin{aligned}
		-\int\varphi(w-w_p)|\D w_p|^q &= \int|\D w_p|^{q-2}\bmetric{\D w_p}{(w-w_p)\D\varphi+\varphi\D w-\varphi\D w_p} \\
		&\leq \int|\D w_p|^{q-1}|w-w_p|\,|\D\varphi|+\varphi|\D w_p|^{q-1}|\D w|-\varphi|\D w_p|^q.
	\end{aligned}\]
	Recall that $|\D w_p|$ is locally uniformly bounded \cite{Kotschwar-Ni_2009}, and $w_p\to w$ in $C^0_{\loc}(D)$. Hence
	\[\limsup_{p\to\infty}\int\varphi|\D w_p|^q\leq\limsup_{p\to\infty}\int\varphi|\D w_p|^{q-1}|\D w|\leq\int\varphi|\D w|.\]
	But the standard lower semi-continuity gives
	\[\int\varphi|\D w|\leq\liminf_{p\to\infty}\int\varphi|\D w_p|^q.\]
	Hence,
	\begin{equation}\label{eq-approx:continuity_of_norm}
		\int\varphi|\D w|=\lim_{p\to\infty}\int\varphi|\D w_p|^q,\qquad\forall\,\varphi\in \Lip_0(D).
	\end{equation}
	
	Combining \eqref{eq-approx:p-imcf} \eqref{eq-approx:convergence_of_vector} \eqref{eq-approx:continuity_of_norm}, we have
	\[\int\metric{\chi_D\nu}{\D\varphi}=-\int\varphi|\D w|\qquad\forall\,\varphi\in\Lip_0(D).\]
	Testing this condition with $\varphi w$ for any $\varphi\in\Lip_0(D)$ gives
	\[\begin{aligned}
		-\int\varphi w|\D w| &= \int\metric{\chi_D\nu}{w\D\varphi+\varphi\D w} \\
		&= \int\varphi\metric{\chi_D\nu}{\D w}+\lim_{p\to\infty}\int\bmetric{|\D w_p|^{q-2}\D w_p}{w_p\D\varphi} \\
		&= \int\varphi\metric{\chi_D\nu}{\D w}+\lim_{p\to\infty}\int\Big[\!-\varphi w_p|\D w_p|^q-\varphi|\D w_p|^q\Big].
	\end{aligned}\]
	Then by \eqref{eq-approx:continuity_of_norm} (applied to both $\varphi$ and $\varphi w$) and the interior gradient estimate for $w_p$, we have
	\[\int\varphi\metric{\chi_D\nu}{\D w}=\int\varphi|\D w|.\]
	This implies $\metric{\chi_D\nu}{\D w}=|\D w|$ a.e.. In summary, we have shown that
	\begin{equation}\label{eq-approx:calibrated_imcf_produced}
		|\nu|\leq1,\qquad \metric{\chi_D\nu}{\D w}=|\D w|,\qquad \div(\chi_D\nu)=|\D w|.
	\end{equation}
	Hence $w|_D$ is a weak IMCF calibrated by $\chi_D\nu$, for each $D\in\cD_{\reg}$.
	
	We form the data $\big(w,\nu,\cD_{\reg},\{\chi_D\}\big)$ with the ambient domain $\Omega_{\reg}$. Recall that $w<\infty$ in $\Omega_{\reg}$, and each $D\in\cD_{\reg}$ is a $C^1$ domain in $\Omega_{\reg}$. It remains to verify that each $\gamma_i$ is a ridge (see Definition \ref{def-cluster:simple}(iv)). Fix a point $z\in\gamma_i$. Let $\rho$ be small enough so that up to a rigid motion, we have $z=0$ and $|\D u-\p_x|<e^{-10}$ in $Q(2\rho)$, and $(\cup\gamma_i)\cap Q(2\rho)=\graphh(f)$ for a function $f$. Let $D_1,D_2\in\cD_{\reg}$ be so that $D_1\cap Q(\rho)=\{y<f(x)\}$ and $D_2\cap Q(\rho)=\{y>f(x)\}$. Note that $|\D^\perp u-\p_y|<e^{-10}$, so $\D v_p$ points upwards, so $\chi_{D_1}=-1$ and $\chi_{D_2}=1$. Therefore, on $\gamma_i\cap Q(\rho)$ we have
	\[\begin{aligned}
		\nu\cdot\nu_{D_1} &= (-e^w\D^\perp u)\cdot\Big(\frac{\D^\perp u}{|\D u|}\Big)=-1=\chi_{D_1}, \\
		\nu\cdot\nu_{D_2} &= (-e^w\D^\perp u)\cdot\Big(\!-\frac{\D^\perp u}{|\D u|}\Big)=1=\chi_{D_2}.
	\end{aligned}\]
	All the conditions in Definition \ref{def-cluster:simple} are met, thus we obtain a simple cluster in $\Omega_{\reg}$.
	
	\vspace{3pt}
	\noindent\textbf{Summary.}
	
	Recall that we constructed a collection of curves $\gamma_i\subset\Omega_0$, and set $\cD$ to be the set of connected components of $\Omega\setminus\bar{\cup\gamma_i}$. We set $\cD_{\Crit},\cD_{\reg},\Omega_{\reg}$ as in \eqref{eq-approx:def_Dreg_Dcrit} \eqref{eq-approx:def_omega_reg_2}. Then, we constructed a continuous function $w:\Omega\to\RR\cup\{+\infty\}$ with $\Omega_{\reg}=\{w<+\infty\}$, by first defining it in $\cD_{\reg}$ then extending it to the entire $\Omega$. Moreover, we showed that $w_p\to w$ in $C^0_{\loc}(D)$, for each $D\in\cD_{\reg}$, along our chosen subsequence. Items \ref{item-approx:gamma_i} \ref{item-approx:Dreg_Dcrit} \ref{item-approx:Omega_reg} \ref{item-approx:simple_cluster} \ref{item-approx:w_on_gamma} \ref{item-approx:z0_on_gamma} all follow from the construction. In particular, (iia) comes from Lemma \ref{lemma-approx:nondeg}(i).
	
	To finish item \ref{item-approx:w_on_path}, observe that the case $z_2\in\bigcup_{D\in\cD_{\reg}}D$ follows from Lemma \ref{lemma-approx:weaker_item_9}. The case $z_2\notin\bigcup_{D\in\cD_{\reg}}D$ follows from
	\[w(z_2)=-\log|\D u(z_2)|\geq-\log\max_\sigma|\D u|\]
	and Lemma \ref{lemma-approx:weaker_item_9_2}.
	
	Item \ref{item-approx:w_on_z0} comes easily from the construction: if $z_0\in\Omega_{\reg}$, then the result follows from items \ref{item-approx:w_on_gamma} \ref{item-approx:z0_on_gamma}. If $z_0\notin\Omega_{\reg}$, then we have $z_0\in\Crit(u)$, hence $w(z_0)=+\infty=-\log|\D u|(z_0)$, recalling that we defined $w\equiv+\infty$ outside $\Omega_{\reg}$.
	
	To verify item \ref{item-approx:w_on_nonconst}, it suffices to show that $w=-\log|\D u|$ on $D\cap\{w\ne\const\}$ for each $D\in\cD_{\reg}$. Recall that $w$ is a weak IMCF in $D$ and is calibrated by the continuous vector field $\nu=-e^w\D^\perp u$. Hence, if $z\in\{w\ne\const\}$, then $|\nu(z)|=1$. This implies $e^{-w(z)}=|\D u(z)|$, namely, $w(z)=-\log|\D u(z)|$.
	
	To prove \ref{item-approx:z0_in_crit}, note that $z_0$ cannot lie on $\bar{\cup\gamma_i}$ (since $\bar{\cup\gamma_i}\subset\bar{\Omega_0}$), and $z_0$ cannot lie in any $D\in\cD_{\reg}$ (otherwise, Lemma \ref{lemma-approx:nondeg}(ii) would violate $v_p(z_0)=0$). So $z_0$ must lie in some $D\in\cD_{\Crit}$.
	
	Finally, let $Y_t,Z_t,\tZ_t,F,F'$ be as in item \ref{item-approx:lvlset}. Since $w\leq-\log|\D u|$, we always have $Y_t\supset Y'_t$, $\tZ_t\supset\tZ'_t$, $Z_t\supset Z'_t$, hence $F\supset F'$. It remains to prove the reverse inclusion. To this end, we may assume $F\ne\emptyset$. We first consider $F\subset Y_t$ or $F\subset Z_t$. For each $z\in F$, we may find a $C^1$ path $\sigma\Subset F$ joining $z_0$ to $z$ (since $F$ is open and connected). Then \eqref{eq-approx:main_w_on_path_1} implies
	\[w(z)\geq-\log\max_\sigma|\D u|\left\{\begin{aligned}
		& >t\qquad\text{(if $F\subset Y_t$)}, \\
		& \geq t\qquad\text{(if $F\subset Z_t$)}. \\
	\end{aligned}\right.\]
	For $F\subset Y_t$, this immediately implies $z\in F'$. For $F\subset Z_t$, this implies $z\in\tZ'_t$ for all $z\in F$, hence $F\subset\Int(\tZ'_t)=Z'_t$. Hence $F\subset F'$, as desired. Finally, consider $F\subset\tZ_t$ and $F'\subset\tZ'_t$. For each $s<t$, let $F_s$ (resp. $F'_s$) be the connected component of $Y_s$ (resp. $Y'_s$) that contains $z_0$. We have proved that $F_s=F'_s$. So for each $s<t$ we have
	\[F\subset\tZ_t\subset Y_s\quad\Rightarrow\quad F\subset F_s\quad\Rightarrow\quad F\subset F'_s\quad\Rightarrow\quad F\subset Y'_s.\]
	Taking $s\nearrow t$, it follows that $F\subset\tZ'_t$, hence $F\subset F'$.
\end{proof}

%\newpage

\section{Simple clusters and heat equation}\label{sec:heat}

Suppose that a smooth family of convex curves $\{\gamma_t\}$ evolve at the speed of $1/\kappa$. Let $h=h(t,\th)$ be the collection of support functions of $\gamma_t$. Then we have
\begin{equation}\label{eq-heat:intro}
	\frac{\p h}{\p t}=\metric{\frac{\p\gamma_t}{\p t}}{\nu_\th}=\kappa^{-1}=\frac{\p^2h}{\p\th^2}+h.
\end{equation}
See Subsection \ref{subsec:imcf}. In Section \ref{sec:ex}, we also built examples of simple clusters where \eqref{eq-heat:intro} holds across the interfaces. The aim of this section is to show that $\p_th=\p_{\th\th}h+h$ holds in a certain sense for weak IMCFs and simple IMCF clusters. Let $w$ be either a weak IMCF or a simple cluster. We fix as usual
\[Y_t=\{w>t\},\qquad \tZ_t=\{w\geq t\},\qquad Z_t=\Int(\tZ_t).\]
If $w$ is a weak IMCF, then $Y_t,Z_t$ are mean concave and $\tZ_t=\bar{Z_t}$ for all $t$. When $w$ is a simple cluster, $\tZ_t$ may differ from $\bar{Z_t}$ by line segments (see Lemma \ref{lemma-cluster:local_model}). Denote
\begin{equation}\label{eq-heat:intro_Gamma}
	\gamma_t=\p Z_t\cap\Omega,\qquad\qquad\Gamma=\bigcup_{t\in\RR}\gamma_t.
\end{equation}
If $w$ is a simple cluster, then it is naturally associated with a continuous vector field $\nu$, with $|\nu|=1$ on $\Gamma$. If $w$ is a weak IMCF, then the collection of inward unit normals of $\p Z_t$ (i.e. outer unit normals of $\p\{w<t\}$), still denoted by $\nu$, is continuous (in fact, $C^{0,1/2}$ by Lemma \ref{lemma-prelim:joint_cont_of_normals}) on $\Gamma$ as well. In either case, we obtain a map
\[\nu:\Gamma\to\SS^1\]
which is continuous and orthogonal to each $\gamma_t$.

The main statement is as below. Recall the notation
\[\nu_\th=(\sin\th,-\cos\th)=e^{i(\th-\pi/2)},\qquad\tau_\th=(\cos\th,\sin\th)=e^{i\th}.\]
Recall our convention that a ``line segment'' means a line segment or ray or line or point, and a ``finite line segment'' means a line segment with finite length, and a ``nontrivial line segment'' means a line segment that is not a point.

\begin{theorem}\label{thm-heat:main}
	Suppose $w$ is a weak IMCF or a simple IMCF cluster in a domain $\Omega$. Let $\gamma_t,\Gamma,\nu$ be as defined above. Assume that the map $(t,\nu^\perp):\Gamma\to\RR\times\SS^1$ lifts to a continuous map
	\[\Th:\Gamma\to\RR\times M,\]
	with $M=\RR$ or $M=\SS^1(2k\pi)$, so that each preimage $\Th^{-1}(t,\th)$ is a line segment. Denote $\cJ=\Th(\Gamma)$. Let $h:\cJ\to\RR$ be the (well-defined) support function such that
	\begin{equation}\label{eq-heat:intro_def_h}
		h(t,\th)=\metric{z}{\nu}=\metric{z}{\nu_\th}\qquad\text{for any }z\in\Th^{-1}(t,\th).
	\end{equation}
	For each $(t,\th)\in\cJ$, denote $\sigma=\Th^{-1}(t,\th)$. Then:
	\begin{enumerate}[label={(\roman*)}, nosep]
		\item\label{item-heat:interior_viscosity} If $\sigma\Subset\Omega$, then for any $\varphi\in C^2(R)$ with $R:=[t-\delta,t]\times[\th-\delta,\th+\delta]$ for some $\delta>0$, so that
		\[\varphi(t,\th)=h(t,\th),\qquad\varphi\geq h\ \ (\text{resp. $\varphi\leq h$})\ \ \text{in}\ \ R\cap\cJ,\]
		it holds
		\[\big(\p_t\varphi-\p_{\th\th}\varphi-\varphi\big)(t,\th)\leq0\ \ (\text{resp. $\geq0$}).\]
		
		\item\label{item-heat:boundary_viscosity} Suppose $\sigma\not\Subset\Omega$. Then $\sigma$ is a nontrivial line segment. Let $\chi_\sigma\in\{\pm1\}$ be the (well-defined) sign so that $-\chi_\sigma\nu$ is the outer unit normal of $Z_t$ on $\Int(\sigma)$. Then for any $\varphi\in C^2(R)$ with $R:=[t-\delta,t]\times[\th-\delta,\th+\delta]$ for some $\delta>0$, so that
		\begin{equation}\label{eq-heat:intro_touching_2}
			\varphi(t,\th)=h(t,\th),\qquad\left\{\begin{aligned}
				& \varphi\geq h\ \ \text{in}\ \ R\cap\cJ\qquad(\text{if $\chi_\sigma=-1$}), \\
				& \varphi\leq h\ \ \text{in}\ \ R\cap\cJ\qquad(\text{if $\chi_\sigma=1$}),
			\end{aligned}\right.
		\end{equation}
		and
		\begin{equation}\label{eq-heat:intro_on_sigma}
			\varphi(t,\th)\nu_\th+\p_\th\varphi(t,\th)\tau_\th\in\sigma,
		\end{equation}
		it holds
		\begin{equation}\label{eq-heat:intro_viscosity_result}
			\big(\p_t\varphi-\p_{\th\th}\varphi-\varphi\big)(t,\th)\left\{\begin{aligned}
				& \leq0\qquad(\text{if $\chi_\sigma=-1$}), \\
				& \geq0\qquad(\text{if $\chi_\sigma=1$}).
			\end{aligned}\right.
		\end{equation}
	\end{enumerate}
\end{theorem}

\begin{figure}[ht]
	\centering
	\includegraphics{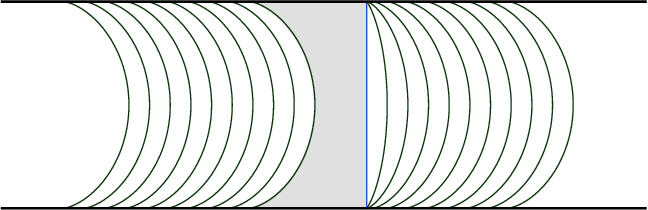}
	\caption{A weak IMCF for which $h(t,\th)$ is not continuous.}\label{fig-heat:noncontinuous}
\end{figure}

Whenever $\cJ$ contains a box $R=(t_1,t_2]\times(\th_1,\th_2)$ so that $\Th^{-1}(R)\Subset\Omega$ and $h$ is continuous in $R$, item \ref{item-heat:interior_viscosity} and the standard viscosity theory \cite{Crandall-Ishii-Lions_1992, Ishii_1995} imply that $h$ is a smooth solution of $\p_th=\p_{\th\th}h+h$ in $R$. However, $h$ need not be continuous in general. Figure \ref{fig-heat:noncontinuous} shows an example where the solution first evolves as a translating soliton, then jumps over the shadowed region at $t=0$, then continues evolving with the outer obstacle (i.e. boundary tangency) condition on $\{y=\pm1\}$. Then $\cJ=\RR\times(0,\pi)$, and $h$ is not continuous across $\{t=0\}$.

\begin{figure}[ht]
	\centering
	\includegraphics{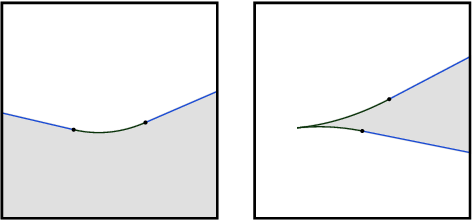}
	\begin{picture}(0,0)
		\put(-227,90){\arrowangle{-90}}
		\put(-220,82){$\nu$}
		\put(-105,90){\arrowangle{-90}}
		\put(-98,82){$\nu$}
		\put(-228,56){$\chi_\sigma=1$}
		\put(-164,64){$\chi_\sigma=1$}
		\put(-2,78){$\chi_\sigma=1$}
		\put(-2,30){$\chi_\sigma=-1$}
		\put(-227,6){$Z_t$}
		\put(-20,50){$Z_t$}
	\end{picture}
	\caption{The orientation $\chi_\sigma$.}\label{fig-heat:chi_sigma}
\end{figure}

Item \ref{item-heat:boundary_viscosity} is used in Theorem \ref{thm-heat:no_interior_segment} below as well as in the $\epsilon$-regularity Theorem \ref{thm-ereg:main_cluster}. We refer to Section \ref{sec:ereg} for related discussions. Figure \ref{fig-heat:chi_sigma} below (in which $\nu\approx-\p_y$) shows how the sign $\chi_\sigma$ is determined.

Theorem \ref{thm-heat:main} is proved using a geometric comparison: if the stated property fails, then the test function $\varphi$ would generate a subsolution or supersolution $\{\tilde\gamma_t\}$ of IMCF. Then $\varphi\leq h$ or $\varphi\geq h$ would allow us to compare the position of $\gamma_t$ and $\tilde\gamma_t$. The maximum principle for weak IMCF finally yields a contradiction. We do not know whether item \ref{item-heat:boundary_viscosity} can be extended to a full comparison, but an obstruction for us proving a full result is explained at around Figure \ref{fig-heat:1b_res}.

\begin{remark}\label{rmk-heat:misc}
	The theorem respects the following symmetries:
	\begin{enumerate}[label={(\roman*)}, nosep]
		\item Rotation and time translation: this corresponds to translation of the objects given by $\cJ\to\cJ+(t_0,\th_0)$ and $h(t,\th)\to h(t-t_0,\th-\th_0)$, $\varphi(t,\th)\to\varphi(t-t_0,\th-\th_0)$.
		\item\label{item-heat:sym_transl} Translation in $\RR^2$: this corresponds to adding $h,\varphi$ by a common factor of $a\sin\th-b\cos\th$, for some $a,b\in\RR$, while leaving $\cJ$ fixed.
		\item Orientation reversing: by Remark \ref{rmk-cluster:simple}(v) this means $\nu'=-\nu$ and $\chi_i'=-\chi_i$ for all regions $D_i$ (where a notation with a prime is associated to the cluster after orientation reversing). The effect on angle parametrization is
		\[\Th'(z)=\Th(z)+(0,\pi),\qquad(\Th')^{-1}(t,\th)=\Th^{-1}(t,\th-\pi),\qquad\cJ'=\cJ+(0,\pi).\]
		By definition we have $h'(t,\th)=\metric{z}{\nu_\th}$ for any $z\in(\Th')^{-1}(t,\th)$, thus
		\[\begin{aligned}
			h'(t,\th) &= \metric{z}{\nu_\th}=-\metric{z}{\nu_{\th-\pi}}\qquad\forall\,z\in\Th^{-1}(t,\th-\pi) \\
			&= -h(t,\th-\pi).
		\end{aligned}\]
		The test function $\varphi$ should be transformed accordingly, so $\varphi'(t,\th)=-\varphi(t,\th-\pi)$. The order between $h$ and $\varphi$ is reversed, so is the sign of $\p_t\varphi-\p_{\th\th}\varphi-\varphi$. In the context of item (ii), note that
		\[\begin{aligned}
			\varphi'(t,\th)\nu_{\th}+\p_\th\varphi'(t,\th)\tau_\th &= \varphi(t,\th-\pi)\nu_{\th-\pi}+\p_\th\varphi(t,\th-\pi)\tau_{\th-\pi} \\
			&\in \Th^{-1}(t,\th-\pi)=(\Th')^{-1}(t,\th),
		\end{aligned}\]
		so \eqref{eq-heat:intro_on_sigma} is preserved. Finally, the orientation becomes $\chi'_\sigma=-\chi_\sigma$, so the compatibility between \eqref{eq-heat:intro_touching_2} and \eqref{eq-heat:intro_viscosity_result} is preserved as well.
		\item\label{item-heat:sym_orient_rot} If we reverse the orientation and then rotate by $180^\circ$: this results in sending $\chi_i\to-\chi_i$, $h\to-h$, $\varphi\to-\varphi$ while not changing $\nu$ and $\cJ$.
		\item\label{item-heat:sym_y_refl} Reflection around the $y$-axis: this fixes $\nu=\p_y$ hence $\th=0$. Therefore, it corresponds to $\th\to-\th$ and $\cJ$ being reflected with respect to $\{\th=0\}$, and $h(t,\th)\to h(t,-\th)$, $\varphi(t,\th)\to\varphi(t,-\th)$.
		\item\label{item-heat:sym_x_refl} Reflection around the $x$-axis: this corresponds to $\th\to\pi-\th$ and $\cJ$ being reflected with respect to $\{\th=\pi/2\}$, and $h(t,\th)\to h(t,\pi-\th)$, $\varphi(t,\th)\to\varphi(t,\pi-\th)$.
		\item\label{item-heat:sym_x_refl_orient} If we reflect around the $x$-axis and then reverse the orientation, then this results in $\chi_i\to-\chi_i$, and $\cJ$ being reflected with respect to $\th=0$, and $h(t,\th)\to -h(t,-\th)$, $\varphi(t,\th)\to-\varphi(t,-\th)$.
	\end{enumerate}
\end{remark}

We prove the following fine regularity for planar weak IMCF, as a consequence of Theorem \ref{thm-heat:main}. In the theorem and its proof, we denote
\[Q_\delta(r)=(-r,r)\times(-\delta r,\delta r).\]

\begin{theorem}\label{thm-heat:no_interior_segment} {\ }
	
	\begin{enumerate}[label={(\roman*)}, nosep]
		\item\label{item-heat:curve_decomp} Let $w$ be a weak IMCF in $\Omega\subset\RR^2$, and $\gamma$ be a connected component of $\p Z_t\cap\Omega$. Then either $\gamma$ is a nontrivial line segment, or we can decompose
		\[\gamma=\gamma_1\cup\gamma_2\cup\gamma_3,\]
		where $\gamma_i$ meets $\gamma_{i+1}$ at endpoints, and $\gamma_1,\gamma_3$ are (possibly empty) line segments, and the interior of $\gamma_2$ is smooth, nonempty, and strictly convex.
		\item\label{item-heat:grad_est} Suppose $\delta\leq1$, and $w$ is a weak IMCF in $Q_\delta(4)$ such that
		\begin{equation}\label{eq-heat:intro_ereg}
			|\nu+\p_y|<e^{-10}\delta\ \ \text{in}\ \ \Gamma,\qquad|w|< e^{-20}\delta^2.
		\end{equation}
		Then
		\begin{equation}\label{eq-heat:intro_curv_est}
			\sup_{Q_\delta(2)}|\D w|\leq e^{-7}\delta.
		\end{equation}
	\end{enumerate}
\end{theorem}

We remark that the gradient estimate for the weak IMCF remains open in general (though weak IMCFs are a priori $\Lip_{\loc}$ by definition), namely, it is unknown whether a weak IMCF $w$ in a domain $(\Omega,g)$ satisfies the uniform bound
\begin{equation}\label{eq-heat:conj_grad_est}
	\sup_K|\D w|\leq C(K,\Omega,g)\qquad\forall\,K\Subset\Omega.
\end{equation}
Such an estimate is only known if $w$ is the limit of $\epsilon$-regularized solutions \cite{Huisken-Ilmanen_2001} or renormalized $q$-harmonic functions \cite{Moser_2007, Kotschwar-Ni_2009, Mari-Rigoli-Setti_2022}. Item \ref{item-heat:grad_est} confirms a special case in $\RR^2$.

\begin{proof}[Proof of Theorem \ref{thm-heat:main}] {\ }
	
	If $w$ is a simple cluster, then let $\{D_i\}, \{\chi_i\}, \{\gamma_j\}$ be the associated regions and orientations and ridges. If $w$ is a weak IMCF, then we set $D_1=\Omega$ and $\chi_1=1$ and $\{\gamma_j\}=\emptyset$.
	
	By rotation and time translation symmetry, it suffices to assume that the basepoint is $(t,\th)=(0,0)$. With a translation in $\RR^2$, see Remark \ref{rmk-heat:misc}\ref{item-heat:sym_transl}, we may assume
	\begin{equation}\label{eq-heat:translation_sym}
		h(0,0)=\varphi(0,0)=\p_\th\varphi(0,0)=0
	\end{equation}
	for any given test function $\varphi$. If the viscosity subsolution (or supersolution) property fails, then there is a test function $\varphi$ defined in a backward neighborhood of $(0,0)$, so that
	\begin{equation}\label{eq-heat:varphi_1}
		h(0,0)=\varphi(0,0)=\p_\th\varphi(0,0)=0,\quad\varphi\geq h\text{\ \ (or $\leq h$, depending on the context)},
	\end{equation}
	but
	\begin{equation}\label{eq-heat:varphi_2}
		\big(\p_t\varphi-\p_{\th\th}\varphi-\varphi\big)(0,0)>0\text{\ \ (or $<0$).}
	\end{equation}
	In the context of item \ref{item-heat:boundary_viscosity}, we also have
	\begin{equation}\label{eq-heat:varphi_3}
		0=\varphi(0,0)\nu_0+\p_\th\varphi(0,0)\tau_0\in\sigma.
	\end{equation}
	Since \eqref{eq-heat:varphi_2} is an open condition, we may replace $\varphi$ by $\varphi\pm\epsilon e^t(\th^2-t)$ and restrict to a smaller neighborhood, thus assume that $\varphi\in C^2\big([-\delta,0]\times[-\delta,\delta]\big)$ for some $\delta>0$, and
	\[\varphi(0,0)=h(0,0),\qquad\varphi>h\text{\ \ (or $<h$)\ \ in\ \ $\big([-\delta,0]\times[-\delta,\delta]\big)\cap \cJ\setminus\{(0,0)\}$,}\]
	but
	\[\p_t\varphi-\p_{\th\th}\varphi-\varphi>0\text{\ \ (or $<0$)\ \ in\ \ $[-\delta,0]\times[-\delta,\delta]$.}\]
	For convenience, we also extend $\varphi$ to positive time: by suitably arranging the extension and further decreasing $\delta$, we may assume that $\varphi\in C^2\big([-\delta,\delta]\times[-\delta,\delta]\big)$, and
	\begin{equation}\label{eq-heat:viscosity_root}
		\left\{\begin{aligned}
			& h(0,0)=\varphi(0,0)=\p_\th\varphi(0,0)=0, \\
			& \varphi>h\text{\ \ (or $<h$)\ \ in\ \ $\big([-\delta,0]\times[-\delta,\delta]\big)\cap \cJ\setminus\{(0,0)\}$,} \\
			& \p_t\varphi-\p_{\th\th}\varphi-\varphi>0\text{\ \ (or $<0$)\ \ in\ \ $[-\delta,\delta]\times[-\delta,\delta]$.}
		\end{aligned}\right.
	\end{equation}
	In the context of item \ref{item-heat:boundary_viscosity}, adding $\pm\epsilon e^t(\th^2-t)$ to $\varphi$ does not affect the validity of \eqref{eq-heat:varphi_3}. Therefore, to prove the theorem, it suffices to show that $\nexists\varphi$ so that \eqref{eq-heat:viscosity_root} holds if in item \ref{item-heat:interior_viscosity}, or \eqref{eq-heat:viscosity_root} \eqref{eq-heat:varphi_3} hold if in item \ref{item-heat:boundary_viscosity} with compatible touching conditions.
	
	\vspace{3pt}
	
	If $w$ is a simple cluster, then $Z_t$ has the local structure given by Lemma \ref{lemma-cluster:local_model}\ref{item-cluster:Yt_Zt_model} (applied with $\delta=1$ there): at each $z\in\cup\gamma_j$, after a rigid motion that sends $z\mapsto0$ and $\nu(z)\mapsto-\p_y$, there is $r>0$ so that
	\[Z_t\cap Q(r)=\Big\{x\in I, g_1(x)<y<g_2(x)\Big\},\]
	with $I=(-r,r)$ or $I=(-r,a)\cup(b,r)$ with $-r\leq a\leq b\leq r$, and $g_1,g_2$ have the properties as in Lemma \ref{lemma-cluster:local_model}\ref{item-cluster:model}. The lifting condition (i.e. each $\Th^{-1}(t,\th)$ is a segment) implies:
	
	\begin{lemma}\label{lemma-heat:no_model_ic}
		It cannot happen that $-r<a=b<r$. Hence, each $Z_t$ is a disjoint union of concave polygons by the definition in Subsection \ref{subsec:curves}.
	\end{lemma}
	\begin{proof}
		Suppose $-r<a=b<r$ occurs for some $z$. We may restrict to a smaller closed square with a different center, thus assume that $a=b=0$. After the recentering, we have
		\[Z_t\cap\bar{Q(r)}=\Big\{x\in[-r,r]\setminus\{0\}, g_1(x)<y<g_2(x)\Big\},\]
		where $g_1$ (resp. $g_2$) is a $e^{-9}$-Lipschitz, concave (resp. convex) function on $[-r,r]$, with $g_1(0)=g_2(0)$ and $g'_1(0)=g'_2(0)$, and $g_1<g_2$ on $[-r,r]\setminus\{0\}$, by Lemma \ref{lemma-cluster:local_model}\ref{item-cluster:g1_g2}.
		
		Denote $\th_0=g'_1(0)=g'_2(0)$. Consider
		\[I_1=\big\{g'_1(x):x\in[-r,r]\big\},\qquad I_2=\{g'_2(x):x\in[-r,r]\}.\]
		So $I_1,I_2$ are both closed intervals containing $\th_0$. The lifting condition implies $I_1\cap I_2=\{\th_0\}$: indeed, if $\th\in(I_1\cap I_2)\setminus\{\th_0\}$, then $\Th^{-1}(t,\arctan\th)\cap\bar{Q(r)}$ is not connected. Also, we have $g'_1(r)<g'_2(r)$ and $g'_1(-r)>g'_2(-r)$, since $g_1(\pm r)<g_2(\pm r)$. This implies either
		\[I_1=\{\th_0\},\qquad I_2=[c,d]\ \ (c<\th_0<d)\]
		or
		\[I_1=[c,d]\ \ (c<\th_0<d),\qquad I_2=\{\th_0\}.\]
		For each $s<t$ close enough to $t$, there is a $C^1$ concave function $g_{1,s}$ and convex function $g_{2,s}$ with $g_{1,s}<g_1\leq g_2<g_{2,s}$, so that
		\[Z_s\cap\bar{Q(r)}=\Big\{x\in[-r,r],\,g_{1,s}(x)<y<g_{2,s}(x)\Big\},\]
		and
		\[g_{1,s}\xrightarrow[s\nearrow t]{C^1([-r,r])}g_1,\qquad g_{2,s}\xrightarrow[s\nearrow t]{C^1([-r,r])}g_2.\]
		Consider the intervals $I_{1,s}=\big\{g'_{1,s}(x):x\in[-r,r]\big\}$, $I_{2,s}=\big\{g'_{2,s}(x):x\in[-r,r]\big\}$, so $I_{1,s}\to I_1$ and $I_{2,s}\to I_2$ as $s\nearrow t$. Hence $I_{1,s}\cap I_{2,s}\ne\emptyset$ for all $s$ close to $t$. Then $\Th^{-1}(s,\arctan\th)$ is not connected whenever $\th\in I_{1,s}\cap I_{2,s}$, contradicting our lifting hypothesis.
	\end{proof}
	
	\begin{figure}[ht]
		\centering
		\includegraphics{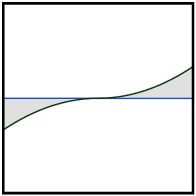}
		\caption{Impossible case (blue segment represents $\sigma$, and $\chi_\sigma$ is not defined).}\label{fig-heat:no_model_ic}
	\end{figure}
	
	Now since $Z_t$ is a disjoint union of concave polygons, the sign $\chi_\sigma$ in item \ref{item-heat:boundary_viscosity} is well-defined (we have ruled out cases such as Figure \ref{fig-heat:no_model_ic}). The following two lemmas are useful:
	
	\begin{lemma}\label{lemma-heat:interior_angle}
		$\Th^{-1}(t,\th)\Subset\Omega$ implies $\th\in\Int(J_t)$, where $J_t=\big((t,\th)\mapsto\th\big)(\Th(\gamma_t))$.
	\end{lemma}
	\begin{proof}
		The shape of $Z_t$ in a neighborhood of $\Th^{-1}(t,\th)$ falls into one of the cases in Figure \ref{fig-heat:int_lemma} ($\Th^{-1}(t,\th)$ is displayed as blue segments). Note that $\th\in\Int(J_t)$ in either case.
		\begin{figure}[ht]
			\centering
			\includegraphics{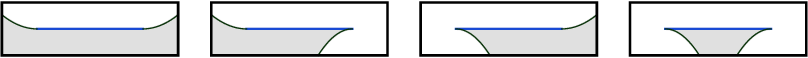}
			\caption{Local models in Lemma \ref{lemma-heat:interior_angle}.}\label{fig-heat:int_lemma}
		\end{figure}
	\end{proof}
	
	\begin{lemma}\label{lemma-heat:0_in_gamma0}
		The existence of test function $\varphi$ implies $0\in\gamma_0$, $\nu(0)=-\p_y$, $\Th(0)=(0,0)$.
	\end{lemma}
	\begin{proof}
		Recall that we have reduced to the case $h(0,0)=\varphi(0,0)=\p_\th\varphi(0,0)=0$. In the context of item \ref{item-heat:boundary_viscosity}, our assumption \eqref{eq-heat:intro_on_sigma} implies
		\[0=\varphi(0,0)\nu_0+\p_\th\varphi(0,0)\tau_0\in\sigma\subset\gamma_0.\]
		Then $\sigma=\Th^{-1}(0,0)$ implies $\Th(0)\in\Th(\sigma)=(0,0)$. For item (i), Lemma \ref{lemma-heat:interior_angle} implies that $h(0,\cdot)$ is defined in a neighborhood of $\th=0$. Then $\varphi(0,0)=h(0,0)$ and $\varphi\leq h$ or $\varphi\geq h$ implies
		\[0=\p_\th\varphi(0,0)\in[\p_\th^- h(0,0),\p_\th^+ h(0,0)],\]
		where $[a,b]=[b,a]$ if $a>b$. On the other hand, \eqref{eq-prelim:recovery_formula_polygon} implies
		\[\begin{aligned}
			\Th^{-1}(0,0) &=\Big\{h(0,0)\nu_0+c\tau_0:c\in[\p_\th^- h(0,0),\p_\th^+ h(0,0)]\Big\} \\
			&= [\p_\th^- h(0,0),\p_\th^+ h(0,0)]\times\{0\}.
		\end{aligned}\]
		Hence $0\in\Th^{-1}(0,0)\subset\gamma_0$. Finally, $\Th(0)=(0,0)$ always implies $\nu(0)=-\p_y$.
	\end{proof}
	
	Depending on whether $0$ lies in a region $D_i$ or on a ridge, and whether $\varphi$ touches $h$ from above or below, the arguments are different. Hence, a case discussion is needed.
	
	\vspace{3pt}
	
	\textbf{Case 1:} $w$ is a weak IMCF, or $0$ lies in some region $D_i$ in the cluster. So $w$ is a weak IMCF in a small neighborhood of $0$. (Warning: we cannot shrink $\Omega$ and assume that $w$ is a weak IMCF in $\Omega$, since the condition $\Th^{-1}(0,0)\Subset\Omega$ in the theorem is not local.) We may reverse the orientation as in Remark \ref{rmk-heat:misc}\ref{item-heat:sym_orient_rot} and assume $\chi_i=1$. Then observe that $\nu$ is the inner unit normal of $\p Z_t$, so in the context of item \ref{item-heat:boundary_viscosity}, we have $\chi_\sigma=1$ and hence the suitable touching condition is $\varphi\leq h$.
	
	Recall $0\in\gamma_0$ and $\nu(0)=-\p_y$, thus the IMCF is downward evolving near $0$. For $r>0$, consider the region
	\[P=[-r,r]\times[-r,r].\]
	Fix $r\ll1$ so that $P\Subset\Omega$, and $w$ is a weak IMCF in a neighborhood of $P$, and
	\begin{equation}\label{eq-heat:almost_horizontal_1}
		|\nu+\p_y|<e^{-10},\qquad|\Th|<e^{-10}\qquad\text{in}\ \ P\cap\Gamma.
	\end{equation}
	Recall that $\gamma_t\to\gamma_0$ graphically in $C^1$ as $t\nearrow0$. Thus, there is $T<0$ so that $\gamma_t\cap P$ is the graph of a $C^{1,1}$ convex function $g_t:[-r,r]\to[-r/2,r/2]$ for each $t\in[T,0]$. Also, we have
	\begin{equation}\label{eq-heat:ft_into}
		g_t\xrightarrow[t\nearrow0]{C^1([-r,r])}g_0,\qquad g_0(0)=g'_0(0)=0.
	\end{equation}
	Denote
	\[\uth_t=\arctan g'_t(-r),\qquad \oth_t=\arctan g'_t(r).\]
	These are the angles of $\gamma_t$ at their intersections with $\p P$. The angle parameter $\Th(\gamma_t\cap P)$ may take range in $[\uth_t+2k_t\pi,\oth_t+2k_t\pi]$ for some $k_t\in\ZZ$, but \eqref{eq-heat:almost_horizontal_1} forces $k_t=0$ for all $t\in[T,0]$. Hence
	\begin{equation}\label{eq-heat:def_uth_oth}
		\Th(\gamma_t\cap P)=[\uth_t,\oth_t]\qquad\forall\,t\in[T,0],
	\end{equation}
	with
	\begin{equation}\label{eq-heat:Th_arctan}
		\Th\big((x,g_t(x))\big)=\arctan g'_t(x),\qquad\forall\,t\in[T,0].
	\end{equation}
	
	The following lemma shows the joint continuity of $g_t$.
	
	% == Latex note ==
	% Multiple Claims in a section. Distinct environment name is defined.
	
	\begin{claimheatmain}\label{claim-heat:1_upper_semi_cont}
		The map $(t,x)\mapsto g_t(x)$ is upper semi-continuous in $[T,0]\times[-r,r]$. If $T'\in[T,0]$ and $\rho\in[0,r]$ satisfy
		\[\Int\big(\{w=t\}\big)\cap\{-\rho\leq x\leq\rho\}=\emptyset\qquad\forall\,t\in[T',0],\]
		then $(t,x)\mapsto g_t(x)$ is continuous in $[T',0]\times[-\rho,\rho]$.
	\end{claimheatmain}
	\begin{proof}
		Note that $t\mapsto g_t(x)$ is decreasing for each $x\in[-r,r]$, and
		\begin{equation}\label{eq-heat:aux3}
			\lim_{s\nearrow t}g_s=g_t,\qquad\lim_{s\searrow t}g_s=g_t^+\qquad\text{in}\ \ C^1([-r,r]),
		\end{equation}
		where $g_t^+$ is such that $\graphh(g_t^+)=\p\{w\leq t\}\cap P$. Since $w$ is downward evolving, we have $g_t^+\leq g_t$ in $[-r,r]$. This shows upper semi-continuity. If there are $T',\rho$ as stated, then we have $g_t=g_t^+$ in $[-\rho,\rho]$ for all $t\in[T',0]$, so the continuity follows from \eqref{eq-heat:aux3}.
	\end{proof}
	
	With the above setups, we are ready to prove the viscosity properties.
	
	\vspace{3pt}
	
	\textbf{Case 1(a).} Supersolution. Suppose otherwise that there is $\delta>0$ and $\varphi\in C^2\big([-\delta,\delta]\times[-\delta,\delta]\big)$, so that
	\begin{equation}\label{eq-heat:1a_viscosity_prep}
		\left\{\begin{aligned}
			& h(0,0)=\varphi(0,0)=\p_\th\varphi(0,0)=0, \\
			& \varphi<h\text{\ \ in\ \ $\big([-\delta,0]\times[-\delta,\delta]\big)\cap \cJ\setminus\{(0,0)\}$,} \\
			& \p_t\varphi-\p_{\th\th}\varphi-\varphi<0\text{\ \ in\ \ $[-\delta,\delta]\times[-\delta,\delta]$.}
		\end{aligned}\right.
	\end{equation}
	As this case appears in both items \ref{item-heat:interior_viscosity} and \ref{item-heat:boundary_viscosity}, we have $0\in\gamma_0$ by Lemma \ref{lemma-heat:0_in_gamma0}, but we are not guaranteed $\Th^{-1}(0,0)\Subset\Omega$. The touching condition implies:
	
	\begin{claimheatmain}\label{claim-heat:1a_nondeg}
		$\p_t\varphi(0,0)\geq0$.
	\end{claimheatmain}
	
	The proof is postponed to the end of this case. Once Claim \ref{claim-heat:1a_nondeg} is shown, we may add a small positive multiple of $t$ to $\varphi$, thus assume that
	\begin{equation}\label{eq-heat:1a_viscosity_root}
		\left\{\begin{aligned}
			& h(0,0)=\varphi(0,0)=\p_\th\varphi(0,0)=0, \\
			& \varphi<h\text{\ \ in\ \ }\big([-\delta,0]\times[-\delta,\delta]\big)\cap\cJ\setminus\{(0,0)\}, \\
			& 0<\p_t\varphi<\p_{\th\th}\varphi+\varphi\text{\ \ in\ \ }[-\delta,\delta]\times[-\delta,\delta].
		\end{aligned}\right.
	\end{equation}
	Consider the comparison curves
	\begin{equation}\label{eq-heat:comp_curves}
		\tilde\gamma_t=\Big\{\varphi(t,\th)\nu_\th+\p_\th\varphi(t,\th)\tau_\th:\th\in(-\delta,\delta)\Big\},\qquad\forall\,t\in(-\delta,\delta),
	\end{equation}
	and the region that they sweep out
	\begin{equation}\label{eq-heat:def_S}
		S=\bigcup_{t\in(-\delta,\delta)}\tilde\gamma_t.
	\end{equation}
	We may further decrease $\delta$, so that $\delta<e^{-10}$ and $S\Subset\Int(P)$. Now $0<\p_t\varphi<\p_{\th\th}\varphi+\varphi$ implies that $\{\tilde\gamma_t\}$ forms a downward evolving smooth supersolution of IMCF in $S$ (note: in a supersolution of IMCF, the normal speed is smaller than $1/\kappa$).
	
	Choose a small $\rho<r$ so that $Q(\rho)\Subset S$. Then find $T'\in[T,0)$ so that $\tilde\gamma_t\cap\bar{Q(\rho)}$, for each $t\in[T',0]$, is the graph of a strictly convex function $\tg_t:[-\rho,\rho]\to[-\rho/10,\rho/10]$. Note that $\tg_0(0)=\tg'_0(0)=0$ since $\varphi(0,0)=\p_\th\varphi(0,0)=0$. By strict convexity and continuity, we may further increase $T'$ and find a constant $\mu>0$, so that
	\begin{equation}\label{eq-heat:1a_2mu_angle}
		\tg'_t(-\rho)<-\tan(2\mu),\quad \tg'_t(\rho)>\tan(2\mu),\qquad\forall\,t\in[T',0].
	\end{equation}
	This means that the angle parameters of $\tilde\gamma_t\cap\bar{Q(\rho)}$ contain $[-2\mu,2\mu]$.
	
	Finally, choose $\rho'<\rho$ so that $|\Th|<\mu$ in $\bar{Q(\rho')}$. Consider $\gamma'_t=\gamma_t\cap Q(\rho')$. We need to compare the position of $\gamma'_t$ and $\tilde\gamma_t\cap\bar{Q(\rho)}=\graphh(\tg_t)$ (the domains of these two graphs are different). A convenient statement is as follows (see Figure \ref{fig-heat:curve_comp} for illustration):
	
	\begin{claimheatmain}\label{claim-heat:1a_curve_comp}
		Let $\tsg_t:\RR\to\RR$ be the convex function obtained by extending $\tg_t$ in the $C^1$ manner by gluing two linear functions on $(-\infty,-\rho]$ and $[\rho,\infty)$. Then
		\[\gamma'_t\subset\big\{y<\tsg_t(x)\big\}\qquad\forall\,t\in[T',0),\]
		and
		\[\gamma'_0\cap\big\{y\geq\tsg_0(x)\big\}=\{0\}.\]
	\end{claimheatmain}
	\begin{proof}
		Clearly, $\gamma'_t$ is a subset of $\graphh(g_t)$ by restricting the domain to $(-\rho',\rho')$. Let us also consider the convex function $\sg_t$ so that $\graphh(\sg_t)$ is the linear extension of $\gamma'_t$ to $\RR$. By construction of $\rho'$, we have
		\[\sg'_t\geq-\tan\mu\ \ \text{at}\ \ -\infty,\qquad\sg'_t\leq\tan\mu\ \ \text{at}\ \ +\infty.\]
		In comparison, \eqref{eq-heat:1a_2mu_angle} implies
		\[\tsg'_t<-\tan(2\mu)\ \ \text{at}\ \ -\infty,\qquad\tsg'_t>\tan(2\mu)\ \ \text{at}\ \ +\infty.\]
		Hence, there is a unique $c\in\RR$ so that $\graphh(\sg_t+c)$ touches $\graphh(\tsg_t)$ from below at some point $z=(x,y)$. Denote $\th=\arctan\sg'_t(x)=\arctan\tsg'_t(x)$, so $\th\in[-\mu,\mu]$. Notice that the support functions of $\graphh(\sg_t)$, $\graphh(\tsg_t)$ are restrictions of $h(t,\cdot)$ and $\varphi(t,\cdot)$. Hence
		\[\varphi(t,\th)=\metric{z}{\nu_\th}=h(t,\th)-c\cos\th>\varphi(t,\th)-c\cos\th\qquad\text{if}\ \ (t,\th)\ne(0,0).\]
		This implies, if $t<0$, that $c>0$ hence $\sg_t<\tsg_t$; and if $t=0$, that $c\geq0$ and only touching point has $\th=0$. The latter, along with $\sg_0(0)=\tsg_0(0)=0$, implies $\{\sg_0\geq\tsg_0\}=\{0\}$.
		\begin{figure}[ht]
			\centering
			\includegraphics{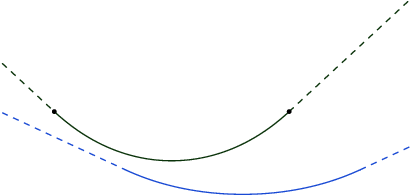}
			\begin{picture}(0,0)
				\put(-124,23){$\tilde\gamma_t$}
				\put(-85,7){$\gamma'_t$}
				\put(-20,10){\arrowangle{26}}
				\put(-12,0){$(\leq\mu)$}
				\put(-200,29){\arrowangle{155}}
				\put(-230,16){$(\geq-\mu)$}
				\put(-195,60){\arrowangle{135}}
				\put(-187,66){$(<-2\mu)$}
				\put(-28,82){\arrowangle{46}}
				\put(-59,91){$(>2\mu)$}
				\put(-1,92){$\graphh(\tsg_t)$}
				\put(-1,21){$\graphh(\sg_t)$}
			\end{picture}
			\caption{Comparing linear extensions of curves from their support functions.}\label{fig-heat:curve_comp}
		\end{figure}
	\end{proof}
	
	Summarizing heuristically, we have curves $\{\gamma'_t\}$ evolving by $1/\kappa$ and curves $\{\tilde\gamma_t\}$ evolving more slowly than $1/\kappa$, so that $\tilde\gamma_t$ lies above $\gamma'_t$ for each $t<0$, but $\tilde\gamma_0$ touches $\gamma_0$ at $0$. This should violate the maximum principle. The precise derivation is as follows.
	
	Choose $T''\in[T',0)$ so that $|\tg_t|<\rho'/2$ in $[-\rho',\rho']$ for all $t\in[T'',0]$. Claim \ref{claim-heat:1a_curve_comp} implies
	\[g_t<\tg_t\leq\rho'/2\text{\ \ in\ \ $[-\rho',\rho']$,}\qquad\forall\,t\in[T'',0),\]
	and
	\[\{g_0\geq\tg_0\}\cap[-\rho',\rho']=\{0\}.\]
	
	Notice that $\tg_t(x)$ is jointly continuous in $t$ and $x$, and $g_t(x)$ is jointly upper semi-continuous by Claim \ref{claim-heat:1_upper_semi_cont}. Hence, there is a small positive constant $y_1>0$ so that
	\begin{equation}\label{eq-heat:1a_curve_comp}
		\begin{aligned}
			& g_t(x)<\tg_t(x)-y_1 \\
			&\hspace{72pt} \forall(t,x)\in\Big([T'',0]\times[-\rho',\rho']\Big)\setminus\Big([T''/2,0]\times[-\rho'/2,\rho'/2]\Big)
		\end{aligned}
	\end{equation}
	Consider the region
	\[\Omega'=\Big\{\!-\rho'<x<\rho',\ -\rho'<y<\tg_{T''}(x)\Big\}.\]
	Then the function
	\[w_1(x,y)=\min\big\{w(x,y-y_1),0\big\}\]
	is a weak IMCF in $\Omega'$, see Remark \ref{rmk-prelim:max_principles}\ref{item-prelim:max_prin_truncate}. Define the function $\tw$ so that $\tw=t$ on $\graphh(\tg_t)$ for all $t\leq0$, and $\tw=0$ in $\{y\leq\tg_0(x)\}$. Then $\tilde w$ is a weak supersolution of IMCF in $\Omega'$, see Remark \ref{rmk-prelim:max_principles}\ref{item-prelim:max_prin_smooth}\ref{item-prelim:max_prin_truncate}. Note that
	\[\{w_1<t\}=\Omega'\cap\{y>g_t(x)+y_1\},\qquad\{\tw< t\}=\Omega'\cap\{y>\tg_t(x)\},\qquad\forall\,t\in[T'',0).\]
	So \eqref{eq-heat:1a_curve_comp} implies that $w_1\leq\tw$ in a neighborhood of $\p\Omega'$ in $\bar{\Omega'}$, namely, $\{w_1>\tilde w\}\Subset\Omega'$. However,
	\[w_1(0,y_1)=\min\{w(0,0),0\}=0>\tilde w(0,y_1).\]
	These violate the maximum principle, see Remark \ref{rmk-prelim:max_principles}\ref{item-prelim:max_principle}, hence proving Case 1(a).
	
	\begin{proof}[Proof of Claim \ref{claim-heat:1a_nondeg}] {\ }
		
		Suppose $\p_t\varphi(0,0)<0$. If $\p_\th^2\varphi(0,0)>0$, then $\varphi$ would be nonnegative in a backward neighborhood of $(0,0)$. On the other hand, for each $t\in[T,0)$, there is a point $z_t\in\gamma_t\cap P$ with the smallest distance from $0$, at which we have
		\[h(t,\Th(z_t))=\metric{z_t}{\nu(z_t)}=-|z_t|<0.\]
		This contradicts $\varphi\leq h$. Hence $\p_\th^2\varphi(0,0)\leq0$. Replacing $\varphi$ by $\varphi-\epsilon\th^2$ for some $\epsilon\ll1$, and further decreasing $\delta$, we may assume
		\[\p_t\varphi<\p_{\th\th}\varphi+\varphi<0\qquad\text{in}\ \ [-\delta,\delta]\times[-\delta,\delta].\]
		Then consider the comparison curve $\tilde\gamma_t$ as in \eqref{eq-heat:comp_curves}. We may decrease $\delta$ and assume that $\tilde\gamma_t\Subset\Int(P)$ for all $t\in[-\delta,0]$. Then $\tilde\gamma_t$ is the graph of a strictly concave function $\tg_t$ over some interval in $[-r,r]$. Extend $g_t,\tg_t$ linearly to $\RR$ in the manner of Claim \ref{claim-heat:1a_curve_comp}, and denote the resulting functions by $\sg_t,\tsg_t$. Note that
		\[\begin{aligned}
			& \hspace{18pt}\sg'_0(-r)=g'_0(-r)\leq0,\qquad\sg'_0(r)=g'_0(r)\geq0; \\
			& \tsg'_t(-r)=\tan\delta,\qquad\tsg'_t(r)=-\tan\delta,\qquad\forall\,t\in[-\delta,0].
		\end{aligned}\]
		Therefore, there exists $t<0$ so that $\sg'_t(-r)<\tsg'_t(-r)$ and $\sg'_t(r)>\tsg'_t(r)$. Then, note that
		\[\sg_t>\sg_0\geq0\geq\tsg_0>\tsg_t\qquad\text{in}\ \ [-r,r],\]
		where the last inequality comes from $\varphi(0,\cdot)<\varphi(t,\cdot)$ and the argument in Claim \ref{claim-heat:1a_curve_comp}. These together imply $\sg_t>\tsg_t$ in $\RR$. Hence, there is $c>0$ so that $\graphh(\sg_t)$ contacts $\graphh(\tsg_t+c)$ tangentially at a point $z=(x,y)$. Denoting $\th=\arctan\sg'_t(x)=\arctan\tsg'_t(x)$, we have
		\[h(t,\th)=\metric{z}{\nu_\th}=\bmetric{z-(0,c)}{\nu_\th}-c\cos\th=\varphi(t,\th)-c\cos\th,\]
		contradicting $\varphi\leq h$.
	\end{proof}
	
	\vspace{3pt}
	
	\textbf{Case 1(b).} Subsolution. Suppose otherwise that there is $\delta>0$ and $\varphi\in C^2\big([-\delta,\delta]\times[-\delta,\delta]\big)$, so that
	\begin{equation}\label{eq-heat:1b_viscosity_prep}
		\left\{\begin{aligned}
			& h(0,0)=\varphi(0,0)=\p_\th\varphi(0,0)=0, \\
			& \varphi>h\text{\ \ in\ \ $\big([-\delta,0]\times[-\delta,\delta]\big)\cap \cJ\setminus\{(0,0)\}$,} \\
			& \p_t\varphi-\p_{\th\th}\varphi-\varphi>0\text{\ \ in\ \ $[-\delta,\delta]\times[-\delta,\delta]$.}
		\end{aligned}\right.
	\end{equation}
	Only item \ref{item-heat:interior_viscosity} is not vacuous for this case, thus $\Th^{-1}(0,0)\Subset\Omega$ by our assumption.
	
	\begin{claimheatmain}\label{claim-heat:1b_nondeg}
		We have $g_0(\pm r)>0$, $g'_0(-r)<0$, $g'_0(r)>0$, $(\p_{\th\th}\varphi+\varphi)(0,0)>0$.
	\end{claimheatmain}
	\begin{proof}
		Denote $\sigma=\Th^{-1}(0,0)\Subset\Omega$. Then Lemma \ref{lemma-heat:0_in_gamma0} implies $0\in\sigma$. Write $\sigma=[x_1,x_2]\times\{0\}$ with $x_1\leq 0\leq x_2$. Recall that we have reduced to the case where $\nu(0)=-\p_y$ and $w$ is positively oriented in a neighborhood of $0$. Hence, the shape of $Z_0$ near $\sigma$ must fall in one of the cases in Figure \ref{fig-heat:1b_nondeg} (if $\sigma$ is a point, then only case 1 can happen since $w$ is an IMCF near $0$). It follows that (see Table \ref{table-prelim:models}) in either case
		\[\p_\th^- h(0,0)=x_1,\qquad\p_\th^+ h(0,0)=x_2.\]
		Since $\varphi$ touches $h$ from above at $(0,0)$, this forces $x_1=x_2=0$. Hence $\sigma$ is a point, and $g_0(\pm r)>0$. Hence $g'_0(-r)<0$ and $g'_0(r)>0$ by convexity. Finally, $(\p_{\th\th}\varphi+\varphi)(0,0)>0$ follows from $\p_{\th\th}h+h>0$ in the barrier sense.
		\begin{figure}[ht]
			\centering
			\includegraphics{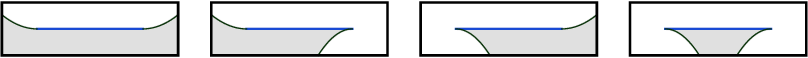}
			\caption{Local models in Claim \ref{claim-heat:1b_nondeg}.}\label{fig-heat:1b_nondeg}
		\end{figure}
	\end{proof}
	
	We may further decrease $\delta$ and assume that
	\begin{equation}\label{eq-heat:1b_viscosity_root}
		\left\{\begin{aligned}
			& h(0,0)=\varphi(0,0)=\p_\th\varphi(0,0)=0, \\
			& \varphi>h\text{\ \ in\ \ $\big([-\delta,0]\times[-\delta,\delta]\big)\cap \cJ\setminus\{(0,0)\}$,} \\
			& 0<\p_{\th\th}\varphi+\varphi<\p_t\varphi\text{\ \ in\ \ $[-\delta,\delta]\times[-\delta,\delta]$.}
		\end{aligned}\right.
	\end{equation}
	By Claim \ref{claim-heat:1b_nondeg} and $g_t\to g_0$ in $C^1([-r,r])$, there is $\epsilon>0$ and $T'\in[T,0)$ so that
	\begin{equation}\label{eq-heat:1b_bd_slope}
		g'_t(-r)<-\tan(2\epsilon),\qquad g'_t(r)>\tan(2\epsilon),\qquad\forall\,t\in[T',0].
	\end{equation}
	
	Consider the comparison curves $\{\tilde\gamma_t\}$ and swept out region $S$ as in \eqref{eq-heat:comp_curves} \eqref{eq-heat:def_S}. We may decrease $\delta$ so that $\delta<\epsilon$ and $S\Subset\Int(P)$. Then $0<\p_{\th\th}\varphi+\varphi<\p_t\varphi$ implies that $\{\tilde\gamma_t\}$ forms a downward evolving smooth subsolution of IMCF (the normal speed is greater than $1/\kappa$). Let $\rho$ be small enough so that $Q(\rho)\Subset S$.
	
	Using the same comparison argument as in Claim \ref{claim-heat:1a_curve_comp}, noticing that $\tilde\gamma_t$ have angle $\delta<\epsilon$ while $\graphh(g_t)$ have angle $>2\epsilon$ at endpoints, we obtain
	\begin{equation}\label{eq-heat:1b_curve_comp_pre}
		\tilde\gamma_t\subset\{y<g_t(x)\}\quad\forall\,t\in[T',0),\qquad\tilde\gamma_0\cap\{y\geq g_0(x)\}=\{0\}.
	\end{equation}
	
	\begin{claimheatmain}\label{claim-heat:1b_continuity}
		There is $T''\in[T',0)$ and $\rho'<\rho$, so that
		\[\left\{\begin{aligned}
			& |g_t|<\rho'/2 \\
			& \Int\big(\{w=t\}\big)\cap\{-\rho'\leq x\leq\rho'\}=\emptyset
		\end{aligned}\right.\qquad\forall\,t\in[T'',0].\]
		Thus $(t,x)\mapsto g_t(x)$ is continuous in $[T'',0]\times[-\rho',\rho']$.
	\end{claimheatmain}
	\begin{proof}
		Let $g_t^+$ be such that $\graphh(g_t^+)=\p\{w\leq t\}\cap P$. Recall that in a weak IMCF, the set $\p\{w\leq t\}\setminus\p\{w<t\}$ is a minimal surface. Then the set $I_t:=\{x:g_t^+(x)<g_t(x)\}$ is either $[-r,r]$ or is of the form $[-r,a)\cup(b,r]$ with $-r\leq a\leq b\leq r$. Moreover, $g_t^+$ is linear in each connected component of $\{g_t^+<g_t\}$.
		
		Recall that $g_0\in C^{1,1}([-r,r])$ and $g_0(0)=g'_0(0)=0$ and $g_0(\pm r)>0$. Thus there exists $\rho'<\rho$ so that
		\[g_0(\pm 2\rho')<\frac{r-2\rho'}{r-\rho'}g_0(\pm \rho')+\frac{\rho'}{r-\rho'}g_0(\pm r).\]
		Since $g_t^+\to g_0$ in $C^1([-r,r])$ as $t\nearrow0$, there is $T''\in[T',0)$ so that
		\begin{equation}\label{eq-heat:1b_aux1}
			g_t^+(\pm 2\rho')<\frac{r-2\rho'}{r-\rho'}g_t^+(\pm \rho')+\frac{\rho'}{r-\rho'}g_t^+(\pm r),\qquad\forall\,t\in[T'',0].
		\end{equation}
		It follows that $I_t\cap[-\rho',\rho']=\emptyset$ for all $t\in[T'',0]$: otherwise, $I_t$ would contain either $[-r,-\rho']$ or $[\rho',r]$, and $g_t^+$ would be linear in these intervals, contradicting \eqref{eq-heat:1b_aux1}. Then $I_t\cap[-\rho',\rho']=\emptyset$ implies the no-jump statement. Finally, the continuity statement follows from Claim \ref{claim-heat:1_upper_semi_cont}, and the condition $|g_t|<\rho'/2$ follows from further increasing $T''$.
	\end{proof}
	
	For each $t\in[T'',0]$, let $\tg_t:[-\rho',\rho']\to[-\rho'/2,\rho'/2]$ be so that $\tilde\gamma_t\cap\bar{Q(\rho')}=\graphh(\tg_t)$. By \eqref{eq-heat:1b_curve_comp_pre} and Claim \ref{claim-heat:1b_continuity}, there is a small $y_1>0$ so that
	\begin{equation}\label{eq-heat:1b_curve_comp}
		\begin{aligned}
			& \tg_t(x)<g_t(x)-y_1, \\
			&\hspace{72pt} \forall\,(t,x)\in\Big([T'',0]\times[-\rho',\rho']\Big)\setminus\Big([T''/2,0]\times[-\rho'/2,\rho'/2]\Big).
		\end{aligned}
	\end{equation}
	To summarize, we have a downward evolving weak IMCF $\{\gamma_t\}$ and a subsolution $\{\tilde\gamma_t\}$, so that each $\tilde\gamma_t$ lies strictly below $\gamma_t$ with the only exception $\tilde\gamma_0\cap\gamma_0=\{0\}$. Using the same maximum principle argument as in the end of Case 1(a), we may derive a contradiction.
	
	\vspace{3pt}
	
	To summarize, we have just proved Theorem \ref{thm-heat:main} for the case where $w$ is a weak IMCF in a neighborhood of $0$.
	
	We are not able to prove (i) when $\Th^{-1}(0,0)\not\Subset\Omega$: the obstruction is that $\varphi>h$ may no longer imply that $\tilde\gamma_t$ lies below $\gamma_t$. See for example Figure \ref{fig-heat:1b_res}: here $\tilde\gamma_0$ lies above $\gamma_0$ in the $x>0$ part, but we cannot rule this out from $\varphi\geq h$ since the angle parameter of $\gamma_0$ does not contain the $\th>0$ portion.
	
	% == Optional ==
	% In Claim \ref{claim-heat:1a_curve_comp}, the angle condition at $\pm\infty$ is essential for this reason. And it is the reason that one needs to truncate $\gamma_t$ at angle $\pm\mu$ in Case 1(a).
	
	\begin{figure}[ht]
		\centering
		\includegraphics{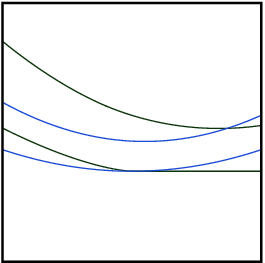}
		\begin{picture}(0,0)
			\put(-4,53){$\tilde\gamma_0$}
			\put(-4,71){$\tilde\gamma_{-\epsilon}$}
			\put(-143,64){$\gamma_0$}
			\put(-148,106){$\gamma_{-\epsilon}$}
		\end{picture}
		\caption{A case where $\varphi\geq h$ but $\tilde\gamma_t$ does not lie below $\gamma_t$.}\label{fig-heat:1b_res}
	\end{figure}
	
	\vspace{3pt}
	
	\textbf{Case 2:} $w$ is a simple cluster, and $0$ lies on a ridge $\gamma$. We still denote
	\[P=[-r,r]\times[-r,r].\]
	Recall that $0\in\gamma_0$ and $\nu(0)=-\p_y$. For $r$ small enough, we may assume
	\[|\nu+\p_y|<e^{-10},\quad|\Th|<e^{-10}\qquad\text{in\ \ }P\cap\Gamma.\]
	Let $f_{\ridge}:[-r,r]\to[-e^{-9}r,e^{-9}r]$ be the $e^{-9}$-Lipschitz function so that
	\[\gamma\cap P=\graphh(f_{\ridge}).\]
	Denote
	\[D_1=P\cap\{y<f_{\ridge}(x)\},\qquad D_2=P\cap\{y>f_{\ridge}(x)\}.\]
	Thus, $w$ consists of an upward evolving weak IMCF in $D_1$ and downward evolving weak IMCF in $D_2$. Note that $D_1$ has negative orientation (i.e. $\chi_1=-1$).
	
	Recall the local model of $Z_0$ given by Lemma \ref{lemma-cluster:local_model}\ref{item-cluster:Yt_Zt_model}, and that Lemma \ref{lemma-heat:no_model_ic} excludes $-r<a=b<r$. Hence, up to a $180^\circ$ rotation and orientation reversing (Remark \ref{rmk-heat:misc}\ref{item-heat:sym_orient_rot}), we fall into one of the two cases below depending on whether $0$ is a vertex or not:
	
	\begin{enumerate}[label={\arabic*.}, topsep=1pt, itemsep=-0.5ex]
		\item after further decreasing $r$, we have
		\begin{equation}\label{eq-heat:model_1}
			Z_0\cap P=\Big\{x\in(0,r],\ g_{0,1}(x)<y<g_{0,2}(x)\Big\},
		\end{equation}
		where $g_{0,1}$ is concave, and $g_{0,2}$ is convex, and $g_{0,1}(0)=g_{0,2}(0)=g'_{0,1}(0)=g'_{0,2}(0)=0$;
		\item after further decreasing $r$, we have
		\begin{equation}\label{eq-heat:model_2}
			Z_0\cap P=\Big\{x\in[-r,r],\ y<g_{0,2}(x)\Big\},
		\end{equation}
		where $g_{0,2}$ is convex with $g_{0,2}(0)=g'_{0,2}(0)=0$ and $g_{0,2}\geq f_{\ridge}$.
	\end{enumerate}
	
	\noindent We first assume that possibility 2 holds. Then the (well-defined) new function
	\[w'=\left\{\begin{aligned}
		& w\qquad\text{in}\ \ \big\{|x|<r,\ g_{0,2}(x)\leq y<r\big\} \\
		& 0\qquad\text{in}\ \ \big\{|x|<r,\ -r<y\leq g_{0,2}(x)\big\}
	\end{aligned}\right.\]
	is a weak IMCF in $\Int(P)$, by Remark \ref{rmk-prelim:imcf_properties}\ref{item-prelim:imcf_extension}. Consider the restricted data in $Q(r)$:
	\[\begin{aligned}
		& \gamma'_t=\p\{w'<t\}\cap Q(r)\subset\gamma_t,\qquad\Gamma'=\bigcup_{t\in\RR}\gamma'_t\subset\Gamma,\qquad\nu'=\nu|_{\Gamma'} \\
		&\hspace{36pt} \Th'=\Th|_{\Gamma'},\qquad \cJ'=\Th'(\Gamma')\subset\cJ,\qquad h'=h|_{\cJ'}.
	\end{aligned}\]
	Clearly, $\Th'$ is a lift of $(t,\nu^\perp)|_{\Gamma'}$ to $\RR\times M$, with each preimage being a line segment. Since $0\in\gamma'_0$ and $\Th'(0)=(0,0)$, we have $h'(0,0)=0$ as well. If $\varphi$ touches $h$ from below, then $\varphi$ also touches $h'$ from below in $\cJ'$. If $(\Th')^{-1}(0,0)\Subset Q(r)$, then by the IMCF case of Theorem \ref{thm-heat:main}(i), which was already proved in Case 1, we have $(\p_t\varphi-\p_{\th\th}\varphi-\varphi)(0,0)\geq0$. If $(\Th')^{-1}(0,0)\not\Subset Q(r)$, then note that
	\[\varphi(0,0)\nu_0+\p_\th\varphi(0,0)\tau_0=0\in(\Th')^{-1}(0,0),\]
	hence \eqref{eq-heat:intro_on_sigma} is satisfied for $w'$. Notice that $\chi_\sigma=1$ for our case here. Hence the IMCF case of Theorem \ref{thm-heat:main}(ii) implies $(\p_t\varphi-\p_{\th\th}\varphi-\varphi)(0,0)\geq0$ as well. If $\varphi$ touches $h$ from above, then we have assumed $\Th^{-1}(0,0)\Subset\Omega$. Similar to Claim \ref{claim-heat:1b_nondeg}, this implies that $\Th^{-1}(0,0)$ is a point, hence $(\Th')^{-1}(0,0)=\{0\}\Subset Q(r)$. So the IMCF case of Theorem \ref{thm-heat:main}(i) gives $(\p_t\varphi-\p_{\th\th}\varphi-\varphi)(0,0)\leq0$, as desired.
	
	It remains to treat the case \eqref{eq-heat:model_1}. Let us summarize the known conditions: we have
	\begin{equation}\label{eq-heat:2a_Z0_model}
		Z_0\cap P=\Big\{x\in(0,r],\ g_{0,1}(x)<y<g_{0,2}(x)\Big\},
	\end{equation}
	where $g_{0,1}:[0,r]\to[-e^{-9}r,0]$ is concave, $g_{0,2}:[0,r]\to[0,e^{-9}r]$ is convex, the graphs of both functions are $\nu$-orthogonal, and $g_{0,1}(0)=g_{0,2}(0)=g'_{0,1}(0)=g'_{0,2}(0)=0$, and $g_{0,1}\leq f_{\ridge}\leq g_{0,2}$ in $[0,r]$, and $g_{0,1}<g_{0,2}$ in $(0,r]$, by Lemma \ref{lemma-cluster:local_model}\ref{item-cluster:Yt_Zt_model}\ref{item-cluster:g1_g2}.
	
	Next, using Lemma \ref{lemma-cluster:local_model}\ref{item-cluster:tZt_model}, we have
	\[\tZ_0\cap P=\Big\{x\in I,\ \tg_{0,1}(x)\leq y\leq\tg_{0,2}(x)\Big\},\]
	where either $I=[-r,r]$ or $I=[-r,a]\cup[b,r]$, and $\tg_{0,1},\tg_{0,2}$ are concave resp. convex in each interval in $I$. Since $Z_0=\Int(\tZ_0)$, it holds $\tg_{0,1}=g_{0,1}$ and $\tg_{0,2}=g_{0,2}$ in $[0,r]$. By \eqref{eq-heat:2a_Z0_model}, we must have $\tg_{0,1}=\tg_{0,2}$ in $I\cap[-r,0]$. Therefore, depending on whether $I=[0,r]$ or not, we fall into one of the following two possibilities by further decreasing $r$:
	\begin{equation}\label{eq-heat:tZ0_case_1}
		\tZ_0\cap P=\Big\{x\in[0,r],\ g_{0,1}(x)\leq y\leq g_{0,2}(x)\Big\},
	\end{equation}
	or
	\begin{equation}\label{eq-heat:tZ0_case_2}
		\tZ_0\cap P=\Big\{x\in[-r,r],\ \tg_{0,1}(x)\leq y\leq\tg_{0,2}(x)\Big\},
	\end{equation}
	where $\tg_{0,i}=g_{0,i}$ in $[0,r]$ and $\tg_{0,i}\equiv0$ in $[-r,0]$. If the second possibility holds, then
	\[w_1=\left\{\begin{aligned}
		& w\qquad\text{in}\ \ \big\{|x|<r,\ -r<y\leq\tg_{0,1}(x)\big\} \\
		& 0\qquad\text{in}\ \ \big\{|x|<r,\ \tg_{0,1}(x)\leq y<r\big\}
	\end{aligned}\right.\]
	and
	\[w_2=\left\{\begin{aligned}
		& w\qquad\text{in}\ \ \big\{|x|<r,\ \tg_{0,2}(x)\leq y<r\big\} \\
		& 0\qquad\text{in}\ \ \big\{|x|<r,\ -r<y\leq\tg_{0,2}(x)\big\}
	\end{aligned}\right.\]
	are respectively an upward evolving and downward evolving weak IMCF in $\Int(P)$, by Remark \ref{rmk-prelim:imcf_properties}\ref{item-prelim:imcf_extension}. Note that $w_1$ is negatively oriented. If $\varphi$ touches $h$ from below, then we may restrict the data to $w_2$ as in the previous page, then use the IMCF case of Theorem \ref{thm-heat:main} to obtain $(\p_t\varphi-\p_{\th\th}\varphi-\varphi)(0,0)\geq0$. If $\varphi$ touches $h$ from above, then we may reflect around the $x$ axis and change the orientation for $w_1$ (see Remark \ref{rmk-heat:misc}\ref{item-heat:sym_x_refl_orient}). Then $\varphi$ would touch $h$ from below, and Theorem \ref{thm-heat:main} implies the desired viscosity property as well.
	
	Now, we may assume that \eqref{eq-heat:tZ0_case_1} holds. We give further setups on $Z_t\cap P$ for $t<0$. By Lemma \ref{lemma-cluster:local_model}\ref{item-cluster:Yt_Zt_model} and Lemma \ref{lemma-heat:no_model_ic} and $Z_t\supset Z_0$, we have
	\[Z_t\cap P=\Big\{x\in I_t,\ g_{t,1}(x)<y<g_{t,2}(x)\Big\},\]
	where $I_t=[-r,r]$ or $I_t=[-r,a_t)\cup(b_t,r]$ with $-r\leq a_t<b_t<0$. By \eqref{eq-heat:tZ0_case_1} and $\tZ_0=\bigcap_{t<0}Z_t$ and that $Z_t$ is decreasing in $t$, there exists $T<0$ so that
	\begin{equation}\label{eq-heat:2a_Zt_model}
		Z_t\cap P=\Big\{x_t<x\leq r,\ g_{t,1}(x)<y<g_{t,2}(x)\},\qquad\forall\,t\in[T,0].
	\end{equation}
	(Namely, $I_t$ does not contain the interval $[-r,a_t)$ for all $t\in[T,0]$.) Also, note that $\lim_{t\nearrow0}x_t=0$. We may further increase $T$ and assume that $x_t\in[-r/2,0]$ for all $t\in[T,0]$.
	
	Lemma \ref{lemma-cluster:local_model}\ref{item-cluster:g1_g2} implies that $g_{t,1}:[x_t,r]\to[-r/2,r/2]$ is concave and $e^{-9}$-Lipschitz, and $g_{t,2}:[x_t,r]\to[-r/2,r/2]$ is convex and $e^{-9}$-Lipschitz, and the graphs of both functions are $\nu$-orthogonal. Regarding the relation with the ridge, we have
	\[g_{t,1}(x)\leq f_{\ridge}(x)\leq g_{t,2}(x)\qquad\forall\,x\in[x_t,r],\]
	and
	\[g_{t,1}(x)<g_{t,2}(x)\qquad\forall\,x\in(x_t,r],\]
	and
	\[g_{t,1}(x_t)=f_{\ridge}(x_t)=g_{t,2}(x_t),\qquad g'_{t,1}(x_t)=f'_{\ridge}(x_t)=g'_{t,2}(x_t).\]
	The vertex $x_t$ satisfies
	\[x_t\in[-r/2,0]\quad \forall\,t\in[T,0],\qquad\text{and}\qquad x_t\nearrow0\ \ \text{as}\ \ t\nearrow0.\]
	By the convergence and $\nu$-orthogonality, we also have
	\[g_{t,1}\xrightarrow[\text{increasing}]{C^1([0,r])}g_{0,1},\qquad g_{t,2}\xrightarrow[\text{decreasing}]{C^1([0,r])}g_{0,2}\qquad\text{as}\ \ t\nearrow0.\]
	
	We are ready to prove the viscosity property. We first reduce all statements to one case. If $\Th^{-1}(0,0)\Subset\Omega$, then concave polygon geometry (Figure \ref{fig-prelim:models}, Table \ref{table-prelim:models}) yields
	\begin{equation}\label{eq-heat:2a_reduction1}
		\exists\,\varphi\text{ touching $h$ from below}\quad\Rightarrow\quad g_{0,1}(x)<0\ \ \ \forall\,x>0,
	\end{equation}
	and
	\begin{equation}\label{eq-heat:2a_reduction2}
		\exists\,\varphi\text{ touching $h$ from above}\quad\Rightarrow\quad g_{0,2}(x)>0\ \ \ \forall\,x>0.
	\end{equation}
	Through a reflection around the $x$-axis and an orientation reversal, see Remark \ref{rmk-heat:misc}\ref{item-heat:sym_x_refl_orient}, proving the first case \eqref{eq-heat:2a_reduction1} implies the second case \eqref{eq-heat:2a_reduction2}. Next, if $\sigma:=\Th^{-1}(0,0)\not\Subset\Omega$, then either $g_{0,1}\equiv0$ or $g_{0,2}\equiv0$. If $g_{0,2}\equiv0$, then $g_{0,1}<0$ in $(0,r]$. Then, the definition of $\chi_\sigma$ in Theorem \ref{thm-heat:main} implies $\chi_\sigma=1$, so $\varphi$ should touch $h$ from below in item \ref{item-heat:boundary_viscosity} of the theorem. If $g_{0,1}\equiv0$, then $g_{0,2}>0$ in $(0,r]$, hence $\chi_\sigma=-1$, hence $\varphi$ should touch $h$ from above. Through a reflection around $x$-axis and orientation reversal, the first case implies the second case as well. Therefore, we are reduced to assuming the condition
	\begin{equation}\label{eq-heat:2a_after_reduction}
		g_{0,1}<0\ \ \text{in}\ \ (0,r],\qquad\varphi\text{ touches $h$ from below,}
	\end{equation}
	and trying to prove
	\begin{equation}\label{eq-heat:2a_to_prove}
		(\p_t\varphi-\p_{\th\th}\varphi-\varphi)(0,0)\geq0.
	\end{equation}
	Recall that $g_{0,1}$ is concave with $g_{0,1}(0)=g'_{0,1}(0)=0$, so $g'_{0,1}(r)<0$. Since $\lim_{t\nearrow0}g_{t,1}=g_{0,1}$ in $C^1([0,r])$, we may further increase $T$ and find an $\epsilon>0$, so that
	\begin{equation}\label{eq-heat:2a_bd_angle}
		g'_{t,1}(r)<-\tan(2\epsilon)\qquad\forall\,t\in[T,0].
	\end{equation}
	It remains to prove the essentially unique case \eqref{eq-heat:2a_after_reduction} \eqref{eq-heat:2a_to_prove}.
	
	\begin{figure}[ht]
		\centering
		\includegraphics{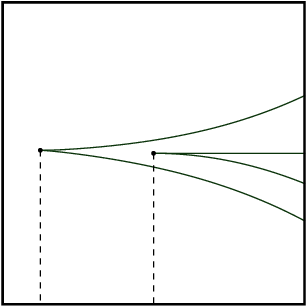}
		\begin{picture}(0,0)
			\put(-2,100){$g_{t,2}$}
			\put(-2,72){$g_{0,2}$}
			\put(-2,57){$g_{0,1}$}
			\put(-2,39){$g_{t,1}:g'_{t,1}(r)<-\tan(2\epsilon)$}
			\put(-130,5){$x_t$}
			\put(-75,5){$0$}
		\end{picture}
		\caption{Shape of $Z_t$.}\label{fig-heat:2a_info}
	\end{figure}
	
	\vspace{3pt}
	
	\textbf{Case 2(a).} If the viscosity property fails, then by our reduction at the beginning of the entire proof, there are $\delta<\min\{\epsilon,e^{-10}\}$ and $\varphi\in C^2\big([-\delta,\delta]\times[-\delta,\delta]\big)$ so that
	\begin{equation}\label{eq-heat:2a_viscosity_prep}
		\left\{\begin{aligned}
			& h(0,0)=\varphi(0,0)=\p_\th\varphi(0,0)=0, \\
			& \varphi<h\text{\ \ in\ \ $\big([-\delta,0]\times[-\delta,\delta]\big)\cap \cJ\setminus\{(0,0)\}$,} \\
			& \p_t\varphi-\p_{\th\th}\varphi-\varphi<0\text{\ \ in\ \ $[-\delta,\delta]\times[-\delta,\delta]$.}
		\end{aligned}\right.
	\end{equation}
	
	\begin{claimheatmain}
		$(\p_{\th\th}\varphi+\varphi)(0,0)\leq0$.
	\end{claimheatmain}
	\begin{proof}
		The local model \eqref{eq-heat:2a_Z0_model} implies $\p_\th^- h(0,0)=0$. So for each $\mu$ with $0<\mu\ll1$, there is a constant $c_\mu$ so that $\varphi(0,\th)-\mu\th+c_\mu$ touches $h$ from below at a point $\th_\mu\in(-\delta,0)$. The concavity of $g_{0,1}$ implies that $(\p_{\th\th}h+h)(0,\th)<0$ in the barrier sense for all $\th\in(-\delta,0)\subset(\arctan g'_{0,1}(r),0)$. By \eqref{eq-heat:2a_viscosity_prep} we have $c_\mu\to0$ and $\th_\mu\to0$ when $\mu\to0$. Hence
		\[\begin{aligned}
			(\p_{\th\th}\varphi+\varphi)(0,0) &= \lim_{\mu\to0}\Big[\p_{\th\th}\varphi(0,\th_\mu)+\varphi(0,\th_\mu)-\mu\th_\mu+c_\mu\Big] \\
			&= \lim_{\mu\to0}\Big[(\p_{\th\th}+\id)\big(\varphi-\mu\th+c_\mu\big)\Big](0,\th_\mu)\leq0.\qedhere
		\end{aligned}\]
	\end{proof}
	
	Thus, by subtracting a small multiple of $\th^2$ from $\varphi$, then further decreasing $\delta$, we may assume
	\begin{equation}\label{eq-heat:2a_viscosity_root}
		\left\{\begin{aligned}
			& h(0,0)=\varphi(0,0)=\p_\th\varphi(0,0)=0,\qquad\delta<\epsilon, \\
			& \varphi<h\text{\ \ in\ \ $\big([-\delta,0]\times[-\delta,\delta]\big)\cap \cJ\setminus\{(0,0)\}$,} \\
			& \p_t\varphi<\p_{\th\th}\varphi+\varphi<0\text{\ \ in\ \ $[-\delta,\delta]\times[-\delta,\delta]$.}
		\end{aligned}\right.
	\end{equation}
	Define the comparison curves $\{\tilde\gamma_t\}$ and region $S$ similarly as the previous case:
	\[\tilde\gamma_t=\Big\{\varphi(t,\th)\nu_\th+\p_\th\varphi(t,\th)\tau_\th:\th\in(-\delta,\delta)\Big\},\qquad S=\bigcup_{t\in(-\delta,\delta)}\tilde\gamma_t.\]
	Now $\{\tilde\gamma_t\}$ is an upward evolving subsolution of IMCF ($\text{speed}>1/\kappa$), with $0\in\tilde\gamma_0$. We may decrease $\delta$ and assume that $S\Subset\Int(P)$.
	
	For convenience, we merge $f_{\ridge}$ and $g_{t,1}$ into a $C^1$ function:
	\[G_t(x)=\left\{\begin{aligned}
		& f_{\ridge}(x)\qquad x\in[-r,x_t], \\
		& g_{t,1}(x)\qquad x\in[x_t,r].
	\end{aligned}\right.\]
	The continuity of $G_t$ and comparison between $\gamma_t,\tilde\gamma_t$ are stated as follows:
	
	\begin{claimheatmain}\label{claim-heat:2a_no_jump}
		There are $T'\in[T,0)$ and $\rho<r$ so that $G_t(x)$ is jointly continuous for $(t,x)\in[T',0]\times[-\rho,\rho]$.
	\end{claimheatmain}
	
	\begin{claimheatmain}\label{claim-heat:2a_curve_comp}
		There is $T''\in[T',0)$ and $\rho'<\rho$, so that
		\[g_{T'',1}<f_{\ridge}\qquad\text{in}\ \ [-\rho',\rho'],\]
		and $\tilde\gamma_t\cap\{|x|\leq\rho'\}=\graphh(\tg_t)$ with
		\[\tg_t(x)>G_t(x)\qquad\forall\,(t,x)\in[T'',0]\times[-\rho',\rho']\setminus\{(0,0)\}.\]
	\end{claimheatmain}
	
	Combining Claims \ref{claim-heat:2a_no_jump} and \ref{claim-heat:2a_curve_comp}, there is $y_1>0$ so that
	\begin{equation}\label{eq-heat:2a_strict_comp}
		\begin{aligned}
			& \tg_t(x)>G_t(x)+y_1 \\
			&\hspace{72pt} \forall(t,x)\in\Big([T'',0]\times[-\rho',\rho']\Big)\setminus\Big([T''/2,0]\times[-\rho'/2,\rho'/2]\Big).
		\end{aligned}
	\end{equation}
	
	Finally, consider the region
	\[\begin{aligned}
		A=\Big\{\!-\rho'<x<\rho',\ g_{T'',1}(x)<y\leq f_{\ridge}(x)\Big\}
	\end{aligned}\]
	together with the outer obstacle
	\[O=\Big\{\!-\rho'<x<\rho',\ y=f_{\ridge}(x)\Big\}.\]
	Let $\tilde w$ be the function in $A$ so that $\tilde w=t$ on $\tilde\gamma_t-(0,y_1)$ for each $t\leq0$, and $\tilde w\equiv0$ elsewhere. Then $\tilde w$ is a weak subsolution of IMCF in $\Int(A)$, see Remark \ref{rmk-prelim:max_principles}\ref{item-prelim:max_prin_truncate}. On the other hand, by the definition of simple clusters, $w|_A$ is a continuously calibrated weak IMCF in $\Int(A)$ with outer obstacle $O$. Finally, note that
	\[\{w<t\}\cap A\subset\big\{y\leq G_t(x)\big\}\cap A,\quad\{\tw<t\}\cap A=\big\{y<\tg_t(x)-y_1\big\}\cap A,\]
	thus \eqref{eq-heat:2a_strict_comp} implies
	\[\{\tilde w>w\}\Subset A.\]
	Moreover, we have
	\[\tilde w(0,-y_1)=0>w(0,-y_1).\]
	These together violate the maximum principle with obstacle, see Remark \ref{rmk-prelim:max_principles}\ref{item-prelim:max_principle_imcfoo}.
	
	\begin{proof}[Proof of Claim \ref{claim-heat:2a_no_jump}] {\ }
		
		Recall $D_1=\{y<f_{\ridge}(x)\}$, and $\{w<t\}\cap D_1=\{y<G_t(x)\}$. Since $\graphh(G_t)$ is orthogonal to the continuous vector field $\nu$, it follows that $\{G_t(x)\}$ is a uniformly $C^1$ family of functions for $t\in[T,0]$. Also, we have
		\begin{equation}\label{eq-heat:aux6}
			\lim_{s\nearrow t}G_s=G_t\leq G_t^+=\lim_{s\searrow t}G_s\qquad\text{in}\ \ [-r,r],
		\end{equation}
		where $G_t^+(x)$ is the function so that $\{w\leq t\}\cap\bar{D_1}=\{y\leq G_t^+(x)\}$.
		
		We claim that $\{G_t<G_t^+\}$ must be of the form $(b,r]$ for some $b\in[x_t,r]$. If this is not true, then there is a component $(b,c)$ of $\{G_t<G_t^+\}$ with $x_t\leq b<c<r$. Then $w$ is constant in the region $\big\{b\leq x\leq c,\ G_t(x)\leq y\leq G_t^+(x)\big\}$, and the divergence theorem gives
		\begin{equation}\label{eq-heat:aux8}
			\big|\graphh(G_t|_{[b,c]})\big|=\big|\graphh(G_t^+|_{[b,c]})\big|.
		\end{equation}
		However, $G_t$ is concave in $(b,c)$ since $G_t<G_t^+\leq f_{\ridge}$. Then \eqref{eq-heat:aux8} is impossible to hold.
		
		\begin{figure}[ht]
			\centering
			\includegraphics{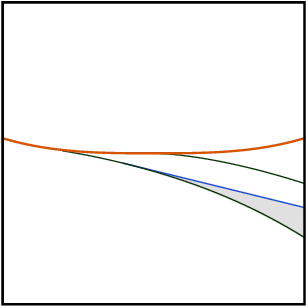}
			\begin{picture}(0,0)
				\put(-2,78){$f_{\ridge}$}
				\put(-2,55){$G_0$}
				\put(-2,41){$G_t^+$}
				\put(-2,27){$G_t$}
				\put(-21,38){\arrowangle{58}}
				\put(-45,28){jump}
			\end{picture}
			\caption{Objects in Claim \ref{claim-heat:2a_no_jump}.}\label{fig-heat:2a_jump}
		\end{figure}
		
		Recall $G_0(0)=G'_0(0)=0$, and $G_0$ is concave in $[0,r]$ with $G_0(r)<0$. Similar to Claim \ref{claim-heat:1b_continuity}, we may find $T'<0$ and $\rho<r$ so that $G_t^+$ cannot be a linear function in $[2\rho,r]$ for any $t\in[T',0]$. This implies $\{G_t^+>G_t\}\subset(\rho,r]$ for all $t\in[T',0]$: indeed, if not, then $G_t<G_t^+<G_0\leq f_{\ridge}$ in $(\rho,r]$, hence $\graphh(G_t^+|_{(\rho,r]})$ is a line segment (since it does not touch the obstacle $\graphh(f_{\ridge})$), contradicting our nonlinearity conclusion above.
	\end{proof}
	
	\begin{proof}[Proof of Claim \ref{claim-heat:2a_curve_comp}] {\ }
		
		Let $\tg_t$ be the (concave) function so that $\tilde\gamma_t=\graphh(\tg_t)$, so $\tg_t$ is defined in an interval that depends continuously on $t$. Recall that for all $t\in[T,0]$, there are concave functions $g_{t,1}:[x_t,r]\to[-r/2,r/2]$ and convex functions $g_{t,2}:[x_t,r]\to[-r/2,r/2]$, so that $\gamma_t\cap P=\graphh(g_{t,1})\cup\graphh(g_{t,2})$. Set
		\[\th_t=\arctan g'_{t,1}(x_t)=\arctan g'_{t,2}(x_t)=\arctan f'_{\ridge}(x_t)\]
		the angles of the vertices. We collect the following known facts:
		\begin{itemize}[nosep]
			\item $x_t\in[-r/2,0]$ for all $t\in[T,0]$, and $x_t\nearrow0$ and $\th_t\to0$ as $t\nearrow0$;
			\item $g_{s,1}(x)<g_{t,1}(x)$ and $g_{s,2}(x)>g_{t,2}(x)$ for all $T\leq s<t\leq0$ and $x\in[x_t,r]$;
			\item $g_{t,1}\nearrow g_{0,1}$ and $g_{t,2}\searrow g_{0,2}$ in $C^1([0,r])$, as $t\nearrow0$;
			\item $g'_{t,1}(r)<-\tan(2\epsilon)$ for all $t\in[T,0]$, by \eqref{eq-heat:2a_bd_angle}, and $\delta<\epsilon$ by our setups;
			\item $\tg_t<0$ for all $t<0$;
			\item $(t,x)\mapsto G_t(x)$ is jointly continuous in $[T',0]\times[-\rho,\rho]$, where $T'\in[T,0]$ and $\rho<r$. Here we made use of Claim \ref{claim-heat:2a_no_jump}.
		\end{itemize}
		We may find $T_1\in[T'/2,0]$ so that
		\begin{itemize}[nosep]
			\item $|\th_t|<\delta/2$ for all $t\in[2T_1,0]$;
			\item $0\in\dom(\tg_t)$ for all $t\in[2T_1,0]$.
		\end{itemize}
		To summarize, we have the angle relations
		\[-2\epsilon<-\delta<-\delta/2<\th_t<\delta/2<\delta,\qquad\forall\,t\in[T_1,0].\]
		Next, we select a $\rho_1<\rho$ so that the following holds:
		\begin{itemize}[nosep]
			\item $\dom(\tg_t)\Supset[-\rho_1,\rho_1]$ for all $t\in[T_1,0]$;
			\item $Q(\rho_1)\Subset Z_{T_1}$; thus $x_{T_1}<-\rho_1$ and $g_{T_1,1}(x)<-\rho_1$, $g_{T_1,2}(x)>\rho_1$ for all $x\in[-\rho_1,\rho_1]$;
			\item $x_{T_1/2}<-\rho_1$ and $g_{T_1/2,1}<f_{\ridge}$ in $[-\rho_1,\rho_1]$.
		\end{itemize}
		To summarize, $g_{t,1},g_{t,2}$ both have domains $[x_t,r]$, while $[-\rho_1,\rho_1]\Subset\dom(\tg_t)\Subset(-r,r)$.
		
		Let us extend $g_{t,1}$, $g_{t,2}$, $\tg_t$ in the $C^1$ manner from their own domains to $\RR$, by attaching two linear functions at their endpoints (same as in Claim \ref{claim-heat:1a_curve_comp}). Let $\sg_{t,1}$, $\sg_{t,2}$, $\tsg_t$ denote the resulting functions. Consider the convex sets
		\[A_t^1=\big\{y\leq \sg_{t,1}(x)\big\},\qquad A_t^2=\big\{y\geq \sg_{t,2}(x)\big\},\qquad \tilde A_t=\big\{y\leq \tsg_t(x)\big\},\]
		See Figure \ref{fig-heat:curve_comp_2} for a visualization. Note that $A_t^1\cap A_t^2=\{x\leq x_t, y=\sg_{t,1}(x)=\sg_{t,2}(x)\}$ is a ray. The support functions of $\p A_t^1,\p A_t^2$ are given by
		\[h(t,\cdot)|_{[\arctan g'_{t,1}(r),\th_t]}\qquad\text{and}\qquad
		h(t,\cdot)|_{[\th_t,\arctan g'_{t,2}(r)]}\]
		respectively, while the support function of $\p\tilde A_t$ coincides with $\varphi(t,\cdot)$. Since $t\mapsto\varphi(t,\th)$ is decreasing for each $\th$, the argument in Claim \ref{claim-heat:1a_curve_comp} implies that $t\mapsto\tsg_t(x)$ is strictly increasing in $t$ for each $x\in\RR$. In $[-\rho_1,\rho_1]$ we clearly have $\sg_{t,1}=g_{t,1}$ and $\sg_{t,2}=g_{t,2}$ and $\tsg_t=\tg_t$, for all $t\in[T_1,0]$, so we may freely switch between notations.
		
		The comparison $\varphi<h$ (with exception $\varphi(0,0)=h(0,0)=0$) implies: $\nexists\,c\geq0$ so that $\p\tilde A_t+(0,c)$ contacts tangentially with $\p A_t^1\cup\p A_t^2$ for $t<0$ or $\p\tilde A_0+(0,c)$ contacts tangentially with $\p A_0^1\cup\p A_0^2$ at any nonzero angle.
		
		\begin{figure}[ht]
			\centering
			\includegraphics{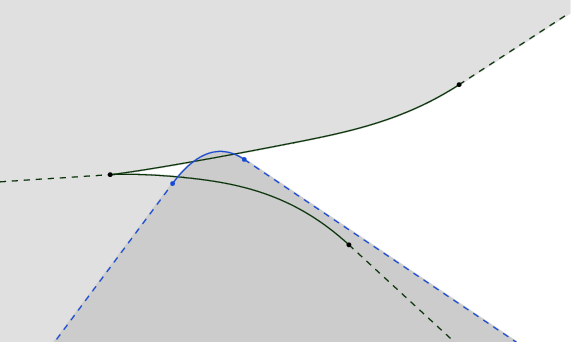}
			\begin{picture}(0,0)
				\put(-182,98){$\tilde\gamma_t$}
				\put(-85,121){$\gamma_t^2$}
				\put(-138,48){$\gamma_t^1$}
				\put(-32,7){$\tsg_t(x)$}
				\put(-69,-11){$\sg_{t,1}(x)$}
				\put(-2,155){$\sg_{t,2}(x)$}
				\put(-263,82){\arrowangle{5}}
				\put(-251,88){$\th_t$}
				\put(-40,128){\arrowangle{35}}
				\put(-27,126){$(\geq\th_t)$}
				\put(-100,18){\arrowangle{-43}}
				\put(-134,7){$(<\!-2\epsilon)$}
				\put(-84,38){\arrowangle{-34}}
				\put(-72,41){$-\delta$}
				\put(-241,32){\arrowangle{51}}
				\put(-239,42){$\delta$}
				\put(-229,88){$x_t$}
				\put(-276,6){$A_t^1$}
				\put(-276,151){$A_t^2$}
				\put(-236,6){$\tilde A_t$}
				\put(-168,102){\rotatebox[origin=c]{-110}{$\xrightarrow{\hspace{24pt}}$}}
				\put(-157,123){$x=\beta_t$}
				\put(-199,98.5){\rotatebox[origin=c]{-70}{$\xrightarrow{\hspace{24pt}}$}}
				\put(-224,120){$x=\alpha_t$}
			\end{picture}
			\vspace{9pt}
			\caption{Objects in Claim \ref{claim-heat:2a_curve_comp}.}\label{fig-heat:curve_comp_2}
		\end{figure}
		
		\begin{claimheatmain}\label{claim-heat:aux_nonempty_int}
			$A_t^2\cap\tilde A_t$ must have a nonempty interior for all $t\in[T_1,0)$.
		\end{claimheatmain}
		\begin{proof}
			Otherwise, we may translate $\p\tilde A_t$ upward so that it meets $\p A_t^2$ tangentially.
		\end{proof}
		
		\begin{claimheatmain}
			$A_t^2\cap\tilde A_t$ is compact for all $t\in[T_1,0]$.
		\end{claimheatmain}
		\begin{proof}
			This follows from the slope comparison $\th_t<\delta$ and $-\delta<\th_t$ at $\pm\infty$.
		\end{proof}
		
		By convexity, there are two intersections of $\p\tilde A_t$ and $\p A_t^2$. Denote their $x$-coordinates by $\alpha_t<\beta_t$, for all $t\in[T_1,0)$. By Claim \ref{claim-heat:aux_nonempty_int} we have $\sg_{t,2}<\tsg_t$ in $(\alpha_t,\beta_t)$.
		
		\begin{claimheatmain}\label{claim-heat:2a_betat_comp}
			$x_t<\beta_t<0$ for all $t\in[T_1,0)$.
		\end{claimheatmain}
		\begin{proof}
			The fact $\beta_t<0$ comes from $\sg_{t,2}>\tsg_t$ in $[0,\infty)$, due to
			\[\sg_{t,2}|_{[0,r]}\geq0,\quad\tsg_t|_{[0,r]}<0,\quad\sg'_{t,2}|_{[r,\infty)}\geq\tan\th_t>-\tan\delta=\tsg'_t|_{[r,\infty)}\qquad\forall\,t\in[T_1,0).\]
			If $\beta_t\leq x_t$ for some $t\in[T_1,0)$, then $\p\tilde A_t$ would intersect $\p A_t^2\cap\{x\leq x_t\}$ twice, hence intersect $\p A_t^1$ at least twice. On the other hand, for a large enough $c>0$, the translated curve $\p\tilde A_t+(0,c)$ only has one intersection with $\p A_t^1$ in $\{x<x_t\}$, due to the slope comparison $-2\epsilon<-\delta$ at $+\infty$. So there exists $c'\geq0$ so that $\p\tilde A_t+(0,c')$ contacts $\p A_t^1$ tangentially, contradicting the support function comparison.
		\end{proof}
		
		\begin{claimheatmain}\label{claim-heat:2a_alphat_comp}
			$\alpha_t<\max\{-\rho_1,x_t\}$ for all $t\in[T_1,0)$.
		\end{claimheatmain}
		\begin{proof}
			Suppose $\alpha_t\geq\max\{-\rho_1,x_t\}$ for some $t\in[T_1,0)$. Then notice that:
			\begin{itemize}[nosep]
				\item $\sg_{t,2}<\tsg_t$ in $(\alpha_t,\beta_t)$, see above Claim \ref{claim-heat:2a_betat_comp};
				\item $s\mapsto\sg_{s,2}(x)$ is decreasing for all $x\in[\alpha_t,\beta_t]$ and $s\leq t$, since $[\alpha_t,\beta_t]\subset[x_t,r]$;
				\item $\tsg_s(x)$ is increasing in $s$ for all $x\in[\alpha_t,\beta_t]$ and $s\leq t$;
				\item $\sg_{T_1,2}(x)>\rho_1$ and $\tsg_{T_1}(x)<0$ for all $x\in[-\rho_1,\rho_1]$, by our choice of $\rho_1$.
			\end{itemize}
			By continuity (Claim \ref{claim-heat:2a_no_jump}), there is $s\in(T_1,t)$ so that $\sg_{s,2}$ contacts tangentially with $\tsg_s$ at some point $x\in(\alpha_t,\beta_t)$. This contradicts the support function comparison.
		\end{proof}
		
		\begin{claimheatmain}\label{claim-heat:2a_right_part}
			$\tsg_t(x)>\sg_{t,1}(x)$ for all $t\in[T_1,0)$ and $x\geq\max\{-\rho_1,x_t\}$.
		\end{claimheatmain}
		\begin{proof}
			If $\tsg_t(x)\leq\sg_{t,1}(x)$, then clearly $\tsg_t(x)\leq\sg_{t,2}(x)$. Hence either $\alpha_t\geq x$ or $\beta_t\leq x$. The former case is ruled out by Claim \ref{claim-heat:2a_alphat_comp}. If the latter case occurs, then note that
			\[\tsg_t>\sg_{t,2}\geq\sg_{t,1}\ \ \text{in}\ \ (\alpha_t,\beta_t),\qquad\tsg_t(x)\leq\sg_{t,1}(x),\qquad\tsg_t(x')>\sg_{t,1}(x')\quad\forall\,x'\gg1.\]
			Then $\p\tilde A_t+(0,c)$ contacts $\p A_t^1$ tangentially for some $c\geq0$, contradiction.
		\end{proof}
		
		\begin{claimheatmain}\label{claim-heat:2a_left_part_weak}
			$\tsg_t(x)\geq f_{\ridge}(x)$ for all $t\in[T_1,0]$ and $x\in[-\rho_1,x_t]$, whenever $x_t\geq-\rho_1$.
		\end{claimheatmain}
		\begin{proof}
			Suppose $\tsg_t(x)<f_{\ridge}(x)$ for some $t\in[T_1,0]$ and $-\rho_1\leq x\leq x_t$. Then $x\leq x_t$ implies $G_t(x)=f_{\ridge}(x)$. Recall that $(t,x)\mapsto G_t(x)$ is continuous in $[T_1,0]\times[-\rho_1,\rho_1]$. Next, recall $x_{T_1}<-\rho_1$ and $g_{T_1,1}<g_{T_1/2,1}<f_{\ridge}$ in $[-\rho_1,\rho_1]$ by our setup, hence $G_{T_1}(x)<f_{\ridge}(x)$. By continuity, there is $t_1\in(T_1,t)$ so that $\tsg_t(x)<G_{t_1}(x)<f_{\ridge}(x)$. The second inequality implies $x>x_{t_1}$, hence $x\geq\max\{-\rho_1,x_{t_1}\}$ and $G_{t_1}(x)=g_{t_1,1}(x)=\sg_{t_1,1}(x)$. By Claim \ref{claim-heat:2a_right_part} we have $\tsg_{t_1}(x)>\sg_{t_1,1}(x)$. Finally, these lead to a contradiction since
			\[G_{t_1}(x)>\tsg_t(x)>\tsg_{t_1}(x)>\sg_{t_1,1}(x)=G_{t_1}(x).\qedhere\]
		\end{proof}
		
		\begin{claimheatmain}\label{claim-heat:2a_left_part}
			$\tsg_t(x)>f_{\ridge}(x)$ for all $t\in[T_1/2,0]$ and $x\in[-\rho_1/2,x_t]$ with $(t,x)\ne(0,0)$.
		\end{claimheatmain}
		\begin{proof}
			Suppose $\tsg_t(x)=f_{\ridge}(x)$ for some $(t,x)$ as stated. Then $x<x_t$ by Claim \ref{claim-heat:2a_right_part} (recall $f_{\ridge}(x_t)=\sg_{t,1}(x_t)$). If $x\leq x_s$ for some $s<t$, then Claim \ref{claim-heat:2a_left_part_weak} implies
			\[\tsg_t(x)>\tsg_s(x)\geq f_{\ridge}(x),\]
			as desired. The only remaining possibility is $\lim_{s\nearrow t}x_s\leq x<x_t$. In this case, note that
			\[\lim_{s\nearrow t}g_{1,s}=\lim_{s\nearrow t}g_{2,s}=f_{\ridge}\qquad\text{on\ \ }[x,x_t],\]
			hence $f_{\ridge}$ is linear on $[x,x_t]$. But $\tsg_t=\tg_t$ is strictly concave in $[-\rho_1,\rho_1]$, and $\tsg_t(x)=f_{\ridge}(x)$, and $\tsg_t\geq f_{\ridge}$ on $[-\rho_1,x_t]$ by Claim \ref{claim-heat:2a_left_part_weak}, which is impossible.
		\end{proof}
		
		\begin{claimheatmain}\label{claim-heat:2a_right_part_2}
			$\tsg_0(x)>\sg_{0,1}(x)$ for all $x\in(0,\rho_1]$.
		\end{claimheatmain}
		\begin{proof}
			Recall $\tsg_0(0)=\sg_{0,1}(0)=\tsg'_0(0)=\sg'_{0,1}(0)=0$ and $\tsg'_0(x)=-\tan\delta>-\tan(2\epsilon)>\sg'_{0,1}(x)$ for $x\gg1$. If $\tg_0(x)\leq g_{0,1}(x)$ for some $x\in(0,\rho_1]$, then $\p\tilde A_0+(0,c)$ would contact $\p A_0^1$ tangentially at a point with nonzero slope for some $c\geq0$, contradiction.
		\end{proof}
		
		Recalling that
		\[\left.\begin{aligned}
			G_t(x) &= f_{\ridge}(x)\qquad\forall\,x\in[-\rho_1,x_t] \\
			G_t(x) &= g_{t,1}(x) = \sg_{t,1}(x)\qquad\forall\,x\in\big[\!\max\{x_t,-\rho_1\},\rho_1\big]\ \ \\
			\tg_t(x) &= \tsg_t(x)\qquad\forall\,x\in[-\rho_1,\rho_1]
		\end{aligned}\right\}\ \forall\,t\in[T_1,0],\]
		Claim \ref{claim-heat:2a_curve_comp} now follows from Claims \ref{claim-heat:2a_right_part}, \ref{claim-heat:2a_left_part}, \ref{claim-heat:2a_right_part_2}, by setting $T''=T_1/2$ and $\rho'=\rho_1/2$. The fact $g_{T'',1}<f_{\ridge}$ in $[-\rho',\rho']$ follows from our setup for $\rho_1$.
	\end{proof}
\end{proof}

\begin{proof}[Proof of Theorem \ref{thm-heat:no_interior_segment}] {\ }
	
	(i) Let $\gamma$ be as stated. We may assume $t=0$, so $\gamma\subset\p Z_0\cap\Omega=\p\{w<0\}\cap\Omega$. Let $\varth$ be an angle parametrization of $\gamma$, and $J_0=\varth(\gamma)$. Fix any $\th\in\Int(J_0)$, denote $\sigma=\varth^{-1}(\th)$.
	
	\begin{claimheatreg}\label{claim-heat:sigma_cpt}
		$\sigma\Subset\Omega$.
	\end{claimheatreg}
	\begin{proof}
		This follows from the connectedness of $\gamma$ and the fact that $\varth:\gamma\to\RR$ is a monotone map (under a length parametrization) and $\th\in\Int(J_0)$.
	\end{proof}
	
	It suffices to show that $\sigma$ is a point, and $\gamma$ is a smooth strictly convex curve near $\sigma$. Once this is done, it follows that
	\[\gamma=\varth^{-1}\big(\Int(J_0)\big)\cup\varth^{-1}\big(\p J_0\big),\]
	where the first set is a $C^\infty$ strictly smooth curve (and is empty iff $\gamma$ is a line segment), and the second set is the union of two (possibly empty) line segments, as desired.
	
	By symmetry, we may assume $\th=0$ and $\sigma=[-l,l]\times\{0\}$ for some $l\geq0$, and $\nu=-\p_y$ on $\sigma$. Denote
	\[P=[-l-r,l+r]\times[-r,r].\]
	Fix $r\ll1$ so that $P\subset\Omega$ and $|\nu+\p_y|<e^{-10}$ in a neighborhood of $P$. Then $\gamma\cap P$ is the graph of a $e^{-9}$-Lipschitz convex function $g_0:[-l-r,l+r]\to[-e^{-9}r,e^{-9}r]$ with $g_0|_{[-l,l]}=0$ and $g_0>0$ elsewhere. Fix
	\[a=\frac12\min\big\{g_0(-l-r),g_0(l+r)\big\},\qquad\delta=\frac12\arctan\min\big\{\!-g'_0(-l-r),g'_0(l+r)\big\}.\]
	Note that $a,\delta>0$. For each $t<0$, let $g_t,g_t^+$ be such that
	\[\p\{w<t\}\cap P=\graphh(g_t),\qquad\p\{w\leq t\}\cap P=\graphh(g_t^+).\]
	Since $\lim_{t\nearrow0}g_t=g_0$ in $C^1([-l-r,l+r])$, there is $\uT<0$ so that
	\begin{equation}\label{eq-heat:no_large_jump}
		g_{2\uT}\leq r/10\ \ \text{in}\ \ [-l-r,l+r],\qquad g_{2\uT}(0)\leq a,
	\end{equation}
	and
	\begin{equation}\label{eq-heat:boundary_slope}
		g'_t(-l-r)<-\tan(1.5\delta),\quad g'_t(l+r)>\tan(1.5\delta)\qquad\forall\,t\in[2\uT,0].
	\end{equation}
	Recall that $g_t,g_t^+$ are convex for all $t\in[\uT,0]$, and satisfy
	\[\lim_{s\nearrow t}g_s=g_t\geq g_t^+,\qquad\lim_{s\searrow t}g_s=\lim_{s\searrow t}g_s^+=g_t^+\qquad\text{in}\ \ C^1([-l-r,l+r]).\]
	Denote
	\[\Gamma_0=\Big(\bigcup_{t\in(\uT,0]}\graphh(g_t)\Big)\cup\Big(\bigcup_{t\in[\uT,0)}\graphh(g_t^+)\Big).\]
	Note that $\Gamma_0$ is compact. Since $|\nu+\p_y|<e^{-10}$, the map $(t,\nu^\perp)$ on $\Gamma_0$ has a unique lift
	\[\Th_0:\Gamma_0\to\RR\times(-e^{-10},e^{-10}).\]
	
	\begin{claimheatreg}\label{claim-heat:preimg_segment}
		Each preimage $\Th_0^{-1}(t,\th)$ in $\Gamma\cap P$ is a line segment.
	\end{claimheatreg}
	\begin{proof}
		Each preimage of $t$ is $\graphh(g_t)\cup\graphh(g_t^+)$. By the property of jumps, the set $\{x:g_t^+(x)<g_t(x)\}$ is either $[-l-r,l+r]$ or is of the form $[-l-r,c)\cup(d,l+r]$ with $-l-r\leq c\leq d\leq l+r$. And $g_t^+$ is linear in each interval where $g_t^+<g_t$. Since
		\[g_t^+(0)<g_{2\uT}(0)\leq a<g_0\big(\!\pm(l+r)\big)<g_t^+\big(\!\pm(l+r)\big),\]
		we have $\{g_t^+<g_t\}=[-l-r,c)\cup(d,l+r]$ with $-l-r\leq c\leq d\leq l+r$, namely, it cannot be $[-l-r,l+r]$. Then note that $g_t=g_t^+$ and $g'_t=(g_t^+)'$ at $c,d$, thus
		\[\Th_0^{-1}(t,\th)=\Big\{(x,g_t(x)):g'_t(x)=\tan\th\Big\}\cup\Big\{(x,g_t^+(x)):(g_t^+)'(x)=\tan\th\Big\}\]
		is a line segment.
	\end{proof}
	
	Denote $\bar\cJ=\Th_0(\Gamma_0)$. Then $\bar\cJ$ is compact since $\Gamma_0$ is compact and $\Th_0$ is continuous. The continuous function $z\mapsto\metric{z}{\nu}$ induces a support function $h:\bar\cJ\to\RR$ by Claim \ref{claim-heat:preimg_segment}, and $h$ is continuous since $\Th_0$ is a quotient map.
	
	To use Theorem \ref{thm-heat:main}, we need to consider the restricted data
	\[\Gamma=\bigcup_{t\in(\uT,0]}\gamma_t\cap\Int(P),\qquad\Th=\Th_0|_\Gamma,\qquad\cJ=\Th(\Gamma).\]
	The set $\Gamma$ here aligns with the one in Theorem \ref{thm-heat:main} with $\Omega:=\Int(P)\cap\{y<g_\uT(x)\}$. We clearly have $\cJ\subset\bar{\cJ}$, and the support function on $\cJ$ is given by $h|_{\cJ}$, which is also continuous in $\cJ$. Consider $R=(\uT,0]\times(-\delta,\delta)$. Then \eqref{eq-heat:boundary_slope} implies $R\subset\cJ$ and $\Th^{-1}(t,\th)\Subset\Omega$ for all $(t,\th)\in R$. Using Theorem \ref{thm-heat:main}\ref{item-heat:interior_viscosity} on the weak IMCF $\min\{w,0\}$ in $\Omega$, and then using the classical viscosity theory, $h$ is a smooth solution of $\p_th=\p_{\th\th}h+h$ in $R$. Hence $h(0,\cdot)$ is smooth in $(-\delta,\delta)$ and satisfies $\p_{\th\th}h+h>0$. Hence $\sigma$ is a point and $\gamma$ is a smooth strictly convex curve near $\sigma$.
	
	\vspace{3pt}
	
	(ii) Denote $\gamma_t=\p\{w<t\}\cap Q_\delta(4)$ and $\gamma_t^+=\p\{w\leq t\}\cap Q_\delta(4)$, so they are $e^{-9}\delta$-Lipschitz convex graphs over $(-4,4)$ whenever they intersect $\bar{Q_\delta(3)}$. Denote $P=[-3,3]\times[-3.5\delta,3.5\delta]$. Let $g_t,g_t^+:[-3,3]\to(-3.5\delta,3.5\delta)$ be so that
	\[\gamma_t\cap P=\graphh(g_t),\qquad\gamma_t^+\cap P=\graphh(g_t^+),\]
	whenever these curves intersect $\bar{Q_\delta(3)}$. Thus $g_t,g_t^+$ are decreasing in $t$, and $g_t\geq g_t^+$, and
	\begin{equation}\label{eq-heat:aux7}
		\lim_{s\nearrow t}g_s=\lim_{s\nearrow t}g_s^+=g_t,\qquad\lim_{s\searrow t}g_s=\lim_{s\searrow t}g_s^+=g_t^+\qquad\text{in}\ \ C^1([-3,3]).
	\end{equation}
	Denote $a_t=\arctan g'_t(-3)$ and $b_t=\arctan g'_t(3)$.
	
	Given any $z_0\in Q_\delta(2)\cap\gamma_{t_0}$ so that $\gamma_{t_0}$ is smooth and strictly convex near $z_0$. We need to show that
	\begin{equation}\label{eq-heat:reg_to_show}
		\text{the curvature of $\gamma_{t_0}$ at $z_0$ is at most $e^{-7}\delta$.}
	\end{equation}
	Once this is done, recalling that $|\D w|$ is the curvature of $\gamma_t$ almost everywhere for almost every $t$, this implies the theorem since $z_0$ is arbitrary.
	
	Let $\th_0=\arctan g'_{t_0}(x_0)$ where $x_0$ is the first coordinate of $z_0$. Consider
	\begin{equation}\label{eq-heat:def_uT}
		\begin{aligned}
			\uT &= \min\Big\{t:\gamma_t^+\cap\bar{Q_\delta(3)}\ne\emptyset,\ \text{$g_s^+$ is not a linear function in} \\
			&\hspace{212pt}
			\text{$[-3,3]$ for all $s\in(t,t_0)$}\Big\}.
		\end{aligned}
	\end{equation}
	By the fact that $z_0\in Q_\delta(2)$ and $g_{t_0}$ is not linear and \eqref{eq-heat:aux7}, this minimum is attained with $\uT<t_0$. The choice of $\uT$ is aimed to avoid large jumps in $(\uT,t_0)$. Denote
	\[\Gamma_0=\Big(\bigcup_{t\in(\uT,t_0]}\graphh(g_t)\Big)\cup\Big(\bigcup_{t\in[\uT,t_0)}\graphh(g_t^+)\Big).\]
	Then $\Gamma_0$ is compact by \eqref{eq-heat:aux7}. Let
	\[\Th_0:\Gamma_0\to(-e^{-20}\delta^2,e^{-20}\delta^2)\times(-e^{-10}\delta,e^{-10}\delta)\]
	be the unique continuous lift of $(t,\nu^\perp)$. The range of $\Th_0$ follows from \eqref{eq-heat:intro_ereg}. Then $\bar\cJ:=\Th_0(\Gamma_0)$ is compact. Arguing as in Claim \ref{claim-heat:preimg_segment} using the nonlinearity condition in \eqref{eq-heat:def_uT}, one may show that each preimage of $\Th_0$ is a line segment in $\Gamma_0$. Hence, the induced support function $h$ is well-defined and is continuous in $\bar\cJ$. We prove a technical lemma:
	
	\begin{claimheatreg}\label{claim-heat:jordan_loop}
		$\p\bar\cJ$ is a Jordan loop.
	\end{claimheatreg}
	\begin{proof}
		For each $t\in[\uT,t_0)$, denote $a_t^+=\arctan (g_t^+)'(-3)$ and $b_t^+=\arctan(g_t^+)'(3)$. Then
		\begin{equation}\label{eq-heat:reg_aux1}
			\lim_{s\searrow t}a_s=\lim_{s\searrow t}a_s^+=a_t^+<b_t^+=\lim_{s\searrow t}b_s^+=\lim_{s\searrow t}b_s,\qquad\forall\,t\in(\uT,t_0),
		\end{equation}
		where $a_t^+<b_t^+$ follows from the nonlinearity in the choice of $\uT$, and
		\begin{equation}\label{eq-heat:reg_aux2}
			\lim_{s\searrow\uT}a_s=\lim_{s\searrow\uT}a_s^+=a_\uT^+\leq b_\uT^+=\lim_{s\searrow\uT}b_s^+=\lim_{s\searrow\uT}b_s,
		\end{equation}
		and $a_t\leq a_t^+$, $b_t\geq b_t^+$ for all $t\in(\uT,t_0)$ due to the geometry of jumps, and
		\begin{equation}\label{eq-heat:reg_aux3}
			\lim_{s\nearrow t}a_s=\lim_{s\nearrow t}a_s^+=a_t<b_t=\lim_{s\nearrow t}b_s^+=\lim_{s\nearrow t}b_s,\qquad\forall\,t\in(\uT,t_0].
		\end{equation}
		Observe that $\p\bar\cJ$ is the union of the two connected compact sets
		\[\begin{aligned}
			\Gamma_1 &= \Big(\bigcup_{t\in(\uT,t_0)}\{t\}\times[a_t,a_t^+]\Big)\cup\big\{(\uT,a_\uT^+),(t_0,a_{t_0})\big\}, \\
			\Gamma_2 &= \Big(\bigcup_{t\in(\uT,t_0)}\{t\}\times[b_t^+,b_t]\Big)\cup\big\{(\uT,b_\uT^+),(t_0,b_{t_0})\big\},
		\end{aligned}\]
		and the two segments $\bar J_{t_0}$, $\bar J_\uT$, where these portions only intersect at endpoints ($\bar\cJ_\uT$ may be a point). Thus, it is sufficient to show that $\Gamma_1,\Gamma_2$ are homeomorphic to $[0,1]$.
		
		We prove the case for $\Gamma_1$. Define an order on $\Gamma_1$, where $(s,\omega)\prec(t,\th)$ if $s<t$ or $s=t$ and $\omega<\th$. Then $\prec$ is a linear order, and $(\uT,a_\uT^+)$ is minimal, and $(t_0,a_{t_0})$ is maximal. By the convergences \eqref{eq-heat:reg_aux1} \eqref{eq-heat:reg_aux2} \eqref{eq-heat:reg_aux3}, $\Gamma_1\cap R$ is open in the order topology for any open box $R$ in $\RR^2$. Let us show that any open interval in the order topology is open in the subspace topology. Given $(t,\th)\in\Gamma_1$ and an open interval $I\ni(t,\th)$.
		\begin{itemize}[nosep]
			\item If $t=\uT$, then $I$ contains $\big((-\infty,t+\epsilon)\times\RR\big)\cap\Gamma_1$ for some $\epsilon>0$.
			\item If $t=t_0$, then $I$ contains $\big((t-\epsilon,\infty)\times\RR\big)\cap\Gamma_1$ for some $\epsilon>0$.
			\item Suppose $t\in(\uT,t_0)$. If $a_t=a_t^+$, then $I\supset\big((t-\epsilon,t+\epsilon)\times\RR\big)\cap\Gamma_1$ for some $\epsilon>0$.
			\item If $t\in(\uT,t_0)$ and $a_t\leq\th<a_t^+$, then $I$ contains $\big((t-\epsilon,t)\times\RR\big)\cap\Gamma_1$ and $\big(\{t\}\times(-\infty,\th+2\epsilon)\big)\cap\Gamma_1$, for some $\epsilon>0$. But \eqref{eq-heat:reg_aux1} implies $\big([t,t+\epsilon')\times(-\infty,\th+\epsilon)\big)\cap\Gamma_1\subset\{t\}\times(-\infty,\th+\epsilon)$ for some $\epsilon'>0$. Hence $I\supset\big((t-\epsilon,t+\epsilon')\times(-\infty,\th+\epsilon)\big)\cap\Gamma_1$.
			\item If $t\in(\uT,t_0)$ and $a_t<\th\leq a_t^+$, then arguing similarly, we may show that $I$ contains $\big((t-\epsilon,t+\epsilon')\times(\th-\epsilon,\infty)\big)\cap\Gamma_1$ for some small $\epsilon,\epsilon'>0$.
		\end{itemize}
		Hence, the order topology coincides with the subspace topology. Therefore, $\Gamma_1$ with the order topology is connected, compact, metrizable, with endpoints. Hence $\Gamma_1\cong[0,1]$, by \cite[Theorem 6.16, 6.17]{Nadler}.
	\end{proof}
	
	To use Theorem \ref{thm-heat:main}, we again need to consider the restricted data
	\[\begin{aligned}
		& \Gamma=\bigcup_{t\in(\uT,t_0]}\gamma_t\cap\Int(P), \qquad \Th=\Th_0|_\Gamma,\qquad \Omega=\Int(P)\cap\{y<g_\uT(x)\} \\
		&\hspace{40pt} J_t=\Th\big(\gamma_t\cap\Int(P)\big),\qquad\cJ=\Th(\Gamma)=\bigcup_{t\in(\uT,t_0]}\{t\}\times J_t.
	\end{aligned}\]
	Then $\cJ\subset\bar\cJ$, and $\Th$ is a lift of $(t,\nu^\perp)$ on $\Gamma$. The $t$-slice of $\bar\cJ$ is equal to $[a_t,b_t]$, which also coincides with the closure of $J_t$. We thus denote it by $\bar J_t$. Note that $\bar J_t$ is not a point for all $t\in(\uT,t_0]$: otherwise, $\gamma_t\cap\Int(P)$ would be a line segment, hence $\gamma_t^+\cap\Int(P)$ is also a line segment, contradicting \eqref{eq-heat:def_uT}.
	
	\begin{claimheatreg}\label{claim-heat:backward_open}
		If $t\in(\uT,t_0]$ and $\th\in\Int(J_t)$, then there is $\epsilon>0$ so that
		\[R\subset\cJ,\qquad\Th^{-1}(t,\th)\Subset\Omega\ \ \forall\,(t,\th)\in R,\qquad\text{where}\qquad R:=(t-\epsilon,t]\times(\th-\epsilon,\th+\epsilon).\]
	\end{claimheatreg}
	\begin{proof}
		The conditions imply $g'_t(-3)<\tan\th<g'_t(3)$ and $\graphh(g_t)\cap Q_\delta(3)\ne\emptyset$. Since $\lim_{s\nearrow t}g_s=g_t$ in $C^1([-3,3])$, there is $\epsilon<t-\uT$ so that $g'_s(-3+\epsilon)<\tan(\th-\epsilon)$ and $g'_s(3-\epsilon)>\tan(\th+\epsilon)$ and $\graphh(g_s)\cap Q_\delta(3)\ne\emptyset$ for all $s\in[t-\epsilon,t]$, as desired.
	\end{proof}
	
	Using Theorem \ref{thm-heat:main}(i) in $\Omega$ with the data $\Gamma,\Th,\cJ$, and using Claim \ref{claim-heat:backward_open}, we have that $h$ is a smooth solution of $\p_th=\p_{\th\th}h+h$ in a backward box near each $t\in(\uT,t_0]$ and $\th\in\Int(J_t)$. Hence, $h(t,\cdot)$ is smooth and satisfies $\p_{\th\th}h+h>0$ in each $\Int(J_t)$. Recall that $\p_{\th\th}h+h$ is the inverse of curvature.
	
	Suppose the curvature bound \eqref{eq-heat:reg_to_show} fails. Recall $(t_0,\th_0)=\Th(z_0)$, so $(\p_{\th\th}h+h)(t_0,\th_0)<e^7\delta^{-1}$. We need to obtain a contradiction. Consider
	\[\phi(t,\th)=h(t_0,\th_0)\cos(\th-\th_0)+\p_\th h(t_0,\th_0)\sin(\th-\th_0)+\frac{e^7e^{t-t_0}}{\delta}\Big((t-t_0)+\frac12(\th-\th_0)^2\Big),\]
	which solves $\p_t\phi=\p_{\th\th}\phi+\phi$. We first obtain some preliminary estimates. The conditions $z_0\in Q_\delta(2)$ and $|t|,|t_0|<e^{-20}\delta^2$ and $|\th|,|\th_0|<e^{-10}\delta$ and $\delta\leq1$ imply
	\begin{align}
		|h(t_0,\th_0)| &= |\metric{z_0}{\nu_{\th_0}}|
			\leq |x_0|\,|\sin\th_0|+|y_0|\,|\cos\th_0|
			\leq 2e^{-10}\delta+2\delta\leq2.1\delta, \label{eq-heat:h_basept}\\
		|\p_\th h(t_0,\th_0)| &= |\metric{z_0}{\tau_{\th_0}}|
			\leq |x_0|\,|\cos\th_0|+|y_0|\,|\sin\th_0|\leq2+2e^{-10}\delta^2\leq2.1, \label{eq-heat:dh_basept}
	\end{align}
	and also
	\begin{equation}\label{eq-heat:pth_varphi}
		\begin{aligned}
			|\p_\th\phi| &\leq |h(t_0,\th_0)|\,|\sin(\th-\th_0)|+|\p_\th h(t_0,\th_0)|+e^7e^{t-t_0}\delta^{-1}|\th-\th_0| \\
			&\leq 2.1\delta\cdot2e^{-10}\delta+2.1+e^7e^{2e^{-20}\delta^2}\delta^{-1}\cdot 2e^{-10}\delta \\
			&\leq 2.3, \qquad\qquad\forall\,(t,\th)\in(-e^{-20}\delta^2,e^{-20}\delta^2)\times(-e^{-10}\delta,e^{-10}\delta).
		\end{aligned}
	\end{equation}
	
	Denoting $\eta=\phi-h$, we have $\eta(t_0,\th_0)=\p_\th\eta(t_0,\th_0)=0$ by construction, and
	\[\p_\th^2\eta(t_0,\th_0)=-\p_\th^2 h(t_0,\th_0)-h(t_0,\th_0)+e^7\delta^{-1}>0.\]
	So there exist $t_1\in(\uT,t_0)$ and $\th_1<\th_2<\th_3\in\Int(J_{t_1})$ with
	\[\eta(t_1,\th_1)>0,\qquad \eta(t_1,\th_2)<0,\qquad \eta(t_1,\th_3)>0.\]
	For any $\epsilon>0$, consider the function
	\[\psi=e^{-\epsilon t-t}\eta.\]
	
	\begin{claimheatreg}\label{claim-heat:no_interior_extremal}
		Suppose $L$ is relatively open in $\bar\cJ$. Let $\p'L$ be the relative boundary of $L$ in $\bar\cJ$.
		\begin{enumerate}[label={(\roman*)}, nosep]
			\item If $\psi\leq0$ in $\bar L$ and $\min_{\bar L}\psi<\min_{\p'L}\psi$, then $\psi$ cannot be minimized in $\bar L$ at a point $(t,\th)$ with $t>\uT$ and $\th\in\Int(\bar{J_t})$.
			\item If $\psi\geq0$ in $\bar L$ and $\max_{\bar L}\psi>\max_{\p' L}\psi$, then $\max_{\bar L}\psi$ cannot be attained in $\{t>\uT\}$.
		\end{enumerate}
	\end{claimheatreg}
	\begin{proof}
		Note that
		\[\p_t\psi-\p_{\th\th}\psi=e^{-\epsilon t-t}\big(\p_t\eta-\p_{\th\th}\eta-\epsilon\eta-\eta\big)=-\epsilon\psi\]
		in each backward box near $(t,\th)$ with $t\in(\uT,t_0]$, $\th\in\Int(\bar J_t)$. Then (i) follows from the standard maximum principle. Suppose that $\psi$ is maximized in $\bar L$ at a point $(t,\th)$ with $t>\uT$. Then again we must have $\th\in\p\bar{J_t}$. Recall that we defined $a_t=\arctan g'_t(-3)$ and $b_t=\arctan g'_t(3)$, so $\bar{J_t}=[a_t,b_t]$. We may assume $\th=a_t$ (same argument for the case $\th=b_t$). The maximizing condition implies
		\[\p_\th\phi(t,a_t)\leq\p_\th^+ h(t,a_t).\]
		The quantity $\p_\th^+ h(t,a_t)$ makes sense since $\bar J_t$ is not a point (see above Claim \ref{claim-heat:backward_open}). Combined with \eqref{eq-heat:pth_varphi}, we have $\p_\th^+ h(t,a_t)\geq-2.3$. Observe that $h(t,\cdot)$ coincides with the support function of $\gamma_t\cap P$, thus \eqref{eq-prelim:recovery_formula_endpt} gives
		\[\Th_0^{-1}(t,a_t)\cap\gamma_t\cap P=\Big\{h(t,a_t)\nu_{a_t}+c\tau_{a_t}:c\leq\p_\th^+ h(t,a_t)\Big\}\cap P.\]
		The facts $\p_\th^+ h(t,a_t)\geq-2.3$ and $|\nu+\p_y|<e^{-10}\delta$ imply that $\Th_0^{-1}(t,a_t)\cap\gamma_t\cap P$ is a nontrivial line segment in $\gamma_t$. Hence $a_t\in J_t$ as well. Denoting $\sigma=\Th^{-1}(t,a_t)\cap\Int(P)\ne\emptyset$, we have
		\begin{equation}\label{eq-heat:reg_preim_of_Th}
			\sigma=\Big\{h(t,a_t)\nu_{a_t}+c\tau_{a_t}:c\leq\p_\th^+h(t,a_t)\Big\}\cap\Int(P).
		\end{equation}
		And clearly, we have $\sigma\not\Subset\Omega$.
		
		We are ready to apply the boundary viscosity principle in Theorem \ref{thm-heat:main}\ref{item-heat:boundary_viscosity}. Denote $\lambda=\psi(t,a_t)=\max_{\bar L}\psi>0$, and set $\varphi=\phi-e^{\epsilon t+t}\lambda$. We clearly have
		\[\varphi(t,a_t)=h(t,a_t).\]
		Since $L$ is relatively open in $\bar\cJ$, and $(t,\th)\notin\p'L$ by our assumption, there is $\delta>0$ so that
		\[\varphi\leq h\qquad\text{in}\ \ \Big([t-\delta,t]\times[\th-\delta,\th+\delta]\Big)\cap\bar\cJ.\]
		Denote $z=\varphi(t,a_t)\nu_{a_t}+\p_\th\varphi(t,a_t)\tau_{a_t}$. We need to show $z\in\sigma$. Since $\varphi(t,a_t)=h(t,a_t)$ and $\p_\th\varphi(t,a_t)=\p_\th\phi(t,a_t)\leq\p_\th^+ h(t,a_t)$, $z$ indeed lies in the first part of \eqref{eq-heat:reg_preim_of_Th}. To show $z\in\Int(P)$, we calculate
		\[\begin{aligned}
			|\varphi(t,a_t)|=|h(t,a_t)| &= \max_{z\in\sigma}|\metric{z}{\nu_{a_t}}| \\
			&\leq \max_{(x,y)\in\sigma}\Big(|x|\,|\sin a_t|+|y|\Big)\leq3\cdot e^{-10}\delta+3.1\delta\leq3.2\delta,
		\end{aligned}\]
		where $|y|\leq3.1\delta$ follows from $\gamma_t\cap\bar{Q_\delta(3)}\ne\emptyset$ and $e^{-9}\delta$-Lipschitzness. We also have
		\[|\p_\th\varphi(t,a_t)|=|\p_\th\phi(t,a_t)|\leq2.3,\]
		these apply to show $z=\varphi(t,a_t)\nu_{a_t}+\p_\th\varphi(t,a_t)\tau_{a_t}\in\Int(P)$ in view that $|\nu_{a_t}+\p_y|\leq e^{-10}\delta$. Also, recall $\chi_\sigma=1$ from our sign conditions. Finally, we may apply Theorem \ref{thm-heat:main}\ref{item-heat:boundary_viscosity} to the solution $\min\{w,t_0\}$ in $\Omega$, and obtain $(\p_t\varphi-\p_{\th\th}\varphi-\varphi)(t,a_t)\geq0$. This contradicts
		\[\p_t\varphi-\p_{\th\th}\varphi-\varphi=\underbrace{\p_t\phi-\p_{\th\th}\phi-\phi}_{=0}-\epsilon e^{\epsilon t+t}\lambda<0. \qedhere\]
	\end{proof}
	
	Let $V_1,V_3\subset\bar\cJ$ be the connected components of $\{\eta>0\}$ containing $(t_1,\th_1)$ and $(t_1,\th_3)$ respectively, and $V_2\subset\bar\cJ$ be the connected component of $\{\eta<0\}$ containing $(t_1,\th_2)$.
	
	\begin{claimheatreg}\label{claim-heat:preservation_of_root}
		$V_i\cap \bar J_t\ne\emptyset$ for all $t\in[\uT,t_1)$ and $i=1,2,3$.
	\end{claimheatreg}
	\begin{proof}
		If $V_1\cap\bar J_T=\emptyset$ for some $T\in[\uT,t_1)$, then by connectivity we have $V_1\subset\{t>T\}$. Since $\eta>0$ in $V_1$ and $\eta=0$ on $\p'V_1$ and on $\bar{V_1}\cap\{t=\uT\}$ (the latter set is nonempty only if $t=\uT$), there is a small $\epsilon$ so that $\max_{\bar{V_1}}\psi$ is attained in $V_1\cap\{t>\uT\}$. Then Claim \ref{claim-heat:no_interior_extremal}(ii) with $L=V_1$ yields a contradiction. The case of $V_3$ may be proved in the same manner.
		
		If $V_2\cap\bar J_T=\emptyset$ for some $T\in[\uT,t_1]$, then for some small $\epsilon>0$, the minimum of $\psi$ in $\bar{V_2}$ is attained at some $(T',\th)$ with $T'>\uT$. By Claim \ref{claim-heat:no_interior_extremal}(i) we have $\th\in\p\bar J_{T'}$. But topologically this cannot happen, since $V_1,V_3$ are connected and they both intersect $\{t=\uT\}$.
	\end{proof}
	
	By \eqref{eq-heat:def_uT}, at least one of the possibilities below occurs at time $\uT$:
	\begin{itemize}[nosep]
		\item $\gamma_t^+\cap\big([-3,3]\times[-3.2\delta,3.2\delta]\big)=\emptyset$ for all $t<\uT$;
		\item $g_{t_i}^+$ is linear in $[-3,3]$ for a sequence $t_i\nearrow\uT$;
		\item $\gamma_{\uT}^+\cap\bar{Q_\delta(3)}\ne\emptyset$ but $\gamma_t^+\cap\bar{Q_\delta(3)}=\emptyset$ for all $t<\uT$.
	\end{itemize}
	If the first case happens, then $g_\uT^+$ is linear on $[-3,3]$, thus $\bar J_\uT$ is a point, contradicting Claim \ref{claim-heat:preservation_of_root}. If the second case happens, then $g_\uT=\lim_{i\to\infty}g_{t_i}^+$ is linear on $[-3,3]$, hence $g_\uT^+$ is linear on $[-3,3]$, which is again impossible. So the first two possibilities cannot happen. Since $g_\uT^+$ is nonlinear, we have $g_\uT=g_\uT^+$ somewhere in $[-3,3]$, hence
	\[\min_{[-3,3]}g_\uT^+\geq\max_{[-3,3]}g_\uT-0.1\delta\]
	by $e^{-9}\delta$-Lipschitzness. Recall that
	\[\lim_{t\nearrow\uT}g_t^+=g_\uT,\qquad\lim_{t\searrow\uT}g_t=g_\uT^+.\]
	This and $\gamma_t^+\cap\bar{Q_\delta(3)}=\emptyset$ for all $t<\uT$ implies $g_\uT\geq3\delta$, hence $g_\uT^+\geq2.9\delta$, hence $g_t>2.8\delta$ for all $t>\uT$ close enough to $\uT$. For these $t$ and all $\th\in\bar{J_t}$, we have
	\[\begin{aligned}
		h(t,\th) &= \metric{z}{\nu_\th}\qquad\text{(for any $z=(x,y)\in\Th^{-1}(t,\th)$)} \\
		&= x\sin\th-y\cos\th
		\leq 3\cdot e^{-10}\delta-2.8\delta(1-e^{-10}\delta)<-2.7\delta.
	\end{aligned}\]
	On the other hand, by \eqref{eq-heat:h_basept} \eqref{eq-heat:dh_basept} we have
	\[\phi(t,\th)\geq-2.1\delta-2.1\cdot2e^{-10}\delta-e^7e^{2e^{-20}\delta^2}\delta^{-1}\big(2e^{-20}\delta^2+2e^{-20}\delta^2\big)>-2.2\delta\]
	for all $(t,\th)\in(-e^{-20}\delta^2,e^{-20}\delta^2)\times(-e^{-10}\delta,e^{-10}\delta)$. Hence $\phi(t,\cdot)>h(t,\cdot)$ in $\bar{J_t}$, contradicting $V_2\cap\bar J_t\ne\emptyset$. This concludes the proof.
\end{proof}

%\newpage

\section{\texorpdfstring{$C^{1,1/3}$}{C\^(1,1/3)} estimates for almost linear solutions}\label{sec:ereg}

For $\delta,r>0$, consider the anisotropic box
\[Q_\delta(r)=(-r,r)\times(-\delta r,\delta r).\]

\begin{theorem}\label{thm-ereg:main_infhar}
	Suppose $\delta\leq1$, and $u$ is $\infty$-harmonic in $Q_\delta(4)$ where it holds
	\begin{equation}\label{eq-ereg:intro_ereg_du}
		\Big|\frac{\D u}{|\D u|}-\p_x\Big|<e^{-15}\delta,\qquad\big|\log|\D u|\big|<e^{-30}\delta^2.
	\end{equation}
	Then
	\begin{equation}\label{eq-ereg:intro_main_1}
		\Big\|\frac{\D u}{|\D u|}-\p_x\Big\|_{C^{0,1/3}(Q_\delta(1))}\leq e^{10}\delta^{2/3},\qquad
		\big\|\log|\D u|\big\|_{C^{0,2/3}(Q_\delta(1))}\leq e^{10}\delta^{4/3},
	\end{equation}
	and hence
	\begin{equation}\label{eq-ereg:intro_main_2}
		\big\|\D u-\p_x\big\|_{C^{0,1/3}(Q_\delta(1))}\leq e^{15}\delta^{2/3}.
	\end{equation}
\end{theorem}

It was already observed by Evans-Savin \cite{Evans-Savin_2008} that boxes of ratio $1\times\delta$ are suitable for an optimal $\epsilon$-regularity theorem for the $\infty$-Laplacian in 2D. Theorem \ref{thm-ereg:main_infhar} and \cite[Theorem 5.2]{Evans-Savin_2008} together yield a clean form:
\begin{theorem}\label{thm-ereg:es_combined}
	There is a universal constant $C$ so that: if $\delta\leq1$, and $u$ is $\infty$-harmonic in $Q_\delta(4)$ in which
	\begin{equation}\label{eq-ereg:intro_ereg_u}
		|u-x|<C^{-1}\delta^2,
	\end{equation}
	then
	\[\Big\|\frac{\D u}{|\D u|}-\p_x\Big\|_{C^{0,1/3}(Q_\delta(1))}\leq e^{15}\delta^{2/3},\qquad
	\big\|\log|\D u|\big\|_{C^{0,2/3}(Q_\delta(1))}\leq e^{15}\delta^{4/3},\]
	and hence
	\[\|u-x\|_{C^{1,1/3}(Q_\delta(1))}\leq e^{20}\delta^{2/3}.\]
\end{theorem}

The order on $\delta$ is sharp in these results: considering $u_0=|x|^{4/3}-|y|^{4/3}$ and its scaling
\[u_\delta(x,y)=\frac34\frac1{\delta^2}\Big[u_0(1+\delta^2x,\delta^2y)-1\Big],\]
it can be verified that
\[\|u_\delta-x\|_{C^0(Q_\delta(1))}=O(\delta^2),\]
and
\[\begin{aligned}
	&\hspace{18pt} \Big\|\frac{\D u_\delta}{|\D u_\delta|}-\p_x\Big\|_{C^0(Q_\delta(1))}=O(\delta),\qquad
	\big\|\log|\D u_\delta|\big\|_{C^0(Q_\delta(1))}=O(\delta^2), \\
	& \Big\|\frac{\D u_\delta}{|\D u_\delta|}-\p_x\Big\|_{C^{0,1/3}(Q_\delta(1))}=O(\delta^{2/3}),\qquad\big\|\log|\D u_\delta|\big\|_{C^{0,2/3}(Q_\delta(1))}=O(\delta^{4/3}).
\end{aligned}\]
Thus, all conditions and conclusions above have the consistent order of $\delta$ as in $u_\delta$.

By Theorem \ref{thm-ereg:es_combined} and the known $C^1$ regularity, we have:

\begin{theorem}\label{thm-ereg:away_from_crit}
	If $\Delta_\infty u=0$ in $\Omega$, then
	\[u\in C^{1,1/3}_{\loc}\big(\Omega\setminus\Crit(u)\big),\qquad|\D u|\in C^{0,2/3}_{\loc}\big(\Omega\setminus\Crit(u)\big).\]
\end{theorem}

Theorem \ref{thm-ereg:main_infhar} is proved by combining $p$-harmonic approximation (Theorem \ref{thm-approx:global}) and the following regularity of simple clusters:

\begin{theorem}\label{thm-ereg:main_cluster}
	Let $\delta\leq1$, and $\big(w,\nu,\{D_i\},\{\chi_i\}\big)$ be a simple IMCF cluster in $Q_\delta(4)$ with
	\[|\nu+\p_y|< e^{-10}\delta,\qquad|w|< e^{-20}\delta^2.\]
	Then
	\begin{equation}\label{eq-ereg:intro_cluster_w}
		\|w\|_{C^{0,2/3}(Q_\delta(1/4))}\leq \delta^{4/3},
	\end{equation}
	and for any ridge $\gamma$ intersecting $Q_\delta(1/8)$, we have
	\begin{equation}\label{eq-ereg:intro_cluster_gamma}
		\|\gamma\|_{C^{1,1/2}(Q_\delta(1/4))}\leq \delta.
	\end{equation}
\end{theorem}

The deduction from Theorem \ref{thm-ereg:main_cluster} to \ref{thm-ereg:main_infhar} is almost immediate. The following analytic lemma is a main ingredient:

\begin{lemma}[$=$\,Lemma \ref{lemma-prelim:two_graphs}]\label{lemma-ereg:two_graphs} {\ }
	
	Suppose $f,g$ are functions on $(-1,1)$, with $f'=g'$ whenever $f=g$, and
	\[\sup_{(-1,1)}|f'|\leq1,\qquad\sup_{(-1,1)}|g'|\leq1,\]
	and
	\[\|f'\|_{C^{0,\alpha}(-1,1)}\leq\mu,\qquad\|g'\|_{C^{0,\alpha}(-1,1)}\leq\mu,\]
	for some $\alpha\in(0,1]$ and $\mu>0$. Let $\nu$ be the (well-defined) upper normal vector field of $\graphh(f)\cup\graphh(g)$. For $s\leq1$, denote $S=\big(\graphh(f)\cup\graphh(g)\big)\cap\{|x|\leq 1-s\}$. Then
	\[\|\nu\|_{C^{0,\alpha/(1+\alpha)}(S)}\leq 16\mu^{\frac1{1+\alpha}}s^{-\alpha}.\]
\end{lemma}

The H\"older exponent $\frac{\alpha}{1+\alpha}$ is sharp, and it serves to turn the order $\frac12$ in \eqref{eq-ereg:intro_cluster_gamma} to the order $\frac13=\frac{1/2}{1+1/2}$ in the first expression of \eqref{eq-ereg:intro_main_1}.

\vspace{3pt}

We now sketch the main ideas for proving Theorem \ref{thm-ereg:main_cluster}. Recall that in a simple cluster, the ``level sets'' $\gamma_t=\p Z_t$ are unions of concave\,/\,convex curves with cusps at endpoints, see Lemma \ref{lemma-cluster:local_model}. The key to showing Theorem \ref{thm-ereg:main_cluster} is to bound the $C^{1,1/2}$ norms of these curves up to the cusps. This is done in Claim \ref{claim-ereg:dth3_h}\,--\,\ref{claim-ereg:reg_gt_E4} below. Once this is shown, we note that $\gamma_t$ are tangential to the ridge at the cusps for each $t$, therefore, the ridge may be viewed as an envelope curve of the family $\{\gamma_t\}$. The $C^{1,1/2}$ regularity of $\gamma_t$ would imply the $C^{1,1/2}$ regularity of the ridge. This is done as a part of Claim \ref{claim-ereg:reg_gt_E4}. The regularity of $w$ in \eqref{eq-ereg:intro_cluster_w} follows in a rather routine way, thus we refer to Claim \ref{claim-ereg:reg_w_E4} for the proof.

\vspace{3pt}

\begin{figure}[ht]
	\centering
	\includegraphics{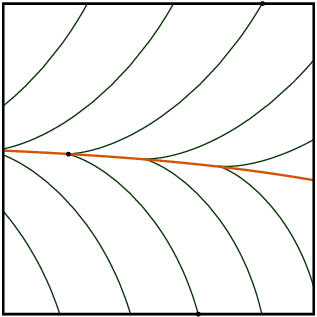}
	\begin{picture}(0,0)
		\put(-59,6){$\uth_t$}
		\put(-29,136){$\oth_t$}
		\put(-122,83){$\rth_t$}
	\end{picture}
	\caption{A smooth cuspidal IMCF.}\label{fig-ereg:intro1}
\end{figure}

To show the idea for proving $C^{1,1/2}$ regularity of $\gamma_t$, let us first consider the model case where the simple cluster is a smooth IMCF by cuspidal curves (Figure \ref{fig-ereg:intro1}). For each time $t$, the cuspidal curve $\gamma_t$ has an angle parameter $\varth$ whose range is denoted by $(\uth_t,\oth_t)$. Let $\rth_t$ be the angle at the cusp in $\gamma_t$.
Then the support function $h(t,\th)$ solves $\p_th=\p_{\th\th}h+h$ in the generalized cylindrical domain $\cJ=\bigcup\,\{t\}\times(\uth_t,\oth_t)$. The quantity
\[\kappa^{-1}=\p_{\th\th}h+h\]
is the signed inverse curvature of $\gamma_t$, with $\kappa^{-1}<0$ in $(\uth_t,\rth_t)$ and $\kappa^{-1}>0$ in $(\rth_t,\oth_t)$ due to our sign conventions. Hence the root $\rth_t$ of $\p_{\th\th}h+h$ is nondegenerate. The parabolic Hopf lemma then implies
\begin{equation}\label{eq-ereg:intro_hopf}
	|\kappa^{-1}|\geq c|\th-\rth_t|
\end{equation}
for some $c>0$. Let $s$ be the length parameter of $\gamma_t$, with $s=0$ corresponding to the cusp. Then we obtain the ODE inequality
\[|d\th/ds|=\frac1{|\kappa^{-1}|}\leq\frac1{c|\th-\rth_t|},\]
which implies
\[\Big|\frac{d(\th-\rth_t)^2}{ds}\Big|\leq2/c\]
and hence
\[|\th(s_1)-\th(s_2)|\leq\sqrt{2/c}\cdot|s_1-s_2|^{1/2}.\]
Hence the edges of $\gamma_t$ are $C^{1,1/2}$ up to the cusp.

There are two main issues in extending this argument to full generality: we need a uniform lower bound on the constant $c$, and to treat the analytic subtleties due to the low regularity of weak IMCF. To resolve the first issue, we observe that the Sturmian principle works as a quantitative Hopf lemma. Similar idea was also employed by Angenent \cite{Angenent_1991}. The following is a toy case for explaining the idea:

\begin{lemma}\label{lemma-ereg:hopf_instance}
	Suppose $f(t,x)$ solves $\p_tf=\p_{xx}f$ smoothly in $[0,1]\times[0,1]$, with $f(t,0)\leq-1$ and $f(t,1)\geq1$ for all $t\in[0,1]$, and $\p_x f(0,x)\geq1$ for all $x\in[0,1]$. Then
	\[\p_x f(t,x)\geq\frac12,\qquad\forall\,(t,x)\text{ with $|f(t,x)|\leq1/2$}.\]
\end{lemma}
\begin{proof}
	If this is not true for some $(t_0,x_0)$, then there exists $c<1/2$ so that $(t_0,x_0)$ is a degenerate root of
	\[g(t,x)=f(t,x)-f(t_0,x_0)-c(x-x_0).\]
	Note that $g$ has no root on $\{x=0\}\cup\{x=1\}$. By the Sturmian theory, $g(0,\cdot)$ must have at least 2 roots in $[0,1]$. But this contradicts $\p_x g(0,x)>0$.
\end{proof}

This argument is applied to show a statement of the form
\[\text{$\p_{\th}\kappa^{-1}\geq c$ for all $(t,\th)$ with $|\kappa^{-1}|(t,\th)\leq c'$,}\]
where $c,c'$ are universal constants. See Claim \ref{claim-ereg:dth3_h}. This replaces the Hopf Lemma \eqref{eq-ereg:intro_hopf}.

\begin{figure}[ht]
	\centering
	\includegraphics{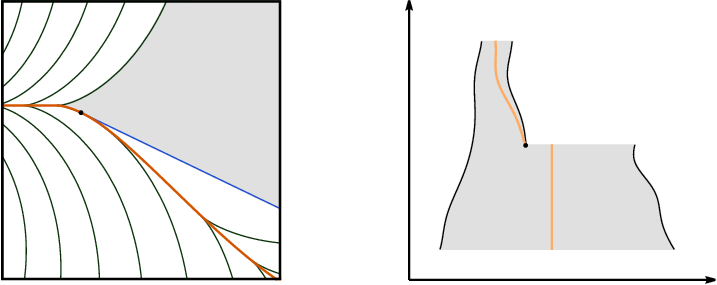}
	\begin{picture}(0,0)
		\put(-305,82){$\th_0$}
		\put(-271,8){$\uth_t$}
		\put(-268,123){$\oth_t$}
		\put(-134,65){$\uth_t$}
		\put(-36,51){$\oth_t$}
		\put(-104,58){$\th_0$}
		\put(-80,39){$\rth_t$}
		\put(-3,5){$\th$}
		\put(-148,136){$t$}
	\end{picture}
	\caption{Discontinuity of the nodal line of $\kappa^{-1}$.}\label{fig-ereg:mixed}
\end{figure}

The second issue (low regularity of the weak IMCF) also causes significant difficulty. Consider Figure \ref{fig-ereg:mixed} (left), where a jump in the weak IMCF is present and partially overlaps with the ridge. Figure \ref{fig-ereg:mixed} (right) describes the image $\cJ$ of the angle parameter. The orange curve represents the angle of the cusps -- namely, the set of points $(t,\th)$ where $\kappa^{-1}(t,\th)=0$. As an effect of the jump, the right boundary of $\cJ$ changes discontinuously in $t$, and the quantity $\kappa^{-1}=\p_{\th\th}h+h$ is not defined at $(t,\th_0)$. In general, $\kappa^{-1}$ may lose continuity near $\p\cJ$ at each jump, and we do not have a strong enough boundary condition for $\kappa^{-1}$ to avoid nodal lines ending on the two sides of $\p\cJ$.

The resolution of this issue is not to compare $\p_{\th\th}h+h$ with a linear function, but to compare $h$ with a degree 3 Taylor polynomial, in the use of Sturmian principle. One may show that $h$ remains continuous up to $\p\cJ$. The rough idea becomes the following: if $\p_{\th}(\p_{\th\th}h+h)$ is too small at a point $(t_0,\th_0)$, then for a suitable Taylor polynomial $\phi$, the function $\phi-h$ would have 3 roots for all $t<t_0$. In the interior of $\cJ$, this is guaranteed by the Sturmian principle. On the two sides of $\cJ$, Theorem \ref{thm-heat:main}\ref{item-heat:boundary_viscosity} would guarantee that a nodal line cannot end; see Claim \ref{claim-ereg:no_interior_extrema} for details. The sign conditions in Theorem \ref{thm-heat:main}\ref{item-heat:boundary_viscosity} are exactly in the favorable direction. Then, $\phi-h$ would have at least 3 roots at the initial time (which is appropriately set up in the proof), contradicting other known facts. See Claim \ref{claim-ereg:dth3_h} for details.

\vspace{3pt}

The same technical difficulties have also arisen in Theorem \ref{thm-heat:no_interior_segment}(ii), and the proof there already implemented the above ideas. The main difference is that we compared $h$ with its second order Taylor expansion in Theorem \ref{thm-heat:no_interior_segment}. Another technical trouble in Theorem \ref{thm-heat:no_interior_segment}\ref{item-heat:grad_est} arises from the fact that $h$ may be discontinuous in $\cJ$ due to large jumps (see around Figure \ref{fig-heat:noncontinuous}). In the case here, $h$ must be continuous due to the presence of cusps, see Claim \ref{claim-ereg:fine_structure_2} and \ref{claim-ereg:preimage_is_segment}.

\begin{proof}[Proof of Theorem \ref{thm-ereg:es_combined} assuming Theorem \ref{thm-ereg:main_infhar}] {\ }
	
	\cite[Theorem 5.2]{Evans-Savin_2008} implies the existence of universal constants $\lambda_0,C_1$ so that: if $\lambda\leq\lambda_0$, and $|u-x|\leq\lambda^2$ in $Q_\lambda(4)$, then
	\begin{equation}\label{eq-ereg:aux2}
		|\p_xu-1|\leq C_1\lambda^2,\qquad|\p_yu|\leq C_1\lambda\qquad\text{in}\ \ Q_\lambda(3).
	\end{equation}
	We remind that the notation ``$Q_\lambda$'' in \cite{Evans-Savin_2008} corresponds to $(-1,1)\times(-\lambda^{1/2},\lambda^{1/2})$. When $\lambda$ is small enough, \eqref{eq-ereg:aux2} implies
	\begin{equation}\label{eq-ereg:aux4}
		\Big|\frac{\D u}{|\D u|}-\p_x\Big|\leq C_2\lambda,\qquad\big|\log|\D u|\big|\leq C_2^2\lambda^2\qquad\text{in}\ \ Q_\lambda(3),
	\end{equation}
	where $C_2$ is another universal constant. Now choose $C>\max\{\lambda_0^{-2},e^{30}C_2^2\}$ and set $\lambda=C^{-1/2}\delta<\lambda_0$. For $u$ satisfying \eqref{eq-ereg:intro_ereg_u} in $Q_\delta(4)$, we have $|u-x|\leq(C^{-1/2}\delta)^2=\lambda^2$. Thus, using \eqref{eq-ereg:aux4} in all vertical translations of $Q_\lambda(4)$ we have
	\[\Big|\frac{\D u}{|\D u|}-\p_x\Big|\leq C_2C^{-1/2}\delta,\qquad\big|\log|\D u|\big|\leq C_2^2C^{-1}\delta^2\qquad\text{in}\ \ Q_\delta(3).\]
	The result follows from the scaled versions of Theorem \ref{thm-ereg:main_infhar}.
\end{proof}

\begin{proof}[Proof of Theorem \ref{thm-ereg:main_infhar} assuming Theorem \ref{thm-ereg:main_cluster}] {\ }
	
	Take any $z_1=(x_1,y_1)$, $z_2=(x_2,y_2)$ in $Q_\delta(1/8)$. Denote
	\[\bw=-\log|\D u|,\qquad\bnu=-\frac{\D^\perp u}{|\D u|}.\]
	Apply Theorem \ref{thm-approx:global} with basepoints $z_1,z_2$ respectively: we obtain two simple IMCF clusters $w_1,w_2$ in $Q_\delta(4)$, with the associated calibrations $\nu_i=-e^{w_i}\D^\perp u$. Notice that $\D u=e^{-w_i}\nu_i^\perp$, and
	\begin{equation}\label{eq-ereg:wi_eq_w}
		\nu_i(z)=\bnu(z)\qquad\text{whenever}\qquad w_i(z)=\bw(z).
	\end{equation}
	
	By Theorem \ref{thm-approx:global}\ref{item-approx:w_leq_-logdu} and \eqref{eq-approx:main_w_on_path_1} and \eqref{eq-ereg:intro_ereg_du}, we have
	\[w_i\leq\bw\leq e^{-30}\delta^2,\qquad w_i\geq-\log\sup_{Q_\delta(4)}|\D u|\geq-e^{-30}\delta^2\]
	in $Q_\delta(4)$, for $i=1,2$. Also,
	\[\begin{aligned}
		|\nu_i+\p_y| &= \big|\!-e^{-w_i}\D^\perp u+\p_y\big|
		= \Big|\!-e^{-w_i}\D^\perp u+\frac{\D^\perp u}{|\D u|}\Big|+\Big|\!-\frac{\D^\perp u}{|\D u|}+\p_y\Big| \\
		&\leq \big|e^{-w_i}|\D u|-1\big|+e^{-15}\delta<e^{-10}\delta.
	\end{aligned}\]
	Hence, $w_1,w_2$ satisfy the condition for Theorem \ref{thm-ereg:main_cluster}.
	
	Also, we fall into the $\epsilon$-regularity realm as appeared below Remark \ref{rmk-approx:misc}. From there it follows that $w_i$ ($i=1,2$) contains only one ridge $\gamma_i\ni z_i$, and $\gamma_1,\gamma_2$ are $e^{-9}\delta$-Lipschitz graphs. By Theorem \ref{thm-approx:global}\ref{item-approx:w_on_gamma}, we have $w_i=\bw$ and $\nu_i=\bnu=-(\gamma'_i)^\perp$ on $\gamma_i$, for each $i=1,2$. Here $\gamma'_i$ denotes the right-pointing unit tangent vector of $\gamma_i$.
	
	Consider the unique point $z'_1\in\gamma_1$ with the same $x$-coordinate as $z_2$. Denote $z'_1=(x_2,y'_1)$. Using Theorem \ref{thm-ereg:main_cluster} on $w_1$, we have
	\begin{equation}\label{eq-ereg:cpnt1}
		\begin{aligned}
			\big|\bw(z_1)-\bw(z'_1)\big| &= \big|w_1(z_1)-w_1(z'_1)\big|\leq \delta^{4/3}|z_1-z'_1|^{2/3}, \\
			\big|\bnu(z_1)-\bnu(z'_1)\big| &= \big|(\gamma'_1)^\perp(z_1)-(\gamma'_1)^\perp(z'_1)\big|\leq \delta|z_1-z'_1|^{1/2}.
		\end{aligned}
	\end{equation}
	It remains to connect $z'_1$ with $z_2$. Recall from Theorem \ref{thm-approx:global}\ref{item-approx:w_on_path}, see also \eqref{eq-approx:ep_reg_w}, that
	\begin{equation}\label{eq-ereg:formula_w1w2}
		\begin{aligned}
			w_1(x_2,y) &= -\log\max\Big\{|\D u(x_2,y')|:y'\in[y,y'_1]\Big\}, \\
			w_2(x_2,y) &= -\log\max\Big\{|\D u(x_2,y')|:y'\in[y,y_2]\Big\}, \\
		\end{aligned}
	\end{equation}
	where we denote $[a,b]=[b,a]$ if $b<a$. Take $y_3\in[y'_1,y_2]$ so that
	\[|\D u(x_2,y_3)|=\max\Big\{|\D u(x_2,y)|:y\in[y'_1,y_2]\Big\}.\]
	Then \eqref{eq-ereg:formula_w1w2} yields
	\[w_1(x_2,y_3)=-\log|\D u(x_2,y_3)|=w_2(x_2,y_3).\]
	Denoting $z_3=(x_2,y_3)$, we have by \eqref{eq-ereg:intro_cluster_w}
	\begin{equation}\label{eq-ereg:cpnt2}
		\begin{aligned}
			\big|\bw(z'_1)-\bw(z_2)\big| &\leq \big|\bw(z'_1)-\bw(z_3)\big|+\big|\bw(z_3)-\bw(z_2)\big| \\
			&=\big|w_1(z'_1)-w_1(z_3)\big|+\big|w_2(z_3)-w_2(z_2)\big|
			\leq 2\delta^{4/3}|z'_1-z_2|^{2/3}.
		\end{aligned}
	\end{equation}
	Then, $\gamma_1,\gamma_2$ have $C^{1,1/2}$ norms bounded by $\delta$ in $Q_\delta(1/4)$, by \eqref{eq-ereg:intro_cluster_gamma}. At each possible point $z\in\gamma_1\cap\gamma_2$ we have
	\[\gamma'_1(z)=\bnu^\perp(z)=\gamma'_2(z),\]
	hence the intersection of $\gamma_1$ and $\gamma_2$ is tangential. Finally, Lemma \ref{lemma-ereg:two_graphs} with $\alpha=1/2$, $\mu=\delta/2$ and $s=1/2$ yields
	\begin{equation}\label{eq-ereg:cpnt3}
		\big|\bnu(z'_1)-\bnu(z_2)\big|=\big|(\gamma'_1)^\perp(z'_1)-(\gamma'_2)^\perp(z_2)\big|\leq 32\delta^{2/3}|z'_1-z_2|^{1/3}.
	\end{equation}
	Combining \eqref{eq-ereg:cpnt1} \eqref{eq-ereg:cpnt2} \eqref{eq-ereg:cpnt3}, we thus have shown
	\[\|\bw\|_{C^{0,2/3}(Q_\delta(1/8))}\leq4\delta^{4/3},\qquad\|\bnu+\p_y\|_{C^{0,1/3}(Q_\delta(1/8))}\leq64\delta^{2/3}.\]
	This implies \eqref{eq-ereg:intro_main_1} by a covering argument. It is elementary that \eqref{eq-ereg:intro_ereg_du} \eqref{eq-ereg:intro_main_1} imply \eqref{eq-ereg:intro_main_2}.
\end{proof}

\begin{proof}[Proof of Theorem \ref{thm-ereg:main_cluster}] {\ }
	
	Let us prove the following: if in addition $0$ lies on a ridge $\gamma$, then
	\begin{equation}\label{eq-ereg:to_prove}
		\|w\|_{C^{0,2/3}(Q_\delta(1))}\leq e^{-1}\delta^{4/3},\qquad\|\gamma\|_{C^{1,1/2}(Q_\delta(1))}\leq e^{-1}\delta.
	\end{equation}
	Indeed, this implies the theorem as follows: if $Q_\delta(3/8)$ intersects no ridges, then \eqref{eq-ereg:intro_cluster_gamma} is vacuous. Next, Theorem \ref{thm-heat:no_interior_segment}\ref{item-heat:grad_est} in $Q_\delta(3/8)$ and $|\nu+\p_y|<e^{-10}\delta$, $\metric{\nu}{\D w}=|\D w|$ imply
	\[|\D w|\leq 16e^{-7}\delta,\qquad\metric{\D w}{\p_x}\leq 16e^{-17}\delta^2\qquad\text{in}\ \ Q_\delta(1/4).\]
	Hence, for two points $z_1=(x_1,y_1)\in Q_\delta(1/4)$ and $z_2=(x_2,y_2)\in Q_\delta(1/4)$, we have
	\[\begin{aligned}
		|w(z_1)-w(z_2)| &\leq 16e^{-17}\delta^2|x_1-x_2|+16e^{-7}\delta|y_1-y_2| \\
		&\leq 16e^{-17}\delta^2|x_1-x_2|^{2/3}+16e^{-7}\delta\cdot\delta^{1/3}|y_1-y_2|^{2/3} \\
		&\leq 32e^{-7}\delta^{4/3}|z_1-z_2|^{2/3},
	\end{aligned}\]
	hence showing \eqref{eq-ereg:intro_cluster_w}. If $Q_\delta(3/8)$ intersects a ridge, then it is the unique ridge appearing in the cluster (Lemma \ref{lemma-cluster:unique_ridge_in_square}). Denote it by $\gamma$. Fixing a basepoint $z_0\in\gamma\cap Q_\delta(3/8)$ and applying \eqref{eq-ereg:to_prove} in the ambient box $Q_\delta(z_0,3)$, we obtain \eqref{eq-ereg:intro_cluster_w} \eqref{eq-ereg:intro_cluster_gamma} as desired.
	
	It remains to prove \eqref{eq-ereg:to_prove} assuming that $0$ lies on a ridge $\gamma$. We encourage the reader to read the proof of Lemma \ref{lemma-cluster:local_model} to gain more familiarity with the setups in this proof.
	
	Let $f_{\ridge}:(-4,4)\to(-4e^{-9}\delta,4e^{-9}\delta)$, with $f_{\ridge}(0)=0$ and $|f'_{\ridge}|\leq e^{-9}\delta$, be such that $\gamma=\graphh(f_{\ridge})$, where $\gamma$ denotes the ridge containing $0$. So $w$ consists of an upward evolving IMCF in $\{y<f_{\ridge}(x)\}$ and a downward evolving IMCF in $\{y>f_{\ridge}(x)\}$, and
	\begin{equation}\label{eq-ereg:monotonicity_in_y}
		y\mapsto w(x,y)\text{\ \ is\ \ }\left\{\begin{aligned}
			& \text{non-decreasing\ \ when\ \ $y\leq f_{\ridge}(x)$}, \\
			& \text{non-increasing\ \ when\ \ $y\geq f_{\ridge}(x)$},
		\end{aligned}\right.\qquad\forall\,x\in(-4,4).
	\end{equation}
	As usual, denote
	\[Y_t=\{w>t\},\qquad\tZ_t=\{w\geq t\},\qquad Z_t=\Int(\tZ_t).\]
	By Lemma \ref{lemma-cluster:local_model}\ref{item-cluster:model}, each $Z_t$ has the form
	\[Z_t=\Big\{x\in I_t,\ g_{1,t}(x)<y<g_{2,t}(x)\Big\},\]
	where $I_t=(-4,4)$ or $I_t=(-4,a_t)\cup(b_t,4)$ with $-4\leq a_t\leq b_t\leq4$, and $g_{1,t},g_{2,t}:\bar I_t\cap(-4,4)\to[-4\delta,4\delta]$ satisfy $g_{1,t}\leq f_{\ridge}\leq g_{2,t}$ and have the convexity properties stated in Lemma \ref{lemma-cluster:local_model}\ref{item-cluster:g1_g2}. We also have
	\[Y_t=\Big\{x\in I_t^+,\ g_{1,t}^+(x)<y<g_{2,t}^+(x)\Big\},\]
	where $I_t^+=(-4,4)$ or $I_t^+=(-4,a_t^+)\cup(b_t^+,4)$ with $-4\leq a_t^+\leq b_t^+\leq 4$, and $g_{1,t}^+,g_{2,t}^+:\bar{I_t^+}\cap(-4,4)\to[-4\delta,4\delta]$ satisfy $g_{1,t}^+<f_{\ridge}<g_{2,t}^+$ in $I_t^+$, and $g_{1,t}^+,g_{2,t}^+$ have the properties in Lemma \ref{lemma-cluster:local_model}\ref{item-cluster:g1_g2}. The fact $Y_t\subset Z_t$ implies $I_t^+\subset I_t$ and $g_{1,t}^+\geq g_{1,t}$, $g_{2,t}^+\leq g_{2,t}$ in $I_t^+$.
	
	\vspace{3pt}
	
	Consider
	\[\uT=\operatorname{min}_{\gamma\cap\bar{Q_\delta(3)}}(w).\]
	Then $I_t,I_t^+\supset[-3,3]$ for all $t<\uT$, but $I_\uT^+\not\supset[-3,3]$.
	
	Therefore, $I_\uT^+=(-4,a_\uT^+)\cup(b_\uT^+,4)$ with $-4\leq a_\uT^+\leq b_\uT^+\leq4$ and $[a_\uT^+,b_\uT^+]\cap[-3,3]\ne\emptyset$, and $w(x,f_{\ridge}(x))\equiv\uT$ in $[a_\uT^+,b_\uT^+]\cap[-3,3]$. Note that $g_{1,\uT}^+,g_{2,\uT}^+$ are concave\,/\,convex in each interval in $I_\uT^+$, by Lemma \ref{lemma-cluster:local_model}\ref{item-cluster:g1_g2}.
	
	Consider the merged functions
	\[\bg_{1,\uT}^+(x)=\left\{\begin{aligned}
		& f_{\ridge}(x)\qquad x\in[a_\uT^+,b_\uT^+] \\
		& g_{1,\uT}^+(x)\qquad\text{otherwise}
	\end{aligned}\right.,\qquad
	\bg_{2,\uT}^+(x)=\left\{\begin{aligned}
		& f_{\ridge}(x)\qquad x\in[a_\uT^+,b_\uT^+] \\
		& g_{2,\uT}^+(x)\qquad\text{otherwise}
	\end{aligned}\right.,\]
	Notice that $\bg_{1,\uT}^+\leq f_{\ridge}$ and $\bg_{2,\uT}^+\geq f_{\ridge}$, and
	\begin{equation}\label{eq-ereg:w_const_uT}
		w\big(x,\bg_{i,\uT}^+(x)\big)=\uT,\qquad\forall\,x\in[-3,3],\ i=1,2.
	\end{equation}
	Then, consider the regions
	\[E_1=\Big\{x\in[-3,3],\ -3\delta\leq y\leq\bg_{1,\uT}^+(x)\Big\},\qquad
	E_2=\Big\{x\in[-3,3],\ \bg_{2,\uT}^+(x)\leq y\leq3\delta\Big\}.\]
	It is clear that $E_1\cup E_2=\bar{Q_\delta(3)}\cap\{w\leq\uT\}$.
	
	\begin{claim}
		The functions
		\[w_1=\left\{\begin{aligned}
			& w\qquad\text{in}\ \ E_1 \\
			& \uT\qquad\text{in}\ \ \bar{Q_\delta(3)}\setminus E_1
		\end{aligned}\right.,\qquad
		w_2=\left\{\begin{aligned}
			& w\qquad\text{in}\ \ E_2 \\
			& \uT\qquad\text{in}\ \ \bar{Q_\delta(3)}\setminus E_2
		\end{aligned}\right.\]
		are both weak IMCFs in $Q_\delta(3)$.
	\end{claim}
	\begin{proof}
		Notice by \eqref{eq-ereg:w_const_uT} that $w_1,w_2$ are Lipschitz in $\bar{Q_\delta(3)}$. Let us show that $w_2$ is a weak IMCF in $Q_\delta(3)$ (same argument for $w_1$). If $w\equiv\uT$ in $(-3,3)\times[\delta,2\delta]$, then $w\equiv\uT$ in the region $\big\{x\in(-3,3),\,\bg_{2,\uT}^+(x)<y<2\delta\big\}$ by \eqref{eq-ereg:monotonicity_in_y}, hence $w_2\equiv\uT$ in $(-3,3)\times(-3\delta,2\delta)$. Then Remark \ref{rmk-prelim:imcf_properties}\ref{item-prelim:imcf_extension} with $r=3$ and $l=3\delta$ and $f\equiv1.5\delta$ implies the result. Otherwise, $w(z_1)<\uT$ for some $z_1\in(-3,3)\times[\delta,2\delta]$. Then for all $t\in(w(z_1),\uT)$, the set $\{w_2<t\}$ is a convex epigraph in $Q_\delta(3)\cap\{y>f_{\ridge}(x)\}$. Taking $t\nearrow\uT$, it follows that $\{w_2<\uT\}$ is the epigraph of a convex function $f\geq f_{\ridge}$. Then Remark \ref{rmk-prelim:imcf_properties}\ref{item-prelim:imcf_extension} with $r=3$, $l=3\delta$ implies the result.
	\end{proof}
	
	Using Theorem \ref{thm-heat:no_interior_segment}\ref{item-heat:grad_est} on $w_1,w_2$ respectively, we obtain
	\begin{equation}\label{eq-ereg:reg_w_E1E2}
		|\D w|\leq2e^{-7}\delta\qquad\text{in}\ \ (E_1\cup E_2)\cap Q_\delta(2).
	\end{equation}
	
	\begin{claim}\label{claim-ereg:reg_g_uT}
		The curvature of $\graphh(\bg_{1,\uT}^+)$ and $\graphh(\bg_{2,\uT}^+)$ in $Q_\delta(2)$ are bounded by $2e^{-7}\delta$. Furthermore,
		\begin{equation}\label{eq-ereg:reg_ridge_uT}
			\|f_{\ridge}\|_{C^{1,1}([a_\uT^+,b_\uT^+]\cap(-2,2))}\leq 2e^{-7}\delta.
		\end{equation}
	\end{claim}
	\begin{proof}
		Since $I_t\supset[-3,3]$ for all $t<\uT$, we have
		\[\tZ_\uT\cap\bar{Q_\delta(3)}=\Big\{x\in[-3,3],\ \tg_{1,\uT}(x)\leq y\leq\tg_{2,\uT}(x)\Big\},\]
		where $\tg_{1,\uT}=\lim_{s\nearrow\uT}g_{1,s}$ and $\tg_{2,\uT}=\lim_{s\nearrow\uT}g_{2,s}$. By \eqref{eq-ereg:reg_w_E1E2}, note that $\graphh(\tg_{1,\uT})\cap Q_\delta(2)$ and $\graphh(\tg_{2,\uT})\cap Q_\delta(2)$ have curvature bounded by $2e^{-7}\delta$.
		
		We claim that $\bg_{1,\uT}^+$ is linear in each open interval in
		\[[-3,3]\cap\big\{\tg_{1,\uT}<\bg_{1,\uT}^+<\tg_{2,\uT}\big\}.\]
		Let $(c,d)$ be such an interval. Then, for every subinterval $[c',d']\subset(c,d)$, there is $\epsilon>0$ so that $w\geq\uT$ in the region
		\[D=\Big\{x\in(c',d'),\ \bg_{1,\uT}^+(x)-\epsilon<y<\bg_{1,\uT}^+(x)+\epsilon\Big\}.\]
		Notice that $\graphh(\bg_{1,\uT}^+)|_{(c',d')}$ is orthogonal to $\nu$, and $w\equiv\uT$, $|\nu|\equiv1$ there. Hence, the vector field $\nu'=e^{\uT-w}\nu$ is a calibration of $\graphh(\bg_{1,\uT}^+)$ in $D$, as desired.
		
		As a result, in $[-3,3]$, the function $\bg_{1,\uT}^+$ either is linear or coincides with $\tg_{1,\uT}$ or $\tg_{2,\uT}$. The curvature bound of $\tg_{1,\uT},\tg_{2,\uT}$ then implies that of $\bg_{1,\uT}^+$, $\bg_{2,\uT}^+$. Finally, noticing that $f_{\ridge}=\bg_{1,\uT}^+|_{[a_\uT^+,b_\uT^+]}$, the curvature bound of $\bg_{1,\uT}^+$ implies \eqref{eq-ereg:reg_ridge_uT}.
	\end{proof}
	
	\begin{figure}[ht]
		\centering
		\includegraphics{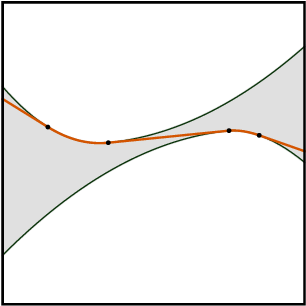}
		\begin{picture}(0,0)
			\put(-3,61){$\tg_{1,\uT}$}
			\put(-3,75){$f_{\ridge}$}
			\put(-3,120){$\tg_{2,\uT}$}
		\end{picture}
		\caption{An extreme case where $Y_\uT=\emptyset$ and $w\equiv\uT$ in the shaded region.}\label{fig-ereg:flat_ridge}
	\end{figure}
	
	Consider the remaining regions
	\[\begin{aligned}
		E_3 &= \Big\{x\in[-3,a_\uT^+],\,g_{1,\uT}^+(x)\leq y\leq g_{2,\uT}^+(x)\Big\}, \\
		E_4 &= \Big\{x\in[b_\uT^+,3],\,g_{1,\uT}^+(x)\leq y\leq g_{2,\uT}^+(x)\Big\}.
	\end{aligned}\]
	Notice $E_3=\emptyset$ if $a_\uT^+<-3$, and $E_4=\emptyset$ if $b_\uT^+>3$. The main regularity claim is as follows:
	
	\begin{figure}[ht]
		\centering
		\includegraphics{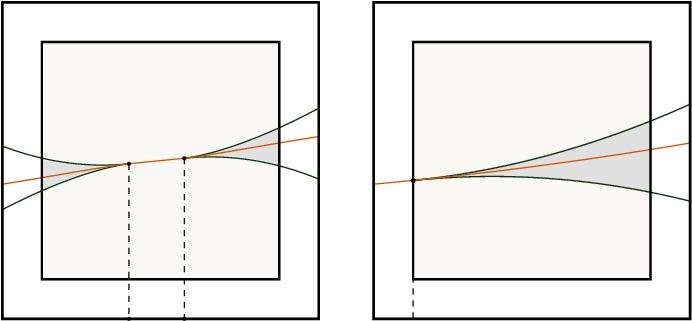}
		\begin{picture}(0,0)
			% Left
			\put(-313,26){$E_1$}
			\put(-313,121){$E_2$}
			\put(-308,80){\arrowangle{-100}}
			\put(-307,90){$E_3$}
			\put(-226,90){\arrowangle{-64}}
			\put(-229,101){$E_4$}
			\put(-290,8){$a_\uT^+$}
			\put(-262,8){$b_\uT^+$}
			% Right
			\put(-134,26){$E_1$}
			\put(-134,121){$E_2$}
			\put(-55,91){\arrowangle{-68}}
			\put(-58,102){$E_4$}
			\put(-151,8){$b_\uT^+$}
		\end{picture}
		\caption{In the second picture, $a_\uT^+=-4$ and $E_3=\emptyset$.}\label{fig-ereg:Ek}
	\end{figure}
	
	\begin{claim}\label{claim-ereg:reg_in_E3E4}
		We have
		\begin{equation}\label{eq-ereg:reg_w_in_E3E4}
			\|w\|_{C^{0,2/3}(E_4\cap Q_\delta(1))}\leq e^{-2}\delta^{4/3}
		\end{equation}
		and
		\begin{equation}\label{eq-ereg:reg_f_in_E3E4}
			\|f_{\ridge}\|_{C^{1,1/2}([b_\uT^+,3]\cap[-1,1])}\leq e^{-2}\delta.
		\end{equation}
		A similar fact holds for $E_3$.
	\end{claim}
	
	Having known this fact, noting that
	\[E_1\cup E_2\cup E_3\cup E_4\supset\bar{Q_\delta(3)},\]
	the main results \eqref{eq-ereg:to_prove} follow from Claim \ref{claim-ereg:reg_g_uT}, \ref{claim-ereg:reg_in_E3E4} and \eqref{eq-ereg:reg_ridge_uT}, by triangle inequalities.
	
	\vspace{3pt}
	
	It remains to prove Claim \ref{claim-ereg:reg_in_E3E4}. Since $E_3,E_4$ are symmetric via a reflection around the $y$-axis, it suffices to treat $E_4$. We may assume $b_\uT^+<1$, otherwise the claim is vacuous. Consider
	\[\oT=\max_{E_4}(w)=\max_{E_4\cap\gamma}(w),\]
	where the second equality follows from \eqref{eq-ereg:monotonicity_in_y}. Notice that $b_t,b_t^+<3$ for all $t\in[\uT,\oT)$. Note: it may be possible that $Z_\oT=\emptyset$ or $b_\oT>3$.
	
	Recall from above that $I_\uT^+\not\supset[-3,3]$. This implies:
	\begin{itemize}[nosep]
		\item $b_\uT^+\geq-3$;
		\item $I_t=(-4,a_t)\cup(b_t,4)$ and $I_t^+=(-4,a_t^+)\cup(b_t^+,4)$ with $b_t,b_t^+>-3$, for all $t>\uT$.
	\end{itemize}
	By Lemma \ref{lemma-cluster:local_model}\ref{item-cluster:g1_g2} we have:
	\begin{itemize}[nosep]
		\item $g_{1,t}$ (resp. $g_{2,t}$) is concave (resp. convex) in $[b_t,4)$, for all $t>\uT$;
		\item $g_{1,t}^+$ (resp. $g_{2,t}^+$) is concave (resp. convex) in $[b_t^+,4)$, for all $t\geq\uT$;
		\item $g_{1,t}\leq f_{\ridge}\leq g_{2,t}$ in $[b_t,4]$ for all $t>\uT$;
		\item $g_{1,t}^+<f_{\ridge}<g_{2,t}^+$ in $(b_t^+,4]$ for all $t\geq\uT$;
		\item $w(x,f_{\ridge}(x))$ is non-decreasing for $x\in[b_\uT^+,3]$;
		\item the graphs of $g_{1,t},g_{2,t},g_{1,t}^+,g_{2,t}^+,f_{\ridge}$ all have $\nu$ as unit normal vector;
		\item the following convergences hold:
		\[\lim_{s\nearrow t}b_s=\lim_{s\nearrow t}b_s^+\leq b_t\leq b_t^+=\lim_{s\searrow t}b_s=\lim_{s\searrow t}b_s^+,\]
		and
		\[\begin{aligned}
			& \lim_{s\nearrow t}g_{1,s}=\lim_{s\nearrow t}g_{1,s}^+=g_{1,t}\qquad\text{in}\ \ C^1([b_t,3]),\qquad\forall\,t\in(\uT,\oT], \\
			& \lim_{s\searrow t}g_{1,s}=\lim_{s\searrow t}g_{1,s}^+=g_{1,t}^+\qquad\text{in}\ \ C^1([b_t^+,3]),\qquad\forall\,t\in[\uT,\oT).
		\end{aligned}\]
		Similar facts hold for $g_{2,t}, g_{2,t}^+$.
	\end{itemize}
	Let us denote
	\[\begin{aligned}
		\gamma_t &= \graphh\big(g_{1,t}|_{[b_t,3]}\big)\cup\graphh\big(g_{2,t}|_{[b_t,3]}\big)=\p Z_t\cap E_4,\qquad\forall\,t\in(\uT,\oT], \\
		\gamma_t^+ &= \graphh\big(g_{1,t}^+|_{[b_t^+,3]}\big)\cup\graphh\big(g_{2,t}^+|_{[b_t^+,3]}\big)=\p Y_t\cap E_4,\qquad\forall\,t\in[\uT,\oT),
	\end{aligned}\]
	(where $\gamma_\oT$ may be empty,) and
	\begin{equation}\label{eq-ereg:def_Gamma0}
		\Gamma_0=E_4\setminus\bigcup_{t>\uT}\Int\big(\{w=t\}\big).
	\end{equation}
	Clearly, $\Gamma_0$ is compact and $\gamma_t,\gamma_t^+\subset\Gamma_0$. It is an elementary topological fact
	\begin{equation}\label{eq-ereg:def_Gamma0_2}
		\Gamma_0=E_4\cap\bigcup_{t\geq\uT}\big(\p\tZ_t\cap\p Y_t\big),
	\end{equation}
	hence $|\nu|=1$ on $\Gamma_0$. As $|w|<e^{-20}\delta^2$ and $|\nu+\p_y|<e^{-10}\delta$, the map $(w,\nu^\perp):\Gamma_0\to\RR\times\SS^1$ lifts to a continuous map
	\[\Th_0:\Gamma_0\to(-e^{-20}\delta^2,e^{-20}\delta^2)\times(-e^{-10}\delta,e^{-10}\delta).\]
	Denote $\bar\cJ=\Th_0(\Gamma_0)$. We aim to show that each preimage of $\Th_0$ is a line segment.
	
	For each $t\in(\uT,\oT]$, consider the quantity
	\begin{equation}\label{eq-ereg:def_btminus}
		b_t^-=\lim_{s\nearrow t}b_s.
	\end{equation}
	Hence $b_\oT^-\leq3$. We mention again that $b_t<3$ for all $t\in(\uT,\oT)$, but possibly $b_\oT>3$.
	
	\begin{claim}\label{claim-ereg:fridge_linear}
		$f_{\ridge}$ is linear and $w(x,f_{\ridge}(x))\equiv t$ in $[b_t^-,b_t]\cap[-3,3]$, for each $t\in(\uT,\oT]$.
	\end{claim}
	\begin{proof}
		Call $I=[b_t^-,b_t]\cap[-3,3]$. For all $s\in(\uT,t)$, note that $g_{1,s}$ is concave and $g_{2,s}$ is convex in $I$. By the definition of $b_t$, the ascending limit $\lim_{s\nearrow t}g_{1,s}$ and descending limit $\lim_{s\nearrow t}g_{2,s}$ must agree (hence must agree with $f_{\ridge}$) on $I$. The linearity of $f_{\ridge}$ follows from convexity. The constancy of $w$ follows from $w=s$ on $\graphh(g_{1,s})\cup\graphh(g_{2,s})$.
	\end{proof}
	
	We in fact have
	\begin{equation}\label{eq-ereg:tZt_description}
		\tZ_t\cap E_4=\Big\{x\in[b_t^-,3]: \tg_{1,t}(x)\leq y\leq\tg_{2,t}(x)\Big\},\qquad\forall\,t\in(\uT,\oT],
	\end{equation}
	where $\tg_{1,t}:=\lim_{s\nearrow t}g_{1,s}$ and $\tg_{2,t}:=\lim_{s\nearrow t}g_{2,s}$. Notice that
	\begin{equation}\label{eq-ereg:tgt_properties}
		\tg_{i,t}=g_{i,t}\text{ on }[b_t,3],\qquad \tg_{1,t}=\tg_{2,t}=f_{\ridge}\text{ is linear on }[b_t^-,b_t].
	\end{equation}
	This follows from $\tZ_t=\bigcap_{s<t}Y_s$ and the proof of Claim \ref{claim-ereg:fridge_linear}, see also Lemma \ref{lemma-cluster:local_model}\ref{item-cluster:tZt_model}.
	
	\begin{claim}\label{claim-ereg:fine_structure_2}
		For all $t\in(\uT,\oT)$, there exists $i\in\{1,2\}$ so that
		\[g_{i,t}(b_t^+)=g_{1,t}^+(b_t^+)=g_{2,t}^+(b_t^+),\qquad\quad g_{i,t}\equiv f_{\ridge},\ \ w(\cdot,f_{\ridge}(\cdot))\equiv t\ \ \text{in}\ \ [b_t,b_t^+].\]
	\end{claim}
	\begin{proof}
		This claim is nontrivial only when $b_t<b_t^+$. Suppose otherwise $g_{1,t}(b_t^+)<g_{1,t}^+(b_t^+)=g_{2,t}^+(b_t^+)<g_{2,t}(b_t^+)$. For $\epsilon>0$, consider
		\[z_1=\big(b_t^+,f_{\ridge}(b_t^+)\big),\quad z_2=\big(b_t^++\epsilon,g_{1,t}^+(b_t^++\epsilon)\big),\quad z_3=\big(b_t^++\epsilon,g_{2,t}^+(b_t^++\epsilon)\big).\]
		When $\epsilon$ is small enough, the triangle with vertices $z_1,z_2,z_3$ has a neighborhood contained in $Z_t$. Recall $w\geq t$ in $Z_t$. Hence, the vector field $e^{t-w}\nu$ calibrates $\graphh\big(g_{1,t}^+|_{[b_t^+,b_t^++\epsilon]}\big)$ and $\graphh\big(g_{2,t}^+|_{[b_t^+,b_t^++\epsilon]}\big)$, showing that they are both line segments, which is impossible.
		
		Hence, either $g_{1,t}(b_t^+)=g_{1,t}^+(b_t^+)=g_{2,t}^+(b_t^+)$ or $g_{2,t}(b_t^+)=g_{1,t}^+(b_t^+)=g_{2,t}^+(b_t^+)$. We may assume that the first possibility occurs. This implies $g_{1,t}(b_t^+)=f_{\ridge}(b_t^+)$ as well. Now we have $w(b_t,f_{\ridge}(b_t))=w(b_t^+,f_{\ridge}(b_t^+))=t$ since $w=t$ on $\p Z_t\cup\p Y_t$. By the monotonicity of $w|_\gamma$, it follows that
		\[w(\cdot,f_{\ridge}(\cdot))\equiv t\qquad\text{in}\ \ [b_t,b_t^+].\]
		The monotonicity in $y$-direction then implies $w\equiv t$ in the set
		\[\Big\{x\in[b_t,b_t^+],\ g_{1,t}(x)\leq y\leq f_{\ridge}(x)\Big\}.\]
		Using $\div(\nu)=e^t\div(e^{-w}\nu)=0$ and the divergence theorem in this region, we obtain
		\[\big|\graphh\big(f_{\ridge}|_{[b_t,b_t^+]}\big)\big|=\big|\graphh\big(g_{1,t}|_{[b_t,b_t^+]}\big)\big|.\]
		As $g_{1,t}$ is concave, this forces $f_{\ridge}\equiv g_{1,t}$ in $[b_t,b_t^+]$, as desired.
	\end{proof}
	
	\begin{claim}\label{claim-ereg:fine_structure_1}
		For all $t\in(\uT,\oT)$, we have $\{g_{1,t}^+>g_{1,t}\}\cap[b_\uT^+,3]=(c,3]$ for some $c\in[b_t^+,3]$, and $g_{1,t}^+$ is linear in this interval. Analogous facts hold for $g_{2,t}^+$.
	\end{claim}
	\begin{proof}
		Note that $g_{1,t}^+$ is linear in any interval in $\big\{g_{1,t}<g_{1,t}^+<f_{\ridge}\big\}$: this follows from the geometry of weak IMCF. Recalling $g_{1,t}^+<f_{\ridge}$ in $(b_t^+,3]$, we have that $g_{1,t}^+$ is linear in each interval in $\{g_{1,t}^+>g_{1,t}\}$. Suppose $\{g_{1,t}^+>g_{1,t}\}$ has a connected component of the form $(c,d)$ with $b_t^+\leq c<d\leq3$. Then
		\[g_{1,t}<g_{1,t}^+<f_{\ridge}\text{\ \ in\ \ }(c,d),\qquad\qquad g_{1,t}(d)=g_{1,t}^+(d).\]
		By concavity, this forces $g_{1,t}(c)<g_{1,t}^+(c)$. By our choice of $(c,d)$ as a connected component, this implies $c=b_t^+$. Then Claim \ref{claim-ereg:fine_structure_2} implies $g_{1,t}^+(b_t^+)=g_{2,t}(b_t^+)$. But this yields a contradiction since
		\[g'_{1,t}(d)=(g_{1,t}^+)'(d)=(g_{1,t}^+)'(b_t^+)=g'_{2,t}(b_t^+)>g'_{1,t}(b_t^+)\geq g'_{1,t}(d). \qedhere\]
	\end{proof}
	
	Combining Claim \ref{claim-ereg:fine_structure_1} and \ref{claim-ereg:fine_structure_2}, we may fully describe the shape of the set $\p Y_t\cap\tZ_t$ for $t\in(\uT,\oT)$. See Figures \ref{fig-ereg:fine_struc}, \ref{fig-ereg:fine_struc2} and Claim \ref{claim-ereg:preimage_is_segment} below. We need one more auxiliary lemma:
	
	\begin{claim}\label{claim-ereg:f_oT}
		If $b_\oT<3$, then there exists $d\in[b_\oT,3]$ and $i\in\{1,2\}$ so that
		\[f_{\ridge}\equiv g_{i,\oT}\ \ \text{in}\ \ [b_\oT,d],\qquad\qquad f_{\ridge}\text{ is linear in }[d,3].\]
	\end{claim}
	\begin{proof}
		Note that $w\equiv\oT$ in the region
		\[\Big\{x\in[b_\oT,3],\ g_{1,\oT}(x)\leq y\leq g_{2,\oT}(x)\Big\}.\]
		The calibration argument implies that $f_{\ridge}$ is linear in any interval where $g_{1,\oT}<f_{\ridge}<g_{2,\oT}$. Arguing similarly as in Claim \ref{claim-ereg:fine_structure_1}, the set $\{g_{1,\oT}<f_{\ridge}<g_{2,\oT}\}$ must take the form $(d,3]$ with $d\in[b_\oT,3]$. The result follows as desired.
	\end{proof}
	
	\begin{figure}[ht]
		\centering
		\includegraphics{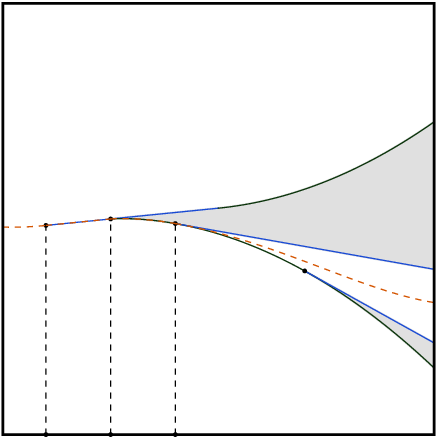}
		\begin{picture}(0,0)
			\put(-206,6){$b_t^-$}
			\put(-173,6){$b_t$}
			\put(-143,6){$b_t^+$}
			\put(-50,105){$w\equiv t$}
			\put(-3,30){$g_{1,t}$}
			\put(-3,43){$g_{1,t}^+$}
			\put(-3,62){$f_{\ridge}$}
			\put(-3,79){$g_{2,t}^+$}
			\put(-3,150){$g_{2,t}$}
		\end{picture}
		\caption{The black and blue solid curves together represent $\p Y_t\cap\p\tZ_t$.}\label{fig-ereg:fine_struc}
	\end{figure}
	
	\begin{figure}[ht]
		\centering
		\includegraphics{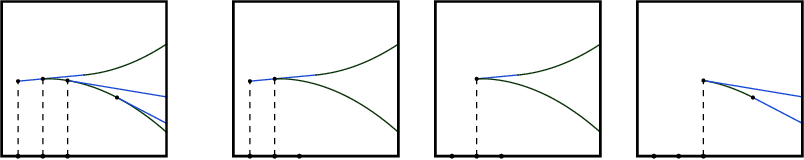}
		\caption{Shapes of $\p Y_t\cap\p\tZ_t$, $\p\tZ_t$, $\p Z_t=\gamma_t$ and $\p Y_t$, respectively.}\label{fig-ereg:fine_struc2}
	\end{figure}
	
	We now proceed to prove that
	
	\begin{claim}\label{claim-ereg:preimage_is_segment}
		$\Th_0^{-1}(t,\th)$ is a line segment for each $(t,\th)\in\bar\cJ$.
	\end{claim}
	\begin{proof}
		The preimage of $\{t=\uT\}$ is $\gamma_\uT^+$, thus the claim holds for $t=\uT$. The preimage of $\{t=\oT\}$ is
		\[\p\tZ_\oT\cap E_4=\Big\{x\in[b_\oT^-,3],\ y=\tg_{1,\oT}(x)\text{ or }y=\tg_{2,\oT}(x)\Big\},\]
		namely, the union of a concave graph and a convex graph with tangential contact at endpoints. This also verifies the claim for $t=\oT$. Assume $t\in(\uT,\oT)$. The preimage of $t$ is
		\[\begin{aligned}
			(\p\tZ_t\cup\p Y_t)\cap E_4 &= \Big\{x\in[b_t^-,3],\ y=\tg_{1,t}(x)\text{ or }y=\tg_{2,t}(x)\Big\} \\
			&\hspace{48pt} \cup\Big\{x\in[b_t^+,3],\ y=g_{1,t}^+(x)\text{ or }y=g_{2,t}^+(x)\Big\},
		\end{aligned}\]
		which by Claim \ref{claim-ereg:fine_structure_2}, \ref{claim-ereg:fine_structure_1} and \eqref{eq-ereg:tgt_properties} is the union of $\gamma_t$ and at most three line segments contacting tangentially (more precisely, with possibly one segment contacting at the cusp and two segments contacting on the edges). The result follows.
	\end{proof}
	
	As a consequence, the support function $h$ defined by
	\[h(t,\th)=\metric{z}{\nu}=\metric{z}{\nu_\th}\quad\text{for any $z\in\Th_0^{-1}(t,\th)$}\]
	is well-defined on $\bar\cJ$, and is continuous since $\Th_0$ is a quotient map. Denote
	\[\begin{aligned}
		& \uth_\uT=\arctan (g_{1,\uT}^+)'(3),\qquad\oth_\uT=\arctan (g_{2,\uT}^+)'(3), \\
		& \uth_\oT=\arctan\tg_{1,\oT}'(3),\qquad\oth_\oT=\arctan\tg_{2,\oT}'(3),
	\end{aligned}\]
	and for each $t\in(\uT,\oT)$, denote
	\[\uth_t=\arctan g'_{1,t}(3),\qquad\oth_t=\arctan g'_{2,t}(3).\]
	Notice that $\oth_t>\uth_t$ for all $t\in[\uT,\oT)$, and from the proof of Claim \ref{claim-ereg:preimage_is_segment}, that
	\[\bar\cJ=\bigcup_{t\in[\uT,\oT]}\{t\}\times[\uth_t,\oth_t].\]
	Also, note the convergence
	\begin{equation}\label{eq-ereg:convergence_th_t}
		\begin{aligned}
			& \lim_{s\nearrow t}\uth_s=\uth_t\leq\arctan(g_{1,t}^+)'(3)=\lim_{s\searrow t}\uth_s, \\
			& \lim_{s\nearrow t}\oth_s=\oth_t\geq\arctan(g_{2,t}^+)'(3)=\lim_{s\searrow t}\oth_s.
		\end{aligned}
	\end{equation}
	
	To invoke Theorem \ref{thm-heat:main}, we need to consider the restricted data
	\[\begin{aligned}
		\Gamma=\bigcup_{t\in(\uT,\oT]}\gamma_t\cap\Int(E_4),\qquad \Th=\Th_0|_\Gamma,\qquad\cJ=\Th(\Gamma)\subset\bar\cJ.
	\end{aligned}\]
	Let $J_t,\bar J_t$ denote the $t$-slice of $\cJ,\bar\cJ$. Thus $\bar J_t=[\uth_t,\oth_t]$ is a nontrivial interval for all $t\in[\uT,\oT)$. Note: $\bar J_\oT$ is not a point if and only if $Z_\oT\cap E_4\ne\emptyset$ (in general, $\tZ_\oT\cap E_4$ could be a line segment and thus $Z_\oT\cap E_4$ could be empty).
	
	\begin{claim}\label{claim-ereg:h_small}
		$|h|\leq e^{-6}\delta$ in $\bar\cJ$.
	\end{claim}
	\begin{proof}
		Due to $g_{1,\uT}^+(b_\uT^+)=g_{2,\uT}^+(b_\uT^+)=f_{\ridge}(b_\uT^+)$ and $f_{\ridge}(0)=0$ and the $e^{-9}\delta$-Lipschitzness of all these functions, we have
		\[E_4\subset[-3,3]\times[-10e^{-9}\delta,10e^{-9}\delta].\]
		Hence, for all $z=(x,y)\in E_4$, we use $|\nu+\p_y|<e^{-10}\delta$ to bound
		\[\big|\metric{z}{\nu}\big|\leq|x|\cdot e^{-10}\delta+|y|\leq 3e^{-10}\delta+10e^{-9}\delta\leq e^{-6}\delta. \qedhere\]
	\end{proof}
	
	\begin{claim}\label{claim-ereg:backward_openness}
		If $t\in(\uT,\oT]$ and $\th\in\Int(J_t)$, then there is $\epsilon>0$ so that
		\[(t-\epsilon,t]\times(\th-\epsilon,\th+\epsilon)\subset\cJ\]
		and
		\[\Th^{-1}(s,\omega)\Subset\Int(E_4)\qquad\forall\,s\in(t-\epsilon,t],\ \omega\in(\th-\epsilon,\th+\epsilon).\]
	\end{claim}
	\begin{proof}
		The condition implies $\uth_t<\th<\oth_t$. Since $\lim_{s\nearrow t}\uth_s=\uth_t$ and $\lim_{s\nearrow t}\oth_s=\oth_t$, there exists $\epsilon<t-\uT$ so that
		\[g'_{1,s}(3)<\tan(\th-\epsilon),\qquad g'_{2,s}(3)>\tan(\th+\epsilon),\qquad\forall\,s\in[t-\epsilon,t].\]
		Hence $[\th-\epsilon,\th+\epsilon]\subset J_s$ and $\Th^{-1}\big(\{s\}\times[\th-\epsilon,\th+\epsilon]\big)\Subset E_4$ for all $s\in[t-\epsilon,t]$.
	\end{proof}
	
	Combining Theorem \ref{thm-heat:main}\ref{item-heat:interior_viscosity} with $\Omega=\Int(E_4)$ and the continuity of $h$ in $\bar\cJ$ (hence in $\cJ$) and the classical viscosity theory, it follows that $h$ is smooth and solves $\p_t h=\p_{\th\th}h+h$ in every parabolic box in Claim \ref{claim-ereg:backward_openness}. Hence $h(t,\cdot)$ is smooth in $\Int(J_t)$ for each $t\in(\uT,\oT]$. Denote
	\[\varpi_t=\arctan g'_{1,t}(b_t)=\arctan g'_{2,t}(b_t)=\arctan f'_{\ridge}(b_t)\]
	the angle at the cusp of $\gamma_t$. For each $t\in(\uT,\oT]$, we may decompose $\gamma_t$ as
	\begin{align}
		\gamma_t &= \Th^{-1}\big(\{\uth_t\}\big)\cup\Th^{-1}\big((\uth_t,\varpi_t]\big)\cup\Th^{-1}\big([\varpi_t,\oth_t)\big)\cup\Th^{-1}\big(\{\oth_t\}\big) \nonumber\\
		&=: \alpha_1\cup\alpha_2\cup\alpha_3\cup\alpha_4, \label{eq-ereg:curve_decomp}
	\end{align}
	where $\alpha_1,\alpha_4$ are (possibly empty or trivial) line segments, and the interior of $\alpha_2,\alpha_3$ are (possibly empty) smooth strictly concave resp. convex curves. Denote
	\[\kappa^{-1}:=\p_{\th\th}h+h,\qquad\forall\,t\in(\uT,\oT],\ \th\in\Int(J_t).\]
	Then $\kappa^{-1}(t,\varpi_t)=0$ whenever $\varpi_t\in\Int(J_t)$, and
	\begin{equation}\label{eq-ereg:sign_of_kappa}
		\kappa^{-1}(t,\cdot)<0\ \ \text{in}\ \ (\uth_t,\varpi_t),\qquad\kappa^{-1}(t,\cdot)>0\ \ \text{in}\ \ (\varpi_t,\oth_t).
	\end{equation}
	By the Sturmian theory, the root $\varpi_t$ is nondegenerate whenever $\varpi_t\in\Int(J_t)$.
	
	\begin{claim}\label{claim-ereg:dth3_h}
		For each $t_0\in(\uT,\oT]$ and $\th_0\in\Int(J_{t_0})$ satisfying
		\begin{equation}\label{eq-ereg:claim_dth3_cond}
			|\p_\th h(t_0,\th_0)|\leq1.6,\qquad |\p_{\th\th}h(t_0,\th_0)|\leq e^5\delta^{-1},
		\end{equation}
		we have
		\[\p_\th^3 h(t_0,\th_0)\geq e^{15}\delta^{-2}.\]
	\end{claim}
	\begin{proof}
		Suppose this is not true at some $(t_0,\th_0)$. Consider
		\[\begin{aligned}
			\phi(t,\th) &= e^{t-t_0}\Big[h(t_0,\th_0)+\p_\th h(t_0,\th_0)(\th-\th_0)+\p_{\th\th} h(t_0,\th_0)\Big((t-t_0)+\frac12(\th-\th_0)^2\Big) \\
			&\hspace{160pt}+\frac{e^{15}}{\delta^2}\Big((t-t_0)(\th-\th_0)+\frac16(\th-\th_0)^3\Big)\Big].
		\end{aligned}\]
		Hence $\p_t\phi=\p_{\th\th}\phi+\phi$. Recall $\bar\cJ\subset(-e^{-20}\delta^2,e^{-20}\delta^2)\times(-e^{-10}\delta,e^{-10}\delta)$. We may estimate
		\begin{align}
			|\p_\th\phi| &= e^{t-t_0}\Big|\p_\th h(t_0,\th_0)+\p_\th^2 h(t_0,\th_0)(\th-\th_0)+\frac{e^{15}}{\delta^2}\Big((t-t_0)+\frac12(\th-\th_0)^2\Big)\Big| \nonumber\\
			&\leq e^{2e^{-20}\delta^2}\Big[1.6+\frac{e^5}\delta\cdot 2e^{-10}\delta+\frac{e^{15}}{\delta^2}\big(2e^{-20}\delta^2+2e^{-20}\delta^2\big)\Big] \nonumber\\
			&\leq 1.7 \label{eq-ereg:p_th_phi}
		\end{align}
		and
		\begin{equation}\label{eq-ereg:p_thth_phi}
			|\p_{\th\th}\phi| = e^{t-t_0}\Big|\p_\th^2 h(t_0,\th_0)+\frac{e^{15}}{\delta^2}(\th-\th_0)\Big|\leq 3.1e^5\delta^{-1},
		\end{equation}
		both inequalities hold in $(-e^{-20}\delta^2,e^{-20}\delta^2)\times(-e^{-10}\delta,e^{-10}\delta)$. Also, by our assumption
		\[\p_\th^3\phi(t_0,\th_0)=e^{15}\delta^{-2}>\p_\th^3h(t_0,\th_0).\]
		
		Denote $\eta=\phi-h$. Thus
		\[\eta(t_0,\th_0)=\p_\th \eta(t_0,\th_0)=\p_\th^2 \eta(t_0,\th_0)=0,\qquad \p_\th^3 \eta(t_0,\th_0)>0.\]
		By the Sturmian theory, for a slightly smaller $t_1<t_0$, there are $\th_1<\th_2<\th_3<\th_4$ so that
		\[\eta(t_1,\th_1)<0,\qquad \eta(t_1,\th_2)>0,\qquad \eta(t_1,\th_3)<0,\qquad \eta(t_1,\th_4)>0.\]
		Similar to the proof of Theorem \ref{thm-heat:no_interior_segment}(ii), we need a boundary maximum principle. For an $\epsilon\ll1$, consider the function and compact set
		\[\psi=e^{-\epsilon t-t}\eta,\qquad\qquad\bar\cJ_1=\bar\cJ\cap\{t\leq t_1\}.\]
		Arguing as in Claim \ref{claim-heat:jordan_loop}, using \eqref{eq-ereg:convergence_th_t}, one may show that $\p\cJ_1$ is a Jordan loop.
		
		\begin{claim}\label{claim-ereg:no_interior_extrema}
			Let $L$ be relatively open in $\bar\cJ_1$. Let $\p'L$ be the relative boundary of $L$ in $\bar\cJ_1$.
			\begin{enumerate}[label={(\roman*)}, nosep]
				\item If $\psi\geq0$ in $\bar L$ and $\max_{\bar L}\psi>\max_{\p'L}\psi$, then $\psi$ must be maximized in $\bar L$ at a point $(t,\th)$ with $t=\uT$ or $\th=\uth_t$.
				\item If $\psi\leq0$ in $\bar L$ and $\min_{\bar L}\psi<\min_{\p'L}\psi$, then $\psi$ must be minimized in $\bar L$ at a point $(t,\th)$ with $t=\uT$ or $\th=\oth_t$.
			\end{enumerate}
		\end{claim}
		\begin{proof}
			We prove (i); then (ii) follows from the same argument with the opposite signs. In each parabolic box $R=(t_2,t_3]\times(\th_5,\th_6)\subset\bar\cJ_1$, we have shown that $h$ is smooth and satisfies $\p_th=\p_{\th\th}h+h$. Hence
			\[\p_t\psi-\p_{\th\th}\psi=e^{-\epsilon t-t}\big(\p_t\eta-\p_{\th\th}\eta-\eta-\epsilon\eta\big)=-\epsilon\psi.\]
			By the standard maximum principle and Claim \ref{claim-ereg:backward_openness}, $\psi$ cannot be maximized in a point $(t,\th)$ with $t>\uT$ and $\th\in\Int(J_t)$. It remains to rule out the case $\th=\oth_t>\uth_t$.
			
			Suppose otherwise that $\psi$ is maximized in $\bar L$ at $(t,\oth_t)$ with $\uT<t\leq t_1$. Denote $\lambda=\psi(t,\oth_t)=\max_{\bar L}\psi>0$. The maximization implies
			\begin{equation}\label{eq-ereg:aux1}
				\p_\th^-\psi(t,\oth_t)\geq0\qquad\Rightarrow\qquad\p_\th^- h(t,\oth_t)\leq\p_\th\phi(t,\oth_t).
			\end{equation}
			The remaining argument is similar to Claim \ref{claim-heat:no_interior_extremal} in Section \ref{sec:heat}. 
			Denote $\varphi=\phi-e^{\epsilon t+t}\lambda$. By the maximization of $\psi$ and the relative openness of $L$, we have $\varphi(t,\oth_t)=h(t,\oth_t)$ and
			\[\varphi\leq h\qquad\text{in}\ \ \bar\cJ\cap\big([t-\delta,t]\times[\oth_t-\delta,\oth_t+\delta]\big)\]
			for some $\delta>0$. 
			
			Denote $z=\varphi(t,\oth_t)\nu_{\oth_t}+\p_\th\varphi(t,\oth_t)\tau_{\oth_t}$. We need to verify $\oth_t\in J_t$ and $z\in\Th^{-1}(t,\oth_t)$. Note that
			\[\Th_0^{-1}(t,\oth_t)\cap\gamma_t=\Big\{h(t,\oth_t)\nu_{\oth_t}+c\tau_{\oth_t}: c\geq\p_\th^- h(t,\oth_t)\Big\}\cap\bar{Q_\delta(3)}.\]
			By \eqref{eq-ereg:aux1} and $\varphi(t,\oth_t)=h(t,\oth_t)$, $\p_\th\varphi=\p_\th\phi$, it is clear that $z$ is contained in the first half of the expression. Using Claim \ref{claim-ereg:h_small} and \eqref{eq-ereg:p_th_phi}, the coordinates of $z$ are estimated as
			\[\begin{aligned}
				|x| &= \big|h(t,\oth_t)\sin\oth_t+\p_\th\varphi(t,\oth_t)\cos\oth_t\big|
				\leq e^{-6}\delta\cdot e^{-10}\delta+1.7<3, \\
				|y| &= \big|\!-h(t,\oth_t)\cos\oth_t+\p_\th\varphi(t,\oth_t)\sin\oth_t\big|\leq e^{-6}\delta+1.7\cdot e^{-10}\delta<3\delta,
			\end{aligned}\] 
			hence $z\in\Th_0^{-1}(t,\oth_t)\cap\gamma_t\cap Q_\delta(3)$. It follows that $z\in\gamma_t\cap\Int(E_4)$, and hence $\oth_t\in J_t$ and $z\in\Th^{-1}(t,\oth_t)$ as desired. Then, observe that $\chi_\sigma=1$ as in Theorem \ref{thm-heat:main}\ref{item-heat:boundary_viscosity}, since the outer unit normal of $Z_t$ points upward on $\sigma$. Now all the conditions for Theorem \ref{thm-heat:main}\ref{item-heat:boundary_viscosity} (with $\Omega=\Int(E_4)$) are met, and it yields $(\p_t\varphi-\p_{\th\th}\varphi-\varphi)(t,\oth_t)\geq0$. But this contradicts
			\[\p_t\varphi-\p_{\th\th}\varphi-\varphi=-\epsilon e^{\epsilon t+t}\lambda<0. \qedhere\]
		\end{proof}
		
		Let $V_1,V_3$ be the connected components of $\{\eta<0\}\cap\bar\cJ_1$ containing $\th_1,\th_3$, and $V_2,V_4$ be the connected components of $\{\eta>0\}\cap\bar\cJ_1$ containing $\th_2,\th_4$, respectively.
		
		\begin{claim}
			We have $V_1\ne V_3$, $V_2\ne V_4$, and $V_k\cap\bar J_t\ne\emptyset$ for all $t\in[\uT,t_1]$ and $k=1,2,3,4$.
		\end{claim}
		\begin{proof}
			If $V_1=V_3$, then it separates $V_2$ from the parabolic boundary $\p\bar\cJ\cap\{t<\oT\}$, and $\psi$ would have a positive maximum $(t,\th)$ in $V_2$ with $t>\uT$ and $\th\in\Int(\bar J_t)$ for some small $\epsilon$, contradicting Claim \ref{claim-ereg:no_interior_extrema}. Similarly, $V_2\ne V_4$. Suppose $V_k\cap\bar J_{t'}=\emptyset$ for some $t'\in[\uT,t_1)$ and $k\in\{1,3\}$. Then the connectedness implies $V_k\subset\{t>t'\}$. By Claim \ref{claim-ereg:no_interior_extrema}, the minimum of $\eta$ in $V_k$ is attained at $(t,\oth_t)$ for some $t\in(t',t_1]$. However, by connectedness, this would imply
			\[V_4\cap\{t=\uT\}=\emptyset,\qquad V_4\cap\{\th=\uth_t\}=\emptyset,\]
			contradicting Claim \ref{claim-ereg:no_interior_extrema}(i) by looking at $\max_{V_4}(\psi)$ for some small $\epsilon$. The same argument works for $V_2,V_4$ as well, by looking at $\min_{V_1}(\psi)$.
		\end{proof}
		
		Hence $V_k\cap\bar J_\uT\ne\emptyset$ for $k=1,2,3,4$. By connectedness, note that
		\[\sup_{V_k\cap\bar J_\uT}(\th)\leq\inf_{V_{k+1}\cap\bar J_\uT}(\th),\]
		hence, there are points $\omega_1<\omega_2<\omega_3<\omega_4\in\bar J_\uT$ so that
		\begin{equation}\label{eq-ereg:uT_has_4_root}
			\eta(\uT,\omega_1)<0,\qquad \eta(\uT,\omega_2)>0,\qquad
			\eta(\uT,\omega_3)<0,\qquad \eta(\uT,\omega_4)>0.
		\end{equation}
		Recall that $\Th_0^{-1}(\bar J_\uT)=\gamma_\uT^+$, which is the union of a concave graph $g_{1,\uT}^+$ and a convex graph $g_{2,\uT}^+$ over $[b_\uT^+,3]$. They have curvature bounded by $2e^{-7}\delta$ in $Q_\delta(2)$, due to Claim \ref{claim-ereg:reg_g_uT}. Note that $h(\uT,\cdot)$ is the support function of $\gamma_\uT^+$. Denoting $\varpi=\arctan f'_{\ridge}(b_\uT^+)$, we have
		\[(\p_{\th\th}h+h)(\uT,\cdot)<0\text{\ \ in\ \ }(\uth_\uT,\varpi),\qquad(\p_{\th\th}h+h)(\uT,\cdot)>0\text{\ \ in\ \ }(\varpi,\oth_\uT)\]
		in the barrier sense, whenever these intervals are nonempty. If $h$ is not differentiable at a point $\th$ in any of these intervals, then $(\p_{\th\th}h+h)(\uT,\cdot)$ forms a negative resp. positive Dirac delta at $(\uT,\th)$, see the proof of \eqref{eq-prelim:viscosity_convexity}. If $h(\uT,\cdot)$ is differentiable at $\th$, then
		\[|\p_\th h(\uT,\th)|\leq1.9\qquad\Rightarrow\qquad h(\uT,\th)\nu_\th+\p_\th h(\uT,\th)\tau_\th\in Q_\delta(2),\]
		by $|h|\leq e^{-6}\delta$ (Claim \ref{claim-ereg:h_small}) and $|\nu+\p_y|<e^{-10}\delta$, hence the curvature bound gives $\pm(\p_{\th\th}h+h)(\uT,\th)\geq e^7/(2\delta)$, see \eqref{eq-prelim:viscosity_convexity}. In summary, we have
		\[\begin{aligned}
			& (\p_{\th\th}h+h)(\uT,\th)\leq-e^7/(2\delta)\qquad\text{whenever}\qquad\th\in(\uth_\uT,\varpi)\ \ \text{and} \\
			&\hspace{168pt} |\p_\th h(\uT,\th)|\leq1.9\text{\ \ or\ \ $\p_\th h(\uT,\th)$ does not exist}
		\end{aligned}\]
		and
		\[\begin{aligned}
			& (\p_{\th\th}h+h)(\uT,\th)\geq e^7/(2\delta)\qquad\text{whenever}\qquad\th\in(\varpi,\oth_\uT)\ \ \text{and} \\
			&\hspace{168pt} |\p_\th h(\uT,\th)|\leq1.9\text{\ \ or\ \ $\p_\th h(\uT,\th)$ does not exist}.
		\end{aligned}\]
		Recalling $|h|\leq e^{-6}\delta$, and $|\p_\th\phi|\leq1.7$ by \eqref{eq-ereg:p_th_phi}, and $|\p_{\th\th}\phi|\leq 3.1e^5\delta^{-1}$ by \eqref{eq-ereg:p_thth_phi}, and the definition $\eta=\phi-h$, we obtain
		\begin{align}
			& \begin{aligned}
				& \p_{\th\th}\eta(\uT,\th)\geq1\qquad\text{whenever}\qquad\th\in(\uth_\uT,\varpi)\ \ \text{and} \\
				&\hspace{168pt} |\p_\th \eta(\uT,\th)|\leq0.1\text{\ \ or\ \ $\p_\th\eta(\uT,\th)$ does not exist},
			\end{aligned} \label{eq-ereg:aux8}\\
			& \begin{aligned}
				& \p_{\th\th}\eta(\uT,\th)\leq-1\qquad\text{whenever}\qquad\th\in(\varpi,\oth_\uT)\ \ \text{and} \\
				&\hspace{168pt} |\p_\th\eta(\uT,\th)|\leq0.1\text{\ \ or\ \ $\p_\th\eta(\uT,\th)$ does not exist},
			\end{aligned}\label{eq-ereg:aux9}
		\end{align}
		both in the barrier sense. On the other hand, recall from \eqref{eq-ereg:uT_has_4_root} that
		\[\eta(\uT,\omega_1)<0,\qquad \eta(\uT,\omega_2)>0,\qquad
		\eta(\uT,\omega_3)<0,\qquad \eta(\uT,\omega_4)>0,\]
		where $\uth_\uT\leq\omega_1<\omega_2<\omega_3<\omega_4\leq\oth_\uT$. Take any $\omega_5\in(\omega_2,\omega_3)$ with $\eta(\uT,\omega_5)=0$. There is a point $\omega_6\in(\omega_1,\omega_5)$ where $\eta(\uT,\cdot)$ is maximized. In view of \eqref{eq-ereg:aux8}, we obtain $\omega_6\geq\varpi$. There is another point $\omega_7\in(\omega_5,\omega_4)$ where $\eta(\uT,\cdot)$ is minimized. In view of \eqref{eq-ereg:aux9}, we obtain $\omega_7\leq\varpi$. Thus $\varpi\geq\omega_7>\omega_5>\omega_6\geq\varpi$, contradiction.
	\end{proof}
	
	\begin{claim}\label{claim-ereg:reg_gt_E4}
		For each $t\in(\uT,\oT]$ it holds
		\begin{equation}\label{eq-ereg:reg_gt_E4}
			\|g_{1,t}\|_{C^{1,1/2}([b_t,3]\cap(-1.5,1.5))}\leq e^{-6}\delta,\qquad \|g_{2,t}\|_{C^{1,1/2}([b_t,3]\cap(-1.5,1.5))}\leq e^{-6}\delta.
		\end{equation}
		For each $t\in[\uT,\oT)$ it holds
		\begin{equation}\label{eq-ereg:reg_gtplus_E4}
			\|g_{1,t}^+\|_{C^{1,1/2}([b_t^+,3]\cap(-1.5,1.5))}\leq e^{-6}\delta,\qquad
			\|g_{2,t}^+\|_{C^{1,1/2}([b_t^+,3]\cap(-1.5,1.5))}\leq e^{-6}\delta.
		\end{equation}
		We have
		\begin{equation}\label{eq-ereg:reg_f_E4}
			\|f_{\ridge}\|_{C^{1,1/2}([b_\uT^+,3]\cap[-1,1])}\leq e^{-2}\delta.
		\end{equation}
		For almost every $z=(x,y)\in E_4\cap\bar{Q_\delta(1)}\setminus\gamma$, it holds
		\begin{equation}\label{eq-ereg:dw_in_E4}
			|\D w(z)|\leq e^{-4}\delta^{4/3}|y-f_{\ridge}(x)|^{-1/3}.
		\end{equation}
	\end{claim}
	\begin{proof}
		For all $t\in(\uT,\oT]$, $\th\in\Int(J_t)$ with $|\p_\th h(t,\th)|\leq1.6$ and $|\p_{\th\th}h|\leq e^5\delta^{-1}$, Claim \ref{claim-ereg:dth3_h} gives
		\[\p_\th(\p_{\th\th}h+h)(t,\th)\geq e^{15}\delta^{-2}-1.6\geq e^{14.9}\delta^{-2}.\]
		Recalling \eqref{eq-ereg:curve_decomp}, we may write $\gamma_t=\alpha_1\cup\alpha_2\cup\alpha_3\cup\alpha_4$, where $\alpha_1,\alpha_4$ are (possibly empty) line segments and $\alpha_2,\alpha_3$ are (possibly empty) smooth strictly convex curves. Denote $\alpha=\alpha_2\cup\alpha_3$. To show \eqref{eq-ereg:reg_gt_E4}, it suffices to bound the $C^{1,1/2}$ norms of $\alpha\cap Q_\delta(1.5)$.
		
		Let $s$ be a length parameter of $\alpha$ so that $\alpha(0)$ is the vertex, and $\alpha_2$ corresponds to $s<0$, and $\alpha_3$ corresponds to $s>0$. The angle $\th$ is thus a function in $s$, with $\th(0)=\varpi_t$. For each $s$ with $\alpha(s)=(x,y)\in Q_\delta(1.5)$, we have
		\[|\p_\th h(t,\th(s))|=\metric{\alpha(s)}{\tau_{\th(s)}}\leq|x|+|y|\,|\sin\th(s)|\leq1.5+1.5\delta\cdot e^{-10}\delta<1.6.\]
		Thus, for these $s$ either
		\[|\p_{\th\th}h(t,\th(s))|>e^5\delta^{-1}\qquad\text{or}\qquad\p_\th(\p_{\th\th}h+h)(t,\th(s))\geq e^{14.9}\delta^{-2}.\]
		In view that $|h|\leq e^{-6}\delta$, for these $s$ either
		\begin{equation}\label{eq-ereg:aux6}
			\big|(\p_{\th\th}h+h)(t,\th(s))\big|>e^{4.9}\delta^{-1}\qquad\text{or}\qquad\p_\th(\p_{\th\th}h+h)(t,\th(s))\geq e^{14.9}\delta^{-2}.
		\end{equation}
		The set $\{s>0:\alpha(s)\in Q_\delta(1.5)\}$ is a (possibly empty) interval $(s_1,s_2)$ with $0\leq s_1$. See Figure \ref{fig-ereg:s1_s4}. For any $s\in(s_1,s_2)$ we argue as follows: if $(\p_{\th\th}h+h)(t,\th(s))<e^{4.9}\delta^{-1}$, then by \eqref{eq-ereg:aux6}, we have $\p_\th(\p_{\th\th}h+h)(t,\th(s'))\geq e^{14.9}\delta^{-2}$ for all $s'\in[s_1,s]$. Hence
		\[(\p_{\th\th}h+h)(t,\th(s))\geq e^{14.9}\delta^{-2}\big(\th(s)-\th(s_1)\big).\]
		If $(\p_{\th\th}h+h)(t,\th(s))\geq e^{4.9}\delta^{-1}$, then trivially
		\[(\p_{\th\th}h+h)(t,\th(s))\geq\frac12e^{14.9}\delta^{-2}\big(\th(s)-\th(s_1)\big)\]
		since $\th(s),\th(s_1)\in(-e^{-10}\delta,e^{-10}\delta)$. Combining both cases, we always have
		\begin{equation}\label{eq-ereg:C112_bound_1}
			(\p_{\th\th}h+h)(t,\th(s))\geq \frac12e^{14.9}\delta^{-2}\big(\th(s)-\th(s_1)\big)\qquad\forall\,s\in(s_1,s_2).
		\end{equation}
		Therefore,
		\[\frac{d\big(\th(s)-\th(s_1)\big)}{ds}=\frac1{(\p_{\th\th}h+h)(t,\th(s))}\leq\frac{2e^{-14.9}\delta^2}{\th(s)-\th(s_1)}.\]
		Integrating this, we obtain
		\begin{equation}\label{eq-ereg:C112_bound_2}
			\th(s')-\th(s)\leq e^{-6.5}\delta(s'-s)^{1/2},\qquad\forall\,(s,s')\subset(s_1,s_2).
		\end{equation}
		Similarly, the set $\{s<0:\alpha(s)\in Q_\delta(1.5)\}$ is an interval $(s_3,s_4)$ with $s_4\leq0$, where it holds
		\begin{equation}\label{eq-ereg:C112_bound_3}
			(\p_{\th\th}h+h)(t,\th(s))\leq -\frac12e^{14.9}\delta^{-2}\big(\th(s_4)-\th(s)\big)
		\end{equation}
		and
		\begin{equation}\label{eq-ereg:C112_bound_4}
			\th(s')-\th(s)\leq e^{-6.5}\delta(s'-s)^{1/2}\qquad\forall\,(s,s')\subset(s_3,s_4).
		\end{equation}
		Then, \eqref{eq-ereg:C112_bound_2} \eqref{eq-ereg:C112_bound_4} already imply \eqref{eq-ereg:reg_gt_E4}. Also, recall the convergences $g_{1,t}^+=\lim_{s\searrow t}g_{1,s}$, $g_{2,t}^+=\lim_{s\searrow t}g_{2,s}$. Then \eqref{eq-ereg:reg_gtplus_E4} follows from \eqref{eq-ereg:reg_gt_E4} and this convergence.
		
		\begin{figure}
			\centering
			\includegraphics{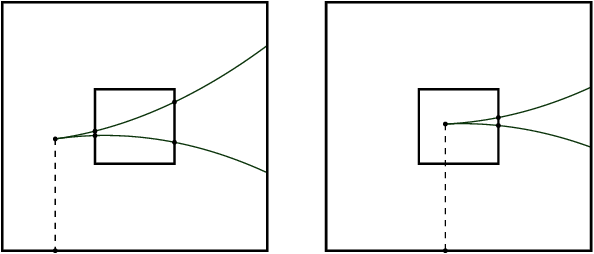}
			\begin{picture}(0,0)
				\put(-255,63){$s_1$}
				\put(-205,81){$s_2$}
				\put(-205,44){$s_3$}
				\put(-255,48){$s_4$}
				\put(-130,67){$s_1=s_4=0$}
				\put(-49,53){$s_3$}
				\put(-49,71){$s_2$}
				\put(-159,33){$\alpha_2$}
				\put(-159,97){$\alpha_3$}
				\put(-3,46){$\alpha_2$}
				\put(-3,75){$\alpha_3$}
				\put(-260,6){$b_t<-1.5$}
				\put(-72,6){$b_t\geq-1.5$}
			\end{picture}
			\caption{The intervals $(s_1,s_2)$ and $(s_3,s_4)$.}\label{fig-ereg:s1_s4}
		\end{figure}
		
		Denote $b=\max\{b_t,-1.5\}$. For our later convenience, let us show that
		\begin{equation}\label{eq-ereg:curv_est_1}
			|\p_{\th\th}h+h|(t,\th(s))\geq e^4\delta^{-1}|x(s)-b|^{1/2},\qquad\forall\,s\in(s_1,s_2)\cup(s_3,s_4),
		\end{equation}
		where $x(s)$ is the $x$-coordinate of $\alpha(s)$. Suppose $s\in(s_1,s_2)$. If $|\p_{\th\th}h+h|(t,\th(s))\geq e^{4.9}\delta^{-1}$, then the result is obvious. If not, then $(\p_{\th\th}h+h)(t,\cdot)$ is increasing in $(s_1,s)$, hence
		\[\begin{aligned}
			x(s)-b &= \int_{\th(s_1)}^{\th(s)}(\p_{\th\th}h+h)(t,\varpi)\cos\varpi\,d\varpi\leq(\p_{\th\th}h+h)(t,\th(s))\cdot\big(\th(s)-\th(s_1)\big) \\
			&\leq e^{-14}\delta^2\cdot(\p_{\th\th}h+h)(t,\th(s))^2\qquad\text{by \eqref{eq-ereg:C112_bound_1}},
		\end{aligned}\]
		as desired. The argument is similar if $s\in(s_3,s_4)$.
		
		\vspace{3pt}
		
		Next, we prove \eqref{eq-ereg:reg_f_E4}. Let us first show that
		\begin{equation}\label{eq-ereg:reg_f_partial}
			\begin{aligned}
				& \big|f'_{\ridge}(b_t^-)-f'_{\ridge}(b_s^+)\big|\leq e^{-3}\delta(b_t^--b_s^+)^{1/2}, \\
				&\hspace{120pt}\text{whenever}\qquad\uT\leq s<t\leq\oT,\ -1.2\leq b_s^+<b_t\leq1.2.
			\end{aligned}
		\end{equation}
		The quantity $b_t^-$ is defined in Claim \ref{claim-ereg:fridge_linear}. It is related to the structure of $\tZ_t$, see \eqref{eq-ereg:tZt_description}. The properties of $b_t^-$ that we need are the following: there is a concave function $\tg_{1,t}:[b_t^-,3]\to(-3\delta,3\delta)$ and a convex function $\tg_{2,t}:[b_t^-,3]\to(-3\delta,3\delta)$, so that
		\begin{equation}\label{eq-ereg:btminus_order}
			g_{1,s}^+<\tg_{1,t}\leq\tg_{2,t}<g_{2,s}^+\qquad\text{in}\ \ [b_t^-,3],
		\end{equation}
		and
		\begin{equation}\label{eq-ereg:btminus_endpt}
			\tg_{1,t}(b_t^-)=\tg_{2,t}(b_t^-)=f_{\ridge}(b_t^-),\qquad \tg_{1,t}'(b_t^-)=\tg_{2,t}'(b_t^-)=f_{\ridge}'(b_t^-),
		\end{equation}
		and
		\begin{equation}\label{eq-ereg:btminus_reg}
			\|\tg_{1,t}\|_{C^{1,1/2}([b_t^-,3]\cap(-1.5,1.5))}\leq e^{-6}\delta,\qquad \|\tg_{2,t}\|_{C^{1,1/2}([b_t^-,3]\cap(-1.5,1.5))}\leq e^{-6}\delta.
		\end{equation}
		The regularity \eqref{eq-ereg:btminus_reg} follows from \eqref{eq-ereg:reg_gt_E4} and the fact that $\tg_{i,t}\equiv g_{i,t}$ in $[b_t,3]$ and $\tg_{i,t}$ is linear in $[b_t^-,b_t]$, see \eqref{eq-ereg:tgt_properties}.
		
		Denote $\Delta x=b_t^--b_s^+$. First, suppose otherwise that $f'_{\ridge}(b_t^-)>f'_{\ridge}(b_s^+)+e^{-3}\delta\Delta x^{1/2}$. Setting $\rho=\Delta x/10$, noting that $b_t^-+\rho<1.5$, we estimate
		\[\begin{aligned}
			\tg_{2,t}(b_t^-+\rho) &\geq \tg_{2,t}(b_t^-)+\tg'_{2,t}(b_t^-)\rho\qquad\text{(since $\tg_{2,t}$ is convex)} \\
			&\geq g_{1,s}^+(b_t^-)+\big(f'_{\ridge}(b_s^+)+e^{-3}\delta\Delta x^{1/2}\big)\rho \qquad\text{(by our hypothesis)}\\
			&\geq g_{1,s}^+(b_s^+)+f'_{\ridge}(b_s^+)\Delta x-\frac23e^{-6}\delta\Delta x^{3/2}+\big(f'_{\ridge}(b_s^+)+e^{-3}\delta\Delta x^{1/2}\big)\rho \\
			&\hspace{310pt}\text{(by \eqref{eq-ereg:reg_gtplus_E4})} \\
			&\geq g_{2,s}^+(b_t^-+\rho)-f'_{\ridge}(b_s^+)(\Delta x+\rho)-\frac23e^{-6}\delta(\Delta x+\rho)^{3/2}\qquad\text{(by \eqref{eq-ereg:reg_gtplus_E4})} \\
			&\qquad +f'_{\ridge}(b_s^+)\Delta x-\frac23e^{-6}\delta\Delta x^{3/2}+\big(f'_{\ridge}(b_s^+)+e^{-3}\delta\Delta x^{1/2}\big)\rho.
		\end{aligned}\]
		Hence
		\[\frac23e^{-6}\delta(\Delta x+\rho)^{3/2}+\frac23e^{-6}\delta\Delta x^{3/2}\geq e^{-3}\delta\Delta x^{1/2}\rho.\]
		Since $\rho=\Delta x/10$, this gives
		\[\frac23e^{-6}\Big(\frac{11}{10}\Big)^{3/2}+\frac23e^{-6}\geq\frac1{10}e^{-3},\]
		which is not true. We thus obtain $f'_{\ridge}(b_t^-)\leq f'_{\ridge}(b_s^+)+e^{-3}\delta\Delta x^{1/2}$. The other inequality $f'_{\ridge}(b_t^-)\geq f'_{\ridge}(b_s^+)-e^{-3}\delta\Delta x^{1/2}$ may be obtained by comparing $\tg_{1,t}(b_t^-+\rho)$ and $g_{1,s}^+(b_t^-+\rho)$ in the same way. This proves \eqref{eq-ereg:reg_f_partial}.
		
		In general, suppose $x_1,x_2\in[b_\uT^+,3]\cap[-1.2,1.2]$ and $x_1<x_2$. Call $t_1=w(x_1,f_{\ridge}(x_1))$ and $t_2=w(x_2,f_{\ridge}(x_2))$. Note that $\uT\leq t_1\leq t_2\leq\oT$, and $\uT<t_1$ if $x_1>b_\uT^+$. Recall that:
		\begin{itemize}[nosep]
			\item $f_{\ridge}$ is linear in $[b_t^-,b_t]\cap[-3,3]$ for each $t\in(\uT,\oT]$ (Claim \ref{claim-ereg:fridge_linear});
			\item $f_{\ridge}$ coincides with either $g_{1,t}$ or $g_{2,t}$ in $[b_t,b_t^+]$ for each $t\in(\uT,\oT)$ (Claim \ref{claim-ereg:fine_structure_2});
			\item $f_{\ridge}$ either coincides with $g_{1,\oT}, g_{2,\oT}$ or is linear in $[b_\oT,3]$ (Claim \ref{claim-ereg:f_oT});
			\item the $C^{1,1/2}$ norms of $g_{1,t},g_{2,t}$ are bounded by $e^{-6}\delta$ in $(-1.5,1.5)$, by \eqref{eq-ereg:reg_gt_E4} \eqref{eq-ereg:reg_gtplus_E4};
			\item \eqref{eq-ereg:reg_f_partial} holds, which connects $b_s^+$ to $b_t^-$.
		\end{itemize}
		The triangle inequality implies
		\[\begin{aligned}
			\big|f'_{\ridge}(x_2)-f'_{\ridge}(x_1)\big| &\leq \big|f'_{\ridge}(x_2)-f'_{\ridge}(b_{t_2}^-)\big|+\big|f'_{\ridge}(b_{t_2}^-)-f'_{\ridge}(b_{t_1}^+)\big| \\
			&\qquad +\big|f'_{\ridge}(b_{t_1}^+)-f'_{\ridge}(x_1)\big| \\
			&\leq \sqrt 3\cdot e^{-3}\delta|x_1-x_2|^{1/2}
			\leq e^{-2}\delta|x_1-x_2|^{1/2}.
		\end{aligned}\]
		Hence, we have shown
		\begin{equation}\label{eq-ereg:reg_f_E4_strong}
			\|f_{\ridge}\|_{C^{1,1/2}([b_\uT^+,3]\cap[-1.2,1.2])}\leq e^{-2}\delta.
		\end{equation}
		This in particular shows \eqref{eq-ereg:reg_f_E4}.
		
		\vspace{3pt}
		
		Finally, we prove \eqref{eq-ereg:dw_in_E4}. It suffices to show that for all $t\in(\uT,\oT)$ and $z=(x,y)\in\alpha\cap\bar{Q_\delta(1)}\setminus\text{(ridge)}$, where $\alpha$ is defined above, it holds:
		\begin{equation}\label{eq-ereg:claim14_to_prove}
			\text{the curvature of $\alpha$ at $z$ is bounded by $e^{-4}\delta^{4/3}|y-f_{\ridge}(x)|^{-1/3}$.}
		\end{equation}
		Let $s$ be such that $z=\alpha(s)$. Note that $s\in(s_1,s_2)\cup(s_3,s_4)$. We may assume $s\in(s_1,s_2)$ (the other case is analogous). In \eqref{eq-ereg:curv_est_1}, if $b_t<-1.2$, then
		\[|\p_{\th\th}h+h|(t,\th(s))\geq e^4\delta^{-1}(0.2)^{1/2}\geq e^4\delta^{-4/3}(6e^{-9}\delta)^{1/3}\geq e^4\delta^{-4/3}|y-f_{\ridge}(x)|^{1/3},\]
		as desired. If $b_t\geq-1.2$, then \eqref{eq-ereg:curv_est_1} gives
		\begin{equation}\label{eq-ereg:curv_est_2}
			|\p_{\th\th}h+h|(t,\th(s))\geq e^4\delta^{-1}\big(x(s)-b_t\big)^{1/2}.
		\end{equation}
		The known $C^{1,1/2}$ bounds in $[-1.2,1.2]$ give
		\begin{equation}\label{eq-ereg:aux12}
			\begin{aligned}
				y(s)-f_{\ridge}(x(s)) &= g_{2,t}(x(s))-f_{\ridge}(x(s))
				\leq \int_{b_t}^{x(s)}|g'_{2,t}-f'_{\ridge}|(x)\,dx \\
				&\leq \int_{b_t}^{x(s)} (e^{-6}+e^{-2})\delta(x-b_t)^{1/2}\,dx
				\leq e^{-2}\delta\big(x(s)-b_t\big)^{3/2}
			\end{aligned}
		\end{equation}
		Combining \eqref{eq-ereg:curv_est_2} \eqref{eq-ereg:aux12}, we also obtain
		\[|\p_{\th\th}h+h|(t,\th(s))\geq e^4\delta^{-1}\Big[e^2\delta^{-1}\big(y(s)-f_{\ridge}(x(s))\big)\Big]^{1/3}. \qedhere\]
	\end{proof}
	
	\begin{claim}\label{claim-ereg:reg_w_E4}
		$\|w\|_{C^{0,2/3}(E_4\cap Q_\delta(1))}\leq e^{-2}\delta^{4/3}$.
	\end{claim}
	\begin{proof}
		\eqref{eq-ereg:dw_in_E4} implies for all $(x,y_1),(x,y_2)\in E_4\cap Q_\delta(1)$ that
		\begin{equation}\label{eq-ereg:reg_w_vertical}
			\big|w(x,y_1)-w(x,y_2)\big|\leq e^{-4}\delta^{4/3}\int_{y_1}^{y_2}|y-f_{\ridge}(x)|^{-1/3}\,dy\leq e^{-3}\delta^{4/3}|y_1-y_2|^{2/3}.
		\end{equation}
		Also, for all $t\in[\uT,\oT)$ and $x$ so that $-1<b_t^+<x<1$, we have
		\[\begin{aligned}
			w(x,f_{\ridge}(x))-w(b_t^+,f_{\ridge}(b_t^+)) &= w(x,f_{\ridge}(x))-w(x,g_{1,t}^+(x)) \\
			&\leq e^{-3}\delta^{4/3}|f_{\ridge}(x)-g_{1,t}^+(x)|^{2/3} \\
			&\leq e^{-4}\delta^2|x-b_t^+|,
		\end{aligned}\]
		where the last line is obtained similarly as \eqref{eq-ereg:aux12}. In view that $w(x,f_{\ridge}(x))$ is non-decreasing for $x\in[b_\uT^+,3]$ and is constant on $[b_t^-,b_t^+]$ for each $t$, it follows that
		\begin{equation}\label{eq-ereg:reg_w_on_gamma}
			\big\|w\big(\cdot,f_{\ridge}(\cdot)\big)\big\|_{\Lip([b_\uT^+,3]\cap(-1,1))}\leq e^{-4}\delta^2.
		\end{equation}
		
		Suppose $z_1=(x_1,y_1)\in E_4\cap Q_\delta(1)$ and $z_2=(x_2,y_2)\in E_4\cap Q_\delta(1)$. If $|y_1-f_{\ridge}(x_1)|\leq 8e^{-9}\delta|x_1-x_2|$, then note that
		\begin{equation}\label{eq-ereg:triangle}
			\big|y_1-f_{\ridge}(x_1)\big|^{2/3}+|x_1-x_2|^{2/3}+|y_2-f_{\ridge}(x_2)|^{2/3}\leq 2|z_1-z_2|^{2/3}.
		\end{equation}
		Hence
		\[\begin{aligned}
			|w(z_1)-w(z_2)| &\leq \big|w(z_1)-w(x_1,f_{\ridge}(x_1))\big|
				+ \big|w(x_1,f_{\ridge}(x_1))-w(x_2,f_{\ridge}(x_2))\big| \\
			&\qquad +\big|w(x_2,f_{\ridge}(x_2))-w(z_2)\big| \\
			&\leq e^{-3}\delta^{4/3}\Big(\big|y_1-f_{\ridge}(x_1)\big|^{2/3}+|y_2-f_{\ridge}(x_2)|^{2/3}\Big)+e^{-4}\delta^2|x_1-x_2| \\
			&\leq e^{-2}\delta^{4/3}|z_1-z_2|^{2/3}.
		\end{aligned}\]
		The argument is the same if $|y_2-f_{\ridge}(x_2)|\leq 8e^{-9}\delta|x_1-x_2|$ or if $y_1-f_{\ridge}(x_1)$ has opposite sign with $y_2-f_{\ridge}(x_2)$, since \eqref{eq-ereg:triangle} remains true in either case. It remains to treat the case
		\[y_1>f_{\ridge}(x_1)+8e^{-9}\delta|x_1-x_2|,\qquad y_2>f_{\ridge}(x_2)+8e^{-9}\delta|x_1-x_2|.\]
		The case with the opposite sign can be proved similarly. We may assume $x_1<x_2$. Find a segment $\sigma\subset E_4$ with slope at most $e^{-9}\delta$, so that it connects $z_1$ with another point $z_3=(x_2,y_3)$. Then note that
		\begin{equation}\label{eq-ereg:triangle2}
			|x_1-x_2|^{2/3}+|y_3-y_2|^{2/3}\leq2|z_1-z_2|^{2/3}.
		\end{equation}
		Notice $|y-f_{\ridge}(x)|\geq 6e^{-9}\delta|x_1-x_2|$ for all $(x,y)\in\sigma$. Then \eqref{eq-ereg:dw_in_E4} and $|\nu+\p_y|<e^{-10}\delta$ yield
		\[|\D w|\leq e^{-1}\delta^{4/3}|x_1-x_2|^{-1/3},\qquad|\metric{\D w}{\p_x}|\leq e^{-10}\delta|\D w|.\]
		Hence
		\begin{equation}\label{eq-ereg:aux13}
			\begin{aligned}
				\big|w(z_1)-w(z_3)\big| &\leq e^{-1}\delta|x_1-x_2|^{-1/3}\big(e^{-10}\delta|x_1-x_2|+|y_1-y_3|\big) \\
				&\leq e^{-9}\delta^2|x_1-x_2|^{2/3}.
			\end{aligned}
		\end{equation}
		Also, by \eqref{eq-ereg:reg_w_vertical},
		\begin{equation}\label{eq-ereg:aux14}
			\big|w(z_3)-w(z_2)\big|\leq e^{-3}\delta^{4/3}|y_3-y_2|^{2/3}.
		\end{equation}
		Combining \eqref{eq-ereg:triangle2} \eqref{eq-ereg:aux13} \eqref{eq-ereg:aux14} finally yields
		\[\big|w(z_1)-w(z_2)\big|\leq e^{-2}\delta^{4/3}|z_1-z_2|^{2/3}.\qedhere\]
	\end{proof}
	
	Combining Claim \ref{claim-ereg:reg_w_E4} and \eqref{eq-ereg:reg_f_E4}, we prove Claim \ref{claim-ereg:reg_in_E3E4}.
\end{proof}

%\newpage

\section{Sublevel sets of \texorpdfstring{$|\D u|$}{|∇u|}}\label{sec:lvlset}

Suppose $u$ is $\infty$-harmonic in a simply-connected domain $\Omega\subset\RR^2$. Denote
\[Y_t=\big\{|\D u|<e^{-t}\big\},\qquad
	\tZ_t=\big\{|\D u|\leq e^{-t}\big\},\qquad Z_t=\Int(\tZ_t).\]

\begin{theorem}\label{thm-lvlset:global}
	Suppose $\Delta_\infty u=0$ in a simply-connected domain $\Omega\subset\RR^2$. Then:
	\begin{enumerate}[label={(\roman*)}, nosep]
		\item Each $Y_t,Z_t$ is a disjoint union of simply-connected concave polygons. Each boundary edge of these polygons is a streamline of $u$.
		\item $\tZ_t\setminus\bar{Z_t}$ is a disjoint union of nontrivial line segments that are streamlines of $u$.
		\item Let $\sigma$ be a connected component of $\tZ_t\setminus\bar{Z_t}$. If $\sigma\Subset\Omega$, then $\bar\sigma\cap\p Z_t\ne\emptyset$.
		\item Each $Z_t\setminus Y_t$ is a disjoint union of nontrivial line segments whose endpoints either lie on $\p Z_t$ or is $\infty$. Each line segment is a streamline of $u$.
	\end{enumerate}
\end{theorem}

\begin{figure}[ht]
	\centering
	\includegraphics{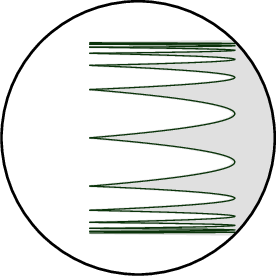}
	\caption{Impossible picture of accumulating edges.}\label{fig-lvlset:impossible}
\end{figure}

Item (i) a priori rules out accumulation of edges of $\p Y_t$ (Figure \ref{fig-lvlset:impossible}), since any compact set can only intersect finitely many edges of a concave polygon. We do not know whether such pathology occurs when $\Omega$ is not simply-connected. Note that $Y_t$ or $Z_t$ may have infinitely many connected components.

\begin{proof}[Proof of Theorem \ref{thm-lvlset:global}] {\ }
	
	Let $\{u_p\}$ be a sequence of $p$-harmonic functions with $p\to\infty$ and $u_p\to u$ in $C^0_{\loc}(\Omega)$. For a basepoint $z_0\in\Omega$ to be determined below, Theorem \ref{thm-approx:global} yields a simple IMCF cluster $\big(w,\nu,\{D_i\},\{\chi_i\}\big)$ in a region $\Omega_{\reg}=\{w<+\infty\}$. For each $t\in\RR$ we denote
	\begin{equation}\label{eq-lvlset:aux1}
		Y'_t=\{w>t\},\qquad\tZ'_t=\{w\geq t\},\qquad Z'_t=\Int(\tZ'_t).
	\end{equation}
	We make implicit their dependence on $z_0$.
	
	\vspace{3pt}
	
	(i) Let $F$ be a connected component of $Y_t$ or $Z_t$. If $F$ is not simply-connected, then there is a domain $D\Subset\Omega$ with $\p D\subset F$ but $D\not\subset F$. If $F\subset Y_t$, then $|\D u|<e^{-t}$ on $\p D$ but $|\D u|\geq e^{-t}$ somewhere in $D$. If $F\subset Z_t$, then $|\D u|\leq e^{-t}$ on $\p D$ but $|\D u|>e^{-t}$ somewhere in $D$. These violate the maximum principle for $|\D u|$ (Lemma \ref{lemma-prelim:max_principle_du}).
	
	Take a basepoint $z_0\in F$. Let $F'$ denote the connected component of either $Y'_t$ or $Z'_t$ that contains $z_0$. By Theorem \ref{thm-approx:global}\ref{item-approx:lvlset}, we have $F=F'$. Recall that $\Omega\setminus\Omega_{\reg}\subset\Crit(u)$ while $|\D u|=e^{-t}$ on $\p F\cap\Omega$, hence $\p F\cap\Omega\subset\Omega_{\reg}$. If $0\in\p F\cap\Omega$ and $\nu(0)=-\p_y$ (the general case is reduced to this via a rigid motion), then by Lemma \ref{lemma-cluster:local_model}(i) (with $\delta=1$, then restricting to the connected component $F$), for a small $r$ we have
	\[F\cap Q(r)=\Big\{x\in I, g_1(x)<y<g_2(x)\Big\},\]
	where $I=(-r,r)$ or $I=(-r,a)\cup(b,r)$ with $-r\leq a\leq b\leq r$, and $g_1,g_2:\bar I\cap(-r,r)\to[-r,r]$ are functions with certain convexity conditions as stated in Lemma \ref{lemma-cluster:local_model}\ref{item-cluster:g1_g2}.
	
	To show that $F$ is a concave polygon, it suffices to rule out the case $a=b=0$. If this happens, then by the connectedness of $F$, we may find a simple loop $\sigma\subset\bar F$ so that $\sigma\setminus\{0\}\subset F$ and $\sigma\cap Q(r)=\graphh((g_1+g_2)/2)$. Let $D$ be the region enclosed by $\sigma$, then $D\not\subset F$. If $F\subset Z_t$, then $|\D u|\leq e^{-t}$ on $\sigma$ but $|\D u|>e^{-t}$ somewhere in $D$, contradicting Lemma \ref{lemma-prelim:max_principle_du}. If $F\subset Y_t$, then $|\D u|<e^{-t}$ on $\sigma\setminus\{0\}$ but $|\D u|\geq e^{-t}$ somewhere in $D$. But Lemma \ref{lemma-prelim:interior_max_gradient} implies that there are at least two points on $\sigma$ on which $|\D u|\geq e^{-t}$, contradiction. Thus, we have shown that $F$ is a concave polygon. Since $\nu=-e^w\D^\perp u$ and $\p F'\perp\nu$, all edges of $\p F$ are streamlines of $u$.
	
	\vspace{3pt}
	
	(ii) For any $z_0\in\tZ_t\setminus\bar{Z_t}$, let $\sigma$ be the connected component of $\tZ_t\setminus\bar{Z_t}$ that contains $z_0$. Let $F$ (resp. $F'$) be the connected component of $\tZ_t$ (resp. $\tZ'_t$) that contains $z_0$. By Theorem \ref{thm-approx:global}\ref{item-approx:lvlset} we have $F=F'\subset\tZ'_t$. Also, we have $Z_t\supset Z'_t$ hence $\bar{Z_t}\supset\bar{Z'_t}$ due to $w\leq-\log|\D u|$. Hence
	\[\sigma\subset F\setminus\bar{Z_t}\subset\tZ'_t\setminus\bar{Z'_t}.\]
	By Lemma \ref{lemma-cluster:local_model}\ref{item-cluster:tZt_setminus_Zt}, the latter is a disjoint union of nontrivial line segments that are contained in the ridges (and recall that ridges are $\nu$-orthogonal). Let $\sigma'$ be the connected component of $\tZ'_t\setminus\bar{Z'_t}$ containing $z_0$ (i.e., it is a maximal line segment). So $\sigma\subset\sigma'$.
	
	Let us further prove that $\bar{\sigma}=\bar{\sigma'}$. Observe that
	\[\tZ'_t\setminus\bar{Z'_t}\subset\{w=t\}\cap\text{(ridges)}\subset\big\{|\D u|=e^{-t}\big\}\subset\tZ_t.\]
	Then clearly $\sigma'\subset\tZ_t$, thus $\sigma'\subset F$. Up to a rigid motion, we may assume that $\sigma'$ is contained in the $x$-axis and $\nu=-\p_y$ on $\sigma'$ (we do not know whether $\sigma'$ is bounded). Take any bounded interval $I$ so that $z_0\in I\times\{0\}\subset\sigma'$. Since $\bar{\sigma'}$ lies in a ridge, for sufficiently small $r$ we have $I\times(-r,r)\Subset\Omega_{\reg}$, and $w$ is an upward evolving weak IMCF in $I\times(-r,0)$ and a downward evolving IMCF in $I\times(0,r)$. Additionally, recall that $w\equiv t$ on $\sigma'$ and that $\sigma'\cap\bar{Z'_t}=\emptyset$. Using the monotonicity of $w$ in the $y$ direction, we have
	\[\tZ'_t\cap\big(I\times(-r,r)\big)=I\times\{0\}.\]
	Since $F=F'$, elementary topology then implies
	\[F\cap\big(I\times(-r,r)\big)=I\times\{0\}.\]
	This shows $I\times\{0\}\subset\tZ_t\setminus\bar{Z_t}$, hence $I\times\{0\}\subset\sigma$. By $\sigma\subset\sigma'$ and the arbitrariness of $I$, it follows that $\bar{\sigma'}=\bar{\sigma}$.
	
	\vspace{3pt}
	
	(iii) Let $\sigma'$ be as in (ii) above, then $\bar\sigma=\bar{\sigma'}$. Then Lemma \ref{lemma-cluster:local_model}\ref{item-cluster:no_interior_segment} gives $\bar{\sigma'}\cap\p Z'_t\ne\emptyset$, hence $\bar\sigma\cap\bar{Z_t}\ne\emptyset$ since $\bar{Z_t}\supset\bar{Z'_t}$, and hence $\bar\sigma\cap\p Z_t\ne\emptyset$ by elementary topology.
	
	\vspace{3pt}
	
	(iv) This is a direct consequence of Lemma \ref{lemma-prelim:interior_max_gradient} with $\text{``$\Omega$''}=Z_t$.
\end{proof}

\section{Precompact sublevel sets}\label{sec:closed}

Let $\Omega\subset\RR^2$ be simply-connected, and $\Delta_\infty u=0$ in $\Omega$. We will implicitly use the fact that $\D u$ is continuous. As usual, we denote
\[Y_t=\big\{|\D u|<e^{-t}\big\},\qquad \tZ_t=\big\{|\D u|\leq e^{-t}\big\},\qquad Z_t=\Int(\tZ_t).\]
We consider the case where a connected component of $Y_t$ is precompact.

\begin{theorem}\label{thm-closed:main}
	Let $\Omega$ be simply-connected and $u$ be $\infty$-harmonic in $\Omega$. Suppose $t_0\in\RR$, and $F$ is a connected component of $Y_{t_0}$ with $F\Subset\Omega$. Then:
	\begin{enumerate}[label={(\roman*)}, nosep]
		\item\label{item-closed:main_finiteness} $\Crit(u)\cap F$ is a nonempty finite set.
		\item\label{item-closed:main_asymp} Let $\{z_k\}$ be the critical points of $u$ in $F$. Then for each $k$ there is $d_k\in\ZZ_{\geq2}$ so that
		\begin{equation}\label{eq-closed:rescaled_sol}
			\lambda^{-\frac{d_k^2}{2d_k-1}}\Big[u\big(\lambda z+z_k\big)-u(z_k)\Big]
		\end{equation}
		converges in $C^1_{\loc}(\RR^2)$ to a quasiradial solution of degree $d_k$, as $\lambda\to0$.
		\item\label{item-closed:main_edges} Let $N$ be the number of edges of $\p F$. Then
		\[N-2=\sum_k\big(2d_k-2\big).\]
		\item\label{item-closed:reg} For each $k$, we have $u\in C^{1,1/3}(B(z_k,r))$ for some $r>0$.
	\end{enumerate}
\end{theorem}

The following is a direct consequence:

\begin{theorem}\label{thm-closed:precompact_crit}
	If $K\Subset\Omega$ and $\Crit(u)\cap\p K=\emptyset$, then $\Crit(u)\cap K$ is a finite set.
\end{theorem}

Recall that each $Y_t,Z_t$ is a disjoint union of simply-connected concave polygons. Thus, each precompact component $E$ of $Y_t,Z_t$ has a well-defined number of edges. We denote it by $N(E)$ in this section. By Lemma \ref{lemma-closed:edge_lower_bound}, we have $N(E)\geq4$ and $N(E)$ is even. We show the following evolution pattern of sublevel sets of $|\D u|$:
\begin{itemize}[nosep]
	\item When $Z_t$ is connected and changes continuously in a time interval $(t_1,t_2)$, the support function of $\p Z_t$ solves $\p_th=\p_{\th\th}h+h$ (Lemma \ref{lemma-closed:1_evolution_lemma}).
	\item When $Z_t$ becomes disconnected or changes noncontinuously at a time $T$, it is cut by finitely many segments into smaller concave polygons (Lemma \ref{lemma-closed:2_split_lemma}).
	\item The quantity $\sum_{E:\,\text{component of $Z_t$}}(N(E)-2)$ is constant in $t$.
	\item There can only be finitely many splitting times. Eventually, the support function solves $\p_th=\p_{\th\th}h+h$ with $t\in[T,\infty)$ for some $T$ in each polygon. The asymptotics at critical points follow by taking $t\to\infty$.
\end{itemize}
All lemmas below apply to any $\infty$-harmonic function $u$ in a simply-connected domain $\Omega$. The following condition plays a central role:
\begin{equation*}\label{eq-closed:condition_Pt}
	Z_t\Subset\Omega,\qquad Z_t=Y_t,\qquad Z_t\text{ is connected}.\tag{\text{$\textrm{P}_t$}}
\end{equation*}

This condition is forward stable:

\begin{lemma}\label{lemma-closed:forward_openness}
	Let $F$ be a connected component of $Y_t$ so that $F\Subset\Omega$. Then there is $T>t$ so that
	$Z_s\cap F$ is connected and is equal to $Y_s\cap F$, for all $s\in(t,T]$.
\end{lemma}

In a time interval where \eqref{eq-closed:condition_Pt} holds, the support function solves $\p_th=\p_{\th\th}h+h$:

\begin{lemma}\label{lemma-closed:1_evolution_lemma}
	Suppose $-\infty\leq t_1<t_2\leq\infty$, and \eqref{eq-closed:condition_Pt} holds for all $t\in(t_1,t_2)$. Then:
	\begin{enumerate}[label={(\roman*)}, nosep]
		\item $t\mapsto N(Y_t)$ is constant in $(t_1,t_2)$;
		\item $\tZ_t=\bar{Z_t}$ for all $t\in(t_1,t_2]$;
		\item $u$ has no critical points in $Y_{t_1}\setminus\tZ_{t_2}$. Denoting the constant in (i) by $N$, the map
		\begin{equation}\label{eq-closed:hodographic_map}
			\Big(\!-\log|\D u|,\,\frac{\D u}{|\D u|}\Big):Y_{t_1}\setminus \tZ_{t_2}\to(t_1,t_2)\times\SS^1
		\end{equation}
		lifts to a degree -1 homeomorphism
		\[\Th:Y_{t_1}\setminus\tZ_{t_2}\to(t_1,t_2)\times\SS^1\big((N-2)\pi\big).\]
		The support function
		\begin{equation}\label{eq-closed:def_h}
			h(t,\th)=\bmetric{\Th^{-1}(t,\th)}{\nu_\th}\qquad\forall\,(t,\th)\in(t_1,t_2)\times\SS^1\big((N-2)\pi\big)
		\end{equation}
		is a smooth solution of $\p_t h=\p_{\th\th}h+h$ in $(t_1,t_2)\times\SS^1\big((N-2)\pi\big)$, and $(\p_{\th\th}h+h)(t,\cdot)$ has $N$ nondegenerate roots for each $t\in(t_1,t_2)$.
	\end{enumerate}
\end{lemma}

In the lemma, if $t_1=-\infty$ (resp. $t_2=\infty$), then our convention is $Y_{t_1}=\Omega$ (resp. $\tZ_{t_2}=\Crit(u)$). By concave polygon geometry, each preimage of $\Th$ will clearly be a line segment (once the quantity $N$ is appropriately set up). The fact that $\Th$ is homeomorphism implies that each preimage is a point. Thus $\p Z_t$ does not contain nontrivial line segments.

At a time when \eqref{eq-closed:condition_Pt} fails to hold, we show that $\p Z_t$ splits into simpler polygons whose number of edges are preserved as in \eqref{eq-closed:intro_edge_preserv}:

\begin{lemma}\label{lemma-closed:2_split_lemma}
	Suppose \eqref{eq-closed:condition_Pt} holds for all $t\in[\uT,T)$, but does not hold for $t=T$. Then:
	\begin{enumerate}[label={(\roman*)}, nosep]
		\item $Z_T\setminus Y_T$ is a finite disjoint union of line segments with endpoints on $\p Z_T$.
		\item Denote $N=N(Z_t)$ for any $t\in(\uT,T)$, which is well-defined by Lemma \ref{lemma-closed:1_evolution_lemma}(i). Let $\{E_k\}_{1\leq k\leq l}$ be the connected components of $Y_T$. Then $2\leq l<\infty$, and $N(E_k)\leq N-2$ for all $k$, and
		\begin{equation}\label{eq-closed:intro_edge_preserv}
			N-2=\sum_{k=1}^l\big(N(E_k)-2\big).
		\end{equation}
	\end{enumerate}
\end{lemma}

We start proving the main results. We need two auxiliary lemmas:

\begin{lemma}\label{lemma-closed:size_lower_bound}
	For each $K\Subset\Omega$ and $T\in\RR$, there is a constant $c=c(K,T)$ so that: if $F$ is a connected component of $Y_t$ or $Z_t$, with $t\leq T$ and $F\subset K$, then $F$ contains a disk with radius $c$.
\end{lemma}
\begin{proof}
	Let $B$ be a largest open disk contained in $F$. Then $\p B\cap\p F$ cannot be contained in an open semicircle in $\p B$. Recalling that edges of $\p F$ are streamlines of $u$ on which $|\D u|=e^{-t}$, we thus find points $z_1,z_2\in\p B\cap\p F$ with
	\[e^T\big|\D u(z_1)-\D u(z_2)\big|\geq\Big|\frac{\D u(z_1)}{|\D u(z_1)|}-\frac{\D u(z_2)}{|\D u(z_2)|}\Big|\geq\sqrt 3.\]
	If the radius of $B$ is too small, this would violate the uniform continuity of $\D u$ in $K$.
\end{proof}

\begin{lemma}\label{lemma-closed:edge_lower_bound}
	For any connected component $F$ of $Y_t,Z_t$ with $F\Subset\Omega$, we have $N(F)\geq4$ and $N(F)$ is even.
\end{lemma}
\begin{proof}
	The vector field $\D u/|\D u|$ is everywhere tangent to the edges of $\p F$. If $N(F)$ is odd, such a vector field does not exist, by Gauss-Bonnet (see Subsection \ref{subsec:curves}). Since there exist no concave 2-gons, we have $N\geq4$.
\end{proof}

\begin{lemma}\label{lemma-closed:Zt_splits_into_Yt}
	If $F$ is a connected component of $Z_t$ with $F\Subset\Omega$, then $F\setminus Y_t$ is a disjoint union of finitely many line segments with endpoints on $\p F$. Letting $\{E_k\}_{1\leq k\leq l}$ be the connected components of $F\cap Y_t$, we have $1\leq l<\infty$ and
	\begin{equation}\label{eq-closed:edges_when_splitting}
		N(F)-2=\sum_k\big(N(E_k)-2\big).
	\end{equation}
\end{lemma}
\begin{proof}
	By Theorem \ref{thm-lvlset:global}(iv), $F\setminus Y_t$ is a disjoint union of line segments whose endpoints lie on $\p F$. Let $\sigma$ be one of the line segments. Then $\sigma$ contacts $\p F$ tangentially since they are both tangential to $\D u$. Then $F\setminus\sigma$ is the disjoint union of two concave polygons $F_1,F_2$ with $N(F)=N(F_1)+N(F_2)-2$: this follows since a new vertex forms at each endpoint of $\sigma$. Applying this argument inductively, we obtain that: for any finite set of segments $\{\sigma_1,\cdots,\sigma_l\}$ in $F\setminus Y_t$, the set $F\setminus\cup\sigma_k$ consists of $l+1$ components $\{F_1,\cdots,F_{l+1}\}$ so that
	\begin{equation}\label{eq-closed:aux1}
		N(F)=\sum_kN(F_k)-2l=\sum_k\big(N(F_k)-2\big)+2.
	\end{equation}
	By Lemma \ref{lemma-closed:edge_lower_bound}, $F\setminus Y_t$ contains at most $N(F)/2$ segments; and \eqref{eq-closed:aux1} implies \eqref{eq-closed:edges_when_splitting}.
\end{proof}

\begin{lemma}\label{lemma-closed:tech_connected}
	If $Z_t\Subset\Omega$ and $Y_t$ is connected, then $Z_t=Y_t$.
\end{lemma}
\begin{proof}
	This is a direct consequence of Lemma \ref{lemma-closed:Zt_splits_into_Yt}.
\end{proof}

\begin{proof}[Proof of Lemma \ref{lemma-closed:forward_openness}] {\ }
	
	Let $c=c(\bar F,t+1)$ be the constant in Lemma \ref{lemma-closed:size_lower_bound}. Since $F$ is connected and open, we may find a connected set $G\Subset F$ with $N(G,c/3)\supset F$. Since $F=\bigcup_{s>t}(Y_s\cap F)$, there is $T\in(t,t+1)$ so that $Y_s\cap F\supset G$ for all $s\in(t,T]$. Then exactly one connected component of $Y_s$ contains $G$. If $Y_s\cap F$ has another connected component $L$, then $L\cap G=\emptyset$ but $L\subset\bar F$. This violates Lemma \ref{lemma-closed:size_lower_bound} since $N(G,c/2)\supset\bar F$. Thus, $Y_s\cap F$ is connected for all $s\in(t,T)$. Next, Lemma \ref{lemma-closed:tech_connected} (with $\text{``$\Omega$''}=F$, note that $F$ is simply-connected by Theorem \ref{thm-lvlset:global}) implies $Z_s\cap F=Y_s\cap F$, as desired.
\end{proof}

\begin{proof}[Proof of Lemma \ref{lemma-closed:1_evolution_lemma}] {\ }
	
	Denote
	\[E=Y_{t_1}\setminus\tZ_{t_2}.\]
	Clearly, $u$ has no critical points in $E$.
	
	Note that $\tZ_{t_2}=\bigcap_{t\in(t_1,t_2)}\bar{Y_t}$, and $\bar{Y_t}$ are connected for all $t\in(t_1,t_2)$ by \eqref{eq-closed:condition_Pt}. Hence $\tZ_{t_2}$ is connected. Also, $Y_{t_1}=\bigcup_{t\in(t_1,t_2)}Y_t$ is connected. Hence, $E$ is homeomorphic to an open annulus. Fix a basepoint $z_0\in\tZ_{t_2}$. Then Theorem \ref{thm-approx:global} gives a simple IMCF cluster $\big(w,\nu,\{D_i\},\{\chi_i\}\big)$ in an open set $\Omega_{\reg}=\{w<\infty\}$. Denote
	\[Y'_t=\{w>t\},\qquad\tZ'_t=\{w\geq t\},\qquad Z'_t=\Int(\tZ'_t).\]
	
	Observe that $Z_t=Z'_t$ for all $t\in(t_1,t_2)$: on one hand, $w\leq-\log|\D u|$ implies $Z_t\supset Z'_t$. On the other hand, letting $F'_t$ be the connected component of $Z'_t$ containing $z_0$, Theorem \ref{thm-approx:global}\ref{item-approx:lvlset} and the connectedness of $Z_t$ imply $Z_t=F'_t\subset Z'_t$.
	
	We claim that
	\begin{equation}\label{eq-closed:lvlset_covering}
		\bar{\bigcup_{t\in(t_1,t_2)}\p Z_t}\supset E.
	\end{equation}
	Otherwise, there is $z\in E$ and $t\in(t_1,t_2)$ and $r>0$, so that
	\begin{equation}\label{eq-closed:aux2}
		B(z,r)\subset Z_s\text{\ \ for all\ \ $s<t$,}\qquad B(z,r)\subset E\setminus Z_s\text{\ \ for all\ \ $s>t$}.
	\end{equation}
	Now $\tZ_t=\bigcap_{s<t}Y_s=\bigcap_{s<t}Z_s$, hence $B(z,r)\subset\tZ_t$, hence $B(z,r)\subset\Int(\tZ_t)=Z_t$. Also, we have $Z_t=Y_t=\bigcup_{s>t}Y_s=\bigcup_{s>t}Z_s$, hence $B(z,r)\subset E\setminus Z_t$, contradiction.
	
	Combining \eqref{eq-closed:lvlset_covering}, the fact $Z_t=Z'_t$, and $|\D u|=e^{-t}$ on each $\p Z'_t$ by Theorem \ref{thm-approx:global}\ref{item-approx:w_on_nonconst}, we obtain $w=-\log|\D u|$ in $E$. Recalling $\nu=-e^w\D^\perp u$ in Theorem \ref{thm-approx:global}\ref{item-approx:simple_cluster}, it follows that
	\begin{equation}\label{eq-closed:nu_is_nuit}
		\nu=-\frac{\D^\perp u}{|\D u|}\qquad\text{in}\ \ E.
	\end{equation}
	In particular, $|\nu|\equiv1$ in $E$, and the map \eqref{eq-closed:hodographic_map} coincides with $(-\log|\D u|,\nu^\perp)$.
	
	Let $-(N'-2)/2$ denote the degree of the map \eqref{eq-closed:hodographic_map}, i.e. the coefficient of the map $\ZZ\cong H_1(E,\ZZ)\to H_1(\SS^1,\ZZ)\cong\ZZ$. Thus, \eqref{eq-closed:hodographic_map} lifts to a degree $-1$ map 
	\[\Th:E\to(t_1,t_2)\times\SS^1\big((N'-2)\pi\big).\]
	Forgetting the time coordinate, it projects to a degree $-1$ map
	\[\varth:E\to\SS^1\big((N'-2)\pi\big).\]
	Since each $\p Z_t$ generates $H_1(E,\ZZ)$, the restriction $\varth|_{\p Z_t}$ also has degree $-1$. But Gauss-Bonnet implies that $\nu:\p Z_t\to\SS^1$ has degree $-\frac{N(Z_t)-2}2$. Comparing the degrees, we have $N'=N(Z_t)$ for all $t\in(t_1,t_2)$. This shows item (i) since $Y_t=Z_t$.
	
	It is now a trivial consequence that $\Th$ is also a lift of the function $(w,\nu^\perp)$ from the set $\bigcup_{t\in\RR}\p Z'_t\cap E$ to $(t_1,t_2)\times\SS^1\big((N'-2)\pi\big)$. Thus, the lifting condition in Theorem \ref{thm-heat:main} (with ambient domain $E$ and $M=\SS^1((N'-2)\pi)$) is satisfied.
	
	Set $N=N'$ (the notation is consistent in item (iii)). Let us show that $\Th^{-1}(t,\th)$ is a line segment for each $(t,\th)\in(t_1,t_2)\times\SS^1\big((N-2)\pi\big)$. First, notice that
	\[\Th^{-1}\big(\{t\}\times\SS^1((N-2)\pi)\big)=\big\{|\D u|=e^{-t}\big\}=\p\tZ_t\]
	since $Z_t=Y_t$ by our assumption \eqref{eq-closed:condition_Pt}. Thus
	\[\Th^{-1}(t,\th)\subset\p\tZ_t=\p Z_t\cup\big(\tZ_t\setminus\bar{Z_t}\big).\]
	Since $\tZ_t\Subset\Omega$, by Theorem \ref{thm-lvlset:global}(ii)(iii), $\tZ_t\setminus\bar{Z_t}$ is a disjoint union of nontrivial line segments $\sigma_i$ which contact $\p Z_t$. However, $\Th^{-1}(t,\th)\cap\p Z_t$ is also a line segment (denoted by $\alpha$), and thus
	\begin{equation}\label{eq-closed:preimage_of_Th}
		\Th^{-1}(t,\th)=\alpha\cup\bigcup\Big\{\sigma_i:\bar{\sigma_i}\cap\alpha\ne\emptyset\Big\},
	\end{equation}
	which is indeed a line segment (recall that $\p Z_t$ and $\sigma_i$ are all tangent to $\D u$).
	
	For each $t_3,t_4$ with $t_1<t_3<t_4<t_2$, the restricted map
	\[\Th:\tZ_{t_3}\setminus Y_{t_4}\to[t_3,t_4]\times\SS^1\big((N-2)\pi\big)\]
	is continuous, surjective, hence is a quotient map. It is shown above that each preimage is a line segment (whose direction is parallel to $\D u$). Hence, the continuous function
	\begin{equation}\label{eq-closed:def_support_1}
		\bar h(z)=\Bmetric{z}{-\frac{\D^\perp u}{|\D u|}(z)}=\metric{z}{\nu(z)}
	\end{equation}
	is constant on each preimage of $\Th$, thus descends to a continuous function
	\begin{equation}\label{eq-closed:def_support_2}
		h:[t_3,t_4]\times\SS^1\big((N-2)\pi\big)\to\RR.
	\end{equation}
	Then Theorem \ref{thm-heat:main}(i) implies that $h$ is a viscosity solution of $\p_t h=\p_{\th\th}h+h$ on $(t_3,t_4]\times\SS^1\big((N-2)\pi\big)$, hence a smooth solution as well. Letting $t_3\searrow t_1$ and $t_4\nearrow t_2$, we obtain a smooth support function $h$ that solves $\p_th=\p_{\th\th}h+h$ in $(t_1,t_2)\times\SS^1\big((N-2)\pi\big)$. So each time slice $h(t,\cdot)$ is smooth, and by counting the number of vertices, $(\p_{\th\th}h+h)(t,\cdot)$ has $N$ roots. The Sturmian theory then implies that all the roots are nondegenerate.
	
	We next show that $\Th$ is a homeomorphism. By the invariance of domain, it suffices to show that $\Th$ is injective, namely, each preimage is a point. For all $t_3\in(t_1,t_2)$ note that
	\[|\p_{\th\th}h+h|\leq C(t_3)\qquad\text{in}\ \ \big[t_3,\min\{t_3+1,t_2\}\big)\times\SS^1\big((N-2)\pi\big).\]
	This implies that all edges of $\p Z_t=\p Y_t$ have curvature at least $1/C(t_3)$. In particular, $\p Z_t$ does not contain nontrivial line segments and hence $\Th|_{\p Z_t}$ is injective.
	
	Suppose that $\tZ_t\setminus\bar{Z_t}$ contains a nontrivial line segment $\sigma$ for some $t\in(t_1,t_2]$. Fix $t_3>\max\{t_1,t-1\}$ and $z\in\Int(\sigma)$. We know that $Z_t=Z'_t$, hence $\tZ_t=\bigcap_{s<t}Z_s=\bigcap_{s<t}\tZ'_s=\tZ'_t$. Hence $\tZ_t\setminus\bar{Z_t}=\tZ'_t\setminus\bar{Z'_t}$. Then apply Lemma \ref{lemma-cluster:local_model}\ref{item-cluster:model}: after a suitable rigid motion, there is a small enough $r$ so that $\tZ_t\cap Q(r)=(-r,r)\times\{0\}$, and for $s<t$ close enough to $t$ we have
	\[Z_s\cap Q(r)=\Big\{x\in(-r,r),\,g_{1,s}(x)<y<g_{2,s}(x)\Big\}\]
	for some concave $g_{1,s}$ and convex $g_{2,s}$ with $g_{1,s}\nearrow0$, $g_{2,s}\searrow0$ as $s\nearrow t$. But this violates the uniform curvature lower bounds of $\graphh(g_{1,s})$ and $\graphh(g_{2,s})$. Hence $\tZ_t=\bar{Z_t}$.
	
	In summary, we have shown that $\tZ_t\setminus\bar{Z_t}=\emptyset$ for all $t\in(t_1,t_2]$, which implies (ii), and that $\Th$ is injective since $\{|\D u|=e^{-t}\}=\tZ_t\setminus Y_t=\p Z_t$ and $\Th|_{\p Z_t}$ is injective. Thus $\Th$ is a homeomorphism. Finally, it is clear that the support function $h$ defined in \eqref{eq-closed:def_h} coincides with the function defined in \eqref{eq-closed:def_support_1} \eqref{eq-closed:def_support_2}.
\end{proof}

\begin{proof}[Proof of Lemma \ref{lemma-closed:2_split_lemma}] {\ }
	
	Note that $Z_T$ is either disconnected or is not $Y_T$, since \hyperref[eq-closed:condition_Pt]{($\text{P}_T$)} does not hold. Denote $E=Y_\uT\setminus\tZ_T$, so $u$ has no critical points in $\bar E$. By Lemma \ref{lemma-closed:1_evolution_lemma}, the map $\phi=\D u/|\D u|:E\to\SS^1$ has degree $-(N-2)/2$, where $N=N(Y_t)=N(Z_t)$ for all $t\in(\uT,T)$. Let $\{F_j\}$ be the connected components of $Z_T$. Then $\sum\p F_j$ is homologous to $\p Y_t$ for any $t\in(\uT,T)$, and by Gauss-Bonnet, the map $\phi$ has degree $-(N(F_j)-2)/2$ on each $\p F_j$. As a result,
	\[N-2=\sum_j\big(N(F_j)-2\big).\]
	Using Lemma \ref{lemma-closed:Zt_splits_into_Yt} on each $F_j$, we obtain (i) and \eqref{eq-closed:intro_edge_preserv}. The fact $N(E_k)\leq N-2$ follows from $N(E_k)\geq4$ for all $k$. The fact $l\geq2$ follows from Lemma \ref{lemma-closed:tech_connected}: if $l=1$, then $Y_T$ is connected, hence $Z_T=Y_T$, hence \hyperref[eq-closed:condition_Pt]{($\text{P}_T$)} holds, contradiction.
\end{proof}

\begin{proof}[Proof of Theorem \ref{thm-closed:main}] {\ }
	
	For each $t\geq t_0$, let $l(t)$ be the number of connected components of $Y_t\cap F$, and denote
	\[N(t)=2+\sum\Big\{N(G)-2: G\text{ is a connected component of $Y_t\cap F$}\Big\}.\]
	Corollary \ref{cor-equi:min_principle} implies that each connected component of $Y_t$ contains a critical point of $u$, hence a connected component of $Y_s$ for all $s>t$. It follows that $l(t)$ is nondecreasing in $t$.
	
	\setcounter{claim}{0}
	
	\begin{claim}\label{claim-closed:split_part}
		For each $t\geq t_0$, there is $T>t$ so that $l(t)$ is constant in $[t,T]$.
	\end{claim}
	\begin{proof}
		Fix $t\geq t_0$. Let $\{F_k\}$ be the connected components of $Y_t\cap F$. By Lemma \ref{lemma-closed:forward_openness} with $\text{``$F$''}=F_k$, there is $T>t$ so that $Y_s\cap F_k$ is connected for all $k$ and all $s\in(t,T]$. So $Y_s\cap F$ has the same number of connected components as $Y_t\cap F$.
	\end{proof}
	
	\begin{claim}\label{claim-closed:continuous_part}
		If $l(t)$ is constant in $[t_1,t_2)$ and $l(t_2^-)<l(t_2)$, then $N(t)$ is constant in $[t_1,t_2]$.
	\end{claim}
	\begin{proof}
		Let $\{F_k\}$ be the set of connected components of $Y_{t_1}\cap F$. Then $Y_t\cap F_k$ is connected for each $t\in[t_1,t_2)$, by the constancy of $l(t)$. Then Lemma \ref{lemma-closed:tech_connected} with $\text{``$\Omega$''}=F_k$ implies $Y_t\cap F_k=Z_t\cap F_k$. Thus, the condition \eqref{eq-closed:condition_Pt} is satisfied for all $t\in(t_1,t_2)$ when restricted to $F_k$. Then Lemma \ref{lemma-closed:1_evolution_lemma}(i) applied inside each $F_k$ implies that $N(Y_t\cap F_k)$ is constant. Hence,
		\[N(t)=2+\sum_k\big(N(Y_t\cap F_k)-2\big)\]
		stays constant in $(t_1,t_2)$ as well.
		
		For each $k$: if $Y_{t_2}\cap F_k$ is connected, then we use Lemma \ref{lemma-closed:forward_openness}, \ref{lemma-closed:tech_connected} and \ref{lemma-closed:1_evolution_lemma}(i) to conclude that $Y_t\cap F_k$ is connected and $N(Y_t\cap F_k)=\text{const}$ for $t\in(t_1,t_2+\epsilon)$, for some $\epsilon>0$. If $Y_{t_2}\cap F_k$ is not connected, then \eqref{eq-closed:condition_Pt} does not hold at $t_2$ when restricted to $F_k$, hence Lemma \ref{lemma-closed:2_split_lemma}(ii) applies. Considering both cases, $N(t)$ is constant in $(t_1,t_2]$.
		
		Finally, notice that $F_k=\bigcup_{t>t_1}(Y_t\cap F_k)$. Denoting
		\[b_k=-\frac{N(Y_t\cap F_k)-2}2\qquad\text{for any}\ \ t\in(t_1,t_2),\]
		by Lemma \ref{lemma-closed:1_evolution_lemma}, the map $\big(\!-\log|\D u|,\frac{\D u}{|\D u|}\big)$ is a degree $b_k$ map from $F_k\setminus\tZ_{t_2}$ to $(t_1,t_2)\times\SS^1$. This map clearly extends to $\bar{F_k}\setminus\tZ_{t_2}$ and has degree $b_k$ as well. Hence, the tangent vector map $\frac{\D u}{|\D u|}:\p F_k\to\SS^1$ has degree $b_k$, which implies $N(F_k)=-2b_k+2=N(Y_t\cap F_k)$ for any $t\in(t_1,t_2)$. Summing this over $k$, it follows that
		\[N(t_1)=2+\sum_k\big(N(F_k)-2\big)=2+\sum_k\big(N(Y_t\cap F_k)-2\big)=N(t)\]
		for any $t\in(t_1,t_2)$. Hence $N(t)$ is constant in $[t_1,t_2]$.
	\end{proof}
	
	Therefore, there is a sequence of times $t_0<t_1<t_2<\cdots$ so that
	\[l(t_0)=l(t_1^-)<l(t_1)=l(t_2^-)<l(t_2)=\cdots.\]
	Claim \ref{claim-closed:continuous_part} implies that $N(t)$ is constant for $t\in[t_0,\lim_{i\to\infty}t_i)$. However, since each connected component of $Y_t\cap F$ has at least 4 edges, we have
	\[N(t_0)=N(t)\geq 2l(t)\qquad\forall\,t\in[t_0,\lim_{i\to\infty}t_i).\]
	Hence $l(t^-)<l(t)$ can only occur for at most $N(t_0)/2$ times. As a result, there are finitely many splitting times, and there is $T\in\RR$ so that $l(t)$ is constant in $[T,\infty)$.
	
	Let $\{F_k\}$ be the connected components of $Y_T\cap F$. Hence, the condition \eqref{eq-closed:condition_Pt} holds for all $t>T$ when restricted to each $F_k$ (see the proof of Claim \ref{claim-closed:continuous_part}). Applying Lemma \ref{lemma-closed:1_evolution_lemma}, the support function $h(t,\th)$ defined in \eqref{eq-closed:def_h} satisfies $\p_t h=\p_{\th\th}h+h$ in $(T,\infty)\times\SS^1\big((N_k-2)\pi\big)$, where $N_k=N(Y_t\cap F_k)$ for all $t>T$. Denote $d_k=N_k/2$ (so $d_k\geq2$). Then we may write
	\[h(t,\th)=\sum_{n\geq0}\exp\Big(t-\frac{n^2}{(d_k-1)^2}t\Big)\Big(a_n\sin\frac{n\th}{d_k-1}+b_n\cos\frac{n\th}{d_k-1}\Big).\]
	Lemma \ref{lemma-closed:1_evolution_lemma}(iii) implies that each $(\p_{\th\th}h+h)(t,\cdot)$ has $N_k=2d_k$ roots for all $t>T$. Hence, all the frequencies with $n<d_k-1$ must vanish, and the $n=d_k$ frequency does not vanish. Hence for some $a,b,c,d\in\RR$ with $(a,b)\ne(0,0)$, we have
	\begin{equation}\label{eq-closed:asymp_of_h}
		\begin{aligned}
			& \exp\Big(\frac{2d_k-1}{(d_k-1)^2}t\Big)\cdot\Big[h(t,\th)-c\sin\th+d\cos\th\Big] \\
			&\hspace{144pt}\xrightarrow[t\to\infty]{\ C^\infty(\SS^1((N_k-2)\pi))\ }a\sin\frac{d_k\th}{d_k-1}+b\cos\frac{d_k\th}{d_k-1}.
		\end{aligned}
	\end{equation}
	Recall that $\p Y_t\cap F_k$ can be described as
	\begin{equation}\label{eq-closed:recovery_formula}
		\p Y_t\cap F_k=\Big\{h(t,\th)\nu_\th+\p_\th h(t,\th)\tau_\th:\th\in\SS^1\big((N_k-2)\pi\big)\Big\},
	\end{equation}
	and in view that $\D u=e^{-t}\tau_\th$ on $\p Y_t\cap F_k$, we have
	\begin{equation}\label{eq-closed:recovery_formula_2}
		\D u\Big(h(t,\th)\nu_\th+\p_\th h(t,\th)\tau_\th\Big)=e^{-t}\tau_\th,\qquad\forall\,(t,\th)\in(T,\infty)\times\SS^1\big((N_k-2)\pi\big).
	\end{equation} 
	Denote $z_k=(c,d)$. Then $h(t,\th)\to c\sin\th-d\cos\th$ in $C^1$ implies $\p Y_t\cap F_k\to\{z_k\}$ in the Hausdorff topology, as $t\to\infty$. Therefore $\Crit(u)\cap F_k=\{z_k\}$. This proves item \ref{item-closed:main_finiteness}.
	
	For \ref{item-closed:main_edges}, the constancy of $N(t)$ implies
	\[N(F)-2=N(t_0)-2=N(T+1)-2=\sum_k(N_k-2)=\sum_k(2d_k-2).\]
	
	We next prove \ref{item-closed:main_asymp}. Fix a $k$. For a scaling factor $\lambda>0$, denote
	\[u_\lambda(z)=\lambda^{-\frac{d_k^2}{2d_k-1}}\Big[u\big(\lambda z+z_k\big)-u(z_k)\Big].\]
	Then
	\[Y_{t,\lambda}:=\big\{|\D u_\lambda|<e^{-t}\big\}=\lambda^{-1}\big(Y_{t'}-z_k\big)\qquad\text{where}\ \ t'=t-\frac{(d_k-1)^2}{2d_k-1}\log\lambda.\]
	Therefore, $\p Y_{t,\lambda}\cap\lambda^{-1}(F_k-z_k)$ has support function
	\[\begin{aligned}
		h_\lambda(t,\th) &= \lambda^{-1}\Big[ h(t',\th)-c\sin\th+d\cos\th\Big] \\
		&= \exp\Big(\!-\frac{2d_k-1}{(d_k-1)^2}t\Big)\exp\Big(\frac{2d_k-1}{(d_k-1)^2}t'\Big)\Big[h(t',\th)-c\sin\th+d\cos\th\Big].
	\end{aligned}\]
	By \eqref{eq-closed:asymp_of_h}, we have
	\begin{equation}\label{eq-closed:asymp_of_h_scaled}
		\lim_{\lambda\to0}h_\lambda(t,\th)=\exp\Big(\!-\frac{2d_k-1}{(d_k-1)^2}t\Big)\Big[a\sin\frac{d_k\th}{d_k-1}+b\cos\frac{d_k\th}{d_k-1}\Big].
	\end{equation}
	Call the right hand side $\bar h$. Let $\bar u$ be the degree $d_k$ quasiradial solution so that $Y_{t,\infty}:=\{|\D\bar u|<e^{-t}\}$ has support function given by $\bar h$ (see Example \ref{ex-ex:quasiradial}). Then \eqref{eq-closed:recovery_formula} implies
	\begin{equation}\label{eq-closed:aux4}
		\p Y_{t,\lambda}\cap\lambda^{-1}(F_k-z_k)\to\p Y_{t,\infty}
	\end{equation}
	in the Hausdorff topology as $\lambda\to0$. By the maximum principle for the gradient, $\{|\D u_\lambda|\}$ is equibounded in $B(R)$ for all $R>0$. By Savin's equicontinuity theorem and Arzel\`a-Ascoli and $u_\lambda(0)=0$, for each sequence $\lambda_n\to0$ we may extract a further subsequential limit $u_\infty$ of $u_{\lambda_n}$ in $C^1_{\loc}(\RR^2)$. Rescaling \eqref{eq-closed:recovery_formula_2} and taking $\lambda\to0$ gives
	\[\D u_\infty\Big(\bar h(t,\th)\nu_\th+\p_\th \bar h(t,\th)\tau_\th\Big)=e^{-t}\tau_\th,\qquad\forall\,(t,\th)\in\RR\times\SS^1\big((N_k-2)\pi\big).\]
	But this shows $\D u_\infty=\D\bar u$, hence $u_\infty=\bar u$ (see Example \ref{ex-ex:quasiradial}). Since the limit object is independent of the sequence $\lambda_n$, it follows that
	\[\lim_{\lambda\to0}u_\lambda=\bar u\qquad\text{in}\ \ C^1_{\loc}(\RR^2),\]
	which shows \ref{item-closed:main_asymp}.
	
	It remains to prove \ref{item-closed:reg}. We may assume $z_k=0$ and $u(0)=0$. By item \ref{item-closed:main_finiteness} and scaling, we may assume $\Crit(u)\cap B(2)=\{0\}$. It remains to show that $u\in C^{1,1/3}(\bar B(1))$. Suppose this does not hold. Then there is a sequence of points $z_i,\zeta_i\in\bar{B(1)}$, so that
	\begin{equation}\label{eq-closed:aux6}
		\frac{|\D u(z_i)-\D u(\zeta_i)|}{|z_i-\zeta_i|^{1/3}}\to\infty.
	\end{equation}
	Then Theorem \ref{thm-ereg:away_from_crit} implies $z_i,\zeta_i\to0$. By \ref{item-closed:main_asymp}, there is $d\in\ZZ_{\geq2}$ and a quasiradial solution $\bar u$ in $\RR^2$, so that
	\begin{equation}\label{eq-closed:limit}
		\lim_{\lambda\to0}\lambda^{-\frac{d^2}{2d-1}}u(\lambda z)=\bar u(z)\qquad\text{in}\ \ C^1_{\loc}(\RR^2).
	\end{equation}
	It follows that
	\[|\D u(z)|\leq C|z|^{1/3}\qquad\forall\,z\in\bar{B(1)}.\]
	Hence, we must have
	\begin{equation}\label{eq-closed:aux5}
		\frac{|z_i-\zeta_i|}{|z_i|+|\zeta_i|}\to0.
	\end{equation}
	Denote $r_i=|z_i|$. Passing to a subsequence, we may assume $r_i^{-1}z_i\to q$. Then \eqref{eq-closed:aux5} implies $r_i^{-1}\zeta_i\to q$ as well, and $r_i^{-1}|z_i-\zeta_i|\to0$. Let $C$ be the constant in Theorem \ref{thm-ereg:es_combined} (with $\delta=1$). Since $\bar u\in C^1$, we may find $r\ll1$ and $A\in\text{SO}(2)$ and $b\in\RR$, $c\ne0$, so that
	\[\sup_{Q(4)}\Big|cr^{-1}\big(\bar u(Arz+q)-b\big)-x\Big|\leq\frac1{4C},\]
	where $x$ is the first coordinate of $z$. Then \eqref{eq-closed:limit} implies
	\[\sup_{Q(4)}\Big|cr^{-1}\Big(\lambda^{-\frac{d^2}{2d-1}}u(\lambda Arz+\lambda q)-b\Big)-x\Big|\leq \frac1{2C},\qquad\forall\,\lambda\ll1.\]
	Theorem \ref{thm-ereg:es_combined} with $\delta=1$ implies
	\[\Big\|c\lambda^{-\frac{(d-1)^2}{2d-1}}A^T\D u(\lambda Arz+\lambda q)-\p_x\Big\|_{C^{0,1/3}(Q(1))}\leq C,\qquad\forall\,\lambda\ll1.\]
	Taking $\lambda=r_i$ and setting $\tilde z_i=r^{-1}A^{-1}(r_i^{-1}z_i-q)$ and $\tilde\zeta_i=r^{-1}A^{-1}(r_i^{-1}\zeta_i-q)$, this gives
	\[\begin{aligned}
		|c|r_i^{-\frac{(d-1)^2}{2d-1}}\big|\D u(z_i)-\D u(\zeta_i)\big| &\leq C\big|\tilde z_i-\tilde\zeta_i\big|^{1/3} \\
		&\leq Cr^{-1/3}r_i^{-1/3}|z_i-\zeta_i|^{1/3},\qquad\forall\,i\gg1.
	\end{aligned}\]
	Since $d\geq2$, we obtain
	\[\big|\D u(z_i)-\D u(\zeta_i)\big|\leq C(r,c,C)|z_i-\zeta_i|^{1/3},\]
	contradicting \eqref{eq-closed:aux6}.
\end{proof}

\begin{proof}[Proof of Theorem \ref{thm-closed:precompact_crit}] {\ }
	
	Set $t=-\log\min_{\p K}|\D u|$. Then $Y_{t+1}\cap K\Subset K$. Lemma \ref{lemma-closed:size_lower_bound} implies that $Y_{t+1}\cap K$ has finitely many connected components, and then Theorem \ref{thm-closed:main}\ref{item-closed:main_finiteness} (applied to each connected component) implies that $\Crit(u)\cap K$ is finite.
\end{proof}

%\newpage

\section{Discreteness of critical set}\label{sec:crit}

In this section we prove that:

\begin{theorem}\label{thm-crit:main}
	Suppose $u$ is nonconstant and $\infty$-harmonic in a domain $\Omega\subset\RR^2$. Then $\Crit(u)$ is discrete in $\Omega$.
\end{theorem}

Theorem \ref{thm-crit:main} is the final ingredient for showing $C^{1,1/3}$ regularity:

\begin{proof}[Proof of Theorem \ref{thm-intro:reg}] {\ }
	
	Suppose $\Delta_\infty u=0$ in a domain $\Omega$ and $u$ is nonconstant. It is known by Theorem \ref{thm-crit:main} that $\Crit(u)$ is discrete. The result follows from Theorem \ref{thm-closed:main}\ref{item-closed:reg} and Theorem \ref{thm-ereg:away_from_crit}.
\end{proof}

\begin{proof}[Proof of Theorem \ref{thm-intro:iso_crit}] {\ }
	
	We first use Theorem \ref{thm-crit:main} to deduce that $\Crit(u)$ is discrete. Then, we apply Theorem \ref{thm-closed:main}\ref{item-closed:main_asymp} in a small disk centered at each critical point.
\end{proof}

A key idea in proving Theorem \ref{thm-crit:main} is a \Phragmenlindelof type principle, which we use to bound the number of edges of the level sets of $|\D u|$. We also mention the work of Granlund-Marola \cite{Granlund-Marola_2015}, where a different \Phragmenlindelof principle is obtained.

\begin{lemma}\label{lemma-crit:phragmen-lindelof2}
	Suppose $0<r<R$, and $D\subset B(R)$ is open with $C^1$ boundary, and $D\cap B(r)\ne\emptyset$. Let $w\in C^0(\bar D\cap B(R))\cap\Lip_{\loc}(D)$ be a continuously calibrated weak IMCF in $D$ with outer obstacle $\p D\cap B(R)$, so that $\inf_D(w)>-\infty$. Then
	\[\inf_{D\cap B(r)}(w)-\inf_D(w)\geq\log\Big(\frac{\sqrt\pi}2\frac{R-r}{|D\setminus B(r)|^{1/2}}\Big).\]
\end{lemma}

\begin{proof}
	
	Denote $t_0=\inf_D(w)$. For any $T>\inf_{D\cap B(r)}(w)$, we set $E_T=\{w<T\}$, so $E_T$ is a $C^1$ open set contained in $D$ with $E_T\cap B(r)\ne\emptyset$. For all $\rho\in(r,R)$, set
	\[A(\rho)=\big|\p E_T\cap B(\rho)\big|,\qquad S(\rho)=\big|\p B(\rho)\cap E_T\big|,\qquad V(\rho)=\big|E_T\cap B(\rho)\big|.\]
	Using \eqref{eq-prelim:excess_ineq} with $F=B(\rho)$ for almost every $\rho\in(r,R)$, we obtain
	\[A(\rho)\leq e^{T-t_0}S(\rho).\]
	The isoperimetric inequality and coarea formula gives
	\[A(\rho)+S(\rho)\geq\sqrt{4\pi V(\rho)},\qquad \frac{dV}{d\rho}=S(\rho),\qquad\forall\,\text{a.e.}\,\rho\in(r,R).\]
	Putting these together, it follows that
	\[\frac{dV}{d\rho}\geq\frac{\sqrt{4\pi V}}{1+e^{T-t_0}}\qquad\Rightarrow\qquad\frac{dV^{1/2}}{d\rho}\geq\frac{\sqrt\pi}{1+e^{T-t_0}}\geq\frac{\sqrt\pi}{2e^{T-t_0}}.\]
	Since $V(r)>0$, we may integrate this to obtain
	\[\frac{\sqrt\pi}{2e^{T-t_0}}(R-r)\leq V(R)^{1/2}-V(r)^{1/2}\leq\big[V(R)-V(r)\big]^{1/2}\leq|D\setminus B(r)|^{1/2}.\]
	The lemma follows by sending $T\searrow\inf_{D\cap B(r)}(w)$.
\end{proof}

When $|D|$ is small, we may upgrade the $\log$ growth rate to a linear rate:

\begin{lemma}\label{lemma-crit:phragmen_lindelof_3}
	Suppose $D\subset B(1)$ is open with $C^1$ boundary, and $D\cap B(1/2)\ne\emptyset$, and $w$ is a continuously calibrated weak IMCF in $D$ with outer obstacle $\p D\cap B(1)$, so that $\inf_D(w)>-\infty$. If $|D|\leq\frac1{32}$, then
	\[\inf_{D\cap B(1/2)}(w)-\inf_D(w)\geq\frac{1}{64|D|}.\]
\end{lemma}
\begin{proof}
	Set $N=\lfloor\frac1{16|D|}\rfloor$, so $N\geq2$ and $\frac1{32|D|}\leq N\leq\frac1{16|D|}$. For $k=0,1\cdots N$, consider
	\[\alpha_k=\inf_{D\cap B(1-k/2N)}(w),\qquad V_k=\big|D\cap B(1-k/2N)\setminus B(1-(k+1)/2N)\big|.\]
	So $\sum V_k\leq|D|$. Lemma \ref{lemma-crit:phragmen-lindelof2} gives
	\[\alpha_{k+1}-\alpha_k\geq\log\Big(\frac{\sqrt\pi}{4NV_k^{1/2}}\Big),\qquad\forall\,k=0,1,\cdots,N-1.\]
	Jensen's inequality then implies
	\[\alpha_N-\alpha_0\geq\sum_{k=0}^{N-1}\log\Big(\frac{\sqrt\pi}{4NV_k^{1/2}}\Big)\geq N\log\Big(\frac{\sqrt\pi}{4N^{1/2}|D|^{1/2}}\Big)\geq\frac{\log\sqrt\pi}{32|D|}\geq\frac1{64|D|}.\qedhere\]
\end{proof}

\begin{proof}[Proof of Theorem \ref{thm-crit:main}] {\ }
	
	Suppose otherwise that $\Crit(u)$ is not discrete and $\Crit(u)\ne\Omega$. By choosing a non-isolated point $z\in\p\Crit(u)\cap\Omega$ and restricting to a small disk centered at $z$, it suffices to rule out the case that $\Delta_\infty u=0$ in $B(1)$, and $\sup_{B(1)}|\D u|<1$, and
	\begin{equation}\label{eq-crit:to_contradict}
		0\in\p\Crit(u),\qquad\text{$0$ is a non-isolated critical point}.
	\end{equation}
	Now assume that there exists $u$ as in \eqref{eq-crit:to_contradict}. Since $0\in\p\Crit(u)$, there is a connected component $V$ of $B(1)\setminus\Crit(u)$ that intersects $B(1/8)$. By the maximum principle for $|\D u|$, we have $V\not\Subset B(1)$. So there is a path $\Gamma\subset V$ joining some $z_1\in V\cap B(1/8)$ and $z_2\in V\cap\p B(7/8)$. Then, we fix $T\geq1$ so that
	\[\min_\Gamma|\D u|>e^{-T}.\]
	
	For the basepoint $z_0=0$, Theorem \ref{thm-approx:global} gives a continuous function $w:B(1)\to\RR\cup\{+\infty\}$, a regular domain $\Omega_{\reg}=\{w<\infty\}$, and a simple IMCF cluster $\big(w,\nu,\{D_i\},\{\chi_i\}\big)$ in $\Omega_{\reg}$. Note that $w(0)=-\log|\D u|(0)=+\infty$. Denote
	\[Y_t=\{w>t\}.\]
	Let $F_t$ resp. $G_t$ denote the connected component of $\{|\D u|<e^{-t}\}$ resp. $Y_t$ that contains $0$. Theorem \ref{thm-approx:global}\ref{item-approx:lvlset} and \ref{thm-lvlset:global}(i) imply that $F_t=G_t$ is a simply-connected concave polygon. Notice that $F_t\not\Subset B(1)$, otherwise $0$ would be an isolated critical point of $u$ by Theorem \ref{thm-closed:main}, contradicting \eqref{eq-crit:to_contradict}.
	
	For each $t\geq T$, let $\cD'$ be the collection of connected components of $D\cap\Int\big(\{w\leq t\}\big)$ with $D$ ranging in $\{D_i\}$ (note: some $D_i$ may be separated into more than one region). We make implicit the dependence on $t$ for notational simplicity. Each element in $\cD'$ is a $C^1$ domain in $B(1)$ by the outer obstacle condition, and each edge of $\p F_t$ lies on the boundary of exactly one $D\in\cD'$. Note that $\bar D\cap B(1)\subset\Omega_{\reg}$, so $w\in C^0\big(\bar D\cap B(1)\big)$, for all $D\in\cD'$.
	
	\setcounter{claim}{0}
	
	\begin{claim}\label{claim-crit:distinct_regions}
		Let $\{\gamma_j\}$ be the edges of $\p F_t$, and $D_j\in\cD'$ be the region adjacent to $\gamma_j$. Then $\{D_j\}$ are all distinct.
	\end{claim}
	\begin{proof}
		If $\gamma_i,\gamma_j$ lie on different connected components of $\p F_t\cap B(1)$, then $D_i,D_j$ lie in different connected components of $B(1)\setminus\bar{F_t}$, as desired. If $\gamma_i$ and $\gamma_j$ are adjacent, i.e. they have a common endpoint $z$, then $z\in\p D_i\cap\p D_j$, hence $D_i,D_j$ have the opposite orientations (as regions in a simple cluster), hence $D_i\ne D_j$.
		
		Assume that $\gamma_i,\gamma_j$ lie on the same connected component of $\p F_t$ but are not adjacent. Take another edge $\gamma_k$ that lies between $\gamma_i$ and $\gamma_j$ on $\p F_t$. We may assume $D_k\ne D_i,D_j$ (otherwise, we may replace $\gamma_i$ or $\gamma_j$ by $\gamma_k$). If $D_i=D_j$, then we may fix two points $z_i\in\gamma_i$ and $z_j\in\gamma_j$, then connect them by a smooth path in $D_i=D_j$ and another smooth path in $F_t$. These two paths form a loop that encloses $\gamma_k$, hence enclosing $D_k$ as well, showing that $D_k\Subset B(1)$. Then, $\div(e^{-w}\nu)=0$ in $D_k$ and the outer obstacle condition yield a contradiction
		\[0=\int_{D_k}\div(e^{-w}\chi_k\nu)=\int_{\p D_k}e^{-w}\metric{\chi_k\nu}{\nu_{D_k}}=\int_{\p D_k}e^{-w}>0. \qedhere\]
	\end{proof}
	
	For each $t\geq T$, there are at most $64\pi t$ elements in $\cD'$ that intersect $B(1/2)$: indeed, otherwise, some $D$ among them would satisfy $|D|<\frac1{64t}$, and Lemma \ref{lemma-crit:phragmen_lindelof_3} implies
	\[\inf_{B(1)}(w)\leq\inf_{D\cap B(1)}(w)\leq\inf_{D\cap B(1/2)}(w)-\frac1{64|D|}\leq t-\frac1{64|D|}<0.\]
	Hence $\p\{w>0\}\cap B(1)\ne\emptyset$. But Theorem \ref{thm-approx:global}\ref{item-approx:w_on_nonconst} implies $|\D u|=e^{-w}=1$ on $\p\{w>0\}\cap B(1)$, which contradicts our hypothesis $\sup_{B(1)}|\D u|<1$, as desired.
	
	Recall $0\in F_t$ and $\bar{F_t}\cap\Gamma\ne\emptyset$ and $F_t\not\Subset B(1)$ and $F_t$ is simply-connected, for all $t\geq T$. These imply that some connected component of $\p F_t\cap B(1/2)$, called $\gamma_t$, must have length at least $1/2$. By Claim \ref{claim-crit:distinct_regions} and the above paragraph, $\gamma_t$ has $\leq64\pi t$ edges, hence $\leq64\pi t$ vertices. Let $E$ be the connected component of $B(1/2)\setminus\gamma_t$ that contains $0$, so $E$ is a concave polygon with boundary $\gamma_t$ in $B(1/2)$. Applying Gauss-Bonnet on $E$, we have
	\[\begin{aligned}
		2\pi &= -\int_{\gamma_t}|\kappa|+\pi\cdot\text{(\# vertices in $\gamma_t$)}+2|\p E\cap\p B(1/2)| \\
		&\hspace{120pt} +\sum\text{(exterior angles at $\bar{\gamma_t}\cap\p B(1/2)$)},
	\end{aligned}\]
	which implies
	\begin{equation}\label{eq-crit:gauss_bonnet}
		\int_{\gamma_t}|\kappa|\leq64\pi^2t+2\pi.
	\end{equation}
	By Cauchy-Schwarz, we have
	\begin{equation}\label{eq-crit:cs}
		\int_{\gamma_t}\frac1{|\kappa|}\geq\frac{|\gamma_t|^2}{\int_{\gamma_t}|\kappa|}\geq\frac{1}{4(64\pi^2t+2\pi)},\qquad\forall\,t\geq T.
	\end{equation}
	Notice that $w\in\BV_{\loc}(\Omega_{\reg})$: this follows from Lemma \ref{lemma-cluster:unique_ridge_in_square} (local finiteness of ridges) and \cite[Lemma 2.24]{Xu_2024_obstacle}. For almost every $t$, we have $\H^1\big(\p Y_t\cap(\cup\gamma_j)\big)=0$ and the curvature of $\p Y_t$ almost everywhere equals $|\D w|$. The coarea formula and \eqref{eq-crit:cs} then gives
	\begin{equation}\label{eq-crit:coarea}
		\frac\pi4\geq\big|\Omega_{\reg}\cap B(1/2)\big|\geq\int_{-\infty}^\infty\int_{\p Y_t\cap B(1/2)}\frac{1}{|\D w|}\,dt\geq\int_T^\infty\int_{\gamma_t}\frac1{|\kappa|}\,dt=\infty,
	\end{equation}
	which is a contradiction.
\end{proof}

%\newpage

\section{Entire solutions with polynomial growth}\label{sec:entire}

The goal of this section is to prove the following:

\begin{theorem}\label{thm-entire:poly_growth}
	Suppose that $u$ is $\infty$-harmonic in $\RR^2$, such that
	\begin{equation}\label{eq-entire:main_growth_cond}
		|u|\leq C\big(1+|z|^n\big)
	\end{equation}
	for some $C,n>0$. Suppose additionally that $u$ is not a linear function. Then
	\begin{enumerate}[label={(\roman*)}, nosep]
		\item\label{item-entire:blow_down} There exists $d\in\ZZ_{\geq2}$ so that the functions
		\[u_\lambda(z)=\lambda^{-\frac{d^2}{2d-1}}u(\lambda z)\]
		converge in $C^1_{\loc}(\RR^2)$ to a degree $d$ quasiradial solution as $\lambda\to\infty$.
		\item\label{item-entire:degree} Let $\{z_1,\cdots,z_l\}$ be the critical points of $u$, and $d_k$ be the degree of $z_k$. Then
		\[\sum_{k=1}^l(d_k-1)=d-1.\]
	\end{enumerate}
\end{theorem}

Here, the degree of $z_k$ is the number $d_k$ appearing in Theorem \ref{thm-closed:main}\ref{item-closed:main_asymp}. As a result, for an entire solution with polynomial growth, the growth at $\infty$ must fall into the discrete spectrum
\[|z|,\quad|z|^{4/3},\quad|z|^{9/5},\quad|z|^{16/7},\quad\cdots.\]

We further have a rigidity result:

\begin{theorem}\label{thm-entire:deg_2}
	Suppose that $u$ is a nonlinear $\infty$-harmonic function in $\RR^2$, such that
	\[|u|\leq C\big(1+|z|^{9/5-\epsilon}\big)\]
	for some $C,\epsilon>0$. Then $u=|x|^{4/3}-|y|^{4/3}$ up to scaling and rigid motions.
\end{theorem}

% == Optional ==
%In addition, one may show that entire solutions with growth $|z|^{9/5}$ at infinity must fall into the functions constructed in Subsection \ref{subsec:entire3}. Indeed, from the proof of Theorem \ref{thm-entire:poly_growth} below and the growth condition, there would exist a trigonometric polynomial of the form
%\[h(t,\th)=\sum_{j=1}^3\exp\Big(t-\frac{j^2t}{(d-1)^2}\Big)\Big(a_j\cos\frac{j\th}3+b_j\sin\frac{j\th}3\Big),\qquad\th\in\SS^1(4\pi),\]
%with $(a_3,b_3)\ne(0,0)$, so that $\big\{|\D u|<e^{-t}\big\}$ is given by support function $h(t,\cdot)$ for all $t\ll-1$.
%
%\begin{theorem}\label{thm-entire:deg_3}
%	If $u$ is $\infty$-harmonic in $\RR^2$, and there is $C>1$ so that
%	\[C^{-1}|z|^{9/5}\leq|u|\leq C|z|^{9/5}\qquad\forall\,|z|\geq1,\]
%	then $u$ is the function given in Subsection \ref{subsec:entire3} for some choice of parameters.
%\end{theorem}

Still, we denote
\[Y_t=\big\{|\D u|<e^{-t}\big\},\qquad \tZ_t=\big\{|\D u|\leq e^{-t}\big\},\qquad Z_t=\Int(\tZ_t).\]
The key step is to show that $Y_t$ is connected, bounded, and has bounded number of edges for all $t\leq T$, for some $T\in\RR$. Then Lemma \ref{lemma-closed:1_evolution_lemma} gives a support function $h$ that is an ancient solution of $\p_t h=\p_{\th\th}h+h$, and $(\p_{\th\th}h+h)(t,\cdot)$ has bounded number of roots. Then, Theorem \ref{thm-entire:ancient_sol} below implies that $h$ has finitely many Fourier frequencies. This eventually leads to the identity
\[h(t,\th)=\sum_{j=1}^d\exp\Big(t-\frac{j^2t}{(d-1)^2}\Big)\Big(a_j\cos\frac{j\th}{d-1}+b_j\sin\frac{j\th}{d-1}\Big),\]
for $\th\in\SS^1\big((2d-2)\pi\big)$ and $t\ll-1$, for some coefficients $a_j,b_j$. By taking $t\to-\infty$ of this identity, we obtain a degree $d$ quasiradial solution as the blow-down limit.

Together with Section \ref{sec:closed}, we may fully describe the behavior of entire solutions with polynomial growth: the curves $\p Y_t$ are asymptotic to those of a quasiradial solution at $t\to-\infty$, then follow the IMCF and splitting pattern as in Section \ref{sec:closed}. It would be interesting to see whether all Fourier states at infinity are realizable when $d\geq4$. For the $d=3$ case, this is answered in Subsection \ref{subsec:entire3}.

\begin{proof}[Proof of Theorem \ref{thm-entire:poly_growth}] {\ }
	
	We have assumed that $u$ is not affine or constant. By changing $u(z)\mapsto a+bu(cz+d)$, we may assume
	\begin{equation}\label{eq-entire:normalization}
		u(0)=0,\qquad|\D u(0)|=1.
	\end{equation}
	For each $t<0$, let $F_t$ be the connected component of $Y_t$ that contains $0$. Notice that $F_t\ne\RR^2$ for all $t$ (otherwise, $u$ is Lipschitz in $\RR^2$ hence linear, by \cite[Theorem 4]{Savin_2005}).
	
	Fix a sequence $p\to\infty$ and $p$-harmonic functions $u_p$, so that $u_p\to u$ in $C^0_{\loc}(\RR^2)$. Apply Theorem \ref{thm-approx:global} with basepoint $z_0=0$. We obtain a simple cluster $\big(w,\nu,\cD_{\reg},\{\chi_D\}\big)$ in an open set $\Omega_{\reg}\subset\RR^2$. Recall that $\RR^2\setminus\Omega_{\reg}\subset\Crit(u)$, and $w$ is continuous from $\RR^2$ to $\RR\cup\{+\infty\}$ with $\Omega_{\reg}=\{w<\infty\}$. Note that $w(0)=-\log|\D u(0)|=0$. By \eqref{eq-approx:main_w_on_path_1} we have
	\begin{equation}\label{eq-entire:no_poisson_pole}
		\inf_{B(R)}(w)\geq-\log\sup_{B(R)}|\D u|>-\infty,\qquad\forall\,R>0.
	\end{equation}
	
	Theorem \ref{thm-approx:global}\ref{item-approx:lvlset} implies that $F_t$ coincides with the connected component of $\{w>t\}$ containing $0$, and Theorem \ref{thm-lvlset:global} implies that $F_t$ is a simply-connected concave polygon.
	
	For each $t<0$, let $\cD'$ be the collection of connected components of $D\cap\Int\big(\{w\leq t\}\big)$ with $D$ ranging in $\cD_{\reg}$. Each element of $\cD'$ is a $C^1$ domain in $\RR^2$, and each edge of $\p F_t$ lies on the boundary of exactly one $D\in\cD'$. Following the proof of Claim \ref{claim-crit:distinct_regions} in Theorem \ref{thm-crit:main}, distinct edges of $\p F_t$ lie in the boundaries of distinct elements of $\cD'$.
	
	\setcounter{claim}{0}
	
	\begin{claim}\label{claim-entire:bound_on_domains}
		Set $N(n)=16^{n+1}$. Then $\#\cD'\leq N(n)$ for all $t<0$.
	\end{claim}
	
	It follows that
	
	\begin{claim}\label{claim-entire:bounded_n_of_edge}
	$\p F_t$ has at most $N(n)$ edges for all $t<0$.
	\end{claim}
	
	\begin{proof}[Proof of Claim \ref{claim-entire:bound_on_domains}]
		Suppose that $\cD'$ contains $M>N(n)$ regions $\{D_1,\cdots,D_M\}$. Fix $R\gg1$ so that $B(R)\cap D_i\ne\emptyset$ for all $i\in\{1,\cdots,M\}$. It follows that
		\begin{equation}\label{eq-entire:aux2}
			\inf_{D_i\cap B(R)}(w)\leq t<0\,\qquad\forall\,i.
		\end{equation}
		For each $i\in\{1,\cdots,M\}$ and $j\geq1$, denote
		\[\kappa_{ij}=\frac{|D_i\cap B(2^jR)|}{|B(2^jR)|},\qquad\osc_{ij}=\inf_{D_i\cap B(2^{j-1}R)}(w)-\inf_{D_i\cap B(2^jR)}(w).\]
		Using Lemma \ref{lemma-crit:phragmen-lindelof2} with $D=D_i\cap B(2^jR)$ and $\text{``$R$''}=2^jR$ and $\text{``$r$''}=2^{j-1}R$, we have
		\begin{equation}\label{eq-entire:use_PL}
			\osc_{ij}\geq\log\Big(\frac{\sqrt\pi}2\frac{2^{j-1}R}{|D_i\cap B(2^jR)\setminus B(2^{j-1}R)|^{1/2}}\Big)\geq\log\Big(\frac1{4\sqrt{\kappa_{ij}}}\Big).
		\end{equation}
		By Jensen's inequality and $\sum_i\kappa_{ij}\leq1$, we have
		\[\sum_{i=1}^M\osc_{ij}\geq M\log\big(\sqrt M/4\big),\qquad\forall\,j\geq1.\]
		Summing this over $j\in\{1,\cdots,J\}$ for $J\geq1$ and using \eqref{eq-entire:aux2} yields
		\[JM\log\big(\sqrt M/4\big)\leq\sum_{i=1}^M\Big(\inf_{D_i\cap B(R)}(w)-\inf_{D_i\cap B(2^JR)}(w)\Big)\leq-\sum_{i=1}^M\inf_{D_i\cap B(2^JR)}(w).\]
		This implies
		\[\inf_{B(2^JR)}(w)\leq-J\log\big(\sqrt M/4\big)<-Jn\log 4,\qquad\forall\,J\geq1.\]
		Since $w=-\log|\D u|$ on each $\p\{w<t\}$, this implies
		\[\Big\{|\D u|=\exp\big(Jn\log 4\big)\Big\}\cap B(2^JR)\ne\emptyset,\qquad\forall\,J\geq1.\]
		But this violates the gradient estimate $|\D u|\leq C\big(1+|z|^{n-1}\big)$ by \eqref{eq-entire:main_growth_cond}.
	\end{proof}
	
	Using a blow-down argument, we show that:
	
	\begin{claim}\label{claim-entire:Ft_bounded}
		$F_t\Subset\RR^2$ for all $t<0$.
	\end{claim}
	\begin{proof}
		Suppose $F_t$ is unbounded for some $t<0$. Let $\gamma_t$ be an unbounded component of $\p F_t$. Fix $z_1\in\gamma_t$. For each $R$, let $\gamma_{t,R}$ be the connected component of $\gamma_t\cap B(R)$ containing $z_1$, and $U$ be the connected component of $B(R)\setminus\gamma_{t,R}$ containing $0$. So $U$ is a concave polygon with boundary $\gamma_{t,R}$ in $B(R)$. Using Gauss-Bonnet on $U$, we have
		\[\begin{aligned}
			2\pi &= -\int_{\gamma_{t,R}}|\kappa|+\pi\cdot\big(\text{\# vertices of $\gamma_{t,R}$}\big) \\
			&\hspace{64pt} +R^{-1}|\p U\cap\p B(0,R)|+\sum\text{(exterior angles at $\bar{\gamma_{t,R}}\cap\p B(R)$)}.
		\end{aligned}\]
		It follows that $\int_{\gamma_{t,R}}|\kappa|\leq\pi N(n)+2\pi$, and taking $R\to\infty$, that
		\[\int_{\gamma_t}|\kappa|\leq \pi N(n)+2\pi.\]
		Hence, $\gamma_t$ must have an unbounded edge that is asymptotic to a ray. Denote this edge by $\sigma_t$, and up to a rotation, we may assume that $\sigma_t$ is asymptotic to the positive $x$-axis.
		
		Consider
		\[M(R)=\max_{\bar B(R)}|u|,\qquad\forall\,R\geq0.\]
		The polynomial growth \eqref{eq-entire:main_growth_cond} implies
		\[\lim_{R\to\infty}\Big[R^{-n-1}M(R)\Big]=0.\]
		Therefore, we may find a sequence $R_i\to\infty$ so that
		\begin{equation}\label{eq-entire:doubling}
			(\lambda R_i)^{-n-1}M(\lambda R_i)\leq R_i^{-n-1}M(R_i),\qquad\forall\,\lambda\in[1,\infty),\ i\in\ZZ_{\geq1}.
		\end{equation}
		Consider the rescaled functions
		\[u_i(z)=\frac{u(R_iz)}{M(R_i)}.\]
		Clearly $\sup_{B(1)}|u_i|=1$, and \eqref{eq-entire:doubling} implies $\sup_{B(\lambda)}|u_i|\leq\lambda^{n+1}$ for each $\lambda\geq1$, $i\in\ZZ_{\geq1}$. By the interior gradient estimate and Arzel\`a-Ascoli, we may take a subsequential limit
		\[u_\infty=\lim_{i\to\infty}u_i\qquad\text{in}\ \ C^1_{\loc}(\RR^2).\]
		Clearly $u_\infty(0)=0$ and $\sup_{B(1)}|u_\infty|=1$, and
		\begin{equation}\label{eq-entire:aux1}
			|\D u_i|=\frac{R_i|\D u(R_iz)|}{M(R_i)}=\frac{R_ie^{-t}}{M(R_i)}\qquad\text{on}\ \ R_i^{-1}\sigma_t.
		\end{equation}
		By comparison with cones, the function $M(R)/R$ is nondecreasing in $R$. Then the superlinear growth implies $\lim_{R\to\infty}\big[M(R)/R\big]=\infty$. So \eqref{eq-entire:aux1} implies that $|\D u_\infty|=0$ on the nonnegative $x$-axis, contradicting Theorem \ref{thm-crit:main}.
	\end{proof}
	
	In summary: each $F_t$ ($t<0$) is the connected component of $Y_t$ that contains $0$, it is bounded and has at most $N(n)$ edges. Hence $u$ has at most $N(n)$ critical points in $\RR^2$: otherwise, $F_t$ would contain more than $N(n)$ critical points for some $t\ll-1$, contradicting Theorem \ref{thm-closed:main}\ref{item-closed:main_edges}. Fix $T<0$ so that $F_T$ contains all the critical points of $u$.
	
	Observe that $Y_t$ is connected for all $t<T$: otherwise, some connected component $V$ of $Y_t$ does not contain critical points of $u$. Theorem \ref{thm-closed:main}\ref{item-closed:main_finiteness} then implies $V\not\Subset\RR^2$. This contradicts $V\subset F_s$ for some $s<t$ (indeed, take $R>0$ so that $V\cap B(R)\ne\emptyset$, then take $s<t$ so that $B(R)\subset F_s$, then $V\subset F_s$ by connectedness).
	
	Next, note that $Z_t\Subset\RR^2$, since $Z_t\subset Y_s$ for all $s<t$. Thus, Lemma \ref{lemma-closed:tech_connected} implies $Y_t=Z_t$ for all $t<T$. Now, the main condition \eqref{eq-closed:condition_Pt} is satisfied for all $t<T$. By Lemma \ref{lemma-closed:1_evolution_lemma}, the map
	\[\Big(\!-\log|\D u|,\frac{\D u}{|\D u|}\Big):\RR^2\setminus\tZ_T\to(-\infty,T)\times\SS^1\]
	lifts to a homeomorphism
	\[\Th:\RR^2\setminus\tZ_T\to(-\infty,T)\times\SS^1\big((N-2)\pi\big),\]
	where $N$ is the number of edges of $\p Y_t$ for all $t<T$. And, the support function
	\[h(t,\th)=\bmetric{\Th^{-1}(t,\th)}{\nu_\th}\]
	is a smooth solution of $\p_th=\p_{\th\th}h+h$, with $(\p_{\th\th}h+h)(t,\cdot)$ having $N$ nondegenerate roots.
	
	Set $d=N/2$ (recall that $N$ is even), so $\th\in\SS^1\big(2(d-1)\pi\big)$. Using Theorem \ref{thm-entire:ancient_sol} below to a suitable parabolic rescaling of $e^{-t}h$, and comparing the number of roots of $(\p_{\th\th}h+h)(t,\cdot)$ with $t\to-\infty$, it follows that
	\[h(t,\th)=\sum_{j=1}^d\exp\Big(t-\frac{j^2t}{(d-1)^2}\Big)\Big(a_j\sin\frac{j\th}{d-1}+b_j\cos\frac{j\th}{d-1}\Big)\qquad\forall\,t\in(-\infty,T),\]
	for some coefficients $a_j,b_j$ with $(a_d,b_d)\ne(0,0)$. Here, note that the constant term must vanish: indeed, the edges of $\p Z_t$ are all streamlines of $u$, thus the alternating length of edges of $\p Z_t$ must be zero, hence $\int h(t,\cdot)=0$ for each $t$, by \eqref{eq-prelim:alternating_length}.
	
	As $t\to-\infty$, the highest frequency term dominates:
	\[\exp\Big(\frac{2d-1}{(d-1)^2}t\Big)h(t,\th)\xrightarrow[t\to-\infty]{C^\infty(\SS^1(2(d-1)\pi))}a_d\sin\frac{d\th}{d-1}+b_d\cos\frac{d\th}{d-1}.\]
	Arguing similarly as in Theorem \ref{thm-closed:main}\ref{item-closed:main_asymp}, one may show that the rescaled solution
	\[u_\lambda(z)=\lambda^{-\frac{d^2}{2d-1}}u(\lambda z)\]
	converges in $C^1_{\loc}(\RR^2)$ to a quasiradial solution as $\lambda\to\infty$. This proves item \ref{item-entire:blow_down}. Using Theorem \ref{thm-closed:main}\ref{item-closed:main_edges} to $Y_t$ for any $t<T$, item \ref{item-entire:degree} follows.
\end{proof}

\begin{proof}[Proof of Theorem \ref{thm-entire:deg_2}] {\ }
	
	The growth condition implies that we must fall into the $d=2$ case of Theorem \ref{thm-entire:poly_growth}, hence from the proof, the set $\p Y_t$ has support function
	\[h(t,\th)=(a_1\sin\th+b_1\cos\th)+e^{-3t}(a_2\sin2\th+b_2\cos2\th)\qquad\forall\,t<T,\ \th\in\SS^1,\]
	for some $T\in\RR$. Hence, after a rigid motion and scaling, we may assume that
	\[h(t,\th)=\frac{27}{128}e^{-3t}\sin2\th,\qquad\forall\,t<T,\ \th\in\SS^1.\]
	Recall the recovery formula
	\[\D u\Big(h(t,\th)\nu_\th+\p_\th h(t,\th)\tau_\th\Big)=e^{-t}\tau_\th,\]
	see \eqref{eq-closed:recovery_formula_2} in the proof of Theorem \ref{thm-closed:main}\ref{item-closed:main_asymp}. For the case here this gives
	\[\D u\Big(\frac{27}{64}e^{-3t}\cos^3\th,-\frac{27}{64}e^{-3t}\sin^3\th\Big)=e^{-t}(\cos\th,\sin\th),\qquad\forall\,t<T,\th\in\SS^1.\]
	This implies
	\[\D u(x,y)=\frac43\big(x^{1/3},-y^{1/3}\big)\]
	for all $(x,y)\in\RR^2\setminus B(R)$, for some $R\gg1$. Thus up to an additive constant, we have
	\[u(x,y)=|x|^{4/3}-|y|^{4/3}\qquad\text{in}\ \ \RR^2\setminus B(R).\]
	By Jensen's maximum principle, this implies $u=|x|^{4/3}-|y|^{4/3}$ in the entire $\RR^2$.
\end{proof}

Finally, we state the finite frequency lemma that we needed in Theorem \ref{thm-entire:poly_growth}:

\begin{theorem}\label{thm-entire:ancient_sol}
	Suppose $f\in C^\infty\big((-\infty,0]\times\SS^1\big)$ solves the heat equation, such that $f(t,\cdot)$ has at most $2N$ roots for all $t\leq0$. Then there is $k\leq N$ and constants $a_j,b_j$ $(0\leq j\leq k)$, such that
	\[f(t,\th)=\sum_{j=0}^ke^{-j^2t}(a_j\cos j\th+b_j\sin j\th).\]
\end{theorem}

The theorem is a consequence of the following lemma:

\begin{lemma}\label{lemma-entire:sobolev_norm}
	Suppose $f\in C^\infty(\SS^1)$ has $2k$ nondegenerate roots. Let $a_j,b_j$ denote the real Fourier coefficients of $f$ (with the convention that $b_0=0$). Then
	\begin{equation}\label{eq-entire:quant_sturm_hurwitz}
		\|f\|_{W^{-2k-2,2}}\leq C(k)\sum_{j=0}^k\big(|a_j|+|b_j|\big).
	\end{equation}
\end{lemma}

\begin{proof}[Proof of Theorem \ref{thm-entire:ancient_sol} assuming \ref{lemma-entire:sobolev_norm}] {\ }
	
	By the Sturmian theory, there are $k\leq N$ and $t_0\leq0$ so that $f(t,\cdot)$ has exactly $2k$ nondegenerate roots in $\SS^1$ for all $t\leq t_0$. Let $a_j,b_j$ denote the Fourier coefficients of $f|_{t=0}$. Then the Fourier coefficients of $f|_{t=-\tau}$ are $e^{j^2\tau}a_j, e^{j^2\tau}b_j$. Lemma \ref{lemma-entire:sobolev_norm} implies
	\[\|f(t,\cdot)\|_{W^{-2k-2,2}(\SS^1)}\leq C(k)e^{k^2|t|}\sum_{j=0}^k\big(|a_j|+|b_j|\big),\qquad\forall\,t\leq t_0.\]
	By the smoothing effect of the heat equation, it follows that
	\[\|f(t,\cdot)\|_{L^2(\SS^1)}\leq C(k)\|f(t-1,\cdot)\|_{W^{-2k-2,2}(\SS^1)}\leq C\big(k,f|_{t=0}\big)e^{k^2|t|},\qquad\forall\,t\leq t_0.\]
	This implies that $f$ has vanishing Fourier frequencies beyond $k$, as desired.
\end{proof}

Lemma \ref{lemma-entire:sobolev_norm} may be viewed as a quantitative version of the Sturm-Hurwitz theorem. The latter states that if a nonzero function $f\in C^\infty(\SS^1)$ has vanishing Fourier coefficients with frequencies $\leq k$, then $f$ has at least $2k+2$ roots. This corresponds to Lemma \ref{lemma-entire:sobolev_norm} with $a_0=a_1=b_1=\cdots=b_k=0$. An elegant proof of this fact (assuming that all roots are nondegenerate) appears in the book of P\'olya-Szeg\"o \cite[II-141]{Polya-Szego}: if otherwise $f$ has $2l$ roots $\th_1,\cdots,\th_{2l}$ ($l\leq k$), then
\[P(\th)=\prod_{j=1}^{2l}\sin\frac{\th-\th_j}{2}\]
is a degree $l$ trigonometric polynomial on $\SS^1$, and satisfies $\int_{\SS^1}fP>0$ (since $f,P$ have the same sign everywhere). These contradict the fact that $f$ has vanishing low frequencies. Our proof is inspired by this argument.

% ==Optional==
% Note that $\|f\|_{L^2}$ cannot be bounded in terms of low frequencies and the number of roots; one counterexample is a family of functions that approximate $\delta_{-\epsilon}-\delta_{\epsilon}$. These functions have small low frequencies (since integrating against $\sin(k\th)$ and $\cos(k\th)$ do not detect nearly canceling $\delta$-functions), but have unbounded $L^2$ norms. In this direction, see also the work \cite{Steinerberger_2020} by Steinerberger.

\begin{proof}[Proof of Lemma \ref{lemma-entire:sobolev_norm}] {\ }
	
	If $k=0$, then $f$ is either positive or negative in $\SS^1$, hence
	\[\Big|\int_{\SS^1}f\varphi\Big|\leq\Big|\int_{\SS^1}f\Big|\cdot\|\varphi\|_{C^0}\leq C|a_0|\cdot\|\varphi\|_{W^{2,2}}.\]
	Now assume $k\geq1$, and $\th_1<\th_2<\cdots<\th_{2k}$ are the roots of $f$. It suffices to show that
	\begin{equation}\label{eq-sturmian:wts}
		\Big|\int_{\SS^1}f\varphi\Big|\leq C(k)\sum_{i=0}^k\big(|a_i|+|b_i|\big)\cdot\|\varphi\|_{C^{2k+1}},\qquad\forall\,\varphi\in C^{2k+1}(\SS^1).
	\end{equation}
	Indeed, \eqref{eq-entire:quant_sturm_hurwitz} then follows from the Sobolev inequality.
	
	It is always possible to find a point $\omega\in\SS^1$, such that $|\omega-\th_j|>\pi/4k$ for each $j$. We fix any valid choice of $\omega$. Consider the (well-defined) trigonometric polynomials
	\[P_j(\th)=\left\{\begin{aligned}
		& \sin\frac{\th-\omega}{2}\cdot\prod_{1\leq s\leq j}\sin\frac{\th-\th_s}2\qquad\text{(if $j$ is odd),} \\
		& \prod_{1\leq s\leq j}\sin\frac{\th-\th_s}2\qquad\text{(if $j$ is even),}
	\end{aligned}\right.\]
	and inductively,
	\[Q_1=\varphi(\th_1),\qquad Q_{j+1}=Q_j+\frac{\varphi(\th_{j+1})-Q_j(\th_{j+1})}{P_j(\th_{j+1})}P_j.\]
	So $P_j,Q_j$ are trigonometric polynomials with degree $\leq\lceil j/2\rceil$, with $Q_j(\th)=\varphi(\th)$ for all $\th\in\{\th_1,\cdots,\th_j\}$ and $j\in\{1,\cdots,2k\}$. In fact, $Q_j$ are the Newton interpolations of $\varphi$.
	
	We claim that the coefficients of $Q_j$ are uniformly bounded by $C(k)\|\varphi\|_{C^{2k+1}}$, for all $1\leq j\leq 2k$. This is clear for $Q_1$. Suppose that our claim holds for $Q_j$. Notice that $\varphi-Q_j$ takes zero value in the set $\{\th_1,\cdots,\th_j\}$. Notice the general fact: if $f$ is a $C^r$ section of a line bundle over $\SS^1$ and satisfies $f(0)=0$, then
	\begin{equation}\label{eq-sturmian:reg_f/x}
		\Big\|\frac{f(\th)}{\sin(\th/2)}\Big\|_{C^{r-1}(\SS^1)}\leq C(r)\|f\|_{C^r(\SS^1)}.
	\end{equation}
	Here, we view $\sin(\th/2)$ as a section of the non-orientable bundle over $\SS^1$. To show \eqref{eq-sturmian:reg_f/x}, we view both sides as periodic functions on $\RR$, and notice that
	\[f(x)=\int_0^x f'(y)\,dy=x\int_0^1 f'(xt)\,dt,\]
	hence
	\[\frac{f(x)}{\sin(x/2)}=\frac{x}{\sin(x/2)}\int_0^1 f'(xt)\,dt,\qquad\forall\,x\in[-\pi,\pi].\]
	Then \eqref{eq-sturmian:reg_f/x} follows since $x/\sin(x/2)$ is real analytic in a neighborhood of $[-\pi,\pi]$. Now if $j$ is even, then applying \eqref{eq-sturmian:reg_f/x} successively with center $\th_1,\cdots,\th_j$, we obtain
	\begin{equation}\label{eq-sturmian:reg_Q/P}
		\Big\|\frac{\varphi-Q_j}{P_j}\Big\|_{C^{2k+1-j}}\leq C(k)\|\varphi-Q_j\|_{C^{2k+1}}\leq C(k)\|\varphi\|_{C^{2k+1}},
	\end{equation}
	where the last inequality uses the induction hypothesis. If $j$ is odd, then applying \eqref{eq-sturmian:reg_f/x} at $\th_1,\cdots,\th_{j}$ gives
	\begin{equation}\label{eq-sturmian:reg_Q/P_2}
		\Big\|\frac{\varphi-Q_j}{\prod_{s\leq j}\sin(\th-\th_s)/2}\Big\|_{C^{2k+1-j}}\leq C(k)\|\varphi-Q_j\|_{C^{2k+1}}\leq C(k)\|\varphi\|_{C^{2k+1}}.
	\end{equation}
	Next, notice that
	\begin{equation}\label{eq-sturmian:bound_of_coef}
		\Big|\frac{\varphi(\th_{j+1})-Q_j(\th_{j+1})}{P_j(\th_{j+1})}\Big|\leq C(k)\|\varphi\|_{C^{2k+1}},\qquad\forall\,j\in\{1,\cdots,2k\}.
	\end{equation}
	Indeed, when $j$ is even, this is a direct consequence of \eqref{eq-sturmian:reg_Q/P}; if $j$ is odd, this follows from \eqref{eq-sturmian:reg_Q/P_2} and $|\omega-\th_{j+1}|>\pi/4k$. Our initial claim then follows from \eqref{eq-sturmian:bound_of_coef} and induction. As a byproduct, we find that \eqref{eq-sturmian:reg_Q/P} holds for $j=2k$.
	
	Finally, we estimate
	\[\begin{aligned}
		\Big|\int f\varphi\Big|
		&\leq \Big|\int fQ_{2k}\Big|+\int |f|\,|\varphi-Q_{2k}| \\
		&\leq \Big|\int fQ_{2k}\Big|+C(k)\|\varphi\|_{C^{2k+1}}\int |f|\,|P_{2k}| \qquad\text{(by \eqref{eq-sturmian:reg_Q/P})}\\
		&\leq \Big|\int fQ_{2k}\Big|+C(k)\|\varphi\|_{C^{2k+1}}\Big|\int fP_{2k}\Big|\qquad\text{($f,P_{2k}$ have the same sign)}.
	\end{aligned}\]
	Recall that the degrees of $P_{2k},Q_{2k}$ are no greater than $k$, and the coefficients of $P_{2k}$, $Q_{2k}$ are bounded respectively by $C(k)$ and $C(k)\|\varphi\|_{C^{2k+1}}$. Then \eqref{eq-sturmian:wts} follows.
\end{proof}

%\newpage
\appendix

\section{Estimates on planar \texorpdfstring{$p$}{p}-harmonic functions}\label{sec:kzz_es}

The aim of this section is to prove Lemma \ref{lemma-equi:estimates}. Let $v$ solve \eqref{eq-equi:reg_eq} with $p\in[3,\infty)$ and $\epsilon\in(0,1]$. Denote $J=(|\D v|^2+\epsilon^2)^{1/2}$. Hence $v$ solves $\div\big(J^{p-2}\D v\big)=0$. It is easy to see that
\begin{equation}\label{eq-kzz:dJ}
	\D J=J^{-1}\D^2 v\cdot\D v,\qquad \D J^2=2\D^2v\cdot\D v.
\end{equation}
The equation for $v$ can now be written as
\begin{equation}\label{eq-kzz:eq_for_v}
	\Delta_\infty v+\frac{J^2}{p-2}\Delta v=0.
\end{equation}

\vspace{3pt}

We prove \eqref{eq-equi:kzz_sobolev} \eqref{eq-equi:kzz_flatness} closely following the arguments in \cite{Koch-Zhang-Zhou_2019} and \cite{Lindgren-Lindqvist_2021}. Recall from \cite[Lemma 2.1]{Koch-Zhang-Zhou_2019} that for any $C^2$ function $v$ it holds
\begin{align}
	-\det\D^2v &= -\frac12\div\big(\Delta v\D v-\D^2 v\cdot\D v\big), \label{eq-kzz:det1}\\
	|\D v|^2(-\det\D^2v) &= |\D^2 v\cdot\D v|^2-\Delta v\Delta_\infty v, \label{eq-kzz:det2}
\end{align}
where all the expressions hold in the classical sense. Combining \eqref{eq-kzz:eq_for_v} \eqref{eq-kzz:det2} we have
\[|\D v|^2(-\det\D^2v)=|\D^2v\cdot\D v|^2+\frac{J^2}{p-2}(\Delta v)^2\geq0,\]
and hence
\[-\det\D^2v\geq\frac{|\D v|^2}{J^2}(-\det\D^2v)=|\D J|^2+\frac{(\Delta v)^2}{p-2}.\]
Combining this with \eqref{eq-kzz:det1}, we obtain
\begin{equation}\label{eq-kzz:test_identity}
	\int\Big(\underbrace{\frac12\Delta v\metric{\D v}{\D\varphi}}_{\termI}\underbrace{-\frac12\D^2 v(\D v,\D\varphi)}_{\termII}\Big)\geq\int\Big(|\D J|^2+\frac{(\Delta v)^2}{p-2}\Big)\varphi
\end{equation}
for any nonnegative test function $\varphi\in C^1_0$.

\begin{proof}[Proof of Lemma \ref{lemma-equi:estimates}(ii)]
	Test \eqref{eq-kzz:test_identity} with $\varphi=J^2\eta^2$ with $\eta$ a cutoff function. We compute
	\[\begin{aligned}
		\termI &= \int\Delta v\bmetric{\D v}{\eta^2\D^2 v\cdot\D v+J^2\eta\D\eta}
		= \int\eta^2\Delta v\Delta_\infty v+J^2\eta\Delta v\metric{\D v}{\D\eta} \\
		&= \int-\frac{J^2\eta^2}{p-2}(\Delta v)^2+J^2\eta\Delta v\metric{\D v}{\D\eta} \\
		&= \int-\frac{J^2\eta^2}{p-2}(\Delta v)^2 \\
		&\qquad\quad
		-\D v\cdot\Big(2\eta\metric{\D v}{\D\eta}\D^2 v\cdot\D v+J^2\metric{\D v}{\D\eta}\D\eta+J^2\eta\D^2 v\cdot\D\eta+J^2\eta\D^2\eta\cdot\D v\Big) \\
		&= \int-\frac{J^2\eta^2}{p-2}(\Delta v)^2+\frac{2J^2\eta}{p-2}\metric{\D v}{\D\eta}\Delta v-J^2\metric{\D v}{\D\eta}^2 \\
		&\qquad\qquad -J^2\eta\D^2 v(\D v,\D\eta)-J^2\eta\D^2\eta(\D v,\D v),
	\end{aligned}\]
	and
	\[\begin{aligned}
		\termII &= \int -\eta^2|\D^2 v\cdot\D v|^2-J^2\eta\D^2 v(\D v,\D\eta).
	\end{aligned}\]
	Hence, \eqref{eq-kzz:test_identity} implies
	\begin{equation}\label{eq-kzz:aux1}
		\begin{aligned}
			\int2\eta^2|\D^2v&\cdot\D v|^2+\frac{2J^2\eta^2}{p-2}(\Delta v)^2 \\
			&\leq \int\frac{2J^2\eta}{p-2}|\D v|\,|\D\eta|\,|\Delta v|
			+ 2J^2\eta|\D^2 v\cdot\D v|\,|\D\eta|
			+ J^2\eta|\D v|^2|\D^2\eta|.
		\end{aligned}
	\end{equation}
	Then using Young's inequality, \eqref{eq-kzz:aux1} implies
	\begin{equation}\label{eq-kzz:sobolev_for_J}
		\int\eta^2\big|\D|\D v|^2\big|^2=4\int\eta^2|\D^2 v\cdot\D v|^2\leq C\int J^4\big(|\D\eta|^2+|\eta|\,|\D^2\eta|\big).
	\end{equation}
	This implies \eqref{eq-equi:kzz_sobolev} by a suitable choice of $\eta$.
\end{proof}

\begin{proof}[Proof of Lemma \ref{lemma-equi:estimates}(iii)]
	Test \eqref{eq-kzz:test_identity} with $\varphi=(v-P)^2\eta^2$ with $\eta$ a cutoff function:
	\[\begin{aligned}
		\termI &= \int
		\Delta v(v-P)\eta^2\metric{\D v}{\D v-\D P}
		+ \Delta v(v-P)^2\eta\metric{\D v}{\D\eta} \\
		&= \int
		-\metric{\D v}{\D v-\D P}^2\eta^2
		- 2(v-P)\eta\metric{\D v}{\D\eta}\metric{\D v}{\D v-\D P} \\
		&\qquad -(v-P)\eta^2\D^2 v(\D v,2\D v-\D P)
		- 2(v-P)\eta\metric{\D v}{\D v-\D P}\metric{\D v}{\D\eta} \\
		&\qquad -(v-P)^2\metric{\D v}{\D\eta}^2
		- (v-P)^2\eta\D^2 v(\D v,\D\eta)
		- (v-P)^2\eta\D^2\eta(\D v,\D v) \\
		&\leq \int
		- \big(|\D v|^2-\metric{\D v}{\D P}\big)^2\eta^2
		+ 4|v-P|\cdot\eta\cdot|\D\eta|\cdot|\D v|^2\cdot\big(|\D v|+|\D P|\big) \\
		&\qquad +2|v-P|\cdot\eta^2\cdot|\D^2 v\cdot\D v|\cdot\big(|\D v|+|\D P|\big) \\
		&\qquad +|v-P|^2\cdot\eta\cdot|\D^2 v\cdot\D v|\cdot|\D\eta|
		+ |v-P|^2\cdot\eta\cdot|\D v|^2\cdot|\D^2\eta|.
	\end{aligned}\]
	By Cauchy-Schwarz, this gives
	\[\begin{aligned}
		\termI &\leq
		- \int\big(|\D v|^2-\metric{\D v}{\D P}\big)^2\eta^2 \\
		&\qquad +C\Big(\int|v-P|^2\big(|\D v|^2+|\D P|^2\big)\eta^2\Big)^{1/2}\Big(\int|\D v|^4|\D\eta|^2\Big)^{1/2} \\
		&\qquad +C\Big(\int|v-P|^2\big(|\D v|^2+|\D P|^2\big)\eta^2\Big)^{1/2}\Big(\int|\D^2 v\cdot\D v|^2\eta^2\Big)^{1/2} \\
		&\qquad +C\Big(\int|v-P|^4|\D\eta|^2\Big)^{1/2}\Big(\int|\D^2 v\cdot\D v|^2\eta^2\Big)^{1/2} \\
		&\qquad + C\Big(\int|v-P|^4\eta|\D^2\eta|\Big)^{1/2}\Big(\int|\D v|^4\eta|\D^2\eta|\Big)^{1/2}.
	\end{aligned}\]
	Moreover,
	\[\begin{aligned}
		\termII &\leq
		\int|\D^2 v\cdot\D v|\cdot\Big(|v-P|\cdot\eta^2\cdot\big(|\D v|+|\D P|\big)+|v-P|^2\cdot\eta\cdot|\D\eta|\Big) \\
		&\leq \Big(\int|\D^2 v\cdot\D v|^2\eta^2\Big)^{1/2}\Big(\int|v-P|^2\big(|\D v|^2+|\D P|^2\big)\eta^2\Big)^{1/2} \\
		&\qquad +\Big(\int|\D^2 v\cdot\D v|^2\eta^2\Big)^{1/2}\Big(\int|v-P|^4|\D\eta|^2\Big)^{1/2}.
	\end{aligned}\]
	Inserting these into \eqref{eq-kzz:test_identity}, then using \eqref{eq-kzz:sobolev_for_J}, then finally choosing an appropriate $\eta$, the inequality \eqref{eq-equi:kzz_flatness} follows.
\end{proof}

Next, we prove \eqref{eq-equi:grad_est} \eqref{eq-equi:es_flatness} closely following the estimates in \cite{Evans-Smart_2011b}. Let $v$ be as stated above, so
\begin{equation}\label{eq-es:p_laplacian_var}
	\Delta_\infty v+\frac{J^2}{p-2}\Delta v=0.
\end{equation}
Consider the associated linearized operator
\begin{equation}\label{eq-es:linearization}
	\L\psi:=\D^2\psi(\D v,\D v)+2\D^2 v(\D\psi,\D v)+\frac{J^2}{p-2}\Delta\psi+\frac{2\Delta v}{p-2}\metric{\D\psi}{\D v}.
\end{equation}
Notice the relation
\begin{equation}\label{eq-es:LXY}
	\L(XY)=X\L Y+Y\L X+2\metric{\D X}{\D v}\metric{\D Y}{\D v}+\frac{2J^2}{p-2}\metric{\D X}{\D Y},
\end{equation}
and $\L\p_xv=\L\p_yv=0$ by translation symmetry or explicit calculation. Then notice that
\begin{equation}\label{eq-es:LA_Ldv2}
	\L|\D v|^2=\L(v_x^2)+\L(v_y^2)=2|\D^2v\cdot\D v|^2+\frac{2J^2}{p-2}|\D^2v|^2.
\end{equation}
We may also calculate
\begin{equation}\label{eq-es:Lv}
	\L v=3\Delta_\infty v+\frac{J^2+2|\D v|^2}{p-2}\Delta v=\Big(3-\frac{J^2+2|\D v|^2}{J^2}\Big)\Delta_\infty v=2\epsilon^2J^{-2}\Delta_\infty v,
\end{equation}
and
\begin{equation}\label{eq-es:Ly}
	\L y = 2\D^2v(\D v,\p_y)+\frac{2\Delta v}{p-2}\p_yv
	= 2\D^2 v\big(\D v,\p_y-J^{-2}(\p_yv)\D v\big).
\end{equation}

\begin{proof}[Proof of Lemma \ref{lemma-equi:estimates}(i)] {\ }
	
	We first assume that $v$ solves \eqref{eq-es:p_laplacian_var}. Let $\eta\in C^2_0(B(z,3.9r))$ be a cutoff function with $\eta|_{B(z,3.5r)}\equiv1$ and
	\[|\D\eta|\leq Mr^{-1},\qquad |\D^2\eta|,|\Delta\eta|\leq Mr^{-2},\]
	where $M$ is a universal constant. We calculate and estimate term by term:
	\[\begin{aligned}
		\L\big(\eta^2|\D v|^2\big) &= 
		\eta^2\Big(2|\D^2v\cdot\D v|^2+\frac{2J^2}{p-2}|\D^2v|^2\Big) \\
		&\qquad +|\D v|^2\Big(\D^2\eta^2(\D v,\D v)
		+ 2\D^2 v(\D\eta^2,\D v)+\frac{J^2\Delta\eta^2}{p-2}
		+ \frac{2\Delta v}{p-2}\metric{\D\eta^2}{\D v}\Big) \\
		&\qquad +2\metric{\D|\D v|^2}{\D v}\metric{\D\eta^2}{\D v}
		+ \frac{2J^2}{p-2}\metric{\D|\D v|^2}{\D\eta^2} \\
		&\geq 2\eta^2|\D^2v\cdot\D v|^2+\frac{2J^2\eta^2}{p-2}|\D^2v|^2
		- 4M^2r^{-2}|\D v|^4\\
		&\qquad -4Mr^{-1}\eta|\D^2v\cdot\D v|\,|\D v|^2
		- \frac{4M^2r^{-2}J^2}{p-2}|\D v|^2
		- 4Mr^{-1}\eta|\D^2v\cdot\D v|\,|\D v|^2 \\
		&\qquad -8Mr^{-1}\eta|\D^2v\cdot\D v|\,|\D v|^2
		- \frac{8Mr^{-1}\eta J^2}{p-2}|\D^2 v\cdot\D v|.
	\end{aligned}\]
	Notice that
	\[16Mr^{-1}\eta|\D^2 v\cdot\D v|\,|\D v|^2\leq\eta^2|\D^2 v\cdot\D v|^2+64M^2r^{-2}|\D v|^4\]
	and
	\[\frac{8Mr^{-1}\eta J^2}{p-2}|\D^2 v\cdot\D v|\leq\frac{J^2\eta^2}{p-2}|\D^2v|^2+\frac{16M^2r^{-2}J^2}{p-2}|\D v|^2.\]
	Hence, we obtain
	\begin{equation}\label{eq-es:grad_aux1}
		\L\big(\eta^2|\D v|^2\big) \geq \eta^2|\D^2v\cdot\D v|^2-C_2r^{-2}|\D v|^4-\frac{C_2r^{-2}J^2}{p-2}|\D v|^2,
	\end{equation}
	for some universal constant $C_2$. Then, by \eqref{eq-es:Lv} \eqref{eq-es:LXY} we note that
	\begin{equation}\label{eq-es:grad_aux2}
		\L v^2\geq-4|v|\,|\D^2v\cdot\D v|\,|\D v|+2|\D v|^4+\frac{2J^2}{p-2}|\D v|^2.
	\end{equation}
	Denote $F=C_2^{-1}r^2\eta^2|\D v|^2+v^2$. Combining \eqref{eq-es:grad_aux1} \eqref{eq-es:grad_aux2} we have
	\[\eta^4\L F\geq C_2^{-1}r^2\eta^6|\D^2v\cdot\D v|^2+\eta^4|\D v|^4-4\eta^4|v|\,|\D^2v\cdot\D v|\,|\D v|.\]
	Notice that
	\[4\eta^4|v|\,|\D^2 v\cdot\D v|\,|\D v|\leq C_2^{-1}r^2\eta^6|\D^2v\cdot\D v|^2+4C_2r^{-2}v^2\eta^2|\D v|^2.\]
	We finally obtain
	\begin{equation}\label{eq-es:LF}
		\eta^4\L F\geq \eta^4|\D v|^4-4C_2r^{-2}v^2\eta^2|\D v|^2.
	\end{equation}
	Suppose $F$ is maximized at a point $z_0\in\bar{B(z,3.9r)}$. If $\eta(z_0)=0$, then for all $z'\in B(z,3.5r)$ we have
	\[C_2^{-1}r^2|\D v(z')|^2\leq F(z')\leq F(z_0)=v(z_0)^2\leq\sup_{B(z,4r)}|v|^2.\]
	If $\eta(z_0)>0$, then $z_0$ is an interior critical point of $F$, and $\L F(z_0)\leq0$. By \eqref{eq-es:LF}, this implies $\eta(z_0)^2|\D v(z_0)|^2\leq4C_2r^{-2}v(z_0)^2$. Hence for all $z'\in B(z,3.5r)$ we have
	\[C_2^{-1}r^2|\D v(z')|^2\leq F(z')\leq F(z_0)\leq C_2^{-1}r^2\cdot 4C_2r^{-2}v(z_0)^2+v(z_0)^2.\]
	In either case, we obtain $|\D v|^2\leq 8C_2r^{-2}\sup_{B(z,4r)}|v|^2$ in $B(z,3.5r)$. This proves the desired estimate for solutions of \eqref{eq-equi:reg_eq}. It is known that $p$-harmonic functions ($3\leq p<\infty$) can be locally approximated by solutions of \eqref{eq-equi:reg_eq}, hence the estimate passes to $p$-harmonic functions as well. Finally, the $\infty$-harmonic case is well-known \cite{Crandall-Evans-Gariepy_2001}. 
\end{proof}

\begin{proof}[Proof of Lemma \ref{lemma-equi:estimates}(iv)] {\ }
	
	Recall that $0<\lambda\leq a$, and $v$ solves \eqref{eq-equi:reg_eq} and satisfies $|v-ay|\leq\lambda r$ in $B(4r)$. By the gradient estimate \eqref{eq-equi:grad_est}, there is a universal constant $C_3$ so that
	\[|\D v|\leq C_3a\qquad\text{in}\ \ B(3.5r).\]
	Denote $\Lambda=C_3a$. Let $\eta\in C^2_0(B(3r))$ be a cutoff function, so that $\eta|_{B(r)}\equiv1$ and
	\[|\D\eta|\leq Mr^{-1},\qquad |\D^2\eta|,|\Delta\eta|\leq Mr^{-2}\qquad\text{in}\ \ B(3r),\]
	for some universal constant $M$.
	
	Combining \eqref{eq-es:Lv} \eqref{eq-es:Ly} we have in $B(3.5r)$
	\[\begin{aligned}
		\big|\L(v-ay)\big| &= \Big|2(\epsilon^2J^{-2})\Delta_\infty v-2a\D^2 v\big(\D v,\p_y-J^{-2}(\p_yv)\D v\big)\Big|
		\leq 4(\Lambda+a)|\D^2 v\cdot\D v|.
	\end{aligned}\]
	Denote $A=|\D v|^2-a\p_yv$. Thus in $B(3.5r)$
	\begin{equation}\label{eq-es:L_vy}
		\begin{aligned}
			\L(v-ay)^2 &= 2(v-ay)\L(v-ay)+2\metric{\D v-a\p_y}{\D v}^2+\frac{2J^2}{p-2}|\D v-a\p_y|^2 \\
			&\geq -8(\Lambda+a)\lambda r|\D^2 v\cdot\D v|+2A^2.
		\end{aligned}
	\end{equation}
	Using \eqref{eq-es:LA_Ldv2} \eqref{eq-es:LXY} and $\L\p_y v=0$, we further compute
	\[\begin{aligned}
		\L(\eta^2A^2) &= \eta^2\Big(2A\L A+2\metric{\D A}{\D v}^2+\frac{2J^2}{p-2}|\D A|^2\Big) \\
		&\qquad +A^2\Big(
				\D^2\eta^2(\D v,\D v)
				+ 2\D^2 v(\D\eta^2,\D v)
				+ \frac{J^2}{p-2}\Delta\eta^2
				+ \frac{2\Delta v}{p-2}\metric{\D\eta^2}{\D v}\Big) \\
		&\qquad +2\metric{\D A^2}{\D v}\metric{\D\eta^2}{\D v}
			+ \frac{2J^2}{p-2}\metric{\D A^2}{\D\eta^2} \\
		&\geq 4\eta^2A|\D^2v\cdot\D v|^2
			+ \frac{4J^2\eta^2A}{p-2}|\D^2v|^2
			+ 2\eta^2\metric{\D A}{\D v}^2
			+ \frac{2J^2\eta^2}{p-2}|\D A|^2 \\
		&\qquad -A^2\Big(
				4M^2r^{-2}\Lambda^2
				+ 4Mr^{-1}\eta|\D^2v\cdot\D v|
				+ \frac{4M^2r^{-2}J^2}{p-2}
				+ 4Mr^{-1}\eta|\D^2v\cdot\D v|\Big) \\
		&\left.\begin{aligned}
			&\qquad -\Big[2\eta^2\metric{\D A}{\D v}^2+8A^2\metric{\D\eta}{\D v}^2\Big]\quad \\
			&\qquad -\Big[\frac{2J^2\eta^2}{p-2}|\D A|^2+\frac{8J^2A^2}{p-2}|\D\eta|^2\Big].
		\end{aligned}\right\}\text{(by Young's inequality)}
	\end{aligned}\]
	If $A>0$, then we use Young's inequality to deduce
	\[\begin{aligned}
		8Mr^{-1}\eta A^2|\D^2v\cdot\D v| &\leq 4\eta^2A|\D^2v\cdot\D v|^2+4M^2r^{-2}A^3 \\
		&\leq 4\eta^2A|\D^2v\cdot\D v|^2+8M^2r^{-2}\Lambda^2A^2,
	\end{aligned}\]
	and notice that
	\[8A^2\metric{\D\eta}{\D v}^2\leq 8M^2r^{-2}\Lambda^2A^2,\qquad \frac{12M^2r^{-2}}{p-2}J^2A^2\leq 36M^2r^{-2}(\epsilon^2p^{-1}+\Lambda^2)A^2.\]
	Inserting these into the main expression, we obtain
	\begin{equation}\label{eq-es:L_et2A2}
		\L(\eta^2A^2) \geq -CM^2r^{-2}(\epsilon^2p^{-1}+\Lambda^2)A^2\qquad\text{ in $B(3.5r)$, if }A>0,
	\end{equation}
	for a universal constant $C$. Denote $C_4=CM^2(\epsilon^2p^{-1}+\Lambda^2)$. Consider the quantity
	\[F=\eta^2\max\{A,0\}^2+C_4r^{-2}(v-ay)^2+C_4^{1/2}\lambda|\D v|^2.\]
	By \eqref{eq-es:LA_Ldv2} \eqref{eq-es:L_vy} \eqref{eq-es:L_et2A2} and Young's inequality, whenever $A>0$ we have
	\begin{equation}\label{eq-es:LF_2}
		\begin{aligned}
			\L F &\geq C_4r^{-2}A^2-8(\Lambda+a)C_4r^{-1}\lambda|\D^2v\cdot\D v|+2C_4^{1/2}\lambda|\D^2v\cdot\D v|^2 \\
			&\geq C_4r^{-2}\big[A^2-8(\Lambda+a)^2C_4^{1/2}\lambda\big].
		\end{aligned}
	\end{equation}
	Suppose $\max_{\bar{B(3.5r)}}(F)$ is attained at a point $z_1\in\bar{B(3.5r)}$. If $A(z_1)\leq0$ or $\eta(z_1)=0$, then for all $z\in B(r)$ with $A(z)>0$, we have
	\[A^2(z)=\eta^2A^2(z)\leq F(z)\leq F(z_1)\leq C_4\lambda^2+C_4^{1/2}\Lambda^2\lambda\leq\big(C_4a+C_4^{1/2}\Lambda^2\big)\lambda.\]
	If $A(z_1)>0$ and $\eta(z_1)>0$ both hold, then $z_1$ must be an interior critical point of $F$, and $\L F(z_1)\leq0$, and \eqref{eq-es:LF_2} is applicable to $z_1$. It follows that $A(z_1)^2\leq8(\Lambda+a)^2C_4^{1/2}\lambda$. Hence for all $z\in B(r)$ with $A(z)>0$, we have
	\[A^2(z)\leq F(z)\leq F(z_1)\leq A(z_1)^2+C_4\lambda^2+C_4^{1/2}\Lambda^2\lambda\leq \big(8(\Lambda+a)^2C_4^{1/2}+C_4a+C_4^{1/2}\Lambda^2\big)\lambda.\]
	Inserting the definitions of $C_4,\Lambda$, we obtain an estimate of the form
	\[A\leq C(\epsilon p^{-1/2}+a)^{3/2}\lambda^{1/2}\qquad\text{in}\ \ B(r),\]
	for a universal constant $C$, as desired.
\end{proof}

%\newpage

\noindent \textit{Department of Mathematics, University of California, Berkeley.}

\vspace{3pt}

\noindent \textit{Email:} \href{mailto:kaixu@berkeley.edu}{kaixu@berkeley.edu}
	
\end{document}